# INTERVAL GROUPOIDS

**W. B. Vasantha Kandasamy**
**Florentin Smarandache**
**Moon Kumar Chetry**

**2010**

# INTERVAL GROUPOIDS

W. B. Vasantha Kandasamy
Florentin Smarandache
Moon Kumar Chetry

# CONTENTS

# PREFACE

This book introduces several new classes of groupoid, like polynomial groupoids, matrix groupoids, interval groupoids, polynomial interval groupoids, matrix interval groupoids and their neutrosophic analogues.

Interval groupoid happens to be the first non-associative structure constructed using intervals built using $Z_n$ or Z or Q or R or $Z^+ \cup \{0\}$ or $Q^+ \cup \{0\}$ and so on.

This book has five chapters. Chapter one is introductory in nature. In chapter two new classes of groupoids and interval groupoids are defined and described.

The analogous neutrosophic study is carried out in chapter three. The applications of this new structure is given in chapter four. The final chapter suggests more than 200 problems.

This book has given 77 new definitions, 426 examples of these new notions and over 150 theorems.

The authors deeply acknowledge Dr. Kandasamy for the proof reading and Meena and Kama for the formatting and designing of the book.

W.B.VASANTHA KANDASAMY
FLORENTIN SMARANDACHE
MOON KUMAR CHETRY

Chapter One

# PRELIMINARY NOTIONS

In this chapter we just give the basic definition of the groupoid which forms the first section. In section two we just recall the properties associated with neutrosophy.

## 1.1 Groupoids

In this section we recall the definition of groupoid and give some examples. For more about groupoids and its properties please refer [17-20].

**DEFINITION 1.1.1:** *Given an arbitrary set P a mapping of $P \times P$ into P is called a binary operation on P. Given such a mapping $\sigma: P \times P \to P$ we use it to define a product $*$ in P by declaring $a * b = c$ if $\sigma(a, b) = c$.*

**DEFINITION 1.1.2:** *A non empty set of elements G is said to form a groupoid if in G is defined a binary operation called the product denoted by * such that a * b ∈ G for all a, b ∈ G.*

It is important to mention here that the binary operation * defined on the set G need not be associative that is (a * b) * c ≠ a * (b * c) in general for all a, b, c ∈ G, so we can say the groupoid (G, *) is a set on which is defined a non associative binary operation which is closed on G.

A groupoid G is said to be a commutative groupoid if for every a, b ∈ G we have a * b = b * a. A groupoid G is said to have an identity element e in G if a * e = e * a = a for all a ∈ G.

We call the order of the groupoid G to be the number of distinct elements in it denoted by o(G) or |G|. If the number of elements in G is finite we say the groupoid G is of finite order or a finite groupoid otherwise we say G is an infinite groupoid.

**DEFINITION 1.1.3:** *Let (G, *) be a groupoid a proper subset H ⊂ G is a subgroupoid if (H, *) is itself a groupoid.*

**DEFINITION 1.1.4:** *A groupoid G is said to be a Moufang groupoid if it satisfies the Moufang identity (xy) (zx) = (x(yz))x for all x, y, z in G.*

**DEFINITION 1.1.5:** *A groupoid G is said to be a Bol groupoid if G satisfies the Bol identity ((xy) z) y = x ((yz) y) for all x, y, z in G.*

**DEFINITION 1.1.6:** *A groupoid G is said to be a P-groupoid if (xy) x = x (yx) for all x, y ∈ G.*

**DEFINITION 1.1.7:** *A groupoid G is said to be right alternative if it satisfies the identity (xy) y = x (yy) for all x, y ∈ G. Similarly we define G to be left alternative if (xx) y = x (xy) for all x, y ∈ G.*

**DEFINITION 1.1.8:** *A groupoid G is alternative if it is both right and left alternative, simultaneously.*

**DEFINITION 1.1.9:** *A groupoid G is said to be an idempotent groupoid if $x^2 = x$ for all $x \in G$.*

**DEFINITION 1.1.10:** *Let G be a groupoid. P a non empty proper subset of G, P is said to be a left ideal of the groupoid G if 1) P is a subgroupoid of G and 2) For all $x \in G$ and $a \in P$, $xa \in P$.*

*One can similarly define right ideal of the groupoid G. P is called an ideal if P is simultaneously a left and a right ideal of the groupoid G.*

**DEFINITION 1.1.11:** *Let G be a groupoid A subgroupoid V of G is said to be a normal subgroupoid of G if*

1. *aV = Va*
2. *(Vx)y = V(xy)*
3. *y(xV) = (yx)V*

*for all x, y, a $\in$ V.*

**DEFINITION 1.1.12:** *A groupoid G is said to be simple if it has no non trivial normal subgroupoids.*

*Example 1.1.1:* The groupoid G given by the following table is simple.

| *     | $a_0$ | $a_1$ | $a_2$ | $a_3$ | $a_4$ | $a_5$ | $a_6$ |
|-------|-------|-------|-------|-------|-------|-------|-------|
| $a_0$ | $a_0$ | $a_4$ | $a_1$ | $a_5$ | $a_2$ | $a_6$ | $a_3$ |
| $a_1$ | $a_3$ | $a_0$ | $a_4$ | $a_1$ | $a_5$ | $a_2$ | $a_6$ |
| $a_2$ | $a_6$ | $a_3$ | $a_0$ | $a_4$ | $a_1$ | $a_5$ | $a_2$ |
| $a_3$ | $a_2$ | $a_6$ | $a_3$ | $a_0$ | $a_4$ | $a_1$ | $a_5$ |
| $a_4$ | $a_5$ | $a_2$ | $a_6$ | $a_3$ | $a_0$ | $a_4$ | $a_1$ |
| $a_5$ | $a_1$ | $a_5$ | $a_2$ | $a_6$ | $a_3$ | $a_0$ | $a_4$ |
| $a_6$ | $a_4$ | $a_1$ | $a_5$ | $a_2$ | $a_6$ | $a_3$ | $a_0$ |

It is left for the reader to verify (G, ∗) = {$a_0, a_1, a_2, \ldots, a_6$, ∗} has no normal subgroupoids. Hence, G is simple.

**DEFINITION 1.1.13:** *A groupoid G is normal if*

1. *xG = Gx*
2. *G(xy) = (Gx)y*
3. *y(xG) = (yx)G*

*for all x, y ∈ G.*

**DEFINITION 1.1.14:** *A Smarandache groupoid G is a groupoid which has a proper subset S, S ⊂ G such that S under the operations of G is a semigroup.*

**DEFINITION 1.1.15:** *Let G be a Smarandache groupoid (SG) if the number of elements in G is finite we say G is a finite SG, otherwise G is said to be an infinite SG.*

**DEFINITION 1.1.16:** *Let G be a Smarandache groupoid. G is said to be a Smarandache commutative groupoid if there is a proper subset, which is a semigroup, is a commutative semigroup.*

For more about groupoids and Smarandache groupoids please refer [20]

## 1.2 Introduction to Neutrosophic Algebraic Structures

In this section we just recall some basic neutrosophic algebraic structures essential to make this book a self contained one. For more refer [8-16].

In this section we assume fields to be of any desired characteristic. We denote the indeterminacy by 'I', as i will make a confusion, as it denotes the imaginary value, viz. $i^2 = -1$ that is $\sqrt{-1} = i$. The indeterminacy I is such that $I \cdot I = I^2 = I$.

Here we recall the notion of neutrosophic groups. neutrosophic groups in general do not have group structure.

**DEFINITION 1.2.1:** *Let (G, \*) be any group, the neutrosophic group is generated by I and G under \* denoted by N(G) = {⟨G ∪ I⟩, \*}.*

***Example 1.2.1:*** Let $Z_7$ = {0, 1, 2, …, 6} be a group under addition modulo 7. N(G) = {⟨$Z_7$ ∪ I⟩, '+' modulo 7} is a neutrosophic group which is in fact a group. For N(G) = {a + bI / a, b ∈ $Z_7$} is a group under '+' modulo 7. Thus this neutrosophic group is also a group.

***Example 1.2.2:*** Consider the set G = $Z_5$ \ {0}, G is a group under multiplication modulo 5. N(G) = {⟨G ∪ I⟩, under the binary operation, multiplication modulo 5}.

N(G) is called the neutrosophic group generated by G ∪ I. Clearly N(G) is not a group, for $I^2$ = I and I is not the identity but only an indeterminate, but N(G) is defined as the neutrosophic group.

Thus based on this we have the following theorem:

**THEOREM 1.2.1:** *Let (G, \*) be a group, N(G) = {⟨G ∪ I⟩, \*} be the neutrosophic group.*
  1.  *N(G) in general is not a group.*
  2.  *N(G) always contains a group.*

*Proof:* To prove N(G) in general is not a group it is sufficient we give an example; consider ⟨$Z_5$ \ {0} ∪ I⟩ = G = {1, 2, 4, 3, I, 2 I, 4 I, 3 I}; G is not a group under multiplication modulo 5. In fact {1, 2, 3, 4} is a group under multiplication modulo 5. N(G) the neutrosophic group will always contain a group because we generate the neutrosophic group N(G) using the group G and I. So G $\subsetneq$ N(G); hence N(G) will always contain a group.

Now we proceed onto define the notion of neutrosophic subgroup of a neutrosophic group.

**DEFINITION 1.2.2:** *Let $N(G) = \langle G \cup I \rangle$ be a neutrosophic group generated by G and I. A proper subset P(G) is said to be a neutrosophic subgroup if P(G) is a neutrosophic group i.e. P(G) must contain a (sub) group of G.*

*Example 1.2.3*: Let $N(Z_2) = \langle Z_2 \cup I \rangle$ be a neutrosophic group under addition. $N(Z_2) = \{0, 1, I, 1 + I\}$. Now we see $\{0, I\}$ is a group under + in fact a neutrosophic group $\{0, 1 + I\}$ is a group under '+' but we call $\{0, I\}$ or $\{0, 1 + I\}$ only as pseudo neutrosophic groups for they do not have a proper subset which is a group. So $\{0, I\}$ and $\{0, 1 + I\}$ will be only called as pseudo neutrosophic groups (subgroups).

We can thus define a pseudo neutrosophic group as a neutrosophic group, which does not contain a proper subset which is a group. Pseudo neutrosophic subgroups can be found as a substructure of neutrosophic groups. Thus a pseudo neutrosophic group though has a group structure is not a neutrosophic group and a neutrosophic group cannot be a pseudo neutrosophic group. Both the concepts are different.

Now we see a neutrosophic group can have substructures which are pseudo neutrosophic groups which is evident from the following example:

*Example 1.2.4:* Let $N(Z_4) = \langle Z_4 \cup I \rangle$ be a neutrosophic group under addition modulo 4. $\langle Z_4 \cup I \rangle = \{0, 1, 2, 3, I, 1 + I, 2I, 3I, 1 + 2I, 1 + 3I, 2 + I, 2 + 2I, 2 + 3I, 3 + I, 3 + 2I, 3 + 3I\}$. $o(\langle Z_4 \cup I \rangle) = 4^2$.

Thus neutrosophic group has both neutrosophic subgroups and pseudo neutrosophic subgroups. For $T = \{0, 2, 2 + 2I, 2I\}$ is a neutrosophic subgroup as $\{0\ 2\}$ is a subgroup of $Z_4$ under addition modulo 4. $P = \{0, 2I\}$ is a pseudo neutrosophic group under '+' modulo 4.

**DEFINITION 1.2.3:** *Let K be the field of reals. We call the field generated by $K \cup I$ to be the neutrosophic field for it involves the indeterminacy factor in it.*

*We define $I^2 = I$, $I + I = 2I$ i.e., $I + ... + I = nI$, and if $k \in K$ then $k.I = kI$, $0I = 0$. We denote the neutrosophic field by $K(I)$ which is generated by $K \cup I$ that is $K(I) = \langle K \cup I \rangle$. $\langle K \cup I \rangle$ denotes the field generated by $K$ and $I$.*

***Example 1.2.5:*** Let R be the field of reals. The neutrosophic field of reals is generated by R and I denoted by $\langle R \cup I \rangle$ i.e. R(I) clearly $R \subset \langle R \cup I \rangle$.

***Example 1.2.6:*** Let Q be the field of rationals. The neutrosophic field of rationals is generated by Q and I denoted by Q(I).

**DEFINITION 1.2.4:** *Let $K(I)$ be a neutrosophic field we say $K(I)$ is a prime neutrosophic field if $K(I)$ has no proper subfield, which is a neutrosophic field.*

***Example 1.2.7:*** Q(I) is a prime neutrosophic field where as R(I) is not a prime neutrosophic field for $Q(I) \subset R(I)$.

**DEFINITION 1.2.5:** *Let $K(I)$ be a neutrosophic field, $P \subset K(I)$ is a neutrosophic subfield of P if P itself is a neutrosophic field. $K(I)$ will also be called as the extension neutrosophic field of the neutrosophic field P.*

We can also define neutrosophic fields of prime characteristic p (p is a prime).

**DEFINITION 1.2.6:** *Let $Z_p = \{0,1, 2, ..., p – 1\}$ be the prime field of characteristic p. $\langle Z_p \cup I \rangle$ is defined to be the neutrosophic field of characteristic p. Infact $\langle Z_p \cup I \rangle$ is generated by $Z_p$ and I and $\langle Z_p \cup I \rangle$ is a prime neutrosophic field of characteristic p.*

***Example 1.2.8:*** $Z_7 = \{0, 1, 2, 3, ..., 6\}$ be the prime field of characteristic 7. $\langle Z_7 \cup I \rangle = \{0, 1, 2, ..., 6, I, 2I, ..., 6I, 1 + I, 1 + 2I, ..., 6 + 6I \}$ is the prime field of characteristic 7.

**DEFINITION 1.2.7:** *Let $G(I)$ by an additive abelian neutrosophic group and K any field. If $G(I)$ is a vector space over K then we call $G(I)$ a neutrosophic vector space over K.*

Elements of these neutrosophic fields will also be known as neutrosophic numbers. For more about neutrosophy please refer [8-15]. We see $Z_nI = \{aI \mid a \in Z_n\}$ is a neutrosophic field called pure neutrosophic field. Likewise QI, RI and $Z_pI$ are neutrosophic fields where p is a prime. Thus $Z_5I = \{0, I, 2I, 3I, 4I\}$ is a pure neutrosophic field.

Pure neutrosophic structures as QI or RI or $Z_nI$ cannot contain any real numbers. However $0.I = 0$, so 0 belongs to them.

**Chapter Two**

# NEW CLASSES OF GROUPOIDS

This chapter introduces seven new classes of groupoids and has seven sections. Section one introduces the new class of matrix groupoids using $Z_n$ or Z or Q or R or C. Polynomial groupoids of five levels are introduced in section two. Special interval groupoids are introduced in section three. Polynomial interval groupoids are introduced in section four. Section five introduces interval matrix groupoids. Smarandache interval groupoids are introduced in section six and in the final section we give the new classes of groupoids built using intervals.

## 2.1 New Classes of Matrix Groupoids

In this section we introduce a new class of matrix groupoids and neutrosophic matrix groupoids and analyse a few of their properties. Throughout this book $(x_1, x_2, \ldots, x_n)$ represents a row matrix with entries $x_i \in Z_n$ or any ring or field; $1 \leq i \leq n$. This row matrix will be known as usual or real row matrix. If $x_i \in Z_nI$ or $N(Z_n)$ or $QI$ or $ZI$ or $RI$ or $N(Q)$ or $N(Z)$ or $N(R)$ then we call the row matrix $(x_1, \ldots, x_n)$ to be a neutrosophic row matrix. Likewise we define

$$\begin{bmatrix} y_1 \\ \vdots \\ y_m \end{bmatrix}$$

if $y_i \in Z_n$ or $Q$ or $R$ or $Z$ to be a usual column matrix and if $y_i \in Z_nI$ or $QI$ or $RI$ or $N(Z_n)$ or $N(Q)$ or $N(R)$ as column neutrosophic matrix $1 \leq i \leq m$.

A matrix $M_{n \times m} = (m_{ij})$ with $m_{ij}$ in $Z_n$ or $Q$ or $R$ or $Z$ will be known as the real matrix and if they belong to $N(Z_n)$ or $ZI$ or $Z_nI$ or $N(Z)$ or $N(Q)$ or $QI$ or $N(R)$ or $RI$ will be known as the neutrosophic matrix. With this understanding we proceed onto describe and define some new classes of groupoids.

**DEFINITION 2.1.1:** *Let $G = \{(x_1, \ldots, x_n) \mid x_i \in Z_m; 1 \leq i \leq n\}$ $m \geq 3$ be a collection of $1 \times n$ row matrices with entries from the modulo integer $Z_m$. Define * a binary operation on G as follows*
$(x_1, \ldots, x_n) * (y_1, \ldots, y_n)$
$\quad = \quad t(x_1, \ldots, x_n) + u(y_1, \ldots, y_n)$
$\quad = \quad (tx_1 + uy_1(mod\ m),\ tx_2 + uy_2(mod\ m),\ \ldots,\ (tx_n + uy_n)$
$\qquad (mod\ m))$
*where $t, u \in Z_n \setminus \{0\}$, $t \neq u$ and $(t, u) = 1$ for all $(x_1, x_2, \ldots, x_n)$, $(y_1, y_2, \ldots, y_n) \in G$. We define $(G, *, (t, u))$ to be a row matrix groupoid using $Z_n$.*

We will illustrate this situation by some simple examples.

***Example 2.1.1:*** Let $G = \{(x_1, x_2, x_3) \mid x_i \in Z_4; 1 \leq i \leq 3\}$. Define * on G by $(x_1, x_2, x_3) * (y_1, y_2, y_3) = 2(x_1, x_2, x_3) + 3(y_1, y_2, y_3)$ where $2, 3 \in Z_4 \setminus \{0\}$ and $(2, 3) = 1$.
Take

$(3, 2, 1) * (1, 0, 3) \quad = \quad 2(3, 2, 1) + 3(1, 0, 3)$
$\qquad \qquad \qquad \quad = \quad (2\ 0\ 2) + (3\ 0\ 1)$
$\qquad \qquad \qquad \quad = \quad (1\ 0\ 3).$

Consider (3 2 1) * [(1 0 3) * (0 2 2)] and [(3 2 1) * (1 0 3)] * (0 2 2). We see

$$[(3\ 2\ 1) * (1\ 0\ 3)] * (0\ 2\ 2) \quad = \quad (1\ 0\ 3) * (0\ 2\ 2)$$
$$= \quad (2\ 0\ 2) + (0\ 2\ 2)$$
$$= \quad (2\ 2\ 0) \quad\quad\quad (1)$$

$$(3\ 2\ 1) * [(1\ 0\ 3) * (0\ 2\ 2)] \quad = \quad (3\ 2\ 1) * [(2\ 0\ 2) + (0\ 2\ 2)]$$
$$= \quad (3\ 2\ 1) * (2\ 2\ 0)$$
$$= \quad (2\ 0\ 2) + (2\ 2\ 0)$$
$$= \quad (0\ 2\ 2) \quad\quad\quad (2)$$

$(2\ 2\ 0) \neq (0\ 2\ 2)$. Thus the operation * is non associative. Thus $(G, *, (2, 3))$ is a $1 \times 3$ row matrix groupoids built using $Z_4$.

We see $1 \times 3$ row matrix groupoids using $Z_4$ with $(t, u) = (1, 2)$ or $(2, 1)$ or $(1, 3)$ or $(3, 1)$ or $(2, 3)$ or $(3, 2)$.

Thus we have 6 distinct $1 \times 3$ row matrix groupoids built using $Z_4$.

***Example 2.1.2:*** Let $G = \{(x_1, x_2, x_3, x_4, x_5, x_6) \mid x_i \in Z_7; 1 \leq i \leq 6\}$. Take $(t, u) = 5, 6)$; $(G, *, (t, u)) = (G, (5, 6), *)$ is a $1 \times 6$ row matrix groupoid built using $Z_7$. We have 22 distinct $1 \times 6$ row matrix groupoids built using $Z_7$. It is left for the reader to find the number of $1 \times m$ matrix groupoids built using $Z_n$.

We can in the definition 2.1.1 replace $Z_n$ by $Z$ or $Q$ or $R$. In these cases we will get infinite number of row matrix groupoids of infinite order. By order of a groupoid G we mean the number of distinct elements in G. If G has finite number of elements we call G a finite groupoid and if G has infinite number of elements then we call G to be an infinite groupoid. We will just illustrate them by some simple examples.

***Example 2.1.3:*** Let $G = \{(x_1, x_2, \ldots, x_{12}) \mid x_i \in Z, 1 \leq i \leq 12\}$; take $(t, u) = (15, -16)$ $\{G, * (15, -16)\}$ is a $1 \times 12$ row matrix groupoid constructed using $Z$. We see the order of G is infinite. Further we have infinite number of pairs $(t, u)$ with $t, u \in Z$, $(t \neq u, (t, u) = 1)$. Thus there are infinite number of $1 \times 12$ row matrix groupoids built using $Z$. It is left for the reader to check whether G is associative.

***Example 2.1.4:*** Let $G = \{(x_1, x_2, \ldots, x_7) / x_i \in R; 1 \le i \le 7\}$. Take $(t, u) = (0.7, -1.52) \in R \times R$. $\{G, *, (t, u)\}$ is a $1 \times 7$ row matrix groupoid built using reals R.

Clearly cardinality of G is infinite and we can construct infinite number of $1 \times 7$ row matrix groupoids using R.

Next we proceed on to build a new class of column matrix groupoids using $Z_n$ or R or Q or Z.

**DEFINITION 2.1.2:** *Let $G = \{(x_1, x_2, \ldots, x_n)^t \mid x_i \in Z_m; 1 \le i \le n\}$ be the collection of all $n \times 1$ column matrices with entries from $Z_m$. Choose $t, u \in Z_m \setminus \{0\}; t \ne u, (t, u) = 1$. For $(x_1, \ldots, x_n)^t$ and $(y_1, \ldots, y_n)^t \in G$. Define*

$$(x_1, \ldots, x_n)^t * (y_1, \ldots, y_n)^t = \begin{bmatrix} x_1 \\ \vdots \\ x_n \end{bmatrix} * \begin{bmatrix} y_1 \\ \vdots \\ y_n \end{bmatrix} = \begin{bmatrix} tx_1 \\ \vdots \\ tx_n \end{bmatrix} + \begin{bmatrix} uy_1 \\ \vdots \\ uy_n \end{bmatrix}$$

$$= \begin{bmatrix} tx_1 + uy_1 (\bmod m) \\ \vdots \\ tx_n + uy_n (\bmod m) \end{bmatrix} = \begin{bmatrix} z_1 \\ \vdots \\ z_n \end{bmatrix} \in G.$$

*Thus $(G, (t, u), *)$ is a groupoid. This groupoid will be known as the $n \times 1$ column matrix groupoid built using $Z_m$.*

We will illustrate this situation by some examples.

***Example 2.1.5:*** Let

$$G = \left\{ \begin{bmatrix} x_1 \\ x_2 \\ x_3 \\ x_4 \\ x_5 \end{bmatrix} \middle| x_i \in Z_{12}; 1 \le i \le 5 \right\}$$

be the collection of all $5 \times 1$ column matrices with entries from $Z_{12}$. Choose $t = 2$ and $u = 3$; $(t, u) = 1$. $(G, (2, 3), *)$ is a $5 \times 1$ column matrix groupoid. Select

$$a = \begin{bmatrix} 10 \\ 2 \\ 0 \\ 0 \\ 1 \end{bmatrix}, b = \begin{bmatrix} 3 \\ 2 \\ 11 \\ 0 \\ 0 \end{bmatrix} \text{ and } c = \begin{bmatrix} 1 \\ 0 \\ 9 \\ 2 \\ 0 \end{bmatrix} \in G.$$

To prove $(a * b) * c \neq a * (b * c)$.

Consider

$$(a * b) * c = \left( \begin{bmatrix} 10 \\ 2 \\ 0 \\ 0 \\ 1 \end{bmatrix} * \begin{bmatrix} 3 \\ 2 \\ 11 \\ 0 \\ 0 \end{bmatrix} \right) * \begin{bmatrix} 1 \\ 0 \\ 9 \\ 2 \\ 0 \end{bmatrix}$$

$$= \left( \begin{bmatrix} 8 \\ 4 \\ 0 \\ 0 \\ 2 \end{bmatrix} + \begin{bmatrix} 9 \\ 6 \\ 9 \\ 0 \\ 0 \end{bmatrix} \right) * \begin{bmatrix} 1 \\ 0 \\ 9 \\ 2 \\ 0 \end{bmatrix} = \begin{bmatrix} 5 \\ 10 \\ 9 \\ 0 \\ 2 \end{bmatrix} * \begin{bmatrix} 1 \\ 0 \\ 9 \\ 2 \\ 0 \end{bmatrix}$$

$$= \begin{bmatrix} 10 \\ 8 \\ 6 \\ 0 \\ 4 \end{bmatrix} + \begin{bmatrix} 3 \\ 0 \\ 3 \\ 6 \\ 0 \end{bmatrix} = \begin{bmatrix} 1 \\ 8 \\ 9 \\ 6 \\ 4 \end{bmatrix} \tag{1}$$

Now $a * (b * c)$

$$= \begin{bmatrix} 10 \\ 2 \\ 0 \\ 0 \\ 1 \end{bmatrix} * \left( \begin{bmatrix} 3 \\ 2 \\ 11 \\ 0 \\ 0 \end{bmatrix} * \begin{bmatrix} 1 \\ 0 \\ 9 \\ 2 \\ 0 \end{bmatrix} \right) = \begin{bmatrix} 10 \\ 2 \\ 0 \\ 0 \\ 1 \end{bmatrix} * \left( \begin{bmatrix} 6 \\ 4 \\ 10 \\ 0 \\ 0 \end{bmatrix} + \begin{bmatrix} 3 \\ 0 \\ 3 \\ 6 \\ 0 \end{bmatrix} \right)$$

$$= \begin{bmatrix} 10 \\ 2 \\ 0 \\ 0 \\ 1 \end{bmatrix} * \begin{bmatrix} 9 \\ 4 \\ 1 \\ 6 \\ 0 \end{bmatrix} = \begin{bmatrix} 8 \\ 4 \\ 0 \\ 0 \\ 2 \end{bmatrix} * \begin{bmatrix} 3 \\ 0 \\ 3 \\ 6 \\ 0 \end{bmatrix} = \begin{bmatrix} 11 \\ 4 \\ 3 \\ 6 \\ 2 \end{bmatrix} \quad (2)$$

Clearly

$$\begin{bmatrix} 1 \\ 8 \\ 9 \\ 6 \\ 4 \end{bmatrix} \neq \begin{bmatrix} 11 \\ 4 \\ 3 \\ 6 \\ 2 \end{bmatrix}.$$

Hence (G, *, (2, 3)) is a groupoid as the operation * in general is non associative.

*Example 2.1.6:* Let

$$G = \left\{ \begin{bmatrix} x_1 \\ x_2 \\ x_3 \end{bmatrix} \middle| x_i \in Z_{19}; 1 \leq i \leq 3 \right\}$$

be the collection of all 3 × 1 column matrices with entries from $Z_{19}$. Choose (t, u) = (8, 9). It is easily verified (G, *, (8, 9)) is a 3 × 1 column matrix groupoid built using $Z_{19}$, integers modulo 19. However in the definition 2.1.2 we can replace $Z_n$ by Q or Z or R and still (G, *, (t, u)) will be the m × 1 column matrix groupoid. The only difference will be those built using $Z_n$ will be column matrix groupoids of finite order where as the groupoids built using Z or R or Q will be infinite groupoids.

Further the number of $n \times 1$ column matrix groupoids using $Z_m$ (for a fixed n and m) will be finite where as the number of $n \times 1$ column matrix groupoids built using Z or Q or R will be infinite in number.

Just we will give only examples of a column matrix groupoid of infinite cardinality before we proceed onto define the notion of $m \times n$ matrix groupoids $m \neq 1$, $n \neq 1$ (m = n or m $\neq$ n can occur).

*Example 2.1.7:* Let

$$G = \left\{ \begin{bmatrix} x_1 \\ x_2 \\ x_3 \\ x_4 \end{bmatrix} \middle| x_i \in Z; 1 \leq i \leq 4 \right\}$$

be the collection of all $4 \times 1$ column matrices with entries from the set Z. Choose (t, u) = (5, –12). Clearly (G, (5, –12), *) is a $4 \times 1$ column matrix groupoid of infinite order.

The reader is left with the task of verifying * on G is non associative in general.

*Example 2.1.8:* Let

$$G = \left\{ \begin{bmatrix} x_1 \\ x_2 \\ x_3 \\ x_4 \\ x_5 \\ x_6 \\ x_7 \\ x_8 \end{bmatrix} \middle| x_i \in Q; 1 \leq i \leq 8 \right\}$$

be the collection of all $8 \times 1$ column matrices with entries from Q. Define * using (t, u) = (–7/8, 5/19) $\in$ Q $\times$ Q. Clearly (G, *, (–7/8, 5/19)) is a 8 × 1 column matrix groupoid. It is easily verified that (G, *, (–7/8, 5/19)) is a groupoid of infinite cardinality.

Now we proceed onto define substructures in these groupoids.

**DEFINITION 2.1.3:** *Let (G, *, (t, u)) be a row (column) matrix groupoid H ⊆ G (H be a proper subset of G); we call (H, (t,u), *) to be a row (column) matrix subgroupoid of G if (H, (t, u), *) is itself a row (column) matrix groupoid of G.*

We will illustrate this be examples.

*Example 2.1.9:* Let $G = \{(x_1, x_2, x_3, x_4, x_5, x_6) \mid x_i \in Z_{12}, 1 \leq i \leq 6\}$ be the collection of all 1 × 6 row matrices. Choose $(t, u) = (5, 8) \in Z_{12} \times Z_{12}$. (G, (t, u), *) is a row matrix groupoid built using $Z_{12}$. Take $H = \{(x, x, x, x, x, x) / x \in Z_{12}\} \subseteq G$; {H, *, (5, 8)} is a row matrix subgroupoid of {G, *, (5, 8)}.

*Example 2.1.10:* Let $G = \{(x_1, x_2, x_3, x_4) \mid x_i \in Q, 1 \leq i \leq 4\}$ be a 1 × 4 row matrices with entries from Q. Take $(t, u) = (7/2, -5)$; (G, (7/2, –5), *) is a 1 × 4 row matrix subgroupoid. Take $H = \{(x, 0, y, z_1) \mid x, y, z_1 \in Q\} \subseteq G$; {H, *, (7/2, – 5)} is a 1 × 4 row matrix subgroupoid of {G, *, (7/2, – 5)}.

*Example 2.1.11:* Let

$$G = \left\{ \begin{bmatrix} a_1 \\ a_2 \\ a_3 \\ a_4 \\ a_5 \\ a_6 \end{bmatrix} \middle| a_i \in Z_{17}; 1 \leq i \leq 6 \right\}$$

be the collection of all 6 × 1 column matrices with entries from $Z_{17}$. {G, *, (9, 8)} is a 6 × 1 column matrix groupoid built using $Z_{17}$. Take

$$H = \left\{ \begin{bmatrix} 0 \\ 0 \\ 0 \\ a_1 \\ a_2 \\ a_3 \end{bmatrix} \middle| a_i \in Z_{17}; 1 \leq i \leq 3 \right\} \subseteq G;$$

{H, *, (9, 8)} is a 6 × 1 column matrix subgroupoid of G.

*Example 2.1.12:* Let

$$G = \left\{ \begin{bmatrix} a_1 \\ a_2 \\ a_3 \\ a_4 \\ a_5 \end{bmatrix} \middle| a_i \in R; 1 \leq i \leq 5 \right\}$$

be the collection of all 5 × 1 column matrices with entries from R. (G, $(-\sqrt{2}, \sqrt{17})$, *) is a 5 × 1 column matrix groupoid built using R. Take

$$H = \left\{ \begin{bmatrix} a \\ a \\ a \\ a \\ 0 \end{bmatrix} \middle| a \in R \right\} \subseteq G.$$

(H, $(-\sqrt{2}, \sqrt{17})$, *) $\subseteq$ (G, $(-\sqrt{2}, \sqrt{17})$, *) is a 5 × 1 column matrix subgroupoid of G built using R.

**DEFINITION 2.1.4:** *Let $G = M_{n \times m} = \{(m_{ij}); 1 \leq i \leq n, 1 \leq j \leq m$ be a $n \times m$ matrix with entries from $Z_n$ or R or Q or Z}. Choose t, u $\in Z_n$ or R or Q or Z such that $t \neq u$, (t, u) = 1. Define * on G; for two matrices $M = (m_{ij})$ and $N = (n_{ij})$ as*
$M * N = (m_{ij}) * (n_{ij})$

$$= (tm_{ij} + un_{ij})$$
$$= P.$$

P is a n × m matrix in G. Thus (G, (t, u), *) is a n × m matrix groupoid built using $Z_n$ (or R or Q or Z).

We will illustrate this situation by examples.

*Example 2.1.13:* Let

$$G = \left\{ \begin{bmatrix} a & b & e \\ c & d & f \end{bmatrix} \middle| a,b,c,d,e,f \in Z_9 \right\}$$

be the collection of all 2 × 3 matrices with entries from $Z_9$. Take (t, u) = (5, 6); we see {G, *, (5, 6)} is a 2 × 3 matrix is a groupoid with entries from $Z_9$.

*Example 2.1.14:* Let

$$G = \left\{ \begin{bmatrix} a_1 & a_2 \\ a_3 & a_4 \\ a_5 & a_6 \\ a_7 & a_8 \end{bmatrix} \middle| a_i \in Q; 1 \le i \le 8 \right\}$$

be the collection of all 4 × 2 matrices with entries from Q. Take (u, v) = (7, –3/2) from Q × Q. (G, (u, v), *) is a 4 × 2 matrix groupoid with entries from Q.

*Example 2.1.15:* Let G = {All 9 × 9 upper triangular matrices with entries from R}. Choose (t, u) = ($\sqrt{3}, -\sqrt{7}/2$). It is easily verified (G, ($\sqrt{3}, -\sqrt{7}/2$), *) is a 9 × 9 matrix groupoid with entries from R.

*Example 2.1.16:* Let

$$G = \left\{ \begin{bmatrix} a_1 & a_2 \\ 0 & a_3 \end{bmatrix} \middle| a_i \in Z_3; 1 \le i \le 3 \right\}.$$

Choose (t, u) = (2, 1), (G, *, (2,1)} is a 2 × 2 matrix groupoid.
The operation carried out in G is in general non associative.

Let
$$A = \begin{pmatrix} 2 & 1 \\ 0 & 1 \end{pmatrix}$$

and
$$B = \begin{pmatrix} 1 & 2 \\ 0 & 1 \end{pmatrix}$$

in G.

$$(A*B) = 2\begin{pmatrix} 2 & 1 \\ 0 & 1 \end{pmatrix} + \begin{pmatrix} 1 & 2 \\ 0 & 1 \end{pmatrix}$$

$$= \begin{pmatrix} 1 & 2 \\ 0 & 2 \end{pmatrix} + \begin{pmatrix} 1 & 2 \\ 0 & 1 \end{pmatrix} = \begin{pmatrix} 2 & 1 \\ 0 & 0 \end{pmatrix} \in G.$$

$$\text{Take } C = \begin{pmatrix} 2 & 2 \\ 0 & 2 \end{pmatrix} \in G$$

$$(A * B) * C = \begin{pmatrix} 2 & 1 \\ 0 & 0 \end{pmatrix} * \begin{pmatrix} 2 & 2 \\ 0 & 2 \end{pmatrix}$$

$$= \begin{pmatrix} 1 & 2 \\ 0 & 0 \end{pmatrix} + \begin{pmatrix} 2 & 2 \\ 0 & 2 \end{pmatrix}$$

$$= \begin{pmatrix} 0 & 1 \\ 0 & 2 \end{pmatrix} \quad (1)$$

$$A* (B*C) = \begin{pmatrix} 2 & 1 \\ 0 & 1 \end{pmatrix} * \left[\begin{pmatrix} 2 & 1 \\ 0 & 2 \end{pmatrix} + \begin{pmatrix} 2 & 2 \\ 0 & 2 \end{pmatrix}\right]$$

$$\begin{pmatrix} 2 & 1 \\ 0 & 1 \end{pmatrix} * \begin{pmatrix} 1 & 0 \\ 0 & 1 \end{pmatrix} = \begin{pmatrix} 1 & 2 \\ 0 & 2 \end{pmatrix} + \begin{pmatrix} 1 & 0 \\ 0 & 1 \end{pmatrix} = \begin{pmatrix} 2 & 2 \\ 0 & 0 \end{pmatrix} \quad (2)$$

$$\begin{pmatrix} 0 & 1 \\ 0 & 2 \end{pmatrix} \neq \begin{pmatrix} 2 & 2 \\ 0 & 0 \end{pmatrix}.$$

Thus the operation * defined in G in general is non associative.

As in case of row (and column) matrix groupoids we can define the notion of matrix subgroupoids. We will only give one example of a matrix subgroupoid.

*Example 2.1.17:* Let

$$G = \left\{ \begin{bmatrix} a & b & g \\ c & d & h \\ e & f & i \end{bmatrix} \middle| a,b,g,c,d,h,e,f,i \in Z \right\}.$$

{G, *, (3, -2)} is a matrix groupoid built over Z.
Let

$$H = \left\{ \begin{bmatrix} a & a & a \\ a & a & a \\ a & a & a \end{bmatrix} \middle| a \in Z \right\} \subseteq G;$$

{H, *, (3, –2)} is a matrix subgroupoid of {G, *, (3, –2)}.

We see all matrix groupoids built using $Z_n$ (n < ∞) modulo integers are finite where as all matrix groupoids built using Z or R or Q are nZ are of infinite order.

We as in case of general groupoids say a matrix groupoid G is a matrix P-groupoid if (AB)A = A (BA), for all A, B ∈ G.

We say a matrix groupoid G is said to be right alternative if (AB)B = A(BB) for all A, B ∈ G.

We will call a matrix groupoid G to be an idempotent groupoid if $A^2 = A$ for all A ∈ G.

As in case of general groupoids we in case of matrix groupoids also define the notion of right ideal, left ideal and ideal.

We also call a matrix subgroupoid V of G to be a normal subgroupoid if

$$aV = Va$$
$$(Vx)\,y = V(xy)$$
$$y\,(xV) = (yx)\,V$$

for all $x, y, a \in V$.

A matrix groupoid G is normal if

$$xG = Gx$$
$$G\,(xy) = (Gx)y$$
$$y\,(xG) = (yx)G$$

for all $x, y \in G$.

Now we define yet another class of matrix groupoids.

**DEFINITION 2.1.5:** *Let $G = \{$row matrix or column matrix or a $m \times n$ matrix with entries from $Z_n$ or $Q$ or $Z$ or $R\}$ 'or' is used in a mutually exclusive sense.*

*Now choose $t, u \in Z_n$ (or $Z$ or $Q$ or $R$) such that $t \neq u$ but $(t, u) \neq 1$. Then if for $x, y \in G$ define $x * y = tx + uy$ then $(G, (t, u), *)$ is a matrix groupoid of type I which is different from those defined earlier.*

We illustrate it by some examples.

*Example 2.1.18:* Let $G = \{$all $2 \times 2$ upper triangular matrices with entries from $Z_8\}$. Choose $(t, u) = (2, 4)$; $2, 4 \in Z_8$. Now $\{G, *, (2, 4)\}$ is a matrix groupoid.

We will just show how * is defined.
Take

$$A = \begin{pmatrix} 3 & 5 \\ 0 & 1 \end{pmatrix} \text{ and } B = \begin{pmatrix} 2 & 3 \\ 0 & 7 \end{pmatrix} \in G.$$

$$A * B = 2\begin{pmatrix} 3 & 5 \\ 0 & 1 \end{pmatrix} + 4\begin{pmatrix} 2 & 3 \\ 0 & 7 \end{pmatrix}$$

$$= \begin{pmatrix} 6 & 2 \\ 0 & 2 \end{pmatrix} + \begin{pmatrix} 0 & 4 \\ 0 & 4 \end{pmatrix} = \begin{pmatrix} 6 & 6 \\ 0 & 6 \end{pmatrix} \in G.$$

We see (G, (2, 4), *) is a finite groupoid. If we do not permit (t, u) = 1 we see this class of matrix groupoids i.e., matrix groupoids of type I form a disjoint class from the matrix groupoids in which (t, u) = 1.

Now as in case of the other matrix groupoids of type I we in case of these groupoids also define all the properties without any modifications.

*Example 2.1.19:* Let

$$G = \left\{ \begin{bmatrix} x_1 \\ x_2 \\ x_3 \\ x_4 \\ x_5 \end{bmatrix} \middle| x_i \in Z; i = 1, 2, 3, 4, 5 \right\}$$

choose (t, u) = (5, 10); 5, 10 ∈ Z. {G, *, (5, 10)} is a matrix groupoid of type I of infinite order.

Now we proceed onto define matrix groupoids of type II.

**DEFINITION 2.1.6:** *Let G = {collection of all row matrices or column matrices or m × n matrices with entries by $Z_n$ or Q or Z or R} 'or' is used in a mutually exclusive sense.*

*Now choose t, u such that u = t. Then (G, *, (t, t)) is another new class of matrix groupoids which we choose to call as matrix groupoids of type II. Clearly class of type I matrix groupoids are disjoint from the class of type II groupoids. Further they are also disjoint from the usual class of matrix groupoids.*

We will give some examples and discuss a few properties about them.

*Example 2.1.20:* Let G = {$(x_1, x_2, x_3, x_4, x_5, x_6, x_7)$ / $x_i \in Z_{21}$, $1 \leq i \leq 7$} be $1 \times 7$ row matrix built using $Z_{21}$.

Take t = 8; {G, *, (8, 8)} is a row matrix groupoid of type II.

Z = (1, 1, 3, 2, 2, 0 1), (3, 2, 0 1, 20, 18, 7) = X and Y = (1, 20, 4, 0, 7, 17, 3) ∈ G.

Z* (X * Y)
  = Z* [(3, 2, 0, 1, 20, 18, 7) * (1, 20, 4, 0, 7, 17, 3)]
  = Z* (3, 16, 0, 8, 13, 18, 14) + (8, 13, 3, 0, 14, 10, 3)]
  = (1, 1, 3, 2, 2, 0, 1) * (11, 8, 3, 8, 6, 7,17)
  = (8, 8, 3, 16, 16, 0, 8) + (11, 3, 3, 3, 6, 14, 10)
  = (19, 11, 6, 19, 1, 14, 18)         (1)

(Z*X) * Y
  = {(1, 1, 3, 2, 2, 0, 1) * (3, 2 0, 1, 20, 18, 7)} *
    (1, 20, 4, 0, 7, 17, 3)
  = [(8, 8, 3, 16, 16, 0, 8) + (3 16 0 8 13 18 14)] *
    (1, 20, 4, 0, 7, 17, 3)
  = (11, 3, 3, 3, 8, 18, 1) + (8, 13, 11, 0, 14, 10, 3)
  = (19, 16, 14, 3, 1, 7, 4)           (2)

We see (19, 11, 6, 19, 1, 14, 18) ≠ (19, 16, 14, 3, 1, 7, 4). Thus * is non associative.

Now we derive some properties of these matrix groupoids of type II.

**THEOREM 2.1.1:** *The matrix groupoids (G, *, (t, t)) are matrix P-groupoids. G is the collection of row (or column or m × n matrix) with entries from $Z_n$.*

*Proof:* Let $A = (m_{ij})$ and $B = (n_{ij})$ be row (or column or m × n) matrices with entries from $Z_n$. Let $t \in Z_n$. To prove G is a P-matrix groupoid, we have to prove (A * B) * A = A * (B * A),
  Now
(A * B) * A   = $(m_{ij}) * (n_{ij}) * (m_{ij})$
           = $(tm_{ij}) + (tn_{ij}) * (m_{ij})$
           = $(t^2 m_{ij}) + (t^2 n_{ij}) + (tm_{ij})$         (1)

Consider
A * (B*A)    = $(m_{ij}) * (tn_{ij}) + (tm_{ij})$
           = $tm_{ij} + t^2 n_{ij} + t^2 m_{ij}$         (2)

We see
$$A * (B*A) = (A*B) * A$$
Thus (G, (t, t), *) is a matrix P-groupoid of type II.

**COROLLARY 2.1.1:** *We see if $G = \{(m_{ij}) / m_{ij} \in Z \text{ or } Q \text{ or } R\}$ {G, *, (t, t)) is a matrix P-groupoid of type II.*

Proof is left as an exercise to the reader.

**THEOREM 2.1.2:** *{G, *, (t, t)}; G is row (or column or m × n) matrix with entries from $Z_p$, p a prime. {G, *, (t, t); 1 < t < n} is not an alternative matrix groupoid.*

*Proof:* To show {G, *, (t, t); 1 < t < n} is not an alternative matrix groupoid we have to show $(x * y) * y \neq x (y * y)$ for some x, y ∈ G. ($x = (m_{ij})$ and $y = (n_{ij})$). Consider

$$
\begin{aligned}
(x * y) * y &= ((m_{ij}) * (n_{ij})) * (n_{ij}) \\
&= (t(m_{ij}) + t(n_{ij})) * (n_{ij}) \\
&= t^2(m_{ij}) + t^2(n_{ij}) + t(n_{ij}) \\
&= t^2(m_{ij}) + (t^2 + t)(n_{ij}) \quad \text{(I)}
\end{aligned}
$$

$$
\begin{aligned}
x * (y * y) &= (m_{ij}) * (t\, n_{ij} + t\, n_{ij}) \\
&= t m_{ij} + t^2 n_{ij} + t^2 n_{ij} \\
&= t m_{ij} + 2t^2 (n_{ij}) \quad \text{(II)}
\end{aligned}
$$

Now I and II are identical if and only if $t^2 \equiv t \pmod{p}$. But this is impossible as p is a prime. Hence I and II are never identical for 1 < t < n.
Thus (G, * (t, t,)) is not a matrix alternative groupoid of type II.

**THEOREM 2.1.3:** *Let (G, * (t, t)) be a matrix groupoid of type two built using $Z_n$ (G, *, (t, t)) is an alternative groupoid if and only if $t^2 \equiv t \pmod{n}$.*

*Proof:* Let $A = (m_{ij})$ and $B = (n_{ij})$ be any two matrices in G entries of A and B are from $Z_n$.

$$(A * B) * B = [(m_{ij}) * (n_{ij})] * (n_{ij})$$
$$= (tm_{ij} + tn_{ij}) \times (n_{ij})$$
$$= (t^2 m_{ij} + t^2 n_{ij}) + tn_{ij}$$
$$= (tm_{ij} + tn_{ij} + t\, n_{ij}) \text{ as } t^2 \equiv t \pmod{p}$$

$$A * (B * B) = (m_{ij}) * (n_{ij} * n_{ij})$$
$$= (m_{ij}) * (tn_{ij} + tn_{ij})$$
$$= t\, m_{ij} + t^2 n_{ij} + t^2 n_{ij}$$
$$= tm_{ij} + tn_{ij} + tn_{ij} \text{ (as } t^2 = t \pmod{p}\text{)}.$$

Thus $A * (B * B) = (A * B) * B$. Hence $(G, (t, t), *)$ is an alternative matrix groupoid of type II. Also all matrix groupoids of type II are commutative.

Next we proceed onto define the notion of matrix groupoids of type III.

**DEFINITION 2.1.7:** *Let G = {set of all column matrices or row matrices or m × n matrices with entries from the field Q or ring Z or R or $Z_n$) 'or' is in the mutually exclusive sense. We define for any A, B ∈ G (A, B can be a row matrix or column matrix or a m × n matrix with entries for $Z_n$ or Q or Z or R) 'or' used in both the places only in the mutually exclusively sense. A * B = tA + uB where u or t is zero u, t ∈ $Z_n$ (or Z or R or Q) we define (G, *, (t, u); u = 0, or t = 0) to a matrix groupoid of type III.*

We will illustrate this by the following examples.

***Example 2.1.21:*** Let $G = \{(a_1, a_2, a_3, a_4, a_5, a_6) \mid a_i \in Z_8; 1 \le i \le 6\}$ be the set of all $1 \times 6$ row matrices with entries from $Z_8$. Define * on G as follows choose $(3, 0) = (t, u)$ and for A, B ∈ G.

Define $A * B = 3A + 0B$ then $(G, (3,0), *)$ is a matrix groupoid of type III. We see for $(5, 3, 2, 1, 6, 7) = A$ and $B = (1, 0, 6, 5, 3, 2)$

$$A * B = 3 (5, 3, 2, 1, 6, 7) + 0 (1\ 0\ 6\ 5\ 3\ 2)$$
$$= (7, 1, 6, 3, 2, 5) \in G.$$

Take $C = (2, 1, 0, 7, 5, 3)$; now

(A * B) * C  = 3 (7, 1, 6, 3, 2, 5) + 0 (2, 1, 0, 7, 5, 3)
  = (5, 3, 2, 1, 6, 7).  I

Consider
A * (B * C)  = A * (3 (1 0 6 5 3 2 ) + 0 (2, 1, 0, 7, 5, 3))
  = (5, 3, 2, 1, 6, 7) * (3, 0, 2, 7, 1, 6)
  = (7, 1, 6, 3, 2, 5)  II

A * (B * C) ≠ (A * B) * C evident from I and II. Thus (G, *, (3, 0)) is a matrix groupoid of type III.

**THEOREM 2.1.4:** *(G, *, (0, t)) where G is row (column or m × n matrix) matrix with entries from $Z_n$; n is not a prime. Then (G, *, (0, t)) is a matrix P-groupoid of type III if and only if $t^2 \equiv t \pmod{n}$*

*Proof:* Consider A, B ∈ G

A * (B * B)  = $(m_{ij}) * [(n_{ij}) * (n_{ij})]$
  = $(m_{ij}) * (t.n_{ij})$ (where A = $(m_{ij})$ and B = $(n_{ij})$)
  = $t^2 n_{ij}$
  = $t n_{ij}$ if and only if $t^2 \equiv t \pmod{n}$.

Consider
(A * B) * B  = (0 + tB) * B
  = 0.t B + tB = tB
  = $tn_{ij}$.

Thus A * (B * B) = (A * B) * B if and only if $t^2 \equiv t \pmod{n}$. Hence (G, *, (0, t)) is a matrix P-groupoid of type III.

**THEOREM 2.1.5:** *Let {G, (t, u), *} be the matrix groupoid with (t, u) = 1 with entries from $Z_n$. This matrix groupoid is a semigroup if and only if $t^2 \equiv t \pmod{n}$ and $u^2 \equiv u \pmod{n}$; t, u ∈ $Z_n \setminus \{0\}$ and (t, u) = 1.*

*Proof:* Given G = {all m × n matrices, 1 ≤ m < ∞ and 1 ≤ n < ∞ with entries from $Z_n$}.

Choose $t, u \in Z_n \setminus \{0\}$; $t^2 \equiv t \pmod{n}$, $u^2 \equiv u \pmod{n}$ and $(t, u) = 1$. To show $(G, (t, u), *)$ is a matrix semigroup.

Let $A = (a_{ij})$, $B = (b_{ij})$ and $C = (c_{ij})$ be three matrices in G. Now

$$\begin{aligned}
(A * B) * C &= ((a_{ij}) * (b_{ij})) * (c_{ij}) \\
&= ((ta_{ij}) + (ub_{ij})) * (c_{ij}) \\
&= (t^2 a_{ij}) + (tu\, b_{ij}) + (uc_{ij}) \\
&= ta_{ij} + (tub_{ij}) + uc_{ij} \qquad \text{I}
\end{aligned}$$

as $t^2 \equiv t \pmod{n}$.

Consider

$$\begin{aligned}
A * (B * C) &= (a_{ij}) * ((b_{ij}) * (c_{ij})) \\
&= (a_{ij}) * ((tb_{ij}) * (uc_{ij})) \\
&= (ta_{ij}) + (tub_{ij}) + (u^2 c_{ij})\ (u^2 \equiv u \pmod n) \\
&= (ta_{ij}) + (tub_{ij}) + (uc_{ij}) \qquad \text{II}
\end{aligned}$$

I and II are the same. Hence $(G, *, (t, u))$ is a semigroup.

*Note:* If n is a prime $(G, *, (t, u))$ is never a semigroup as $u^2 \equiv u \pmod{n}$ and $t^2 \equiv t \pmod{n}$ can never occur.

**THEOREM 2.1.6:** *The matrix groupoid $(G, *, (t, u))$ is an idempotent matrix groupoid if and only if $t + u \equiv 1 \pmod{n}$ (G is the set of all $m \times n$ matrices; $1 \leq m < \infty$ and $1 \leq n < \infty$ with entries from $Z_n$}.*

*Proof:* Given $(G, *, (t, u))$ is matrix groupoid row or column or $m \times n$ matrix; or in the mutually exclusive sense from $Z_n$. Let $A = (a_{ij}) \in G$

$$\begin{aligned}
A * A &= (a_{ij}) * (a_{ij}) \\
&= ta_{ij} + ua_{ij} \\
&= (t + u)(a_{ij}).
\end{aligned}$$

Now $A * A = A$ implies $(a_{ij}) = ((t + u)(a_{ij}))$. That is $((t + u - 1)(a_{ij})) = (0)$. This is possible if and only if $t + u \equiv 1 \pmod{n}$.

Hence the claim.

It is left as exercises for the reader to prove that zero matrix is not an ideal of the matrix groupoid $(G, *, (t, u))$. G is the collection of matrices with entries from $Z_n$ and theorem 2.1.7.

**THEOREM 2.1.7:** *Let $(G, *, (t, u))$ and $(G, *, (u, t))$ be a matrix groupoids with entries from $Z_n$. P is a left ideal of $(G, *, (t, u))$ if and only if P is a right ideal of $(G, *, (u, t))$.*

*Recall a matrix groupoid $(G, *, (t, u))$ is simple if and only if $(G, *, (t, u))$ has no normal subgroupoids.*

It is left as an exercise for the reader to prove the following two theorems.

**THEOREM 2.1.8:** *Let $(G, *, (t, u))$ be a matrix groupoid with entries from $Z_n$. If $n = t + u$ where both t and u are primes then $(G, *, (t, u))$ is a simple matrix groupoid.*

**THEOREM 2.1.9:** *Let $(G, *, (t, u))$ be a matrix groupoid with entries from $Z_n$, n even and $t + u = n$ with $(t, u) = t$. Then $(G, *, (t, u))$ has only one matrix subgroupoid of order $n / t$ and it is a normal matrix subgroupoid of $(G, *, (t, u))$.*

From the above theorem the following conclusion is obvious.

**THEOREM 2.1.10:** *Let $(G, *, (t, u))$ be a matrix groupoid with entries from $Z_n$, n even with $(t, u) = t$ and $t + u = n$. Then $(G, *, (t, u))$ is not a simple matrix groupoid.*

*We say as in case of usual groupoid a matrix groupoid $(G, *, (t, u))$ is a Smarandache matrix groupoid if $(G, *, (t, u))$ has a matrix subgroupoid $(H, *, (t, u))$ which is a matrix semigroup.*

***Example 2.1.22:*** Let $G = \{(x_1, x_2, x_3) \mid x_i \in Z_{10}; 1 \leq i \leq 3\}$ be a collection of row matrices with entries from $Z_{10}$.

Define $*$ on G by $(a_{ij}) * (b_{ij}) = (a_{ij}) + (5.b_{ij})) \pmod{n}$ for all $(a_{ij}), (b_{ij}) \in G$. $(G, *, (1, 5))$ is a matrix groupoid built using $Z_{10}$. $H = \{(0\ 0\ 0), (5\ 5\ 5)\} *, (1, 5))$ is a matrix semigroup. Hence $(G, *, (1, 5))$ is a Smarandache matrix groupoids.

As in case of groupoids we can define in case of Smarandache matrix groupoid $(G, *, (t, u))$ Smarandache matrix subgroupoid $(H, *, (t, u)$ if $(H, *, (t, u))$ is a matrix subgroupoid of (G, *, (t, u)) and H has a proper subset (K, *, (t, u)) such that (K, *, (t, u)) is a matrix semigroup.

The following theorem can be proved by any interested reader.

**THEOREM 2.1.11:** *Every matrix subgroupoid of a matrix groupoid (G, *, (t, u)) need not in general be a Smarandache subgroupoid of (G, *, (t, u)).*

The reader is expected construct examples of the above claim.
    Almost all properties enjoyed by groupoids can be derived in case of matrix groupoids more so the Smarandache properties.
    Also these matrix groupoids are both finite and infinite. We can also study about properties like isomorphism and homomorphism of matrix groupoids. Interested reader is expected to do this regular exercise. However the final chapter contains problems for the reader.

## 2.2 Polynomial Groupoids

In this section we introduce for the first time, the notion of polynomial groupoids built using $Z_n[x]$, $n < \infty$, or $Z[x]$ or $R[x]$ or $Q[x]$.
    Some operation is defined on these sets so that the set together with the operation becomes a groupoid which we call as polynomial groupoids.
    We follow the notation if $a_0 + a_1x + \ldots + a_nx^n$ is a polynomial of degree n then it is represented by the $n + 1$ tuple $(a_0, a_1, \ldots, a_n)$, $a_i \in Z_n$ or R or Q or Z, $0 \leq i \leq n$ 'or' used in the mutually exclusive sense. $(0, 0, \ldots, 0)$ denotes the zero polynomial $0 + 0x + \ldots + 0x^n$. Thus $x^3 + x + 1$ is represented by $(1\ 1\ 0\ 1)$. $x^6 + 2x + 1$ is represented by $(1, 2, 0, 0, 0, 0, 1)$ and so on. Now we proceed onto define polynomial groupoids using $Z_n$.

**DEFINITION 2.2.1:** *Let $f(x) = a_0 + \ldots + a_n x_n$ and $b_0 + \ldots + b_n x_n = g(x)$ be two polynomials in $Z_m^n [x]$, where $Z_m^n [x]$ denotes the* *collection of all polynomials of degree less than or equal to n with coefficients from $Z_m$ ($m < \infty$; $n < \infty$).*

*We define $f(x) * g(x)$ as follows $f(x) = (a_0, a_1, ..., a_n)$ and $g(x) = (b_0, b_1, ..., b_n)$*

$$\begin{aligned} f(x) * g(x) &= (a_0, a_1, ..., a_n) * (b_0, b_1, ..., b_n) \\ &= (a_0 b_1, a_1 b_2 ..., a_{n-1} b_n, a_n) \\ &= a_0 b_1 + a_1 b_2 x + ... + a_{n-1} b_n x^{n-1} + a_n x^n \\ &= h(x) \in Z_m^n [x] \text{ (multiplication of } a_i b_j \text{ is modulo} \end{aligned}$$

*m).*

*Thus the coefficients are reshuffled in this way. ($Z_m^n [x]$, *) is defined as the polynomial groupoid of degree n with entries from $Z_m$. If*

$$\begin{aligned} g(x) &= a_0 + a_1 x + ... + a_n x^n \\ h(x) &= b_0 + b_1 x + ... + b_n x^n \\ t(x) &= t_0 + t_1 x + ... + t_n x^n \end{aligned}$$

$a_i, b_i, t_i \in Z_m$; $1 \leq i \leq n$.

$$\begin{aligned} (g(x) * h(x)) * t(x) &= (a_0 \, a_1 \, ... \, a_n) * (b_0, ..., b_n)) * (t_0, ..., t_n) \\ &= (a_0 b_1 \, a_1 b_2 ... \, a_{n-1} b_n, a_n) * (t_0, t_1, ..., t_n) \\ &= (a_0 b_1 t_1, a_2 b_2 t_2, ..., a_{n-1} b \, t_n, a_n) \quad \text{I} \end{aligned}$$

$$\begin{aligned} g(x) * h(x) * t(x)) &= g(x) [(b_0, b_1, ..., b_n) * (t_0, t_1, ..., t_n)) \\ &= (a_0, a_1, ..., a_n) * (b_0 t_1, b_1 t_2, ..., b_{n-1} t_n, b_n) \\ &= (a_0 b_1 t_2, a_1 b_2 t_3, ..., a_{n-1} b_n, a_n) \quad \text{II} \end{aligned}$$

Clearly the polynomials given by I and II are different. Thus the operation * in general is non associative.

We will first illustrate this by some simple examples.

***Example 2.2.1:*** Let $G = \{(x_1, x_2, x_3, x_4, x_5) / x_i \in Z_5, 1 \leq i \leq 5\}$ all polynomials of degree less than or equal to four with coefficients from $Z_5$. Define * on G as follows.

Let
$$f(x) = 3x^2 + 4x + 1 = (1, 4, 3, 0, 0)$$
and
$$g(x) = 4x^4 + x^3 + 4 = (4, 0, 0, 1, 4)$$

$$\begin{aligned}
f(x) * g(x) &= (1, 4, 3, 0, 0) * (4, 0, 0, 1, 4) \\
&= (0, 0, 3, 0, 0) \\
&= 3x^2
\end{aligned}$$

$$\begin{aligned}
g(x) * f(x) &= (4\ 0\ 0\ 1\ 4) * (1\ 4\ 3\ 0\ 0) \\
&= (1\ 0\ 0\ 0\ 4) \\
&= 1 + 4x^4.
\end{aligned}$$

We see in the first place in general $f(x) * g(x) \neq g(x) * f(x)$. Now we show '*' is non associative.
Choose
$$\begin{aligned}
h(x) &= 3x^4 + 2x^3 + x^2 + 4x + 1 \\
&= (1, 4, 1, 2, 3).
\end{aligned}$$
Now
$$\begin{aligned}
(f(x) * g(x)) * (h(x)) &= (0\ 0\ 3\ 0\ 0) * (1\ 4\ 1\ 2\ 3) \\
&= (0\ 0\ 1\ 0\ 0) \\
&= x^2.
\end{aligned}$$

$$[f(x) * g(x)] * h(x) = x^2 \qquad\qquad \text{I}$$

Consider
$$\begin{aligned}
f(x) * [g(x) * h(x)] &= f(x)[(4\ 0\ 0\ 1\ 4) * (1\ 4\ 1\ 2\ 3)] \\
&= f(x) * [1\ 0\ 0\ 3\ 4] \\
&= (1\ 4\ 3\ 0\ 0)(1\ 0\ 0\ 3\ 4) \\
&= (0\ 0\ 4\ 0\ 0) \\
&= 4x^2 \qquad\qquad \text{II}
\end{aligned}$$

We see I and II are different.

It is important to mention here that one need not restrain one's study to polynomials with coefficients from $Z_m$ ($m < \infty$ modulo integers). We can define polynomial groupoids using Z or Q or R. It is pertinent to record here that authors have solved the open problem proposed in [20] that we can have an infinite class of groupoids which are Smarandache groupoids and the cardinality of these polynomial matrix groupoids are infinite when built over Z or Q or R.

The following example answers the open problem.

*Example 2.2.2:* Let G = {all polynomials of degree less than or equal to n} (G, *) is a polynomial groupoid. Take P = {$mx^n$ | m

∈ Z or R or Q}. (P, *) is a semigroup. Hence (G, *) is a polynomial Smarandache groupoid. Take $p(x) = 20x^n$, $q(x) = 4x^n$ and $r(x) = -5x^n$.

$$(p(x) * q(x) * r(x)) = [(0, \ldots, 0, 20) * (0, 0, \ldots, 0, 4)] *$$
$$(0, \ldots, 0, -5)$$
$$= (0, 0, \ldots, 0, 20) * (0, 0, \ldots, -5)$$
$$= (0, 0, \ldots, 20)$$
$$= p(x).$$

$$p(x) * [q(x) * r(x)] = [(0, \ldots, 0, 20) * (0, 0, \ldots, 0, 4)] *$$
$$(0, \ldots, 0, -5)$$
$$= (0, 0, \ldots, 0, 20) * (0, 0, \ldots, 0, 4)$$
$$= (0, 0, \ldots, 20)$$
$$= p(x)$$

Thus * on P is associative. Hence (P, *) is a semigroup in (G, *). Thus (G, *) is a polynomial Smarandache groupoid which is of infinite order.

Now we have mainly constructed this class of polynomials to show that a solution to the open problem in [20] exists.

Before we proceed onto define classes of polynomial groupoids in a very usual way similar to the one done in matrix groupoids, we extend the notion of polynomial groupoids to any polynomial ring which has polynomials of all degrees.

Thus if

$$G = \left\{ \sum_{i=0}^{\infty} a_i x^i \,\middle|\, a_i \in Z; 0 \leq i \leq \infty \right\}$$

all polynomials of any degree with coefficients from Z in the variable x. Any polynomial

$$p(x) = \sum_{i=0}^{\infty} a_i x^i = (a_0, a_1, \ldots, a_\infty);$$

$a_i \in Z$. If

$$q(x) = \sum_{i=0}^{\infty} b_i x^i = (b_0, b_1, \ldots, b_\infty),$$

$$p(x) * q(x) = (a_0 b_1, a_1 b_2, a_2 b_3, \ldots, a_\infty).$$

Thus (G, *) is again an infinite groupoid.

Now we proceed onto define a different type of operation on G.

**DEFINITION 2.2.2:** *Let $G = \{$all polynomials in the variable $x$ with coefficients form $Z_m$, $m < \infty$ of degree say less than or equal to $n\}$. Define for any $p(x) = p_0 + p_1 x + \ldots + p_n x^n$ and $q(x) = q_0 + q_1 x + \ldots + q_n x^n$ a binary operation * as follows*

$$p(x) * q(x) = (tp_0 + uq_0) + (tp_1 + uq_1) x + \ldots + (tp_n + uq_n)x^n$$

*where $t, u \in Z_m \setminus \{0\}$ $(t, u) = 1$ and $t \neq u$; $p_i, q_i \in Z_m$; $0 \leq i \leq n$. It is easily verified $(G, *, (t, u))$ is a groupoid. This groupoid is defined as polynomial groupoid of type I.*

We will illustrate this situation by some examples.

***Example 2.2.3:*** Let

$$G = \left\{ \sum_{i=0}^{5} a_i x^i \middle| a_i \in Z_8; 0 \leq i \leq 5 \text{ x an indeterminate} \right\}.$$

$(G, *, (3, 2))$ is a polynomial groupoid of finite order of type I.

***Example 2.2.4:*** Let

$$G = \left\{ \sum_{i=0}^{n} a_i x^i \middle| a_i \in Z; n < \infty, 0 \leq i \leq n \right\}$$

all polynomial of degree less than or equal to n with coefficients from Z. Choose $(22, -81) = (t, u)$. It is easily verified $(G, *, (22 - 81))$ is polynomial groupoid of infinite order and is type I.

If we consider polynomials P of the form

$$\sum_{i=0}^{n} a_i x^i; a \in Z;$$

$\{P, *, (22, -81)\}$ is a polynomial subgroupoid of $\{G, *, (22, -81)\}$.

The reader is given the task to prove or disprove $(P, *, (22, -81)$ is a semigroup.

***Example 2.2.5:*** Let

$$G = \left\{ \sum_{i=0}^{3} a_i x^i \,\middle|\, a_i \in Z_3; 0 \leq i \leq 3 \right\}$$

all polynomials in the variable x of degree less than or equal to 3 with coefficients from $Z_3$.

We see we have only two polynomial groupoids of type I constructed in this way. They are {G, *, (2, 1)} and {G, *, (1, 2)}. Here it is important to note that we can have infinite number of polynomial groupoids built using $Z_3$. This is done by varying the degree of the polynomials from 1 to n; n ≤ ∞. However for each degree fixed as m; m ≤ ∞ we can have only two distinct polynomial groupoids of type I built using $Z_3$.

We consider the following examples to give subgroupoids of type I.

***Example 2.2.6:*** Let

$$G = \left\{ \sum_{i=0}^{2} a_i x^i \,\middle|\, a_i \in Z_{10}; 0 \leq i \leq 3 \right\}$$

be all polynomials in the variable x with coefficients from $Z_{10}$ of degree less than or equal to two. (G, *, (1, 5)) is a polynomial groupoid of type I built using $Z_{10}$.

Consider

$$P = \left\{ \sum_{i=0}^{2} a_i x^i \,\middle|\, a = 0 \text{ or } 5 \right\};$$

(P, *, (1, 5)) is a polynomial subgroupoid of G.

Now we proceed onto define type II polynomial groupoids.

**DEFINITION 2.2.3:** *Let*

$$G = \left\{ \sum_{i=0}^{n} a_i x^i \,\middle|\, a_i \in Z_n; 0 \leq i \leq n \right\}$$

be all polynomials in the variable x with coefficients from $Z_n$ of all possible degrees including infinity. Define a binary operation * on G as follows.
   If
$$p(x) = \sum_{i=0}^{n} p_i x^i \text{ and } q(x) = \sum_{i=0}^{n} q_i x^i$$

then

$$p(x) * q(x) = \sum_{i=0}^{n} (tp_i + uq_i) x^i$$

'+' modulo n where $(t,u) \in Z_n \setminus \{0\}$ but $(t,u) \neq 1$. Then $(G, *, (t, u))$ is defined as the polynomial groupoid of type II built using $Z_n$.

*Note:* $Z_n$ can be replaced by Q or Z or R and still polynomial groupoid will continue to be of type II groupoid.

We will illustrate this situation by some examples.

***Example 2.2.7:*** Let

$$G = \left\{ \sum_{i=0}^{12} a_i x^i \,\middle|\, a_i \in Z_{18}; 0 \leq i \leq 12 \right\};$$

all polynomials in the variable x with coefficients from $Z_{18}$ of degree less than or equal to 12. Define * on G as follows

$$p(x) * q(x) = \sum_{i=0}^{12} (tp_i + uq_i)(\mod 18) x^i$$

where $(t, u) = (4, 6)$; $4, 6 \in Z_{18}$ $(G, *, (4, 6))$ is a polynomial groupoid of type II built using $Z_{18}$.

   We have a very large class of polynomial groupoids of type II. We will give some properties enjoyed by them once we define two more types of polynomial groupoids.

**DEFINITION 2.2.4:** *Let*

$$G = \left\{ \sum_{i=0}^{n} a_i x^i \,\middle|\, a_i \in Z_m; 0 \leq i \leq n \right\};$$

all polynomials of degree less than or equal to n with coefficients from $Z_m$}. Define the operation * on G as follows: for any polynomial

$$p(x) = \sum_{i=0}^{n} p_i x^i$$

and

$$q(x) = \sum_{i=0}^{n} q_i x^i .$$

*Define*

$$p(x) * q(x) = \sum_{i=0}^{n} (tp_i + uq_i) \ (mod \ m) \ x^i$$

where $t, u \in Z_n \setminus \{0\}$ and $t = u$.

It is easily verified (G, *, (t, t)) is a polynomial groupoid and this is defined as a polynomial groupoid of type II.
    We will illustrate this by a few examples.

***Example 2.2.8:*** Let

$$G = \left\{ \sum_{i=0}^{7} a_i x^i \ \middle| \ a_i \in Z_{24}; 0 \leq i \leq 7 \right\}$$

all polynomials in the variable x with coefficients from $Z_{24}$ of degree less than or equal to 7.
    Define * on G by

$$p(x) * q(x) = \left( \sum_{i=0}^{7} p_i x^i \right) * \left( \sum_{i=0}^{7} q_i x^i \right)$$

$$= \sum_{i=0}^{7} (4p_i + 4q_i) \ (mod \ 24) \ x^i$$

It is easily verified (G, *, (4, 4)) is a polynomial groupoid of type III using $Z_{24}$.
    We show * in general in non associative.
    Take $p(x) = 5x^3 + 7x + 3$, $q(x) = 8x^5 + 4x^2 + 2x + 1$ and $r(x) = x^6 + 5x^5 + 7x^2 + 13$ in G.
    Consider

(p(x) * q(x)) * r(x)
= {($5x^3 + 7x + 3$) * ($8x^5 + 4x^2 + 2x + 1$)} *
  ($x^6 + 5x^5 + 7x^2 + 13$)
= ($20x^3 + 8x^5 + 16 + 16x^2 + 12x$) * ($x^6 + 5x^5 + 7x^2 + 13$)
= ($4x^6 + 4x^5 + 3x^3 + + 20x^2 + 20$)    I

Consider
p(x) * (q(x) * r(x))
= ($5x^3 + 7x + 3$) * [($8x^5 + 4x^2 + 2x + 1$) *
  ($x^6 + 5x^5 + 7x^2 + 13$)]
= ($5x^3 + 7x + 3$) * [$4x^6 + 8x + 8 + 20x^2 + 4x^5$]
= [$16x^6 + 16x^5 + 20 + 12x + 8x^2 + 20x^3$].    II

We see I and II are not the same polynomial so the operation * in general is non associative. Thus (G, *, (4, 4)) is a polynomial groupoid of type III built using $Z_{24}$.

*Example 2.2.9:* Let

$$G = \left\{ \sum_{i=0}^{4} a_i x^i \,\middle|\, a_i \in Z_7; 0 \le i \le 4 \right\}$$

all polynomials in the variable x of degree less than or equal to 4 with coefficients from $Z_7$.

Define * on G by a * b = 4a + 4b(mod 7) i.e. (G, *, (4, 4)) is a polynomial groupoid of type III. Take
$$a = (3x + 2)$$
$$b = (4x^2 + 2x + 5)$$
and
$$c = 6x^4 + 4x^3 + 5x + 1.$$

(a * b) * c = [($3x + 2$) * ($4x^2 + 2x + 5$)] *
  ($6x^4 + 4x^3 + 5x + 1$)
= ($6x + 2x^2$) * ($6x^4 + 4x^3 + 5x + 1$)
= $3x^4 + 4 + 2x + x^2 + 2x^3$    I

a * (b * c) = ($3x + 2$) * [($4x^2 + 2x + 5$)*($6x^4 + 4x^3 + 5x + 1$)]
= ($3x + 2$) * [$3x^4 + 2x^3 + 2x^2 + 3$]
= [$5x + 6 + 5x^4 + x^3 + x^2$]    II

I and II are different. Hence (G, *, (4, 4)) is a polynomial groupoid with entries from $Z_7$ of type III.

Now we proceed onto define type IV polynomial groupoids.

**DEFINITION 2.2.5:** *Let $G = \left\{ \sum_{i=0}^{n} a_i x^i \; ; \; 0 \leq i \leq n; \; a_i \in Z_m \right\}$ be all polynomials with coefficients from $Z_m$. Define * on G as $g(x) * h(x) = [tg(x) + 0\, h(x)]; \; t \in Z_m$. $\{G, *, (t, 0)\}$ is a polynomial groupoid of type IV with entries from $Z_m$.*

We will illustrate this situation by some examples.

*Example 2.2.10:* Let $Z_8 = \{0, 1, 2, \ldots, 7\}$ ring of modulo integers eight. Choose
$$G = \left\{ \sum_{i=0}^{26} a_i x^i \;\middle|\; a_i \in Z_8; 0 \leq i \leq 26 \right\};$$

all polynomials of degree less than or equal to 26 with coefficients from $Z_8$ in the variable x.

Let $t = 5 \in Z_8$. Define for any two polynomials $g(x), h(x)$ in G. $g(x) * h(x) = tg(x) + 0\, h(x)$. $\{G, *, (5, 0)\}$ is a polynomial groupoid of type IV with entries from $Z_8$.

*Example 2.2.11:* Let
$$G = \left\{ \sum_{i=0}^{5} a_i x^i \;\middle|\; a_i \in Z_{23}; 0 \leq i \leq 5 \right\}$$

be the collection of all polynomials with entries from $Z_{23}$, in the variable x of degree less than or equal to 5. For $p(x), q(x) \in G$ define $p(x) * q(x) = 4p(x) + 0q(x) = (G, *, (4, 0)\}$ is a polynomial groupoid of type IV with entries from $Z_5$.

Suppose $Z_p[x]$ contains polynomials in the variable x with coefficients from $Z_p$. Define operation * on $Z_p[x]$ by $p(x) * q(x) = tp(x) + 0q(x)$, $p(x), q(x) \in Z_p[x]$ and $t \in Z_p$ $\{Z_p[x], *, (t, 0)\}$ is an infinite polynomial groupoid of type IV of infinite order.

When t = 1 we get $(Z_p[x], *, (1, 0))$ to be a semigroup where p is a prime if in $Z_n[x]$; n is not a prime and $1 \neq t \in Z_n \setminus \{0\}$ is such that $t^2 = t \pmod{n}$ then we see $[Z_n[x], *, (t, 0)]$ is a semigroup.

We will illustrate this situation by an example.

***Example 2.2.12:*** Let $Z_6[x]$ be a polynomial ring. Let $[Z_6[x], *, (3, 0)]$ be the polynomial groupoid of type IV. Take $3 \in Z_6$ we see $p(x) * q(x) = 3p(x)$.
Now
$$
\begin{aligned}
(p(x) * q(x)) * r(x) &= (3p(x) + 0) * r(x) \\
&= 3p(x) * r(x) \\
&= 3p(x) \quad (\because 3^2 = 3 \pmod{6}). \\
p(x) * (q(x) * r(x)) &= p(x) * [3q(x)] \\
&= 3p(x) + 0 \\
&= 3p(x).
\end{aligned}
$$

Thus * on $Z_6[x]$ is associative when $t = 3$; $t \in Z_6$.

Motivated by this example we see we have for the class of polynomial groupoids $\{Z_n[x], *, (t, 0)\}$ for $t \in Z_n$, we have some of them to be polynomial semigroups, for if $t = 1$ $\{Z_n[x], *, (1, 0)\}$ is a polynomial semigroup. In view of this we define the Smarandache special class of groupoids.

**DEFINITION 2.2.6:** *Let G(S) = {Class of groupoids defined using the same set S}. If G(S) has atleast one semigroup then we call G(S) to be a Smarandache Special Class of groupoid. (SSC-groupoids).*

We will illustrate this by an example.

***Example 2.2.13:*** Let $G(Z_6[x]) = \{Z_6[x], *, (t, 0); t \in Z_6\}$ be class of polynomial groupoids built using $Z_6[x]$. We see when n = 3, n = 4 and n = 1 we get polynomial semigroups; hence $G(Z_6[x])$ is a Smarandache Special Class of groupoid.

**THEOREM 2.2.1:** *Let $G(Z_n[x]) = \{Z_n[x]; *; (t, 0); t \in Z_n\}$ be a class of polynomial groupoids. $G(Z_n[x])$ is a SSC – groupoid.*

*Proof:* Since $1 \in Z_n$ we see $\{Z_n[x], *, (1, 0)\} \in G(Z_n[x])$ is a polynomial semigroup, hence $G(Z[x])$ is a SSC-groupoid.

*Example 2.2.14:* Let $G(Z_{12}[x])$ be the class of groupoids built using $Z_{12}[x]$. Take $H_1 = \{Z_{12}[x], *, (4, 0), 4 \in Z_{12}\}$, $H_2 = \{Z_{12}[x], *, (9, 0)\}$ and $H_3 = \{Z_{12}[x], *, (1, 0)\}$ in $G(Z_{12}[x])$. $H_1$, $H_2$ and $H_3$ are polynomial semigroups. Hence $G(Z_{12}[x])$ is a SSC-groupoid.

*Example 2.2.15:* Let $G(Z_{30}[x])$ be the class of polynomial groupoids built using $Z_{30}[x]$. Take $H_1 = \{Z_{30}[x], *, (15, 0)\}$, $H_2 = \{Z_{30}[x], *, (1, 0)\}$, $H_3 = \{Z_{30}[x], *, (10, 0)\}$ and $H_4 = \{Z_{30}[x], *, (6, 0)\}$ are polynomial semigroups. So $G(Z_{30}[x])$ is a SSC-groupoid.

Now we will give classes of SSC-groupoid.

**THEOREM 2.2.2:** *The class of groupoids $Z(n) = \{Z_n, *, (t, u)\}$ n, not a prime; is a SSC-groupoid.*

*Proof:* We see $Z(n)$ has semigroups. For if $(t, u) = 1$ with $t, u \in Z_n \setminus \{0\}$ with $t^2 = t \pmod{n}$ and $u^2 = u \pmod{n}$ then $(Z_n, *, (t, u))$ is a semigroup. Hence $Z(n)$ is a SSC-groupoid.

We will illustrate this situation by an example.

*Example 2.2.16:* Take $\{Z_{12}, *, (4, 9)\} \in Z(12)$. $(Z_{12}, *, (4, 9))$ is a semigroup as $4^2 = 4 \pmod{12}$ and $9^2 = 9 \pmod{12}$.

Now we will construct a class of infinite groupoids.

**DEFINITION 2.2.7:** *Let $G = \{Q[x], *, (m, n) \:/\: m, n \in Z^+; (m, n) = 1\}$. It is easily verified G is a groupoid of infinite order. This groupoid will be known as polynomial groupoid with rational coefficients.*

**DEFINITION 2.2.8:** *Let $G^* = \{Q[x], *, (m, n) \,/\, m, n \in Z^+, (m, n) = d\}$. This groupoid is known as the polynomial groupoid with rational coefficients ($d \neq 1$. If $d = 1$ $G = G^*$).*

***Example 2.2.17:*** Let $G = \{Q[x], *, (5, 7)\}$ be a polynomial groupoid with rational coefficients. Let
$$p(x) = x + 1,$$
$$q(x) = 3x^2 + \frac{7}{2}x + \frac{5}{3}$$
and
$$r(x) = 5x^3 - \frac{3}{2}x + 9 \in Q[x].$$

$p(x) * (q(x) * r(x))$
$$= p(x) * [15x^2 + \frac{35x}{2} + \frac{25}{3} + 35x^3 - \frac{21}{2}x + 63]$$
$$= p(x) * [35x^3 + 15x^2 + 7x + \left(\frac{25}{3} + 63\right)]$$
$$= 5(x + 1) + 35 \times 7x^3 + 105x^2 + 49x + \left(\frac{25 + 189}{3}\right) \quad \text{I}$$

Consider
$(p(x) * q(x)) * r(x)$
$$= (5(x+1) + 21x^2 + \frac{49x}{2} + \frac{35}{3}) * (5x^3 - \frac{3x}{2} + 9)$$
$$= 25x + 25 + 105x^2 + \frac{49 \times 5}{2}x + \frac{35 \times 5}{3} + 35x^3 - \frac{21}{2}x + 63$$
$$= 35x^3 + 105x^2 + (25 + \frac{49 \times 5}{2} - \frac{21}{2})x + (25 + 63 + \frac{35 \times 5}{3})$$
$$\quad \text{II}$$

Clearly I and II are not equal so the operation * in general defined on Q[x] is non associative.

Next we proceed onto give an example in which $(m, n) = d \neq 1$.

***Example 2.2.18:*** Let $G^* = \{Q[x], *, (6, 3) = 3\}$ be a polynomial groupoid. Take
$$p(x) = x^7 + 1,$$
$$q(x) = 8x^4$$
and
$$r(x) = 3x^2 + 5x + 1.$$

$$
\begin{aligned}
(p(x) * q(x)) * r(x) &= [6(x^7 + 1) + 3 \times 8x^4] * (3x^2 + 5x + 1) \\
&= (6x^7 + 6 + 24x^4] * (3x^2 + 5x + 1) \\
&= (36x^7 + 36 + 144x^4) + 9x^2 + 15x + 3 \\
&= 36x^7 + 144x^4 + 9x^2 + 15x + 39 \qquad \text{I}
\end{aligned}
$$

$$
\begin{aligned}
p(x) * (q(x) * q(x)) &= p(x) * [48x^4 + 9x^2 + 15x + 3] \\
&= 6x^7 + 6 + 144x^4 + 27x^2 + 45x + 9 \\
&= 6x^7 + 144x^4 + 27x^2 + 45x + 15 \qquad \text{II}
\end{aligned}
$$

I and II are not equal. Hence $G^*$ is a groupoid and is not a semigroup.

**DEFINITION 2.2.9:** *Let $G^{**} = \{Q[x], *, (t, 0) / t \in Z^+ \setminus \{0\}\}$. $G^{**}$ is a polynomial groupoid with polynomials from $Q[x]$.*

We will illustrate this situation by an example.

***Example 2.2.19:*** Let $G^{**} = \{Q[x], *, (5, 0)\}$ be a polynomial groupoid. Take $p(x) = 8x^3$, $q(x) = 5x^2 + 1$ and $r(x) = x^5$ in $Q[x]$.

$$
\begin{aligned}
(p(x) * q(x)) * r(x) &= (40x_3 + 0) * r(x) \\
&= 200x^3 \qquad \text{I} \\
p(x) * (q(x) * r(x)) &= 8x^3[25x^2 + 5 + 0] \\
&= 40x^3 \qquad \text{II}
\end{aligned}
$$

It is pertinent to mention here that we can build groupoids using $Z[x]$ or $R[x]$ or $C[x]$ in a similar way. We can derive many interesting properties about these polynomial groupoids. Infact all these groupoids are of infinite order. Interested reader is expected to derive interesting properties about these groupoids.

## 2.3 Special Interval Groupoids

In this section we introduce yet another special class of groupoids called interval groupoids. Here we construct several interesting properties about them.

*Notations:* Let
$$Z_I^+ = \{[a, b] \mid a, b \in Z^+ \cup \{0\}, a \leq b\}$$
$$Q_I^+ = \{[a, b] \mid a, b \in Q^+ \cup \{0\}, a \leq b\}$$
$$R_I^+ = \{[a, b] \mid a, b \in R^+ \cup \{0\}, a \leq b\}.$$

Clearly $Z_I^+ \subseteq Q_I^+ \subseteq R_I^+$. Consider $Z_n^I = \{[0, r] \mid r \in Z_n\}$ is the set of intervals in $Z_n$.

Now using these four classes of intervals we build classes of groupoids.

**DEFINITION 2.3.1:** *Let $G * (Z_n^I) = \{[0, r], *, (m, p) / r, m, p \in Z_n$ with $(m, p) = 1$, $m$ and $p$ primes$\}$ where $[0, r] * [0, s] = [0, mr + ps \pmod{n}]$. $G * (Z_n^I)$ is a class of groupoids, defined as modulo integer interval groupoids of level one.*

*Clearly $G * (Z_n^I)$ has only finite number of elements.*

We will illustrate this by some examples.

***Example 2.3.1:*** Let $G * (Z_{12}^I) = \{[0, s], *, (5, 7); 5, 7, s \in Z_{12}\}$ be a modulo integer interval groupoid of level one. $G * (Z_{12}^I)$ has 12 elements.

***Example 2.3.2:*** $G * (Z_{19}^I) = \{[0, s], *, (11, 7), 11, 7, s \in Z_{19}\}$ is a modulo integer interval groupoid of level one.

Now we proceed onto define interval groupoids of level one with entries from $Z^+ \cup \{0\}$ or $Q^+ \cup \{0\}$ or $R^+ \cup \{0\}$ or $C^+ \cup \{0\}$.

**DEFINITION 2.3.2:** *G(S) = {[0, a], *, (p, q) = 1; a, p, q ∈ Z$^+$ ∪ {0}} (or C$^+$, Q$^+$ or R$^+$) p and q primes} where * is defined as [0, a] * [0, b] = [0, pa + qb] ∈ G(S) and S = $C_I^+$ or $Z_I^+$ or $R_I^+$. G(S) is a class of interval groupoids of level one.*

*As p, q ∈ Z$^+$ \ {0}, vary over Z$^+$ \ {0} we get infinite number of interval groupoids. S can be $R_I^+$ $C_I^+$ or $Q_I^+$.*

We will illustrate this by some examples.

***Example 2.3.3:*** Let G*($Z_I^+$) = {[0, a], *, (5, 19); a, 5, 19 ∈ Z$^+$} be an interval integer groupoid of level one. It is easily verified * in general is non associative. For take [0, 1], [0, 3] and [0, 2] in G*($Z_I^+$).

$$
\begin{aligned}
\text{[0 1]} * (\text{[0 3]} * \text{[0 2]}) &= \text{[0 1]} * ([0, 15] + [0, 38]) \\
&= \text{[0 1]} * [0, 53] \\
&= [0, 5] + [0, 53 \times 19] \\
&= [0, 5 + 53 \times 19] \quad\quad\quad\quad \text{I}
\end{aligned}
$$

$$
\begin{aligned}
([0, 1] * \text{[0 3]}) * \text{[0 2]} &= ([0, 5] + [0, 57]) * [0, 2] \\
&= [0, 62] * [0, 2] \\
&= [0, 310] + [0, 38] \\
&= [0, 348] \quad\quad\quad\quad\quad\quad \text{II}
\end{aligned}
$$

Clearly I and II are not equal; hence G * ($Z_I^+$) is a integer interval groupoid of level one.

***Example 2.3.4:*** Let G * ($Q_I^+$) = {[0, a], *, (2, 3); a, 2, 3, ∈ Q$^+$} be the rational interval groupoid of level one.

***Example 2.3.5:*** Let G* ($R_I^+$) = {[0, a], *, [3, 23], a, 3, 23 ∈ R$^+$} be a real interval groupoid of level one.

***Example 2.3.6:*** Let G* ($C_I^+$) = {[0, a], *, (3, 2) | a, 3, 2 ∈ C$^+$} be a complex interval groupoid of level one.

***Example 2.3.7:*** Let $G^*(R_I^+) = \{[0, a], *, (29, 53); a, 29, 53 \in R^+\}$ be a real interval groupoid of level one.

Now we proceed onto define the notion of real or complex or integer or rational interval subgroupoid of level one.

**DEFINITION 2.3.3:** *Let $T \in G^*(S)$ be a real or complex or integer or modulo integer or rational interval groupoid of level one.*

*Let $P \subseteq T$ (P a proper subset of T). If P itself is a real or complex or integer or modulo integer or rational interval groupoid of level one then we call P to be a real or complex or integer or modulo integer or rational interval subgroupoid of T of level one.*

We will illustrate this situation by some examples.

***Example 2.3.8:*** Let $\{[0, a], *, (3, 2); a, 3, 2 \in Z_{12}\} \in G^*(Z_{12}^I)$ be a modulo integer groupoid of order 12 of level one. Take $P = \{[0, 2], [0, 0], [0, 4], [0, 6], [0, 8], [0, 10], *, (3, 2)\} \subseteq \{[0, a], *, (3, 2), a, 3, 2, \in Z_{12}\} \in G^*(Z_{12}^I)$; P is a modulo integer interval subgroupoid of order 6.

***Example 2.3.9:*** Let $T = \{[0, a], *, (5, 3) \mid 5, 3, a \in Z_{12}\}$ be a modulo integer interval groupoid of level one.

The reader is requested to find modulo integer interval subgroupoids of level one.

***Example 2.3.10:*** Let $R = \{[0, a], *, (3, 2); a \in Z^+ \cup \{0\}\}$ be an integer interval groupoid of level one.
Take $P = \{[0, a], *, (3, 2); a \in 2Z^+ \cup \{0\}\} \subseteq R$. P is a integer interval subgroupoid of R of level one.

***Example 2.3.11:*** Let $T = \{[0, a], *, (5, 7) / a \in R^+ \cup \{0\}\}$ be a real interval groupoid of level one. Take $P = \{[0, a], *, (5, 7) / a \in Q^+ \cup \{0\}\} \subseteq T$. P is a real interval subgroupoid of level one.

*Remark:* It is important and interesting to note that P is not a real interval groupoid but only rational interval groupoid so how much are we justified in calling it as a real interval subgroupod. But by the rule of convention we call so.

*Example 2.3.12:* Let T = {[0, x], *, (7, 3) / x, 7, 3 ∈ $Q^+ \cup \{0\}$}} be a rational interval groupoid of level one. Take P = {[0, a], *, (7, 3) / a, 7, 3 ∈ $Z^+ \cup \{0\}$}} ⊆ T; P is a rational interval subgroupoid of T of level one.

*Example 2.3.13:* Let T = {[0, a], *, (17, 2) / a, 17, 2 ∈ $C^+ \cup \{0\}$}} be a complex interval groupoid. Let S = {[0, a], *, (17, 2) / a, 17, 2 ∈ $Z^+ \cup \{0\}$}} ⊆ T; S is a complex interval subgroupoid of T of level one.

Now we proceed onto define interval groupoid of level two.

**DEFINITION 2.3.4:** *Let $G^{**}(Z_n^I)$ = {[0 a], *, (p, q) = 1, a, p and q are in $Z_n$;} where [0, a] * [0, b] = [0, pa + qb (mod n)]. $G^{**}(Z_n^I)$ is defined as the class of modulo integer interval groupoid of level two.*

*Example 2.3.14:* Let T = {[0, a], *, (5, 8), a, 5, 8 ∈ $Z_9$} ⊆ $G^{**}(Z_9^I)$ be the modulo integer interval groupoid of level two.

*Example 2.3.15:* Let P = {[0, a], *, (9, 10), a, 9, 10 ∈ $Z_{11}$} ⊆ $G^{**}(Z_{11}^I)$ be the modulo integer interval groupoid of level two.

Like wise we can define class of integer interval groupoid of level two; $G^{**}(Z_I^+) = G^{**}((Z^+ \cup \{0\})^I)$, class of real interval groupoid of level two; $G^{**}(R_I^+) = G^{**}((R^+ \cup \{0\})^I)$, class of rational interval groupoid of level two; $G^{**}(Q_I^+) = G^{**}(Q^+ \cup \{0\})^I)$ and class of complex interval groupoid of level two $G^{**}(C_I^+) = G^{**}((C^+ \cup \{0\}^I)$.

We make an assumption ($C_I^+$) or $Z_I^+$ or $Q_I^+$ or $R_I^+$ contains {0}.

*Example 2.3.16:* Let T = {[0, a], *, (8, 15), 15, 8, a ∈ $Z^+ \cup \{0\}$} be an integer interval groupoid of level two.

*Example 2.3.17:* Let S = {[0, a], *, (19, 16), a ∈ $Z^+ \cup \{0\}$} be the integer interval groupoid of level two.

*Example 2.3.18:* Let P = {[0, a], *, (27, 64), a, 27, 64 ∈ $Q^+ \cup \{0\}$} be the rational interval groupoid of level two.

*Example 2.3.19:* Let T = {[0, x], *, (27, 43), x, 27, 43 ∈ $R^+ \cup \{0\}$} be the real interval groupoid of level two.

*Example 2.3.20:* Let W = {[0, y], *, (43, 7); y, 7, 43 ∈ $C^+ \cup \{0\}$} be the complex interval groupoid of level two.

It is both interesting and important to note the following.

The class of integer interval groupoids of level two ⊆ The class of rational interval groupoids of level two ⊆ The class of real interval groupoids of level two ⊆ The class of complex interval groupoids of level two; that is $G^{**}(Z_I^+) \subseteq G^{**}(Q_I^+) \subseteq G^{**}(R_I^+) \subseteq G^{**}(C_I^+)$.

Same type of containment is true in case of interval groupoids of level one.

Now we will call a level one interval groupoid to be simple if it has no proper interval subgroupoids of level one.

We see $G^*(Z_n^I)$; when n is a prime is a simple modulo integer interval groupoid of level one.

*Example 2.3.21:* Let T = {[0, a], a ∈ $Z_5$, *, (3, 2)} be a modulo integer interval groupoid of level one. T = {[0, 0], [0, 1], [0, 2], [0, 3], [0, 4], * (3, 2)}.

We assume [0, 0] is not an interval groupoid as [0, 0] = [0] the degenerate interval. Now T is a simple modulo integer interval groupoid of level one. Thus we have a class of simple modulo integer interval groupoids of level one.

Now we can define interval subgroupoids of level two.

**DEFINITION 2.3.5:** *Let $T = \{[0, a], *, (p, q), a, p, q \in Z_n\}$ be a modulo integer interval groupoid of level two. If $P \subset T$ is such that $P = \{[0, a], *, (p, q)\} \subseteq T$ is a proper subset of $T$ and $P$ itself is a modulo integer interval groupoid of level two; then we call $P$ to be a modulo integer interval subgroupoid of level two of $T$. When $T$ has no interval subgroupoid of level two; then we call $T$ to be a simple modulo integer interval groupoid of level two.*

We will illustrate both the situation by some examples.

*Example 2.3.22:* Let $T = \{[0, a], *, (3, 8), a, 3, 8 \in Q^+ \cup \{0\}\}$ be a rational interval groupoid of level two. Take $P = \{[0, a], *, (3, 8), 3, 8, a \in Z^+ \cup \{0\}\} \subseteq T$; $P$ is a rational interval subgroupoid of level two.

*Example 2.3.23:* Let $T = [0, a], *, (9, 8), 9, 8, a \in R^+ \cup \{0\}\}$ be a real interval groupoid of level two. Take $P = \{[0, a], *, (9, 8), a, 9, 8 \in Q^+ \cup \{0\}\} \subseteq T$, $P$ is a real interval subgroupoid of level two.

*Example 2.3.24:* Let $T = \{[0, a]; *, (15, 8); 15, 8, a \in Z_{30}\}$ be a modulo integer interval groupoid of level two. Take $P = \{[0, 0], [0, 10], [0, 20], *, (15, 8), 15, 8 \in Z_{30}\} \subseteq T$. $P$ is a modulo integer interval subgroupoid of level two of $T$.

*Example 2.3.25:* Let $S = \{[0, a], *, (9, 8); a, 9, 8 \in Z_{11}\}$ be a modulo integer interval groupoid of level two.
  Clearly $S$ has no modulo integer interval subgroupoid of level two. Hence $S$ is a simple modulo integer interval groupoid.

The following theorem is left as an exercise for the reader to prove.

**THEOREM 2.3.1:** *Let $G^{**}(Z_p^I)$ be the collection of modulo integer interval groupoids of level two. Every modulo integer interval groupoid in $G^{**}(Z_p^I)$ is simple; $p$ a prime.*

Now we proceed onto define interval groupoids of level three.

**DEFINITION 2.3.6:** *Let $G^{**}(Z_n^I) = \{[0, a]; *, (p, q) = d \neq 1\}$ where $*$ is such that $[0, a] * [0, b] = [0, pa+qb \pmod{n}]$ we define $G^{***}(Z_n^I)$ to be the class of modulo integer interval groupoid of level three.*

*If we replace $Z_n^I$ by $Z_I^+$ we call $G^{***}(Z_I^+)$ to be the class of integer interval groupoid of level three.*

*If $Z_n^I$ in the definition is replaced by $Q_I^+$ we call $G^{***}(Q_I^+)$ to be the class of rational interval groupoids of level three. If $Z_n^I$ in the definition is replaced by $R_I^+$ we call $G^{***}(R_I^+)$ to be the class of real interval groupoids of level three. If in the definition $Z_n^I$ is replaced by $C_I^+$ we call $G^{***}(C_I^+)$ to be the class of complex interval groupoid of level three.*

We will illustrate each of the situation by some examples.

*Example 2.3.26:* Let $T = \{[0, a], *, (8, 24); a, 8, 24 \in Z_{30}\}$ be a modulo integer interval groupoid of level three.

*Example 2.3.27:* Let $T = \{[0, a], *, (9, 3); a, 9, 3 \in Z_{11}\}$ be a modulo integer interval groupoid of level three.

*Example 2.3.28:* Let $P = \{[0, a], *, (27, 30); a, 27, 30 \in Z^+ \cup \{0\}\}$ be the integer interval groupoid of level three.

*Example 2.3.29:* Let $P = \{[0, a], *, (11, 66); a, 11, 66 \in Z^+ \cup \{0\}\}$ be the integer interval groupoid of level three.

*Example 2.3.30:* Let $T = \{[0, a], *, (12, 26); a, 12, 26 \in Q^+ \cup \{0\}\}$ be a rational interval groupoid of level three.

*Example 2.3.31:* Let $F = \{[0, a], *, (28, 35); a, 28, 35 \in R^+ \cup \{0\}\}$ be the real interval groupoid of level three.

***Example 2.3.32:*** Let E = {[0, b], *, (7, 497); b, 7, 497 ∈ $C^+$ ∪ {0}} be the complex interval groupoid of level three.

Now we proceed onto define the notion of interval subgroupoid of level three.

**DEFINITION 2.3.7:** *Let T = {[0, a], *, (p, q) = d ≠ 1 p, q, d, a ∈ $Z_n$} be a modulo integer interval groupoid of level three. Let P ⊆ T; if P is a modulo integer interval groupoid of level three and P is a proper subset of T we call P to be the modulo integer interval subgroupoid of T of level three.*

*If T has no proper modulo integer interval subgroupoid then we call T to be a simple modulo integer interval groupoid.*

Analogous definitions hold good in case of integer or real or rational or complex interval groupoid of level three.
We will illustrate this situation by some simple examples.

***Example 2.3.33:*** Let T = {[0, a], *, [15, 10], a, 10, 15, ∈ $Z_{60}$} be a modulo integer interval groupoid of level three.
Let P = {[0, a], [0, 5], [0, 10], [0, 15], [0, 20], [0, 25], [0, 30], [0, 35], [0, 40], [0, 45], [0, 50], [0, 55], *, (15, 10)} ⊆ T be the modulo integer interval subgroupoid of T of level three.

***Example 2.3.34:*** Let T = {[0, a], *, (8, 24), a, 8, 24 ∈ $Z^+$ ∪ {0}} be the integer interval groupoid of level three. Take P = {[0, 2a], *, 2a, 8, 24 ∈ $Z^+$ ∪ {0}} ⊆ T, P is an integer interval subgroupoid of T of level three.

***Example 2.3.35:*** Let T = {[0, a], *, (27, 18), a, 27, 18 ∈ $Q^+$ ∪ {0}} be a rational interval groupoid of level three. Take P = {[0, a], *, (27, 18), a, 27, 18 ∈ $Z^+$ ∪ {0}} ⊆ T. P is a rational subgroupoid of level three.

***Example 2.3.36:*** Let X = {[0, a], *, (20, 22), a, 20, 22 ∈ $R^+$ ∪ {0}} be a real interval groupoid of level three.
Take P = {[0, a], *, (20, 22), a ∈ $Q^+$ ∪ {0}} ⊆ X, P is a real interval subgroupoid of level three of X.

*Example 2.3.37:* Let Y = {[0, x], *, (28, 42), x, 28, 42 ∈ $C^+ \cup \{0\}$} be a complex interval groupoid of level three. X = {[0, x], *, (28, 42), x ∈ $R^+ \cup \{0\}$} ⊆ Y, is the complex interval subgroupoid of level three.

*Example 2.3.38:* Let B = {[0, a], *, (8, 12), a, 8, 12 ∈ $Q^+ \cup Z_{13}$} be a modulo integer interval groupoid of level three. B is a simple modulo integer interval groupoid of level three.

The reader is expected to answer the following problem.
  Let X = {[0, a], *, (r, s) = t ≠ 1, a, r, s, t ∈ $Z_p$; p is a prime} be a modulo interval groupoid of level three. Is X simple?

Before we proceed onto describe more properties about these interval groupoid we give the definition of level four interval groupoids.

**DEFINITION 2.3.8:** *Let G\*\*\*\*($Z_n^I$) = {[0, a], *, (p, 0), 0 ≠ p, a ∈ $Z_n$} where [0, a] * [0, b] = [0, pa + 0b (mod n)] = [0, pa (mod n)]; be a modulo integer interval groupoid defined as the modulo integer interval groupoid of level four.*

If we replace ($Z_n^I$) in the definition 2.3.8 by $Z_I^+$ we get the integer interval groupoid of level four. By replacing $Z_n^I$ by $Q_I^+$ in the definition 2.3.8 we get the rational interval groupoid of level four. If $Z_n^I$ is replaced by $R_I^+$ in definition 2.3.8 we get the real interval groupoid of level four and is denoted by G\*\*\*\*($R_I^+$). Similarly by replacing $Z_n^I$ by $C_I^+$ in definition 2.3.8 we get the complex interval groupoid of level four and is denoted by G\*\*\*\*($C_I^+$). We see
  G\*\*\*\*($Z_I^+$) ⊆ G\*\*\*\*($Q_I^+$) ⊆ G\*\*\*\*($R_I^+$) ⊆ G\*\*\*\*($C_I^+$).

  We will illustrate this situation by some examples.

*Example 2.3.39:* Let B = {[0, a], *, (5, 0), a, 5 ∈ $Z_{15}$} be a modulo integer interval groupoid of level four.

*Example 2.3.40:* Let C = {[0, x], *, (6, 0), x, 6 ∈ $Z_{19}$} be a modulo integer interval groupoid of level four.

*Example 2.3.41:* Let D = {[0, b], *, (8, 0), b, 8 ∈ $Z_I^+$} be a integer interval groupoid of level four.

*Example 2.3.42:* Let X = {[0, a], *, (7, 0), a, 7 ∈ $Q_I^+$} be a rational interval groupoid of level four.

*Example 2.3.43:* Let Y = {[0, b], *, (141, 0), b, 141 ∈ $R_I^+$} be a real interval groupoid of level four.

*Example 2.3.44:* Let P = {[0, t], *, (15, 0), t, 15 ∈ $C_I^+$} be a complex interval groupoid of level four.

We will now proceed onto define interval subgroupoids of level one.

**DEFINITION 2.3.9:** *Let S = {[0, a], *, (p, 0), a, p ∈ $Z_n$} be a modular integer interval groupoid of level four. T = {[0, a], *, (p, 0), a ∈ X ⊆ $Z_n$} ⊆ S be a modulo integer interval groupoid of level four; then we call T to be a modulo integer interval subgroupoid of level four.*

*If S has no proper modulo integer interval subgroupoids then we call S to be a simple modulo integer interval groupoid of level four.*

We can analogously define these concepts in case of other interval groupoids of level four.
We will illustrate this by some examples.

*Example 2.3.45:* Let X = {[0, a], *, (8, 0), a, 8 ∈ $Z_{12}$} be a modulo integer interval groupoid of level four. Choose Y = {[0, 0], [0, 2], [0, 4], [0 6], [0, 8], [0, 10], *, (8, 0)} ⊆ X, Y is a modulo integer interval subgroupoid of level four of X.

***Example 2.3.46:*** Let W = {[0, a], *, (5, 0); a, 5 ∈ $Z_7$} be a modulo integer interval groupoid of level four. We see W is a simple modulo integer interval groupoid of level four.

In view of this we have the following theorem.

**THEOREM 2.3.2:** *Let P = {[0, a], *, (t, 0), a, t ∈ $Z_p$, p a prime} be a modulo integer interval groupoid of level four. P is simple.*

The proof of the theorem is left as an exercise for the reader.

***Example 2.3.47:*** Let T = {[0, a], *, (9, 0), a, 9 ∈ $Z^+ \cup \{0\}$} be an integer interval groupoid of level four. Take W = {[0, a], *, (9, 0), a ∈ $8Z^+ \cup \{0\}$} ⊆ T is an integer interval subgroupoid of level four of T.

***Example 2.3.48:*** Let W = {[0, a], *, (12, 0), a ∈ $Q^+ \cup \{0\}$} be an rational interval groupoid of level four. S = {[0, a], *, (12, 0), a ∈ $3Z^+ \cup \{0\}$} is an integer interval subgroupoid of level four of W.

***Example 2.3.49:*** Let T = {[0, a], *, (19, 0), a ∈ $R^+ \cup \{0\}$} be a real interval groupoid of level four. S = {[0, a], *, (19, 0), a ∈ $Q^+ \cup \{0\}$} ⊆ T; S is a real interval subgroupoid of level four.

***Example 2.3.50:*** Let W = {[0, a], *, (23, 0), a ∈ $C^+ \cup \{0\}$} be a complex interval groupoid of level four. P = {[0, a], *, (23, 0), a ∈ $R^+ \cup \{0\}$} ⊆ W; P is a complex interval subgroupoid of level four.

Now having seen examples of interval subgroupoid of level four. We now proceed onto study the properties of all these interval groupoids in the four levels.

**THEOREM 2.3.3:** *The modulo integer interval groupoids T in $G^{**}(Z_n^I)$ such that T = {[0, a]; *, (t, u), a, t, u ∈ $Z_n$}(mod n) is a modulo integer idempotent interval groupoid of level two if and only if t + u = 1(mod n).*

*Proof:* Given $T = \{[0, a]; *, (t, u), a, t, u \in Z_n\} \in G^{**}(Z_n^I)$ and $t + u \equiv 1 \pmod{n}$. To show T is an idempotent interval groupoid we have to prove $[0, a]*, [0, a] = [0, a]$ for all $[0, a] \in Z_n^I$.

Consider $[0, a] * [0, a] = [0, ta + ua \pmod{n}]$. If T is to be an idempotent interval groupoid we need $[0, a] * [0, a] = [0, a]$. So that $[0, ta + ua \pmod{n}] = [0, a]$. That is $ta + ua = a \pmod{n}$ thus $(t + u - 1)a \equiv 0 \pmod{n}$.

This is possible if and only if $t + u \equiv 1 \pmod{n}$. Hence the claim.

We will illustrate this by a simple example.

***Example 2.3.51:*** Let $T = \{[0, a], *, [4, 5], a, 4, 5 \in Z_8\} \in G^{**}(Z_8^I)$. T is a modulo integer idempotent interval groupoid. For consider $[0, a] * [0, a] = [0, 4a + 5a \pmod{8}] = [0, a]$ as $9a = a \pmod{8}$. This is true for all $[0, a] \in Z_8^I$. Hence T is a modulo integer idempotent interval groupoid.

**COROLLARY 2.3.1:** *The above theorem is also true in case of the class of interval groupoids of level one.*

We now proceed onto define normal subgroupoids and ideals of interval groupoids.

**DEFINITION 2.3.10:** *Let $G^*(Z_n^I)$ or $G^*(Z_I^+)$ or $G^*(Q_I^+)$ or $G^*(R_I^+)$ or $G^*(C_I^+))$ be a class of interval groupoids. A interval subgroupoid V of $T \subseteq G^*(Z_I^+)$ is said to be a normal interval groupoid or interval normal groupoid $T \subseteq G^*(Z_n^I)$ if*

  i. *[0, a] V = V [0, a]*
 ii. *(V [0, a]) [0, b] = V ([0, a] [0, b])*
iii. *[0, a] ([0, b] V) = ([0, a] [0, b])V*

*for all [0, a], [0, b] $\in$ V.*

*Note:* An interval groupoid $T \subseteq G^*(Z_n^I)$ is normal if

    (a) $[0, x] T = T [0, x]$
    (b) $T ([0, x] [0, y]) = (T [0, x]) [0, y]$
    (c) $[0, y] ([0, x] T) = ([0, y] [0, x]) T$

for all $[0, x], [0, y] \in T \subseteq G^*(Z_n^I)$.

This same definition holds good for interval groupoids built using $Z_I^+$, $Q_I^+$, $R_I^+$ or $C_I^+$. Further the same definition is true for all the four levels of interval groupoids.

Thus from here onwards by interval groupoid T we may mean any interval groupoid built using $Z_n^I$ or $Z_I^+$, $Q_I^+$, $R_I^+$ or $C_I^+$ and from the context it will be easily understood to which level they belong to and built using which set.

Now we proceed onto define ideal of interval groupoids.

**DEFINITION 2.3.11:** *Let T be any interval groupoid. P a non empty subset of T. P is said to be a left ideal of T if*

    *(1) P is an interval subgroupid of T*
    *(2) For all $x \in T$ and $a \in P$, $xa \in P$.*

*P is called a right ideal if P is an interval subgroupoid and for all $x \in T$ and $a \in P$ $ax \in P$. If P is both left and right ideal of T then we call P to be the interval ideal of T or just an ideal of T. As in case of usual groupoids we call an interval groupoid to be simple if it has nontrivial interval subgroupoids. The interval groupoids of level one, two, three and four does not have {0} as an ideal.*

We have the following interesting theorem.

**THEOREM 2.3.4:** *Let P be a left ideal of $T \subseteq G^*(Z_n^I)$ where T = {[0, a], \*, (t, u), a, t, u $\in Z_n$} is the interval groupoid then P is a right ideal of $T' = \{[0, a], *, (u, t); a, t, u \in Z_n\}$.*

The proof is left as an exercise for the reader.

**THEOREM 2.3.5:** *Let $T \subseteq G^*(Z_n^I)$ where $T = \{[0, a], *, (t, u), a, t, u \in Z_n$ with $t$ and $u$ primes and $t + u = n\}$ be an interval groupoid then $T$ is simple.*

The proof uses simple number theoretic properties and is left for the reader to prove.

*Note:* Even if n is a prime and $t + u = p$, u and t primes then also $T = \{[0, a], *, (t, u); t + u = p; t, a, u \in Z_p$, p a prime$\}$, the interval groupoid is simple.

**THEOREM 2.3.6:** *The interval groupoid $T = \{[0\ a], *, (0, t); a, t \in Z_n\}$ of level four is a interval P-groupoid and alternative interval groupoid if and only if $t^2 \equiv t \pmod{n}$.*

The proof is left as an exercise for the reader.
    None of the above results hold good in case of interval groupoids built using $Z_I^+$, $Q_I^+$, $R_I^+$ or $C_I^+$. These properties are only valid for the interval groupoids built using $Z_n^I$.

## 2.4 New Classes of Polynomial Interval Groupoids

In this section we introduce the new notion of polynomial interval groupoids. We will first briefly describe the essential notations. $Z_n^I[x] = \{$collection of all polynomials in the variable x with coefficients from the interval set $Z_n^I\}$

$$= \left\{\sum_{i=0}^{\infty}[0, a_i]x^i \,\middle|\, a_i \in Z_n\right\}.$$

This collection will be known as modulo integer interval polynomials.

*Example 2.4.1:* $[0, 1]x^7 + [0, 2]x^6 + [0, 1]x^5 + [0, 2]x^3 + [0, 1]$ is a modulo integer interval polynomials in the variable x with coefficients from $Z_3^I$. Similarly $Z_I^+[x] = \{$collection of all polynomials in the variable x with interval coefficients of the form [0, a] from the interval set $Z_I^+$ } =

$$\left\{ \sum_{i=0}^{\infty} [0,a_i]x^i \,\middle|\, a_i \in Z^+ \cup \{0\} \right\}.$$

This collection will be known as the positive integer interval polynomials or integer interval polynomials. $Q_I^+ [x]$ = {collection of all polynomials in the variable x with coefficients from the intervals in $Q_I^+$ of the form [0, $a_i$];

$$a_i \in Q_I^+ \cup \{0\}\} = \left\{ \sum_{i=0}^{\infty} [0,a_i]x^i \,\middle|\, a_i \in Q^+ \cup \{0\} \right\}.$$

Thus

$$R_I^+ [x] = \left\{ \sum_{i=0}^{\infty} [0,a_i]x^i \,\middle|\, a_i \in R^+ \cup \{0\} \right\}$$

and

$$C_I^+ [x] = \left\{ \sum_{i=0}^{\infty} [0,a_i]x^i \,\middle|\, a_i \in C^+ \cup \{0\} \right\}$$

are defined as real interval polynomials and complex interval polynomials respectively.

Now using these 5 types of interval polynomials we can define 4 levels of polynomial interval groupoids or interval polynomial groupoids.

**DEFINITION 2.4.1:** *Let*

$$T = \left\{ \sum_{i=0}^{\infty} [0,a_i]x^i \,\middle|\, a_i \in Z_n \right\}$$

*be the modulo integer interval polynomials. Define * on T as follows; let*

$$p(x) = \sum_{i=0}^{\infty} [0,a_i]x^i$$

*and*

$$q(x) = \sum_{j=0}^{\infty}[0,b_j]x^j$$

*belong to T;*

$$p(x) * q(x) = \left(\sum_{i=0}^{\infty}[0,a_i]x^i\right) * \left(\sum_{j=0}^{\infty}[0,b_j]x^j\right)$$

$$= \sum_{k=0}^{\infty}[0, ta_i + ub_j (mod\, n)]x^{i+j}\,;$$

*(t, u $\in Z_n \setminus \{0\}$ such that both t and u are primes and (t, u) = 1).*

$$= \sum_{k=0}^{\infty}[0,c_k]x^k$$

*where k = i + j and $c_k$ = $ta_i$ + $ub_j$ (mod n). [T, \*, (t, u)] is defined as the modulo integer interval polynomial groupoid of level one. We by varying t and u get a class of modulo integer interval polynomial groupoid of level one.*

We will illustrate this by some simple examples.

**Example 2.4.2:** Let $\{T, *, (u, t) \mid T = \sum_{i=0}^{\infty}[0,a_i]x^i\,;\, a_i \in Z_5$ and u = 3 and t = 2$\}$ be the modulo integer interval polynomial groupoid of level one.

**Example 2.4.3:** Let $\{T, *, (u, t) \mid T = \sum_{i=0}^{\infty}[0,a_i]x^i\,;\, a_i \in Z_{12}$ and u = 5 and t = 7$\}$ be the modulo integer interval polynomial groupoid of level one.

We can built rational interval polynomial groupoid of level one using the rational intervals $\{[0, a] \,/\, a \in Q^+ \cup \{0\}\}$, real interval polynomial groupoid of level one using real intervals of the form [0, a]; a $\in R^+ \cup \{0\}$, integer interval polynomial groupoid of level one using integer intervals of the form $\{[0, b] \,/\, b \in Z^+ \cup \{0\}\}$ and complex interval polynomial groupoid of level one using the complex numbers.

We will illustrate each of them by one example.

*Example 2.4.4:* Let
$$T = \left\{ \sum_{i=0}^{\infty}[0,a_i] \,;\, *, (5, 7); a_i, 5, 7 \in Z^+ \cup \{0\} \right\}$$
be the integer interval polynomial groupoid of level one using the integers $Z^+ \cup \{0\}$.

*Example 2.4.5:* Let
$$W = \left\{ \sum_{i=0}^{\infty}[0,a_i]x^i \,,\, *, (11, 43); a_i, 11, 43 \in Q^+ \cup \{0\} \right\}$$
be the rational interval polynomial groupoid of level one.

*Example 2.4.6:* Let
$$S = \left\{ \sum_{i=0}^{\infty}[0,a_i]x^i \,,\, *, (47, 19); a_i, 19, 47 \in R^+ \cup \{0\} \right\}$$
be the real interval polynomial groupoid of level one using reals.

*Example 2.4.7:* Let
$$V = \left\{ \sum_{i=0}^{\infty}[0,a_i]x^i \,;\, *, (23, 2); 2, 23, a_i \in C^+ \cup \{0\} \right\}$$
be the complex interval polynomial groupoid of level one using complex numbers.

We see each of the classes of interval polynomial groupoids of level one have infinite cardinality, since pairs of primes is an infinite collection. We can define interval polynomial subgroupoids of level one built using modulo integers or reals or rationals or integers or complex intervals.

**DEFINITION 2.4.2:** *Let T be a interval polynomial groupoid of level one. Let $P \subseteq T$ be a proper subset of T; if P itself is a interval polynomial groupoid of level one under the operations of T then we define P to be a interval polynomial subgroupoid of T of level one.*

*If T has no interval polynomial subgroupoids of level one then we define T to be a simple interval polynomial groupoid of level one.*

We will illustrate this situation by some simple examples.

*Example 2.4.8:* Let
$$T = \left\{ \sum_{i=0}^{\infty} [0, a_i] x^i \, ; \, *, (5, 7), a_i, 5, 7 \in Z_{10} \right\}$$
be a modulo integer interval polynomial groupoid of level one.
   Take
$$P = \left\{ \sum_{i=0}^{\infty} [0, a_i] x^i \, ; \, *, (5, 7), 5, 7 \in Z_{10} \, ; \, a_i \in \{0, 5\} \subseteq Z_{10} \right\} \subseteq T;$$
P is a modulo integer interval polynomial subgroupoid of level one.

*Example 2.4.9:* Let
$$S = \left\{ \sum_{i=0}^{\infty} [0, a_i] x^i \, ; \, *, (19, 11), a_i \in Z^+ \cup \{0\} \right\}$$
be a integer interval polynomial groupoid. Take
$$P = \left\{ \sum_{i=0}^{\infty} [0, a_i] x^i \, ; \, *, 19, 11), a_i \in 3Z^+ \cup \{0\} \right\} \subseteq S;$$
P is a integer interval polynomial subgroupoid of S.

*Example 2.4.10:* Let
$$T = \left\{ \sum_{i=0}^{\infty} [0, a_i] x^i \, ; \, *, (2, 3), a_i \in Q^+ \cup \{0\} \right\}$$
be a rational integer polynomial groupoid.
Let
$$S = \left\{ \sum_{i=0}^{\infty} [0, a_i] x^i \, ; *, (2, 3), a_i \in Z^+ \cup \{0\} \right\} \subseteq T;$$
S is a rational integer polynomial subgroupoid of T.

*Example 2.4.11:* Let
$$M = \left\{ \sum_{i=0}^{\infty} [0, a_i] x^i \, ; *, (29, 3), a_i \in R^+ \cup \{0\} \right\}$$
be a real interval polynomial groupoid.
   Take

$$W = \left\{ \sum_{i=0}^{\infty} [0, a_i] x^i ; *, (29, 3), a_i \in 3Z^+ \cup \{0\} \right\} \subseteq M.$$

W is a real interval polynomial subgroupoid of M.

*Example 2.4.12:* Let
$$V = \left\{ \sum_{i=0}^{\infty} [0, a_i] x^i ; *, (3, 5), a_i \in C^+ \cup \{0\} \right\}$$
be the complex interval polynomial groupoid.
Take
$$T = \left\{ \sum_{i=0}^{\infty} [0, a_i] x^i ; *, (3, 5), a_i \in R^+ \cup \{0\} \right\} \subseteq V,$$
T is a complex interval polynomial subgroupoid of V.

*Example 2.4.13:* Let
$$P = \left\{ \sum_{i=0}^{\infty} [0, a_i] x^i ; *, (3, 2), a_i \in Z_7 \right\}$$
be the modulo integer interval polynomial groupoid. It is not a simple modulo integer interval polynomial groupoid even though $Z_7$ is a prime field for
$$S = \left\{ \sum_{i=0}^{\infty} [0, a_i] x^{2i} \,\middle|\, a_i \in Z_7, *, (3, 2) \right\} \subseteq P$$
is a modulo integer interval polynomial subgroupoid of P of level one.

Now we proceed onto define level two interval polynomial groupoids using $Z_n$, $Z^+ \cup \{0\}$, $R^+ \cup \{0\}$, $Q^+ \cup \{0\}$ and $C^+ \cup \{0\}$.

**DEFINITION 2.4.3:** *Let*
$$S = \left\{ \sum_{i=0}^{\infty} [0, a_i] x^i ; \right.$$
*\*, (p, q), p and q are not primes but (p, q) = 1, p, q, $a_i \in Z_n$}. S under \* is a groupoid defined as the modulo integer interval polynomial groupoids of level two. Thus we can define interval*

*polynomial groupoids of level two using $Z^+ \cup \{0\}$ or $Q^+ \cup \{0\}$ or $R^+ \cup \{0\}$ or $C^+ \cup \{0\}$.*

We will illustrate this situation by some examples.

**Example 2.4.14:** Let
$$S = \left\{ \sum_{i=0}^{\infty} [0, a_i] x^i ; *, (2, 9), a_i, 2, 9 \in Z_{10} \right\}$$
be the modulo integer interval polynomial groupoid of level two.

**Example 2.4.15:** Let
$$P = \left\{ \sum_{i=0}^{\infty} [0, a_i] x^i, *, (27, 32), a_i, 27, 32 \in Z^+ \cup \{0\} \right\}$$
be the integer interval polynomial groupoid of level two.

**Example 2.4.16:** Let
$$S = \left\{ \sum_{i=0}^{\infty} [0, a_i] x^i ; a_i \in R^+ \cup \{0\}, *, (64, 81) \right\}$$
be the real interval polynomial groupoid of level two.

**Example 2.4.17:** Let
$$W = \left\{ \sum_{i=0}^{\infty} [0, a_i] x^i ; *, (14, 27), a_i, 14, 27 \in Q^+ \cup \{0\} \right\}$$
be the rational interval polynomial groupoid of level two.

**Example 2.4.18:** Let
$$B = \left\{ \sum_{i=0}^{\infty} [0, a_i] x^i, *, (36, 13), a_i, 13, 36 \in C^+ \cup \{0\} \right\}$$
be the complex interval polynomial groupoid of level two.

We can as in case of level one interval polynomial groupoids define interval polynomial subgroupoids of level two.

We will illustrate by one or two examples and then proceed on to define level three interval polynomial groupoids.

*Example 2.4.19:* Let
$$P = \{ \sum_{i=0}^{\infty}[0, a_i] x^i, *, (9, 7), a_i, 9, 7 \in Z_{12}\}$$
be the modulo integer interval polynomial groupoid of level two.
Take
$$W = \{ \sum_{i=0}^{\infty}[0, a_i] x^i, *, (9, 7), a_i \in \{0, 2, 4, 6, 8, 10\} \subseteq Z_{12}\} \subseteq P;$$
T is a modulo integer interval polynomial subgroupoid of level two.

*Example 2.4.20:* Let
$$S = \{ \sum_{i=0}^{\infty}[0, a_i] x^i, *, (3, 4), 3, 4, a_i, \in Z_{11}\}$$
be the modulo integer interval polynomial groupoid of level two.
Consider
$$M = \{ \sum_{i=0}^{\infty}[0, a_i] x^{2i} \mid a_i \in Z_{11}, *, (3, 4)\} \subseteq S,$$
M is a modulo integer interval polynomial subgroupoid of S of level two.

*Example 2.4.21:* Let
$$T = \{ \sum_{i=0}^{\infty}[0, a_i] x^i, *, (27, 43), a_i, 27, 43 \in Z^+ \cup \{0\}\}$$
be a integer interval polynomial groupoid of level two.
Take
$$W = \{ \sum_{i=0}^{\infty}[0, a_i] x^i, *, (27, 43), a_i \in 2Z^+\} \subseteq T,$$
W is an integer interval polynomial subgroupoid of level two.

*Example 2.4.22:* Let
$$V = \{ \sum_{i=0}^{\infty}[0, a_i] x^i, *, (27, 8), a_i, 27, 8 \in Q^+\}$$
be a rational polynomial groupoid of level two.
Choose

$$W = \{ \sum_{i=0}^{\infty}[0, a_i]x^i, *, (27, 8), a_i \in Z^+\} \subseteq V$$

is the rational interval polynomial subgroupoid of level two.

*Example 2.4.23:* Let

$$S = \{ \sum_{i=0}^{\infty}[0, a_i]x^i, *, (47, 81), a_i \in R^+ \cup \{0\}\}$$

be the real interval polynomial groupoid of level two.
Choose

$$W = \{ \sum_{i=0}^{\infty}[0, a_i]x^{2i}, *, (47, 81), a_i \in R^+ \cup \{0\}\} \subseteq S,$$

S is a real interval polynomial subgroupoid of S of level two.

*Example 2.4.24:* Let

$$S = \{ \sum_{i=0}^{\infty}[0, a_i]x^i, *, (27, 2), a_i \in C^+ \cup \{0\}\}$$

be a complex interval polynomial groupoid of S of level two.

$$T = \{ \sum_{i=0}^{\infty}[0, a_i]x^i, *, a_i \in R^+ \cup \{0\}\} \subseteq S,$$

T is a complex interval polynomial subgroupoid of S of level two.

We now proceed onto describe and define level three interval polynomial groupoid.

**DEFINITION 2.4.4:** *Let*

$$T = \{ \sum_{i=0}^{\infty}[0, a_i]x^i, *, (p, q) = d \neq 1, p, q, a_i \in Z^+ \cup \{0\}\}$$

*where * is such that*

$$(\sum_{i=0}^{\infty}[0, a_i]x^i) * (\sum_{i=0}^{\infty}[0, b_j]x^j) = \sum_{j+i=0}^{\infty}[0, a_i p + b_j q]x^{i+j},$$

*T is defined to be the integer interval polynomial groupoid of level three built using $Z_I^+$.*

We can on similar lines define polynomial interval groupoid of level three using modulo integers, rationals, reals and complex numbers.

We will only give examples of them and the reader can easily understand which type it belongs to by the context.

*Example 2.4.25:* Let

$$S = \{ \sum_{i=0}^{\infty}[0, a_i] x^i , *, (21, 6), 21, 6, a_i \in Z_{27}\}$$

be the modulo integer interval polynomial groupoid of level three.

*Example 2.4.26:* Let

$$V = \{ \sum_{i=0}^{\infty}[0, a_i] x^i , *, (24, 32), a_i, 24, 32, \in Q^+ \cup \{0\}\}$$

be the rational interval polynomial groupoid of level three.

*Example 2.4.27:* Let

$$B = \{ \sum_{i=0}^{\infty}[0, a_i] x^i , *, (27, 42), a_i, 42, 27, \in R^+ \cup \{0\}\}$$

be the real interval polynomial groupoid of level three.

*Example 2.4.28:* Let

$$C = \{ \sum_{i=0}^{\infty}[0, a_i] x^i , *, (2, 128), a_i, 2, 128, \in C^+ \cup \{0\}\}$$

be the complex interval groupoid of level three.

As in case of other types of interval groupoids we can in case of level three interval polynomial groupoids define subgroupoids.

However we will illustrate the situation by some examples.

*Example 2.4.29:* Let

$$T = \{ \sum_{i=0}^{\infty}[0, a_i] x^i , *, (8, 12), a_i, 8, 12, \in Z^+ \cup \{0\}\}$$

be the integer interval polynomial groupoid of level three.

Consider

$$W = \{ \sum_{i=0}^{\infty}[0, a_i]x^i, *, (8, 12), a_i \in 2Z^+ \cup \{0\}\} \subseteq T$$

is integer interval polynomial subgroupoid of T of level three.

*Example 2.4.30:* Let

$$S = \{ \sum_{i=0}^{\infty}[0, a_i]x^i, *, (4, 20), a_i, 4, 20, \in Q^+ \cup \{0\}\}$$

be the rational interval polynomial groupoid.
Take

$$M = \sum_{i=0}^{\infty}[0, a_i]x^{2i}, *, (4, 20), a_i, 4, 20 \in Q^+ \cup \{0\}\} \subseteq S,$$

M is a rational interval polynomial subgroupoid of S of level three.

*Example 2.4.31:* Let

$$S = \{ \sum_{i=0}^{\infty}[0, a_i]x^i, *, (5, 20), a_i, 5, 20 \in R^+ \cup \{0\}\}$$

be the real interval polynomial groupoid of level three.
Consider

$$W = \{ \sum_{i=0}^{\infty}[0\, a_i]x^{2i}; *, (5, 20), a_i, 5, 20 \in R^+ \cup \{0\}\} \subseteq S,$$

W is a real interval polynomial subgroupoid of level three.

*Example 2.4.32:* Let

$$M = \{ \sum_{i=0}^{\infty}[0, a_i]x^i, *, (27, 15), a_i, 27, 15 \in C^+ \cup \{0\}\}$$

be the complex interval polynomial groupoid of level three. Choose

$$N = \{ \sum_{i=0}^{\infty}[0\, a_i]x^{2i}; *, (27, 15), a_i, 27, 15 \in C^+ \cup \{0\}\} \subseteq M.$$

N is a complex interval polynomial subgroupoid of level three.

Now we will define level four interval polynomial groupoids.

**DEFINITION 2.4.5:** *Let*

$$T = \{ \sum_{i=0}^{\infty}[0, a_i]x^i, *, (t, 0), 0 \neq t, a_i \in R^+ \cup \{0\}\}$$

*where * is defined by*

$$(\sum_{i=0}^{\infty}[0, a_i]x^i) * (\sum_{i=0}^{\infty}[0, b_j]x^j) = \sum_{i+j=0}^{\infty}[0, ta_i]x^{i+j}.$$

*T is a real interval polynomial groupoid of level four.*
*As we vary t in $R^+$ we get a class of real interval polynomial groupoid of level four.*

Now we can define level four interval polynomial groupoids using $Z_n$ or $Z^+ \cup \{0\}$ or $Q^+ \cup \{0\}$ or $C^+ \cup \{0\}$.

We will illustrate all these situations by some examples.

***Example 2.4.33:*** Let

$$T = \{ \sum_{i=0}^{\infty}[0, a_i]x^i, *, (9, 0), 9, a_i \in Z_{12}\},$$

be a modulo integer interval polynomial groupoid of level four.

***Example 2.4.34:*** Let

$$T = \{ \sum_{i=0}^{\infty}[0, a_i]x^i, *, (11, 0), a_i, 11 \in Z_{13}\}$$

be the modulo integer interval polynomial groupoid of level four.

***Example 2.4.35:*** Let

$$S = \{ \sum_{i=0}^{\infty}[0, a_i]x^i, *, (12, 0), a_i, 12 \in Z^+ \cup \{0\}\}$$

be the integer interval polynomial groupoid of level four.

***Example 2.4.36:*** Let

$$N = \{ \sum_{i=0}^{\infty}[0, a_i]x^i, *, (0, 913), a_i, 913 \in R^+ \cup \{0\}\}$$

be the real interval polynomial groupoid of level four.

*Example 2.4.37:* Let
$$Z = \{\sum_{i=0}^{\infty}[0, a_i]x^i, *, (31, 0), a_i, 31 \in Q^+ \cup \{0\}\}$$
be the rational interval polynomial groupoid of level four.

*Example 2.4.38:* Let
$$S = \{\sum_{i=0}^{\infty}[0, a_i]x^i, *, (22, 0), a_i, 22 \in C^+ \cup \{0\}\}$$
be the complex interval polynomial groupoid of level four.

Now having seen examples of level four interval polynomial groupoids, we can as in case of other level groupoids define interval polynomial subgroupoids.

We only illustrate this situation by some examples.

*Example 2.4.39:* Let
$$T = \{\sum_{i=0}^{\infty}[0, a_i]x^i, *, (12, 0), 12, a_i \in Z_{24}\}$$
be the modulo integer interval polynomial groupoid of level four.
Take
$$P = \{\sum_{i=0}^{\infty}[0, a_i]x^i, *, (12, 0), a_i, 12 \in$$
$\{0, 2, 4, 6, 8, 10, 12, 16, 18, 20, 22\}\} \subseteq T$; P is a modulo integer interval polynomial subgroupoid of T of level four.

*Example 2.4.40:* Let
$$P = \{\sum_{i=0}^{\infty}[0, a_i]x^i, *, (3, 0), a_i, 3 \in Z_{17}\}$$
be the modulo integer interval polynomial groupoid of level four. Take
$$W = \{\sum_{i=0}^{\infty}[0, a_i]x^{2i}, *, (3, 0), a_i \in Z_{17}\} \subseteq P$$

be the modulo integer polynomial subgroupoid of P of level four.

*Example 2.4.41:* Let
$$W = \{ \sum_{i=0}^{\infty}[0, a_i] x^i , *, (9, 0), a_i, 9 \in Z^+ \cup \{0\}\}$$
be the integer interval polynomial groupoid of level four.

$$P = \{ \sum_{i=0}^{\infty}[0, a_i] x^{2i} , *, (9, 0), a_i \in 2Z^+ \cup \{0\}\} \subseteq W;$$

P is a integer interval polynomial subgroupoid of W of level four.

*Example 2.4.42:* Let
$$S = \{ \sum_{i=0}^{\infty}[0, a_i] x^i , *, (13, 0), a_i \; 13 \in Q^+ \cup \{0\}\}$$
be the rational interval polynomial groupoid of level four. Consider
$$X = \{ \sum_{i=0}^{\infty}[0, a_i] x^i , *, (13, 0), a_i, 13 \in Z^+ \cup \{0\}\} \subseteq S,$$
X is a rational interval polynomial subgroupoid of level four of S.

*Example 2.4.43:* Let
$$S = \{ \sum_{i=0}^{\infty}[0, a_i] x^i , *, (22, 0), a_i , 22 \in R^+ \cup \{0\}\}$$
be the real interval polynomial groupoid of level four.
Consider
$$Y = \{ \sum_{i=0}^{\infty}[0, a_i] x^{2i} , *, (22, 0), a_i , 22 \in R^+ \cup \{0\}\} \subseteq S;$$
Y is a real interval polynomial subgroupoid of S of level four.

*Example 2.4.44:* Let
$$Z = \{ \sum_{i=0}^{\infty}[0, a_i] x^i , *, (241, 0), a_i , 241 \in C^+ \cup \{0\}\}$$
be the complex interval polynomial groupoid of level four.

Take
$$W = \{ \sum_{i=0}^{\infty}[0, a_i] x^{2i}, *, (241, 0), a_i, 241 \in C^+ \cup \{0\}\} \subseteq Z,$$

W is a complex interval polynomial subgroupoid of Z of level four.

Now having defined four levels of interval polynomial groupoids we can define for them normal interval polynomial subgroupoid, normal interval polynomial groupoid and ideals in interval polynomial groupoid.

We now give some of the properties enjoyed by them. We see when the interval polynomial groupoids are built using intervals of the form [0, a] from $Z^+ \cup \{0\}$ or $Q^+ \cup \{0\}$ or $R^+ \cup \{0\}$ or $C^+ \cup \{0\}$. We see we cannot get many nice properties. Only for modulo integer interval polynomial groupoids which are built using $Z_n$ enjoy several interesting properties.

**THEOREM 2.4.1:** *Let $\{T, *, (u, t), u, t \in Z_n \setminus \{0\}, (u, t) = 1\}$ where*

$$T = \sum_{i=0}^{\infty}[0, a_i]x^i ;$$

*$a_i \in Z_n$ under * is a modulo integer polynomial groupoid of level two.*

$$\{T, *, (u, t); T = \sum_{i=0}^{\infty}[0, a_i]x^i ; a_i \in Z_n\}$$

*is a polynomial semigroup if and only if $t^2 \equiv t \pmod{n}$ and $u^2 \equiv u \pmod{n}$.*

The proof is left as an exercise for the reader.

No modulo integer polynomial groupoid has {0} as an ideal.

**THEOREM 2.4.2:** *Let*

$$X = \{T = \sum_{i=0}^{\infty}[0, a_i]x^i, *, (t, u), a_i, t, u \in Z_n\}$$

*and*

$$Y = \{S = \sum_{i=0}^{\infty}[0, a_i]x^i, *, (u, t), a_i, t, u \in Z_n\}$$

*be modulo integer interval polynomial groupoid of level one. P $\subseteq$ X is a left ideal of X if and only if P is a right ideal of Y.*

This proof is also left as an exercise for the reader.

Suppose we construct yet a new class of interval polynomial groupoids using $Z_n$ say

$$P = \{\sum_{i=0}^{\infty}[0, a_i]x^i; *, (t, t), t \in Z_n \setminus \{0\} \text{ and } a_i \in Z_n\}$$

be the modulo integer interval polynomial groupoid and if $\wp$ = {collection of all P's for varying t; $t \in Z_n$} then each P is a P-interval groupoid.

The following theorem is left as an exercise for the reader to prove.

**THEOREM 2.4.3:** *Let*

$$T = \{\sum_{i=0}^{\infty}[0, a_i]x^i; *, (t, t), t, a_i \in Z_n\}$$

*be the modulo integer interval polynomial groupoid. Clearly T is a P-groupoid.*

*Remark 2.4.1:* We have a class of n – 1 modulo integer interval polynomial groupoids built using $Z_n$. All these (n – 1) groupoids are P-groupoids. Further as n varies in $Z^+$ we get infinite classes of P-groupoids.

Next we will show that when n is a prime and $t \neq 1$ then the class of modulo integer interval polynomial groupoids are not alternative groupoids.

**THEOREM 2.4.4:** *Let*

$$T = \{\sum_{i=0}^{\infty}[0, a_i]x^i, *, (t, t), a_i \in Z_p, 1 < t < p\}$$

*p a prime be a modulo groupoid. Clearly T is not an alternative groupoid.*

*Proof:* Given T is a modulo integer interval polynomial groupoid. To show T is an alternative groupoid we have to prove $(x * y) * y = x * (y * y)$ for all $x, y \in T$.

Now

$$(x * y) * y = \left( \sum_{i=0}^{\infty} [0, a_i] x^i * \sum_{j=0}^{\infty} [0, b_j] x^j \right) * \left( \sum_{j=0}^{\infty} [0, b_j] x^j \right)$$

$$= \left( \sum_{j+i=0}^{\infty} [0, (a_i t + b_j t) \bmod p] x^{i+j} \right) * \sum_{j=0}^{\infty} [0\, b_j] x^j$$

$$= \sum_{j+i=0}^{\infty} [0, (a_i t^2 + b_j t^2 + b_j t)(\bmod p)] x^{i+j+j} \qquad \text{I}$$

$$x * (y * y) = \left( \sum_{i=0}^{\infty} [0, a_i] x^i \right) * \left( \sum_{j=0}^{\infty} [0, b_j] x^j * \sum_{j=0}^{\infty} [0, b_j] x^j \right)$$

$$= \left( \sum_{i=0}^{\infty} [0, a_i] x^i \right) \left[ \sum_{j=0}^{\infty} [0, [tb_j + tb_j](\bmod p)] x^{2j} \right]$$

$$= \sum_{j+i=0}^{\infty} [0, (ta_i + t^2 b_j + t^2 b_j)(\bmod p)] x^{2j+i} \qquad \text{II}$$

I and II can be equal if and only if $t^2 = t \pmod p$ but when p is a prime and $1 < t < p$, $t^2 \equiv t \pmod p$ is impossible.

Hence T is not an alternative interval polynomial groupoid.

**Remark 2.4.2:** *If*

$$T = \{ \sum_{i=0}^{\infty} [0, a_i] x^i , *, (t, t);\ 1 \neq t \in Z^+ \},$$

*T is not a alternative interval polynomial groupoid for $a_i \in R^+$ or $a_i \in Q^+$ or $a_i \in C^+$ or $a_i \in Z^+$. Hence we have a larger classes of interval polynomial groupoids which are not alternative.*

Proof is obvious as $t^2 = t$ is impossible for any $t \in Z^+ \setminus \{0\}$.

The following result can be easily proved by the reader.

**THEOREM 2.4.5:** *Let*

$$T = \{ \sum_{i=0}^{\infty} [0, a_i] x^i, *, (t, t), a_i \in Z_n, n \text{ not a prime number}\},$$

*T is an alternative modulo integer interval polynomial groupoid if and only if $t^2 \equiv t \pmod{n}$. We have as in case of groupoids some of the modulo integer interval polynomial groupoids to be Smaradache groupoids.*

*Example 2.4.45 :* Let

$$T = \{ \sum_{i=0}^{\infty} [0, a_i] x^i, *, (1, 5), ; a_i, 1, 5 \in Z_5\}$$

be a modulo integer interval polynomial groupoid of level two.
Take

$$S = \{ \sum_{i=0}^{\infty} [0, a_i] x^i, *, (1, 5), \{0, 5\} \subseteq Z_5\}.$$

It is easily verified S is a semigroup of interval polynomial. Thus T is a Smarandache modulo integer interval polynomial groupoid of level two.

Several properties enjoyed by general groupoids and Smarandache groupoids can be derived as a mater of routine. This is left as exercise for the reader.

Next we proceed onto define the notion of interval matrix groupoids using modulo integers, positive integers, rational integers, reals and complex.

## 2.5 Interval Matrix Groupoids

Throughout this section if $A = (a_{ij})$ is a matrix then $a_{ij}$'s are intervals from $Z_n$ or Q or R or C or $Q^+ \cup \{0\}$ or $R^+ \cup \{0\}$ or $Z^+ \cup \{0\}$ or Z. Now we proceed onto define different types of interval matrix groupoids.

**DEFINITION 2.5.1:** *Let $A = \{[a_1, ..., a_n]$ be such that each $a_i = [0\ x_i]$, $[0, x_i]$ is an interval in $Z_m^I$; $x_i \in Z_m$, $i=1, 2, ..., n\}$ ($Z_m$ –* *modulo integers). We call A a 1 ×n row modulo integer interval matrix.*

We will illustrate this by some example.

***Example 2.5.1:*** Let A = [[0, 2], [0 1], [0, 5], [0, 3], [0, 11], [0, 4], [0, 2], [0, 7], [0, 3]]; A is a 1 × 9 row modulo integer 12 matrix with intervals from $Z_{12}^I$.

***Example 2.5.2:*** Let A = [[0, 20], [0 2], [0 19] [0 11], [0 12], [0 15], [0 21], [0 1], [0 0], [0 3]] be a 1 × 10 row modulo integer interval matrix with entries from $Z_{23}^I$.

**DEFINITION 2.5.2:** *A = {[a$_1$, ..., a$_n$] is defined to be a 1 × n row integer interval matrix if a$_i$ ∈ Z$^I$}. If in A = [a$_1$, ..., a$_n$], a$_i$ ∈ Q$^I$, 1 ≤ i ≤ n then we define A to be a 1 × n row rational interval matrix.*

Likewise one can define a 1 × n row real interval matrix if the entries are from R$^I$ and 1 × n row complex interval matrix if the entries are from C$^I$.

We will illustrate each of these by examples.

***Example 2.5.3:*** Let A = [[0, 4], [–3, 2], [–7, 14], [15, 21], [11, 17] [15, 22] [–25, –2], [0, 24], [–2, 70]] be a 1 × 9 row integer interval matrix with entries from Z$^I$.

***Example 2.5.4:*** Let B = [[5, 9], [3, 7], [13, 18], [11, 12] [0, 47], [0, 1], [1, 21], [47, 48], [51, 112], [100, 109], [300, 308], [1000, 4001], [0, 1] [0, 21], [13, 18]] be a 1 × 15 row integer interval matrix or 1 × 15 integer interval row matrix with entries from (Z$^+$ ∪ {0})$_I$.

***Example 2.5.5:*** Let T = [[–5/2, 1], [7/25, 19], [19/3, 25], [–21/19, 19], [–43/3, –7], [3, 12], [0, 5], [–21/4, 0]] be a 1 × 8 row rational interval matrix 1 × 8 rational interval row matrix with entries from Q$^I$.

***Example 2.5.6:*** Let W = [[0, 2], [7/5, 20], [5/14, 27/4], [2/5, 21/13], [19/3, 42]] be a 1 × 5 row rational interval matrix with entries from $[Q^+ \cup \{0\}]_I = Q_I^+ \cup \{0\} = \{[a, b] / a, b \in Q^+ \cup \{0\}\}$.

***Example 2.5.7:*** Let P = [[0, $\sqrt{7}/8$], [$9/\sqrt{11}, \sqrt{29}$], [0, $\sqrt{5}$]] be a 1 × 3 row real interval matrix with entries from $R_I^+ \cup \{0\}$.

***Example 2.5.8:*** Let T = [[$-\sqrt{3}, \sqrt{11}$], [$\sqrt{2}, \sqrt{29}$], [3, $\sqrt{43}$], [$\sqrt{51}$ 53], [$-\sqrt{43}$, 0], [7, 15], [− 2, 0], [5/7, 9]] be a 1 × 8 row real interval matrix or 1 × 8 real interval row matrix with entries from $R^I$.

***Example 2.5.9:*** Let S = [[0, i], [− I, 1 + 4i], [2 + 4i, 26 + 47i], [−7i, 0], [−8 + I, 14 + 3i]] be a 1 × 5 row complex interval matrix with entries from $C^I$.

Now having seen 1 × n row interval matrix we now proceed onto define m × 1 column interval matrices. However from the context one can easily understand whether the interval row matrix is real or rational or complex or so on.

**DEFINITION 2.5.3:** *Let*

$$S = \begin{bmatrix} a_1 \\ a_2 \\ \vdots \\ a_m \end{bmatrix}$$

*be a m × 1 column matrix where $a_i$'s are intervals from Q or Z or $Z_n$ or R or C.*

*Then we define S to be a m × 1 column interval matrix or S is a interval m × 1 column matrix.*

We will illustrate this by some simple examples.

***Example 2.5.10:*** Let

$$S = \begin{bmatrix} [0,1] \\ [0,12] \\ [0,9] \\ [0,1] \\ [0,3] \\ [0,9] \end{bmatrix}$$

be a $6 \times 1$ modulo integer interval column matrix with entries from $Z_{13}^I$.

***Example 2.5.11:*** Let

$$V = \begin{bmatrix} [-1,0] \\ [-7,2] \\ [9,12] \\ [-11,-2] \\ [-1,0] \\ [10,16] \\ [25,41] \\ [-42,47] \end{bmatrix}$$

be a $8 \times 1$ column integer interval matrix with entries from $Z^I$.

***Example 2.5.12:*** Let

$$P = \begin{bmatrix} [0,\sqrt{7}] \\ [\sqrt{3},\sqrt{41}] \\ [-\sqrt{31},48] \\ [\sqrt{2},-17] \end{bmatrix}$$

be a $4 \times 1$ column real interval matrix with entries from $R^I$.

*Example 2.5.13:* Let
$$S = \begin{bmatrix} [\sqrt{3}, \sqrt{19}] \\ [20, 25] \\ [\sqrt{14}, \sqrt{247}] \end{bmatrix}$$
be a 3 × 1 real column interval matrix with entries from $R_I^+$.

*Example 2.5.14:* Let
$$W = \begin{bmatrix} [i, \ 2+3i] \\ [0, \ 4i] \\ [-4+i, 20+7i] \\ [-5-7i, \ 0] \\ [21, 14i+42] \\ [2i, 17i] \\ [3 \ \ 4] \\ [0, \ 5i] \\ [-1, \ 0] \\ [27, 27+i] \end{bmatrix}$$
be a 10 × 1 column complex interval matrix with entries from $C^I$.

Now we will define a usual interval matrix.

**DEFINITION 2.5.4:** *Let V = ($v_{ij}$) where 1 ≤ i ≤ n, 1 ≤ j ≤ m be a n × m matrix whose entries $v_{ij}$ are intervals from $Z_n^I$ or $Z^I$ or $Q^I$ or $R^I$ or $C^I$. We define V to be a n × m interval matrix.*

We will illustrate this situation by some examples.

*Example 2.5.15:* Let
$$A = \begin{bmatrix} [0,1] & [0,9] \\ [0,3] & [0,1] \\ [0,9] & [0,6] \end{bmatrix},$$
A is a 3×2 modulo integer interval matrix with entries from $Z_{11}^I$.

*Example 2.5.16:* Let

$$A = \begin{bmatrix} [-3,1] & [2,5] & [0,7] \\ [-7,0] & [-11,0] & [1,5] \\ [-3,1] & [0,0] & [1,1] \\ [2,5] & [-5,3] & [37,101] \\ [-7\ 7] & [1\ 6] & [-7\ 1] \end{bmatrix}$$

be a 5 × 3 integer interval matrix with entries from $Z^I$.

*Example 2.5.17:* Let

$$W = \begin{bmatrix} \left[\sqrt{7},14\right] & \left[\sqrt{3},\sqrt{11}\right] & [-7,0] \\ \left[-3,\sqrt{92}\right] & [0,2] & [-5,2] \\ \left[3,\sqrt{41}\right] & [-5/7,1] & [9,91] \end{bmatrix}$$

be a 3 × 3 real interval square matrix with entries from $R^I$.

*Example 2.5.18:* Let

$$W = \begin{bmatrix} [0,7] & [-3,0] \\ [2,4] & [1,10] \\ [3,7] & [-3,11] \\ [14,19] & [-16,\sqrt{247}] \\ [\sqrt{3},\sqrt{31}] & [\sqrt{5},\sqrt{53}] \\ [-1,2] & [7/9,3/2] \\ [3,3] & [0,1] \\ [0,7] & [14,19] \\ [1,10] & [-1,2] \end{bmatrix}$$

be a 9 × 2 real interval matrix. We say an interval n × m matrix A = ($a_{ij}$) is square if n = m and rectangular if n ≠ m.

***Example 2.5.19:*** Let
$$P = \begin{bmatrix} [2,9] & [3,8] & [7,9] \\ [0,7] & [-3,2] & [-7,0] \end{bmatrix}$$
be a 2 × 3 interval matrix.

Now having defined the notion of interval matrices now we proceed onto make them groupoids.
Let
$P(Z_n^I) = \{[0\ a] / a \in Z_n]$
$P(Z_I^+ \cup \{0\}) = \{[a, b] / b \neq 0, a \leq b\ a \in Z^+ \cup \{0\}\}.$
$P(Q_I^+ \cup \{0\}) = \{[a, b] / a, b \neq 0 \in Q^+ \cup \{0\}\ a \leq b\}.$
$P(R_I^+ \cup \{0\}) = \{[a, b] / a, b \neq 0 \in R^+ \cup \{0\}\ a \leq b\}$ and
$P(C_I^+ \cup \{0\}) = \{[a, b] = [x + iy, a_1 + ib_1]; x \leq a_1 \text{ and } y \leq b_1$
$\quad x, y, a, b \in R^+\}.$

We will now define operations on the interval matrices which takes its entries from $P(Z_n^I)$, $P(Q_I^+ \cup \{0\})$, $P(R_I^+ \cup \{0\})$, $P(Z_I^+ \cup \{0\})$ and $P(C_I^+ \cup \{0\})$.

**DEFINITION 2.5.5:** *Let $G_{Z_n}^{1 \times m} = \{[0, a_1], ..., [0, a_m]] / a_i \in Z_n, 1 \leq i \leq m\}$ denote the collection of all $1 \times m$ interval row matrices. Define a binary operation $*$ on $G_{Z_n}^{1 \times m}$ as follows.*

*Let $A = \{[0, a_1], ..., [0, a_m]]$ and $B = [[0, b_1], ..., [0, b_m]] \in G_{Z_n}^{1 \times m}$ where $a_i, b_i \in Z_n; 1 \leq i \leq m$.*

$A * B = \{[0, a_1], ..., [0, a_m]] * [[0, b_1], ..., [0, b_m]]$
$\quad = [[0, a_1]*[0, b_1], ..., [0, a_m][0, b_m]].$
$\quad = [[0, ta_1 + ub_1 (mod\ n)], ..., [0, ta_m + ub_m (mod\ n)]]$

*where $t, u \in Z_n \setminus \{0\}$. Clearly $(G_{Z_n}^{1 \times m}, *, (t, u))$ is defined as the row interval $1 \times m$ matrix groupoid using $Z_n$.*

We will illustrate this situation by some examples.

***Example 2.5.20:*** Let $G_{Z_5}^{1\times 3} = \{[[0, a_1], [0, a_2], [0, a_3]]/ a_i \in Z_5, 0 \leq a_i \leq 3\}$. Choose $(2, 3) = (t, u)$; where $2, 3 \in Z_5$. $\{G_{Z_5}^{1\times 3}, *, (t, u))$ is a groupoid, that is it is a row interval $1 \times 3$ matrix groupoid. Choose

$$A = [[0, 1], [0, 3], [0, 2]]$$
and $$B = [[0, 4], [0, 1], [0, 2]]$$
in $G_{Z_5}^{1\times 3}$.

$$\begin{aligned}A * B &= [[0, 1], [0, 3], [0, 2]] * [[0, 4], [0, 1], [0, 2]] \\ &= [[0, 1] * [0, 4], [0, 3] * [0, 1], [0, 2] * [0, 2]] \\ &= [[0, 2+12 \text{ (mod 5)}], [0, (6+3) \text{ mod } 5], \\ &\quad [0, 4 + 6 \text{ (mod 5)}]] \\ &= [[0, 4], [0, 4], [0, 0]].\end{aligned}$$

It is easily verified that '*' on $G_{Z_5}^{1\times 3}$ is non associative.

***Example 2.5.21:*** Let $\{G_{Z_3}^{1\times 4}, *, (1, 2)\}$ be a row matrix interval groupoid.

Next we can give examples of groupoids using $Z^+$ or $R^+$ or $Q^+$ or $C^+$.

***Example 2.5.22:*** Let $G_{Z_I^+}^{1\times 5} = \{[a_1, a_2], [a_3, a_4], [a_5, a_6], [a_7, a_8], [a_9, a_{10}] \mid a_i \in Z^+ \cup \{0\}, a_i \leq a_{i+1} \ 1 \leq i \leq 9\}$. $\{G_{Z_I^+}^{1\times 5}, *, (3, 7)\}$ is a row matrix interval groupoids using $Z^+ \cup \{0\}$.

***Example 2.5.23:*** Let $G_{Z_I^+}^{1\times 27} = \{[[a_1, a_2], \ldots, [a_{53}, a_{54}]]; a_i < a_{i+1}; 1 \leq i \leq 53; a_i \in Z^+ \cup \{0\}\}$, $\{G_{Z_I^+}^{1\times 27}, *, (2, 53)\}$ is a row matrix interval groupoid of row interval matrix groupoid using $Z^+ \cup \{0\}$.

***Example 2.5.24:*** Let $G_{Q_I^+}^{1\times 2} = \{[a_1, a_2] [a_3, a_4]] \mid a_i \in Q^+ \cup \{0\}; 1 \leq i \leq 4; a_i < a_{i+1}\}$, $\{G_{Q_I^+}^{1\times 2}, *, (7, 13)\}$ is a row interval matrix groupoid using $Q^+ \cup \{0\}$.

***Example 2.5.25:*** Let $G_{R_I^+}^{1\times 8} = \{[a_1, a_2]\ [a_3, a_4], \ldots, [a_{15}, a_{16}]\ a_i < a_{i+1}; 1 \le i \le 15\ |a_i \in R^+ \cup \{0\}\}$; $\{G_{R_I^+}^{1\times 8}, *, (7, 2)\}$ is a row interval matrix groupoid using $R^+ \cup \{0\}$.

Now we proceed onto define column matrix interval groupoid.

**DEFINITION 2.5.6:** *Let*

$$G_{Q_I^+}^{n\times 1} = \left\{ \begin{bmatrix} [a_1, a_2] \\ [a_3, a_4] \\ \vdots \\ [a_{2n-1}, a_{2n}] \end{bmatrix} \middle| \begin{array}{l} a_i \le a_{i+1} \\ 1 \le i \le 2n \\ a_i \in Q^+ \cup \{0\} \end{array} \right\}$$

*be a column interval matrix using $Q^+ \cup \{0\}$.*
   *Define * on $G_{Q_I^+}^{n\times 1}$ as*

$A * B$ = *{[$a_i\ a_{i+1}$]} * {[$b_i, b_{i+1}$]}*
   = *{[$ta_i + ub_i, ta_{i+1} + ub_{i+1}$]}*
   = *$tA + uB$;*

$t \ne u$, $(t, u) = 1$, $t, u \in Z^+$. *{$G_{Q_I^+}^{n\times 1}$, *, $(t, u)$} is a column interval matrix groupoid (column matrix interval groupoid).*

We will illustrate this situation by some examples.

***Example 2.5.26:*** Let

$$G_{Z_I^+}^{6\times 1} = \left\{ \begin{bmatrix} [a_1, a_2] \\ [a_3, a_4] \\ [a_5, a_6] \\ [a_7, a_8] \\ [a_9, a_{10}] \\ [a_{11}, a_{12}] \end{bmatrix} \middle| \begin{array}{l} a_i \le a_{i+1} \\ 1 \le i \le 11 \\ a_i \in Z^+ \cup \{0\} \end{array} \right\}$$

be a column interval matrix groupoid (column interval matrix groupoid).

$\{ G_{Z_I^+}^{6\times 1}, *, (5, 19)\}$ with binary operation *, that is for A, B ∈ $G_{Z_I^+}^{6\times 1}$ A* B = 5A + 19B.

**Example 2.5.27:** Let ($G_{C_I^+}^{4\times 1}$, *, (7, 1)) where

$$G_{C_I^+}^{4\times 1} = \left\{ \begin{bmatrix} [a_1,a_2] \\ [a_3,a_4] \\ [a_5,a_6] \\ [a_7,a_8] \end{bmatrix} \middle| \begin{array}{l} a_i \le a_{i+1} \\ 1 \le i \le 8 \\ a_i \in C^+ \cup \{0\} \end{array} \right\}$$

be a column interval matrix groupoid built using $C^+ \cup \{0\}$.

**Example 2.5.28:** Let

$$G_{Z_3}^{3\times 1} = \left\{ \begin{bmatrix} [a_1,a_2] \\ [a_3,a_4] \\ [a_5,a_6] \end{bmatrix} \middle| \begin{array}{l} a_i \le a_{i+1} \\ 1 \le i \le 5 \\ a_i \in Z_3 \end{array} \right\}$$

be the collection of all matrix column groupoids with * on $G_{Z_3}^{3\times 1}$; A * B = 2A + B (mod 3).

**DEFINITION 2.5.7:** *Let $G_{Q_I^+}^{m\times n} = \{A = \{[a_{ij}^1, a_{ij}^2]\}_{m\times n}$ be a m × n interval matrix; $a_{ij}^1, a_{ij}^2 \in Q^+ \cup \{0\}$; $a_{ij}^1 \le a_{ij}^2$; $1 \le i \le m$; $1 \le j \le n\}$.*

*Define * on $G_{Q_I^+}^{m\times n}$ by $A*B = tA + uB$; $t, u \in Q^+$; then $G_{Q_I^+}^{m\times n}$, *, (t, u); t ≠ u, (t, u) = 1} is a m × n interval groupoid built using $Q_I^+$.*

We can construct m×n interval groupoid using $Z_n^I$, $Z_I^+$, $R_I^+$, $Q_I^+$ and $C_I^+$.

We will illustrate this situation by some examples.

*Example 2.5.29:* Let

$$G_{Z_{10}^I}^{3\times 2} = \left\{ \begin{bmatrix} [a_1,a_2] & [a_3,a_4] \\ [a_5,a_6] & [a_7,a_8] \\ [a_9,a_{10}] & [a_{11},a_{12}] \end{bmatrix} \middle| \begin{array}{l} a_i \leq a_{i+1} \\ 1 \leq a_i \leq 8 \\ a_i \in Z_{10} \end{array} \right\}$$

be a $3 \times 2$ interval matrix using the groupoid $Z_{10}$. $\{Z_{10}^I, *, (5, 7)\}$ is a $3 \times 2$ interval matrix groupoid using $Z_{10}$.

*Example 2.5.30:* Let

$$G_{Z_I^+}^{5\times 3} = \left\{ \begin{bmatrix} [a_1^1,a_2^1] & [a_3^1,a_4^1] & [a_5^1,a_6^1] \\ [a_1^2,a_2^2] & [a_3^2,a_4^2] & [a_5^2,a_6^2] \\ [a_1^3,a_2^3] & [a_3^3,a_4^3] & [a_5^3,a_6^3] \end{bmatrix} \middle| \begin{array}{l} a_j^i \in Z_I^+ \cup \{0\} \\ 1 \leq i \leq 3 \\ 1 \leq j \leq 6 \end{array} \right\}$$

$a_j^i \leq a_{j+1}^i$, $j = 1, 2, \ldots, 6\}$ be a $5 \times 3$ interval matrix with entries from $Z^+ \cup \{0\}$ or $Z_I^+$. $\{G_{Z_I^+}^{5\times 3}, *, (3, 19)\}$ is a $5 \times 3$ interval matrix groupoid using entries from $Z^+ \cup \{0\}$.

*Example 2.5.31:* Let

$$G_{R_I^+}^{2\times 5} = \left\{ \begin{bmatrix} [a_1^1,a_2^1] & [a_3^1,a_4^1] & [a_5^1,a_6^1] & [a_7^1,a_8^1] & [a_9^1,a_{10}^1] \\ [a_{11}^2,a_{12}^2] & [a_{13}^2,a_{14}^2] & [a_{15}^2,a_{16}^2] & [a_{17}^2,a_{18}^2] & [a_{19}^2,a_{20}^2] \end{bmatrix} \middle| a_j^i \in R_I^+ \cup \{0\}; a_j^i \leq a_{j+1}^i; 1 \leq i \leq 2; 1 \leq j \leq 19 \right\}$$

be a $2 \times 5$ interval matrix built using $R^+ \cup \{0\}$. $\{G_{R_I^+}^{2\times 5}, *, (3, 2)\}$ is a $2 \times 5$ interval matrix groupoid using $R_I^+$.

*Example 2.5.32:* Let

$$G_{Q_I^+}^{4\times 3} = \left\{ \begin{bmatrix} [a_1^1,a_2^1] & [a_3^1,a_4^1] & [a_5^1,a_6^1] \\ [a_1^2,a_2^2] & [a_3^2,a_4^2] & [a_5^2,a_6^2] \\ [a_1^3,a_2^3] & [a_3^3,a_4^3] & [a_5^3,a_6^3] \\ [a_1^4,a_2^4] & [a_3^4,a_4^4] & [a_5^4,a_6^4] \end{bmatrix} \middle| \begin{array}{l} a_j^i \in Q^+ \cup \{0\} \\ a_j^i \leq a_{j+1}^i \\ 1 \leq i \leq 4 \\ 1 \leq j \leq 5 \end{array} \right\}$$

be a 4 × 3 interval matrix with entries from $Q_I^+$. $\{G_{Q_I^+}^{4\times 3}, *, (3, 19)\}$ is a interval matrix groupoid with entries from $Q_I^+$.

If m = n then we see the interval groupoid is a square interval groupoid.

We will illustrate this by some examples.

*Example 2.5.33:* Let

$$G_{Z_I^+}^{3\times 3} = \left\{ \begin{bmatrix} [a_1^1,a_2^1] & [a_3^1,a_4^1] & [a_5^1,a_6^1] \\ [a_1^2,a_2^2] & [a_3^2,a_4^2] & [a_5^2,a_6^2] \\ [a_1^3,a_2^3] & [a_3^3,a_4^3] & [a_5^3,a_6^3] \end{bmatrix} \middle| \begin{array}{l} a_j^i \in Z^+ \cup \{0\} \\ a_j^i \leq a_{j+1}^i \\ 1 \leq i \leq 3 \\ 1 \leq j \leq 5 \end{array} \right\}$$

be a 3 × 3 interval matrix with entries from $Z_I^+$. $\{G_{Z_I^+}^{3\times 3}, *, (3, 13)\}$ is a interval groupoid matrix with entries from $Z_I^+$.

*Example 2.5.34:* Let

$$G_{Z_8^I}^{5\times 5} = \left\{ \begin{bmatrix} [a_1^1,a_2^1] & [a_3^1,a_4^1] & [a_5^1,a_6^1] & [a_7^1,a_8^1] & [a_9^1,a_{10}^1] \\ [a_1^2,a_2^2] & [a_3^2,a_4^2] & [a_5^2,a_6^2] & [a_7^2,a_8^2] & [a_9^2,a_{10}^2] \\ [a_1^3,a_2^3] & [a_3^3,a_4^3] & [a_5^3,a_6^3] & [a_7^3,a_8^3] & [a_9^3,a_{10}^3] \\ [a_1^4,a_2^4] & [a_3^4,a_4^4] & [a_5^4,a_6^4] & [a_7^4,a_8^4] & [a_9^4,a_{10}^4] \\ [a_1^5,a_2^5] & [a_3^5,a_4^5] & [a_5^5,a_6^5] & [a_7^5,a_8^5] & [a_9^5,a_{10}^5] \end{bmatrix} \middle| a_j^i \in Z_8; a_j^i \leq a_{j+1}^i; 1 \leq i \leq 5; 1 \leq j \leq 9 \right\}$$

be the collection of intervals from $Z_8^I$. $\{G_{Z_8^I}^{5\times 5}, *, (3, 5)\}$ is a 5×5 square interval groupoid with entries from $Z_8$.

*Example 2.5.35:* Let

$$G_{R_I^+}^{2\times 2} = \left\{ \begin{bmatrix} \left[a_1^1,a_2^1\right] & \left[a_3^1,a_4^1\right] \\ \left[a_1^2,a_2^2\right] & \left[a_3^2,a_4^2\right] \end{bmatrix} \middle| \begin{array}{l} a_j^i \in R_I^+ \\ a_j^i \leq a_{j+1}^i \\ 1 = 1,2. \\ 1 \leq j \leq 3 \end{array} \right\}$$

be a 2 × 2 interval matrix with entries from $R_I^+$. { $G_{R_I^+}^{2\times 2}$, *, (3, 23)} is a 2 × 2 square interval groupoid from $R_I^+$.

*Example 2.5.36:* Let

$$G_{Q_I^+}^{16\times 16} = \left\{ \left[a_{ij}^t, a_{i+1j+1}^t\right] \middle| \begin{array}{l} 1 \leq t \leq 16, a_{ij}^t \leq a_{i+1j+1}^t, \\ 1 \leq i,j \leq 15, a_{ij}^t \in Q^+ \cup \{0\} \end{array} \right\}$$

be the set of all 16 × 16 interval matrices with entries from $Q_I^+$. { $G_{Q_I^+}^{16\times 16}$, *, (7, 19)} is a square matrix interval groupoid with entries from $Q_I^+$.

*Example 2.5.37:* Let

$$G_{Z_{29}^I}^{8\times 8} = \left\{ \left[a_{ij}^t, a_{i+1j+1}^t\right] \middle| \begin{array}{l} a_{ij}^t \in Z_{29}, 1 \leq t \leq 8, \\ a_{ij}^t \leq a_{i+1j+1}^t, 1 \leq i,j \leq 15 \end{array} \right\}$$

be a collection 8 × 8 interval matrices with entries from $Z_{29}^I$. { $G_{Z_{29}^I}^{8\times 8}$, *, (17, 3)} be a square matrix interval groupoid with entries from $Z_{29}^I$.

Now we proceed onto define some more properties of these matrix interval groupoid.

**DEFINITION 2.5.8:** *Let { $G_{Q_I^+}^{1 \times n}$, *, (t, u), t, u ∈ Q} be a row interval matrix groupoid built over $Q_I^+$.*

*Let $P \subseteq G_{Q_I^+}^{1 \times n}$, if {P, *, (t, u)} is itself a row interval matrix groupoid then we say {P, *, (t, u)} to be a row interval matrix subgroupoid of { $G_{Q_I^+}^{1 \times n}$, *, (t, u)} built using $Q_I^+$.*

We can have subgroupoid of row matrix interval groupoids constructed using $Z_I^+$, $Z_n^I$, $R_I^+$, $C_I^+$ and $Q_I^+$.

We will illustrate this situation by some examples.

***Example 2.5.38:*** Let $G_{Z_{12}^I}^{1 \times 5}$ = {([$a_1$, $a_2$], [$a_3$, $a_4$], [$a_5$, $a_6$], [$a_7$, $a_8$], [$a_9$, $a_{10}$]), *, (3, 5), $a_i \in Z_{12}$; $1 \le i \le 10$; $a_j \le a_{j+1}$} be a row matrix interval groupoid built using $Z_{12}^I$.

Take P = {([$a_1$, $a_2$], [$a_3$, $a_4$], [$a_5$, $a_6$], [$a_7$, $a_8$], [$a_9$, $a_{10}$]) | $a_i \in$ {0, 2, 4, 6, 8, 10}, *, (3, 5), $a_i \le a_{i+1}$} $\subseteq G_{Z_{12}^I}^{1 \times 5}$. P is a row matrix interval subgroupoid of $G_{Z_{12}^I}^{1 \times 5}$.

S = {([$a_1$, $a_2$] [$a_3$, $a_4$] [$a_5$, $a_6$], [$a_7$, $a_8$] [$a_9$, $a_{10}$]) | $a_i \in$ {0, 3, 6, 9}, *, (3, 5), $a_i \le a_{i+1}$} $\subseteq G_{Z_{12}^I}^{1 \times 5}$ is a row matrix interval subgroupoid of $G_{Z_{12}^I}^{1 \times 5}$.

***Example 2.5.39:*** Let $G_{Z_I^+}^{1 \times 8}$ = {([$a_1$, $a_2$], …, [$a_{15}$, $a_{16}$] / $a_i \in Z_I^+$; $a_i \le a_{i+1}$, i = 1, 2, …, 15; *, (17, 2)} be a row matrix interval groupoid built using $Z_I^+$.

Consider P = {([$a_1$, $a_2$], …, [$a_{15}$, $a_{16}$]) | $a_i \in 6 Z_I^+$, i = 1, 2, …, 15, *, (17, 2)} $\subseteq G_{Z_I^+}^{1 \times 8}$; P is a row matrix interval subgroupoid of $G_{Z_I^+}^{1 \times 8}$, built using $Z_I^+$.

Take L = {[$a_1$, $a_2$], [$a_3$, $a_4$], [0 0], [0 0], [0 0], [0 0], [0 0], [0 0]] | $a_1$, $a_2$, $a_3$, $a_4 \in Z_I^+$; $a_i \le a_{i+1}$, *, (17, 2)} $\subseteq G_{Z_I^+}^{1 \times 8}$; L is also a row matrix interval subgroupoid of $G_{Z_I^+}^{1 \times 8}$.

*Example 2.5.40:* Let $G_{R_I^+}^{1\times 2} = \{[a_1, a_2], [a_3, a_4]\ a_i \in R^+ \cup \{0\}; a_i \leq a_{i+1}; 1 \leq i \leq 3, *, (31, 43)\}$ be a row matrix interval groupoid built using $R_I^+$.

Consider $P = \{([a_1, a_2], [a_3, a_4]); a_i \leq a_{i+1}, 1 \leq i \leq 3, *, (31, 43); a_i \in Q^+ \cup \{0\}\} \subseteq G_{R_I^+}^{1\times 2}$; P is a row matrix interval subgroupoid of $G_{R_I^+}^{1\times 2}$.

Consider $T = \{([a_1, a_2], [a_3, a_4]); a_i \in Z^+ \cup \{0\}; a_i \leq a_{i+1}, i = 1, 2, 3; *, (31, 43)\} \subseteq G_{R_I^+}^{1\times 2}$; T is also a row matrix interval subgroupoid of $G_{R_I^+}^{1\times 2}$. Infact it is easy to verify $T \subseteq P \subseteq G_{R_I^+}^{1\times 2}$.

On similar lines we can define interval matrix subgroupoids in case of column matrix interval groupoid and m×n interval matrix groupoid.

We will leave this easy concept to the reader, however substantiate this by some examples.

*Example 2.5.41:* Let

$$G_{Z_{20}^I}^{3\times 1} = \left\{ \begin{bmatrix} [a_1, a_2] \\ [a_3, a_4] \\ [a_5, a_6] \end{bmatrix} \middle| a_i \in Z_{20}; a_i \leq a_{i+1}; 1 \leq i \leq 5 \right\}$$

be a column interval matrix using $Z_{20}^I$. $G = \{G_{Z_{20}^I}^{3\times 1}, *, (5, 7)\}$ is a column interval matrix groupoid.
Take

$$P = \left\{ \begin{bmatrix} [a_1, a_2] \\ [a_3, a_4] \\ [0, 0] \end{bmatrix} \middle| \begin{array}{l} a_i \in Z_{20} \\ 1 \leq i \leq 3, \\ a_i \leq a_{i+1} \end{array} \right\}, *, (5, 7)\}$$

$\subseteq \{G_{Z_{20}^I}^{3\times 1}, *, (5, 7)\}$. P is a column interval matrix subgroupoid. If we take

$$T = \left\{ \begin{bmatrix} [a_1, a_2] \\ [a_3, a_4] \\ [a_5, a_6] \end{bmatrix} \middle| a_i \in \{0, 2, 4, 6, ..., 18\} \subseteq Z_{20}, *, (5, 7), \right.$$

$a_i \leq a_{i+1}; \ 1 \leq i \leq 5 \} \subseteq \{G^{3\times1}_{Z^I_{20}}, *, (5, 7)\} = G$; T is a column interval matrix subgroupoid of G.

*Example 2.5.42:* Consider $G = \{G^{6\times1}_{Q^+_I}, *, (3, 17)\}$ where

$$G^{6\times1}_{Q^+_I} = \left\{ \begin{bmatrix} [a_1, a_2] \\ [a_3, a_4] \\ [a_5, a_6] \\ [a_7, a_8] \\ [a_9, a_{10}] \\ [a_{11}, a_{12}] \end{bmatrix} \middle| \begin{array}{l} a_i \in Q^+ \cup \{0\}; \\ a_i \leq a_{i+1}; \\ 1 \leq i \leq 11 \end{array} \right\}.$$

G is a column interval matrix groupoid built using $Q^+_I$.

Take

$$P = \left\{ \begin{bmatrix} [a_1, a_2] \\ [0, 0] \\ [a_3, a_4] \\ [0, 0] \\ [a_5, a_6] \\ [0, 0] \end{bmatrix} \middle| \begin{array}{l} a_i \in Z^+ \cup \{0\}; \\ a_i \leq a_{i+1}; \\ 1 \leq i \leq 5 \end{array} \right\} \subseteq G;$$

P is a column interval matrix subgroupoid of G.

Take

$$L = \left\{ \begin{bmatrix} [a_1, a_2] \\ [a_3, a_4] \\ [a_5, a_6] \\ [a_7, a_8] \\ [a_9, a_{10}] \\ [a_{11}, a_{12}] \end{bmatrix} \middle| \begin{array}{l} a_i \in 2Z^+ \cup \{0\}; \\ 1 \leq i \leq 11; \\ a_i \leq a_{i+1} \end{array} \right\} \subseteq G;$$

L is a column interval matrix subgroupoid of G.

*Example 2.5.43:* Consider

$$G = G_{R_I^+}^{10 \times 1} = \left\{ \begin{bmatrix} [a_1, a_2] \\ [a_3, a_4] \\ \vdots \\ [a_{19}, a_{20}] \end{bmatrix} \middle| \begin{array}{l} a_i \leq a_{i+1}; \\ 1 \leq i \leq 19; \\ a_i \in R^+ \cup \{0\} \end{array} \right\}, *, (13, 2)\}$$

be a column interval matrix groupoid.
Let

$$P = \left\{ \begin{bmatrix} [a_1, a_2] \\ [a_3, a_4] \\ \vdots \\ [a_{19}, a_{20}] \end{bmatrix} \middle| \begin{array}{l} a_i \leq a_{i+1}; \\ 1 \leq i \leq 19; \\ a_i \in Q^+ \cup \{0\} \end{array} \right\}; *, (13, 2)\} \subseteq G;$$

P is a column interval matrix subgroupoid of G.

*Example 2.5.44:* Let

$$G = \left\{ \begin{bmatrix} [a_1, a_2] \\ [a_3, a_4] \end{bmatrix} \middle| a_i \in C; i = 1, 2, 3, 4 \right.$$

are proper intervals, *, (3, 5)} be a column interval matrix groupoid.
Take

$$T = \left\{ \begin{bmatrix} [a_1, a_2] \\ [0, 0] \end{bmatrix} \middle| a_1, a_2 \in C; *, (3, 5) \right\} \subseteq G;$$

T is a column interval matrix subgroupoid of G.

*Example 2.5.45:* Let

$$G = G_{Z_{36}^I}^{3 \times 2} = \left\{ \begin{bmatrix} [a_1, a_2] & [a_3, a_4] \\ [a_5, a_6] & [a_7, a_8] \\ [a_9, a_{10}] & [a_{11}, a_{12}] \end{bmatrix} \middle| \begin{array}{l} a_i \in Z_{36;} \\ a_i \leq a_{i+1}; *, (3, 11)\} \\ 1 \leq i \leq 11 \end{array} \right.$$

be a 3 × 2 interval matrix groupoid built using $Z_{36}^I$.
Consider

$$P = \left\{ \begin{bmatrix} [a_1, a_2] & [a_3, a_4] \\ [0,0] & [a_7, a_8] \\ [a_9, a_{10}] & [0,0] \end{bmatrix} \middle| \begin{array}{l} a_i \in Z_{36}; \\ a_i \leq a_{i+1}; \\ 1 \leq i \leq 7; \end{array} \right.$$

*, (3, 11)} ⊆ G}, P is a 3 × 2 interval matrix subgroupoid of G.

*Example 2.5.46:* Let

$$G = \left\{ \begin{bmatrix} [a_1, a_2] & [a_3, a_4] \\ [a_5, a_6] & [a_7, a_8] \end{bmatrix} \right.$$

$a_i \in Q^+ \cup \{0\}; a_i \leq a_{i+1}; 1 \leq i \leq 7; *, (19, 7)\}$ be a 2 × 2 interval matrix groupoid built using $Q_I^+$.
Take

$$P = \left\{ \begin{bmatrix} [a_1, a_2] & [a_3, a_4] \\ [a_5, a_6] & [a_7, a_8] \end{bmatrix} \right.$$

$a_i \in Q^+ \cup \{0\}; a_i \leq a_{i+1}; 1 \leq i \leq 7$, *, (19, 7)} ⊆ G, P is a square interval matrix subgroupoid of G.

*Example 2.5.47:* Let

$$G = \left\{ \begin{bmatrix} [a_1, a_2] & [a_7, a_8] & [a_{13}, a_{14}] & [a_{19}, a_{20}] \\ [a_3, a_4] & [a_9, a_{10}] & [a_{15}, a_{16}] & [a_{21}, a_{22}] \\ [a_5, a_6] & [a_{11}, a_{12}] & [a_{17}, a_{18}] & [a_{23}, a_{24}] \end{bmatrix} \right.$$

$a_i \in Z^+ \cup \{0\}; a_i \leq a_{i+1}; 1 \leq i \leq 23; *, (23, 2)\}$ be a 3 × 4 interval matrix groupoid using $Z_I^+$.
    Consider

$$S = \left\{ \begin{bmatrix} [a_1,a_2] & [a_5,a_6] & [0,0] & [0,0] \\ [0,0] & [a_7,a_8] & [a_9,a_{10}] & [a_{11},a_{12}] \\ [a_3,a_4] & [0,0] & [0,0] & [a_{13},a_{14}] \end{bmatrix} \right.$$

$a_i \in Z^+ \cup \{0\}; a_i \leq a_{i+1}; 1 \leq i \leq 13$, *, (23, 2)$\} \subseteq G$; S is a 3×4 interval matrix subgroupoid of G.

Now we define some of the properties enjoyed by these interval matrix groupoids.

**DEFINITION 2.5.9:** *Let G be a interval matrix groupoid (built using $Z_n$ or $Z_I^+$ or $Q_I^+$ or $R_I^+$ or $C_I^+$). A interval matrix subgroupoid V of G is said to be a normal interval matrix subgroupoid or interval matrix normal subgroupoid of G if*
    *(a) aV = Va*
    *(b) (Vx)y = V(xy)*
    *(c) y (xV) = (yx) V*
*for all x, y, a ∈ V.*

    *If a interval matrix groupoid G has no interval matrix normal subgroupoid then we define G to be a simple interval matrix groupoid; we say an interval matrix groupoid G is itself normal if*
    *(a) x G = Gx*
    *(b) G (xy) = (Gx)y*
    *(c) y (xG) = (yx)G*
*for all x, y ∈ G.*

*Example 2.5.48:* Let

$$G = \left\{ \begin{bmatrix} [a,b] \\ [c,d] \end{bmatrix} \middle| \, a, b, c, d \in Z_{12}, *, (4, 8) \right\},$$

G is a column interval matrix groupoid built using $Z_{12}^I$.
Consider

$$P = \left\{ \begin{bmatrix} [a,b] \\ [c,d] \end{bmatrix} \middle| a, b, c, d \in \{0, 4, 8\}, *, (4, 8) \right\}.$$

P is a normal interval matrix subgroupoid of G.

Inview of this we have the following theorem.

**THEOREM 2.5.1:** *Let $G = G_{Z_n^I}^{m \times p} = \{$all $m \times p$ interval matrices with entries from $Z_n^I\}$, \*, $(t, u)$ such that $t + u = n$, $(t, u) = t$; $t/n$ and $n$ even$\}$ be a $m \times p$ matrix interval groupoid with $1 \leq m, p < \infty\}$.*

*Then G has matrix interval normal subgroupoid.*

The proof uses simple number theoretic techniques.

***Example 2.5.49:*** Let $G = \{([a, b], [c, d], [e, f]) \mid a, b, c, d, e, f \in Z_{10}, *, (8, 4)\}$ be a interval matrix groupoid. Take $P = \{([a, b], [c, d], [e, f]) \mid a, b, c, d, e, f \in \{0, 4, 2, 8, 6\}\}, *, (8, 4)\} \subseteq G$ is a matrix interval subgroupoid of G but need not be a interval matrix normal subgroupoid of G as 4/10 and $4 + 8 \neq 10$.

**DEFINITION 2.5.10:** *Let G be a interval matrix groupoid built using $Z_n^I$ or $Z_I^+$ or $Q_I^+$ or $R_I^+$. G is said to be a interval matrix P-groupoid if $(AB) A = A (BA)$ for all $A, B \in G$.*

***Example 2.5.50:*** Let $G = \{([x, y] [a, b], [c, d]) \mid x, y, a, b, c, d, \in Z_n; *, (t, t)\}$ be a $1 \times 3$ row matrix interval groupoid; G is a row interval matrix P-groupoid.
  Take
$$A = ([a_1, b_1], [a_2, b_2] [a_3, b_3])$$
and
$$B = ([x_1 \; y_1] [x_2 \; y_2] [x_3 \; y_3])$$
$A, B \in G$.

Consider
$(A* B)*A = ([a_1, b_1], [a_2, b_2] [a_3, b_3]), *, ([x_1, y_1], [x_2, y_2], [x_3, y_3]) * \{([a_1, b_1], [a_2, b_2] [a_3, b_3])\}$

$$\begin{aligned}
= &\ ([tx_1 + ta_1, ty_1 + tb_1], [tx_2 + ta_2, ty_2 + tb_2], [tx_3 + ta_3, \\
&\ ty_3 + tb_3]) * ([a_1\ b_1], [a_2\ b_2], [a_3, b_3]) \\
= &\ ([t^2x_1 + t^2a_1 + ta_1, t^2y_1, t^2b_1 + tb_1], \\
&\ [t^2x_2 + t^2a_2 + ta_2, t^2y_2 + t^2b_2 + tb_2], \\
&\ [t^2 x_3 + t^2a_3 + ta_3, t^2y_3 + t^2b_3 + tb_3]) \quad\quad (I)
\end{aligned}$$

Consider

$$\begin{aligned}
A * (B * A) = &\ ([a_1, b_1]\ [a_2, b_2], [a_3, b_3]) * \\
&\ ([x_1, y_1], [x_2, y_2], [x_3, y_3]) * \\
&\ ([a_1, b_1], [a_2, b_2], [a_3, b_3]) \\
= &\ ([a_1, b_1]\ [a_2, b_2], [a_3, b_3]) * ([tx_1 + ta_1, ty_1 + tb_1], \\
&\ [tx_2 + ta_2, ty_2 + tb_2], [tx_3 + ta_3, ty_3 + tb_3]) \\
= &\ ([ta_1 + t^2x_1 + t^2a_1, tb_1 + t^2y_1 + t^2b_1], \\
&\ [ta_2 + tx_2 + t^2a_2, tb_2 + t^2y_2 + t^2b_2], \\
&\ [ta_3 + t^2x_3 + t^2a_3, tb_3 + t^2y_3 + t^2b_3]) \quad\quad II
\end{aligned}$$

It is easily verified that $(A*B)*A = A*(B*A)$ for all $A, B \in G$.

Thus G is a interval matrix P-groupoid.

In view of this we have the following theorem.

**THEOREM 2.5.2:** *Let $G = G_{Z_n^I}^{m \times p} = \{\{$all $m \times p$ interval matrices with entries from $Z_n^I\}, *, (t, t), t \in Z_n \setminus \{0\}\}$ be a $m \times p$ interval matrix groupoid built using $Z_n^I$. Clearly G is a $m \times p$ matrix interval P-groupoid.*

**DEFINITION 2.5.11:** *Let G be a $m \times p$ interval matrix groupoid. If $(A * B) * B = A * (B * B)$ for all $A, B \in G$; we define G to be a $m \times p$ interval matrix right alternative groupoid.*

*If $(A * A) * B = A * (A * B)$ for all $A, B \in G$ then we define G to be a $m \times p$ interval matrix left alternative groupoid.*

*G is said to be a $m \times p$ interval matrix alternative groupoid if G is both $m \times p$ interval matrix right alternative groupoid and left alternative groupoid simultaneously.*

*Example 2.5.51:* Let

$$G = \left\{ \begin{bmatrix} [a,b] \\ [c,d] \\ [e,f] \\ [g,h] \end{bmatrix} \middle| a, b, c, d, e, f, g, h \in Z_5, *, (3, 3) \right\}$$

be a $4 \times 1$ row interval matrix groupoid built using $Z_5^I$. G is not an alternative $3 \times 1$ row interval matrix groupoid. Let

$$A = \begin{bmatrix} [a_1,b_1] \\ [a_2,b_2] \\ [a_3,b_3] \\ [a_4,b_4] \end{bmatrix}, B = \begin{bmatrix} [x_1,y_1] \\ [x_2,y_2] \\ [x_3,y_3] \\ [x_4,y_4] \end{bmatrix}$$

be in G.
    Consider

$$(A * B) * B = \begin{bmatrix} [3a_1 + 3x_1, 3b_1 + 3y_1] \\ [3a_2 + 3x_2, 3b_2 + 3y_2] \\ [3a_3 + 3x_3, 3b_3 + 3y_3] \\ [3a_4 + 3x_4, 3b_4 + 3y_4] \end{bmatrix} * \begin{bmatrix} [x_1,y_1] \\ [x_2,y_2] \\ [x_3,y_3] \\ [x_4,y_4] \end{bmatrix}$$

$$= \begin{bmatrix} [9a_1 + 9x_1 + 3x_1, 9b_1 + 9y_1 + 3y_1] \\ [9a_2 + 9x_2 + 3x_2, 9b_2 + 9y_2 + 3y_2] \\ [9a_3 + 9x_3 + 3x_3, 9b_3 + 9y_3 + 3y_3] \\ [9a_4 + 9x_4 + 3x_4, 9b_4 + 9y_4 + 3y_4] \end{bmatrix}$$

$$= \begin{bmatrix} [4a_1 + 2x_1, 4b_1 + 2y_1] \\ [4a_2 + 2x_2, 4b_2 + 2y_2] \\ [4a_3 + 2x_3, 4b_3 + 2y_3] \\ [4a_4 + 2x_4, 4b_4 + 2y_4] \end{bmatrix} \quad (I)$$

Now consider

$$A * (B * B) = \begin{bmatrix} [a_1, b_1] \\ [a_2, b_2] \\ [a_3, b_3] \\ [a_4, b_4] \end{bmatrix} * \begin{bmatrix} [6x_1, 6y_1] \\ [6x_2, 6y_2] \\ [6x_3, 6y_3] \\ [6x_4, 6y_4] \end{bmatrix}$$

$$= \begin{bmatrix} [3a_1 + 18x_1, 3b_1 + 18y_1] \\ [3a_2 + 18x_2, 3b_2 + 18y_2] \\ [3a_3 + 18x_3, 3b_3 + 18y_3] \\ [3a_4 + 18x_4, 3b_4 + 18y_4] \end{bmatrix}$$

$$= \begin{bmatrix} [3a_1 + 3x_1, 3b_1 + 3y_1] \\ [3a_2 + 3x_2, 3b_2 + 3y_2] \\ [3a_3 + 3x_3, 3b_3 + 3y_3] \\ [3a_4 + 3x_4, 3b_4 + 3y_4] \end{bmatrix} \qquad \text{II}$$

I and II are not equal. Thus G is not a $4 \times 1$ interval matrix alternative groupoid.

In view of this we have the following theorem which can be easily proved using simple number theoretic techniques.

**THEOREM 2.5.3:** *Let $G = G_{Z_p^I}^{m \times p}$ = {all collection of $m \times n$ interval matrices with entries from $Z_p^I$, \*, (t, t), t < p} be a m×n interval matrix groupoid built using $Z_p^I$; G is not a $m \times n$ interval matrix alternative groupoid if p is a prime.*

**THEOREM 2.5.4:** *Let $G = G_{Z_n^I}^{m \times p}$ = {all $m \times p$ interval matrices with entries from $Z_n^I$ n-not a prime, \*, (t, t), ; $t \in Z_n \setminus \{0\}$} be a $m \times p$ matrix interval groupoid. G is an alternative m×p interval matrix groupoid if and only if $t^2 \equiv t \pmod{n}$.*

***Example 2.5.52:*** Let $G = G_{Z_6^I}^{7\times 3}$ = {all 7 × 3 interval matrices with entries from $Z_6^I$, *, (3, 3)} be a 7 × 3 interval matrix alternative groupoid as $3^2 \equiv 3 \pmod 6$.

**THEOREM 2.5.5:** *Let $G = G_{Z_n^I}^{m\times p}$ = {all m×p interval matrix with entries from $Z_n^I$, * (t, 0), $t \in Z_n \setminus \{0\}$} be a m × p matrix interval groupoid. G is a P-groupoid and alternative groupoid if and only if $t^2 \equiv t \pmod n$.*

***Example 2.5.53:*** Let G = {3 × 8 interval matrices with entries from $Z_{12}^I$, *, (4, 0)} be a 3 × 8 interval matrix groupoid. G is clearly a interval matrix P-groupoid and interval matrix alternative groupoid.

## 2.6 Smarandache Interval Groupoid

As in case of usual groupoids we can in case of matrix interval groupoids also define the concept of Smarandache groupoids.
We will illustrate this concept.

***Example 2.6.1:*** Let $G = \{G_{Z_{10}^I}^{3\times 8}, *, (1, 5)\}$ be a 3 × 8 matrix interval groupoid. Take P = {3 × 8 interval matrices from the set {0, 5}, *, (1, 5)} $\subseteq$ G is easily verified to be a 3 × 8 matrix interval semigroup. Thus G is a Smarandache 3 × 8 matrix interval groupoid or 3 × 8 interval matrix Smarandache groupoid.

**THEOREM 2.6.1:** *Let $G = \{G_{Z_n^I}^{m\times p}, *, (t, u)$ such that $t, u \in Z_n \setminus \{0\}$; $(n > 5)$; $t + u \equiv 1, \pmod n$ and $(t, u) = 1\}$ be a m × p interval matrix groupoid using $Z_n^I$. Then G is a Smarandache m × p interval matrix groupoid.*

(We give only hint of the proof).

*Proof:* For every $m \in Z_n$ is such that $m * m = mt + mu = m$. Thus every singleton is a semigroup. Hence the claim.

We have yet another new result about Smarandache interval matrix groupoid using $Z_n^I$.

**THEOREM 2.6.2:** *Let $G = \{ G_{Z_n^I}^{m \times p} ; *, (t, u)$ where $t + u = 1 \pmod{n} \}$ be a $m \times p$ interval matrix groupoid. $G$ is a $m \times p$ matrix interval Smarandache P-groupoid if and only if $t^2 = t \pmod{n}$ and $u^2 = u \pmod{n}$.*

*Proof:* It is already proved $G$ is a $m \times p$ matrix interval Smarandache groupoid if $t + u \equiv 1 \pmod{n}$. Now to show $G$ is a $m \times p$ matrix interval Smarandache P-groupoid if and only if $t^2 = t \pmod{n}$ and $u^2 = u \pmod{n}$.

Let $A_{m \times p} = ([a_{ij}, a_{i+1j+1}])$ and $B_{m \times p} = ([b_{ij}, b_{i+1j+1}])$ with $1 \leq i \leq m - 1$ and $1 \leq j \leq p - 1$. To show $G$ is a $m \times p$ matrix interval S-P-groupoid (Smarandache-P-groupoid) we have to show $(A * B) * A = A * (B * A)$ for all $A, B \in G$.
Consider
$(A * B) * A = [([a_{ij}, a_{i+ij+1}]) * ([b_{ij}, b_{i+ij+1}]) * ([a_{ij}, a_{i+ij+1}])$
$= [(ta_{ij} + ub_{ij}) \pmod{n}, (ta_{i+1j+1} + ub_{i+1j+1}) \pmod{}] * ([a_{ij}, a_{i+ij+1}])$
$= [(t^2 a_{ij} + tub_{ij} + ua_{ij}) \pmod{n}, (t^2 a_{i+1j+1} + tub_{i+1j+1} + ta_{i+ij+1})]$
                                                                                                 I

Consider
$A * (B * A) = ([a_{ij}, a_{i+ij+1}]) * ([(tb_{ij} + ua_{ij}) \bmod n,$
                  $tb_{i+1j+1} + ua_{i+1j+1}) \bmod n])$
             $= [(ta_{ij} + utb_{ij} + u^2 a_{ij}) \bmod n,$
                  $(ta_{i+1j+1} + utb_{i+1j+1} + u^2 a_{i+1j+1}) \bmod n]$     II

Using the fact $t^2 \equiv t \pmod{n}$ and $u^2 \equiv u \pmod{n}$ we see I and II are equal. Hence $G$ is a $m \times p$ interval matrix S-P-groupoid.

We will illustrate this situation by some simple examples.

*Example 2.6.2:* Let $G = \{ G_{Z_{12}^I}^{3 \times 4}, *, (4, 9) \}$ be a $3 \times 4$ interval matrix groupoid. It is easily verified that $G$ is a $3 \times 4$ interval matrix S-P-groupoid.

**THEOREM 2.6.3:** *Let $G = \{ G^{m \times p}_{Z_n^I}, *, (t, u); t, u \in Z_n \setminus \{0\}$ with $t + u \equiv 1 \pmod{n}\}$ be a $m \times p$ matrix interval groupoid. $G$ is a $m \times p$ matrix interval Smarandache alternative groupoid ($m \times p$ matrix interval S-alternative groupoid) if and only if $t^2 \equiv t \pmod{n}$ and $u^2 \equiv u \pmod{n}$.*

The proof is simple and can be obtained by using number theoretic techniques. For definition of different types of Smarandache groupoid using interval matrices please refer [20-21]. Thus using those definition in [20] we will prove results related with them.

***Example 2.6.3:*** Let $G = \{ G^{3 \times 1}_{Z_{12}^I}, *, (3, 9)\}$ be a $3 \times 1$ column interval matrix groupoid built using $Z_{12}^I$. $G$ is a column interval matrix Smarandache Moufang groupoid.

It can be easily verified for all $A, B, C \in G$; $(A * B) * (C * A) = (A * (B*C)) * A$ is true.

Likewise we give an example of a matrix interval Smarandache Bol groupoid.

***Example 2.6.4:*** Let $G = \{ G^{5 \times 7}_{Z_4^I}, *, (2, 3)\}$ be a $5 \times 7$ matrix interval groupoid. It is easily verified $G$ is a Smarandache Bol groupoid. For $A, B, C$ in $P$ we have $((A*B)*C) * B = A* [(B*C)*B]$ where $P = \{$All $5 \times 7$ interval matrices with entries from the subset $\{0, 2\} \subseteq Z_4, *, (2, 3)\} \subseteq G$. However the identity is not true for all elements in $G$.

***Example 2.6.5:*** Let $G = \{ G^{12 \times 5}_{Z_{10}^I}, * (5, 6)\}$ be a interval $12 \times 5$ matrix groupoid built using the intervals of $Z_{10}^I$.

It is easily verified for every $12 \times 5$ interval matrices $A, B, C$ in $G$; $(A * B) * (C * A) = [A * (B * C)] * A$. Thus $G$ is a interval $12 \times 5$ matrix Smarandache strong Moufang groupoid.

**THEOREM 2.6.4:** *Let $G = \{ G_{Z_n^I}^{q \times p}, *, (m, m); m \in Z_n$ with $m + m \equiv 1 \pmod{n}$ and $m^2 \equiv m \pmod{n}\}$ be a $q \times p$ matrix interval groupoid};*

1. *$G$ is a $q \times p$ matrix interval Smarandache strong P-groupoid.*
2. *$G$ is $q \times p$ matrix interval Smarandache idempotent groupoid.*
3. *$G$ is a $q \times p$ matrix interval Smarandache strong Bol groupoid.*
4. *$G$ is a $q \times p$ matrix interval Smarandache strong Maufang groupoid.*
5. *$G$ is a $q \times p$ matrix interval Smarandache strong alternative groupoid.*

*Note:* q and p are positive finite integers $1 < p, q < \infty$; p = q can also occur.

The proof is left as an exercise for the reader.

**THEOREM 2.6.5:** *Let $G = \{ G_{Z_{2m}^I}^{p \times q}, *, (2, 0), 2 \in Z_{2m}\}$ be a $p \times q$ interval matrix groupoid built using intervals from $Z_{2m}^I$. $G$ is a $p \times q$ interval matrix S-groupoid.*

*Proof:* Follows from the fact when intervals $\{[0, 0], [m, m], [0, m]]\} \subseteq Z_{2m}^I$ are taken as entries of the p × q interval matrix the collection is a p × q interval matrix semigroup under the operation *. Hence the claim.

**THEOREM 2.6.6:** *Let $p \times q$ matrix interval groupoid $G = \{G_{Z_n^I}^{p \times q}, *, (m, 0)\}$ built using $Z_n^I$, n = 2m is a Smarandache $p \times q$ matrix interval groupoid.*

Proof is left as an exercise for the reader.

**COROLLARY 2.6.1:** $G = \{ G_n^{t \times s}, *, (p,0)$ *where p is a prime and p / n}, the $t \times s$ interval matrix groupoid is a $t \times s$ interval matrix S-groupoid.*

Several other results in this direction can be derived. Now we proceed onto define classes of groupoids built using intervals.

## 2.7 Classes of Groupoids Built Using Intervals

In this section we indicate the classes of groupoids built using intervals from $Z_n^I$ or $(Z^+ \cup \{0\})_I =$

$$Z_I^+ = \{[a, b] | a \leq b, a, b \in Z^+ \cup \{0\}\} \text{ or}$$
$$Q_I^+ = \{[a, b] | a \leq b, a, b \in Q^+ \cup \{0\}\}$$
$$R_I^+ = \{[a, b] | a \leq b, a, b \in R^+ \cup \{0\}\} \text{ and}$$
$$C_I^+ = \{[a, b] | a \leq b, a, b \in C^+ \cup \{0\}\}$$
$$C^+ = \{x + iy \, / \, x, y \text{ are positive real}\}\}.$$

$C(Z_n^I, p \times q) = \{\{G_{Z_I^n}^{p \times q} = \{\text{all } p \times q \text{ matrices with entries from } Z_I^n\}, *, (t, u) \, / \, t, u \in Z_n\} \, / \, t$ and u vary over $Z_n$ with all possibilities $(t, u) = 1$ or $(t, u) = d$ or $(t, u) = t$ or $(t, u) = u$ or $(t, 0)$ or $(0, t)$ or $(t, t)\}$. Thus $C(Z_n^I, p \times q)$ is defined as the class of $p \times q$ matrix interval groupoids built using $Z_I^n$.

Clearly $C(Z_n^I, p \times q)$ has only finite number of such $p \times q$ matrix interval groupoids for a fixed p and q (p and q are positive integers). If $p = 1$ we call $1 \times q$ interval matrix groupoid as interval row matrix groupoid built using $Z_n^I$ and is denoted by $(Z_I^n, 1 \times q)$ for a fixed q. Similarly $(Z_I^n, p \times 1)$ denotes the class of interval column matrix groupoid for a fixed p.

Likewise $(Z_I^n, m \times m)$ denotes the class of interval $m \times m$ square matrix groupoid built using $Z_I^n$. In $Z_n$, n is finite but n can be prime. Thus we have 4 types of class of matrix interval groupoids. Further the collection of interval $1 \times 9$ row matrix groupoids for varying q is infinite so we have an infinite number of classes $(Z_I^n, 1 \times q)$ where $q \in Z^+$. This is true of any $p \times q$ or $p \times 1$ or $m \times n$ matrix interval groupoids p, q, m $\in Z^+$.

We will describe some of the properties enjoyed by them.

**THEOREM 2.7.1:** *The groupoids in this class ($Z_I^n$, $t \times u$, $(0, p)$) are Smarandache matrix interval groupoids (n-even).*

**COROLLARY 2.7.1:** *No matrix interval groupoid in this class of groupoids $C(Z_I^n, t \times u, (0, p))$ is a n-even matrix internal S idempotent groupoid.*

Both the results hold good for the class of groupoids $C(Z_I^n, t \times u, (p, 0))$ n-even.

**THEOREM 2.7.2:** *Let $C(Z_I^n, t \times u (m, m))$ where $m + m = 1$ (mod n) and $m^2 = m$ (mod n) be the class of matrix interval groupoids for the specific m's then*
1. *$C(Z_I^n, t \times u, (m, m))$ are matrix interval S-idempotent groupoids.*
2. *$C(Z_I^n, t \times u, (m, m))$ for those special m are matrix interval Smarandache strong P-groupoids.*
3. *$C(Z_I^n, t \times u, (m, m))$ for that specific m's are matrix interval Smarandache strong Bol-groupoids.*
4. *$C(Z_I^n, t \times u, (m, m))$ are matrix interval Smarandache strong Moufang groupoids.*
5. *$C(Z_I^n, t \times u, (m, m))$ are Smarandache strong matrix interval strong alternative groupoid.*

**THEOREM 2.7.3:** *$C(Z_I^n, m \times s, (1, p)) / p$ is a prime and $p /n) \subseteq C(Z_I^n, m \times s, (1, p))$ contained in the class of interval m×s matrix groupoids are interval m×s matrix S-groupoids.*

**THEOREM 2.7.4:** *The subclass of groupoids $C(Z_I^n, m \times s, (t, u) / t + u \equiv 1 (mod\ n)) \subseteq C(Z_I^n, m \times s (t, u))$ are $m \times s$ matrix interval S-idempotent groupoids.*

**THEOREM 2.7.5:** *The subclass of groupoids $C(Z_I^n, m \times s, (t, u) / t + u \equiv 1 (mod\ n)$ and $t^2 = t (mod\ n)$ and $u^2 = u (mod\ n)) \subseteq C(Z_I^n, m \times s, (t, u))$ are matrix interval Smarandache P-groupoids.*

**THEOREM 2.7.6:** *The subclass of groupoids $C((Z_I^n, m{\times}s, (t, u))$ / $t+u \equiv 1 \pmod n$, $u^2 = u \pmod n$, $t^2 = t \pmod n) \subseteq C(Z_I^n, m{\times}s, (t, u))$ are matrix interval Smarandache alternative groupoids.*

**THEOREM 2.7.7:** *The subclass of groupoids $C(Z_I^n, m{\times}s, (t, u))$ / $t + u \equiv 1 \pmod n$, $u^2 = u \pmod n$, $t^2 = t \pmod n) \subseteq C(Z_I^n, m{\times}s, (t, u))$ are $m{\times}s$ matrix interval Smarandache strong Bol groupoids.*

**THEOREM 2.7.8:** *The subclass of groupoids $C(Z_I^n, m \times s, (t,u))$ / $t + u \equiv 1 \pmod n$, $u^2 = u \pmod n$, $t^2 = t \pmod n) \subseteq C(Z_I^n, m \times s, (t, u))$ are $m \times s$ matrix interval Smarandache strong Maufang groupoid.*

The above theorems can be easily derived using the definition and simple number theoretic techniques.

Now we proceed onto define classes of groupoids using $Z^+ \cup \{0\}$ etc.

$C(Z_I^n, m \times n, (p, q)) = \{$(collection of all $m \times n$ interval matrices with entries from $Z_I^+$ together with $*$ binary operation on $G_{Z_I^+}^{m \times n}$ such that for $A, B \in G_{Z_I^+}^{m \times n}$, $A * B = pA + qB$, where $p, q \in Z^+$ and $p, q$ vary over $Z^+$ to give the class of $m{\times}n$ matrix interval groupoids using $Z_I^+\}$. Clearly this is an infinite collection. Likewise $C(Q_I^+, m \times n, (p, q))$ $p, q \in Q^+$ is again a class of $m \times n$ matrix interval groupoid of infinite cardinality.
Similarly $C(R_I^+, m \times n, (p, q))$ $p, q \in R^+$ and $C(C_I^+, m \times n, p, q \in R^+)$ are classes of $m \times n$ matrix interval groupoids of infinite cardinality. Here $m$ and $n$ are fixed for that class as $m$ and $n$ varies in $Z^+$ we get different classes of matrix interval groupoids using $Q_I^+$ or $Z_I^+$ or $R_I^+$ or $C_I^+$.

**Chapter Three**

# ON SOME NEW CLASSES OF NEUTROSOPHIC GROUPOIDS

In this chapter we for the first time introduce the new concept of neutrosophic groupoids. The algebraic structures like neutrosophic groups, neutrosophic semigroups, neutrosophic rings, neutrosophic fields and neutrosophic vector spaces have been introduced by the authors. It is pertinent to mention that all those algebraic structures were associative structures. We assume by a neutrosophic set a non empty set S with the element I, the indeterminate. So if S is a neutrosophic set I ∈ S. This chapter has seven sections. Section one introduces the concept of neutrosophic groupoids. Section two introduces new classes of neutrosophic groupoids using $Z_n$. Neutrosophic polynomial groupoids are introduced in section three. Section four introduces the notion of neutrosophic matrix groupoids. Neutrosophic interval groupoids are introduced in section five.

Neutrosophic interval matrix groupoids are introduced in section six and neutrosophic interval polynomial groupoids in section seven are described.

## 3.1 Neutrosophic Groupoids

We now proceed on to define a neutrosophic groupoid and illustrate it by some examples.

**DEFINITION 3.1.1:** *Let S be a neutrosophic set which is non empty. Let \* be a closed binary operation on S such that a\*(b\*c) ≠ (a \* b) \* c for some a, b, c ∈ S. We call (S, \*) to be a neutrosophic groupoid.*

It is interesting to note that all neutrosophic semigroups are neutrosophic groupoids but a neutrosophic groupoid in general is not a neutrosophic semigroup. Thus the class of all neutrosophic semigroups is contained in the class of neutrosophic groupoids.

We will illustrate this new structure by some examples.

*Example 3.1.1:* Let $\{Z \cup I\} = N(Z) = \{a + bI \mid a, b \in Z\}$. Define the binary operation subtraction '–' on N(Z). {N(Z), –} is a neutrosophic groupoid.
    It is easily verified '–' operation on N(Z) is non associative. Thus (N(Z), –) is not a neutrosophic semigroup as '–' operation on N(Z) is not associative.

*Example 3.1.2:* Let $N(Z_3) = \{a + bI / a, b \in Z_3\} = \{0, I, 2I, 1, 2, 1 + I, 2 + I, 2I + 1, 2I + 2\}$. Define a binary operation \* on $N(Z_3)$ as follows. For $x, y \in N(Z_3)$ define $x * y = x + 2y \pmod 3$ '+' usual addition.

It is easily verified $(N(Z_3), *)$ is a neutrosophic groupoid.

*Example 3.1.3:* Let $G = \{M_{2 \times 3} = (m_{ij}) / m_{ij} \in N(R), 1 \le i \le 2$ and $1 \le j \le 3\}$. Define a operation \* on G as follows. If $A = (a_{ij})$ and B = ($b_{ij}$) are in G. A * B = ($2a_{ij} + 3b_{ij}$). It is easily verified (G, *) is a neutrosophic groupoid as '*' is a closed non associative binary operation on G.

Now having seen some examples of neutrosophic groupoids we now proceed onto define the cardinality of these algebraic structures.

**DEFINITION 3.1.2:** *Let (G, *) be a neutrosophic groupoid. If the number of distinct elements in G is finite, then we say G has finite cardinality or finite order and we denote it by |G| < ∞. If G has infinite number of elements then we say G has infinite cardinality and is denoted by |G| = ∞ or o(G) = ∞.*

We see neutrosophic groupoids given in examples 3.1.1 and 3.1.3 are of infinite order and the neutrosophic groupoid given in example 3.1.2 is of finite order.

**DEFINITION 3.1.3:** *Let (S, *) be a neutrosophic groupoid. Suppose H ⊆ S be a proper subset of S and if (H, *) is itself a neutrosophic groupoid then we call (H, *) to be a neutrosophic subgroupoid of (S, *).*

We will illustrate this situation by some examples.

*Example 3.1.4:* Let S = {N(Z), *} be a neutrosophic groupoid where for any x, y ∈ N(Z); x * y = 5x + 3y. Take H = {a + bI | a, b ∈ $Z^+$}; (H, *) is a neutrosophic subgroupoid of S.

*Example 3.1.5:* Let S = {N($Z_5$), * where a * b = 2a + b} be a neutrosophic groupoid.
    Take H = {$Z_5$I}, (H, *) is a neutrosophic subgroupoid of S.
In example 3.1.4 if we take $H_1$ = {$Z^+$} ⊆ N(Z); we see ($H_1$, *) is a groupoid but ($H_1$, *) is not a neutrosophic groupoid. Likewise in example 3.1.5. $H_2$ = ($Z_5$, *) is a subgroupoid of (N($Z_5$), *) but $H_2$ is also not a neutrosophic groupoid.

Thus in view of this we give a definition.

**DEFINITION 3.1.4:** *Let (S, *) be a neutrosophic groupoid. If P = H such that (P, *) is a groupoid but is not a neutrosophic groupoid then we call (P, *) to be a pseudo neutrosophic subgroupoid of (S, *).*

We will first give some examples of this definition.

***Example 3.1.6:*** Let $\{N(Z_{12}), *\} = G$ be a neutrosophic groupoid where '*' is such that $x * y = 3x + 7y$ (where '+' is addition modulo 12); for all $x, y \in N(Z_{12})$. Take $H = \{Z_{12}, *\}$ we see H is a groupoid and is not a neutrosophic groupoid. So H is a pseudo neutrosophic subgroupoid of G.

***Example 3.1.7:*** Let $V = \{N(Z), *$ where * is such that $a * b = 8a + 9b$ for $a, b \in N(Z)\}$. V is a neutrosophic groupoid.
Take $W = \{Z, *\} \subseteq V = \{N(Z), *\}$; W is a pseudo neutrosophic subgroupoid of V.

**DEFINITION 3.1.5:** *Let (V, *) be a neutrosophic groupoid. If V has no neutrosophic subgroupoid; then we call V to be a neutrosophic simple groupoid or simple neutrosophic groupoid.*

***Example 3.1.8:*** Let $V = \{Z_7I, *; a * b = 2a + 3b(\mod 7)\}$ be a neutrosophic groupoid. Clearly V is a simple neutrosophic groupoid.

***Example 3.1.9:*** Let $V = \{Z_3I, *;$ where $a * b = 2a + b(\mod 3)\}$.

| * | 0 | I | 2I |
|---|---|---|----|
| 0 | 0 | I | 2I |
| I | 2I | 0 | I |
| 2I | I | 2I | 0 |

be a neutrosophic groupoid. V is also given the table representation. V is a simple neutrosophic groupoid.
Thus we have a class of simple neutrosophic groupoids.

Now we proceed onto define pseudo simple neutrosophic groupoids.

**DEFINITION 3.1.6:** *Let (V, \*) be a neutrosophic groupoid. If V has no nontrivial pseudo neutrosophic subgroupoids then we define (V, \*) to be a pseudo simple neutrosophic groupoid.*

We will illustrate this by some examples.

*Example 3.1.10:* Let V = {QI, \*} be a neutrosophic groupoid. Clearly V is a pseudo simple neutrosophic groupoid as V has no nontrivial pseudo neutrosophic subgroupoids.

*Note:* If $\{0\} \subseteq V$ are call {0} to be a trivial pseudo neutrosophic subgroupoid.

*Example 3.1.11:* Let V = {$Z_{11}I$, \*, where \* is defined for any a, b $\in Z_{11}I$ as a \* b = 3a + 5b(mod 11)}, V is a neutrosophic groupoid which is a simple pseudo neutrosophic groupoid as V has no nontrivial pseudo neutrosophic subgroupoids.

**DEFINITION 3.1.7:** *Let (V, \*) be a neutrosophic groupoid if V is both a simple neutrosophic groupoid as well as pseudo simple neutrosophic groupoid then we call (V, \*) to be a doubly simple neutrosophic groupoid.*

We will illustrate this by some simple examples.

*Example 3.1.12:* Let V = {$Z_7I$, \*} where for a, b $\in Z_7I$; a \* b = a + 2b (mod 7). V is a doubly simple neutrosophic groupoid.

*Example 3.1.13:* Let V = {ZI, \* where a \* b = 8a + 11b for a, b $\in$ ZI} be a pseudo neutrosophic simple groupoid and is not a doubly simple neutrosophic groupoid.
    Thus we see a simple neutrosophic groupoid need not in general be a doubly simple neutrosophic groupoid. In fact by appropriately defining '\*' on $Z_pI$; p a prime we can get a class of doubly simple neutrosophic groupoids.
    Also we have got a class of groupoids which are not doubly simple neutrosophic groupoids. For take {$N(Z_p)$, \*} = V (p a prime) forms a class of neutrosophic groupoids which are not

doubly simple neutrosophic groupoids as $\{Z_pI, *\} \subseteq \{N(Z_p), *\}$ is a nontrivial neutrosophic subgroupoid of V and $\{Z_p, *\} \subseteq \{N(Z_p) *\}$ is a nontrivial pseudo neutrosophic subgroupoid of V. Hence the claim.

Having defined substructures of a neutrosophic groupoid we now proceed onto define the notion of neutrosophic groupoid homomorphisms.

**DEFINITION 3.1.8:** *Let (G, \*) and (H, o) be two neutrosophic groupoids. A map $\eta : G \to H$ is said to be a neutrosophic groupoid homomorphism if the following conditions hold good.*

*(i) $\eta(I) = I$*
*(ii) $\eta(a * b) = \eta(a) \, o \, \eta(b)$ for all a, b $\in$ G.*

The interested reader is expected to give examples of them. It is important to mention that I should be mapped only onto I and no other element can be substituted as I stands for the indeterminacy.

Another factor which is to be noted in case of these algebraic structures is that neutrosophic groupoids may or may not have identity. An element e in a neutrosophic groupoid such that e \* x = x \* e = x for all x $\in$ G is called the identity. Even if 0 $\in$ G it does not in general imply 0 \* x = x \* 0 = x for all x $\in$ G. Thus with these special properties it is not always possible to define kernel of the homomorphism.

However we can define as in case of other algebraic structures isomorphism. A homomorphism which is one to one and onto is defined as an isomorphism.

*Example 3.1.14:* Let (G, *) = $\{Z_7I, *$ is such that a * b = a + 2b$\}$ and (H, o) = $\{N(Z_7)$ 'o' is such that a o b = a + 2b$\}$ be any two neutrosophic groupoids. Define $\eta : (G, *) \to (H, o)$ by $\eta(a * b) = a \, o \, b$ for all a, b $\in$ G. It is easily verified $\eta$ is a homomorphism of G to H.

The interested reader can construct any number of neutrosophic groupoid homomorphism.

Having defined the notion of neutrosophic groupoid homomorphism we can now proceed onto define the another substructure.

**DEFINITION 3.1.9:** *Let (G, \*) be a neutrosophic groupoid. Let H $\subseteq$ G be a proper subset of G. We say H is a neutrosophic ideal of G if the following conditions are satisfied;*

   i. *(H, \*) is a subgroupoid*
   ii. *For all g $\in$ G and h $\in$ H, g \* h and h \* g are in H.*

We will illustrate this situation by some examples.

*Example 3.1.15:* Let G = {0, I, 2I, …, 15I} = $Z_{16}$I. Define '\*' on G by a \* b = 2a + 3b (mod 16I) for a, b $\in$ G. (G, \*) is a neutrosophic groupoid. Take H = {4I, 8I, 0, 12I}, {H, \*} is a neutrosophic subgroupoid of (G, \*).

It is easily verified that H is not a neutrosophic ideal of the neutrosophic groupoid G. It is easily verified G has no neutrosophic ideals.

*Example 3.1.16:* Let G = $Z_8$I = {0, I, 2I, 3I, 4I, 5I, 6I, 7I}, define \* on G by a \* b = 2a + 4b (mod 8I) for a, b $\in$ G. (G, \*) is a neutrosophic groupoid. Take {0, 2I, 4I, 6I} = H $\subseteq$ G (H, \*) is a neutrosophic subgroupoid of (G, \*).

It is easily verified that (H, \*) is a neutrosophic ideal of (G, \*).
   Thus we see some neutrosophic groupoids contain neutrosophic ideals and others do not contain any neutrosophic ideals.

**DEFINITION 3.1.10:** *Let (G, \*) be a neutrosophic groupoid. Suppose G has no proper neutrosophic ideals then we call (G, \*) to be ideally simple neutrosophic groupoid.*

Can we have non trivial ideally simple neutrosophic groupoids. We explicitly show a class of such neutrosophic groupoids.

**THEOREM 3.1.1:** *Let $G = Z_{2^n}I = \{0, I, 2I, \ldots, 2^nI - I\}$ (modulo integers $2^nI$, $n \geq 2$). Define a operation $*$ on $G$ by $a * b = ma + nb \pmod{2^nI}$ for all $a, b \in G$ and $m + n = $ a prime in $G$ and $m, n \in Z_{2^n}$, $(G, *)$ is a ideally simple neutrosophic groupoid.*

*Proof:* Given G is a neutrosophic groupoid of a special order $2^n$; $n \geq 2$. Now the operation on G is also such that $a * b = ma + nb \pmod{2^n}$ for all $a, b \in G$ with $m + n = $ a prime in G and $m, n \in G$.

Let H be any proper neutrosophic subgroupoid of G. Any element in H is of the form $\{2^sI \mid s = 1, 2, \ldots, 2^n\}$. Let $m, n \in G$ such that $m + n = $ prime.

Now for any odd prime $a \in G$ and $b \in H$ we must have $a * b$ to belong to H.

Thus $ma + nb \notin H$ as $ma + mb \neq 2^sI$ for any $s$, $s = 1, 2, \ldots, 2^n$. For if $m + n = $ prime at least one of m or n must be of the form $p.q$ where q or p is an odd neutrosophic prime such that ultimately $m + n = $ a neutrosophic prime, $\omega$ and $\omega < 2^n$.

This is not possible. Hence $G = Z_{2^n}I$ is ideally simple neutrosophic groupoid.

*Note:* When we say n is a neutrosophic prime; $n = pI$ where p is a prime. Likewise when we say pure neutrosophic modulo integers we mean $Z_nI = \{0, I, 2I, \ldots, (nI - I) = (n - 1)I\}$. Thus $Z_3I = \{0, I, 2I\}$; $2I + I = 3I = 0 \pmod{3I}$.

$N(Z_n) = \{a + bI / a, b \in Z_n\}$. Clearly $N(Z_n)$ is neutrosophic modulo integers (but $N(Z_n) \neq Z_nI$. $Z_nI \subseteq N(Z_n)$) and is not a pure neutrosophic modulo integers. Further $Z_nI \not\subseteq Z_n$ both are different as one is pure neutrosophic where as other is pure modulo integers; only $|Z_nI| = |Z_n| = n$ and $Z_nI \cap Z_n = \{0\}$.

Now we proceed onto define new classes of neutrosophic groupoids using $Z_n$, Z, Q and R.

## 3.2 New Classes of Neutrosophic Groupoids Using $Z_n$

In this section we present some new classes of groupoids constructed using modulo integers.

**DEFINITION 3.2.1:** *Let $Z_nI = \{0, I, 2I, …, (n – 1)I\}$; $n \geq 3$. For a, $b \in Z_nI \setminus \{0\}$ define a binary operation * on $Z_nI$ as follows $a * b = ta + ub \pmod{nI}$ t, u are two distinct elements in $Z_nI \setminus \{0\}$ or $Z_n \setminus \{0\}$. + is the usual addition of two neutrosophic integers modulo nI. $\{Z_nI, (t, u), *\}$ is a neutrosophic groupoid in short denoted by $Z_nI (t, u)$.*

It is interesting to note for a given $Z_nI$ we can construct several neutrosophic groupoids by varying t and u. ($t \neq u$; $t, u \in Z_nI$). Clearly every neutrosophic groupoid $Z_nI (t, u)$ is of order n we denote the class of neutrosophic groupoids built in this way using $Z_nI$ by $Z(n)I$.

*Example 3.2.1:* Let $Z_3I = \{0, I, 2I\}$. The neutrosophic groupoid $\{Z_3I, (I, 2I), *\} = Z_3I (I, 2I)$ is given by the following table.

| *  | 0  | I  | 2I |
|----|----|----|----|
| 0  | 0  | 2I | I  |
| I  | I  | 0  | 2I |
| 2I | 2I | I  | 0  |

Clearly this neutrosophic groupoid is non associative and non commutative and its order is 3.

Likewise $\{Z_3I, (2I, I), *) = Z_3I (2I, I)$ is also a neutrosophic groupoid of order three given by the table.

| *  | 0  | I  | 2I |
|----|----|----|----|
| 0  | 0  | I  | 2I |
| I  | 2I | 0  | I  |
| 2I | I  | 2I | 0  |

We see from the tables that $Z_3I$ (I, 2I) is not the same as $Z_3I$ (2I, I). However it is important to note that we can have only two neutrosophic groupoids built over $Z_3I$.

Now consider the mixed neutrosophic modulo integers $N(Z_3) = \{a + bI \,/\, a, b \in Z_3\} = \{0, 1, 2, I, 2I, 1 + I, 1 + 2I, 2 + 2I, 2 + I\}$. We will study how many neutrosophic groupoids can be constructed using $N(Z_3)$.

*Example 3.2.2:* Let $N(Z_3)$ be the mixed neutrosophic integers modulo $Z_3$.

Define * on $N(Z_3)$ as follows:
$$a * b = 1.a + (1 + I)b$$
for all $a, b \in N(Z_{23})$.

The table for this neutrosophic groupoid is given in the following.

| *   | 0   | 1   | 2    | I    | 2I   | 1+I  | 2+I  | 1+2I | 2+2I |
|-----|-----|-----|------|------|------|------|------|------|------|
| 0   | 0   | 1+I | 2+2I | 2I   | I    | 1    | 2+I  | 1+2I | 2    |
| 1   | 1   | 2+I | 2I   | 1+2I | 1+I  | 2    | I    | 2+2I | 0    |
| 2   | 2   | I   | 1+2I | 2+2I | 2+I  | 0    | 1+I  | 2I   | 1    |
| I   | I   | 1+2I | 2   | 0    | 2I   | 1+I  | 2+2I | 1    | 2+I  |
| 2I  | 2I  | 1   | 2+I  | I    | 0    | 2I+1 | 2    | 1+I  | 2+2I |
| 1+I | 1+I | 2+2I | 0   | 1    | 1+2I | 2+I  | 2I   | 2    | I    |
| 2+I | 2+I | 2I  | 1    | 2    | 2+2I | I    | 1+2I | 0    | 1+I  |
| 1+2I | 1+2I | 2  | I    | 1+I  | 1    | 2+2I | 0    | 2+I  | 2I   |
| 2+2I | 2+2I | 0  | 1+I  | 2+I  | 2    | 2I   | 1    | I    | 1+2I |

It is easily verified $P = (N(Z_3), (1, 1 + I), *)$ is a non associative structure. Thus P is a neutrosophic groupoid we see using $N(Z_3)$ we can construct 56 number of neutrosophic groupoids. Thus we have a very large number of groupoids using $N(Z_3)$.

*Example 3.2.3:* Let $N(Z_4) = \{0, 1, I, 2I, 3I, 2, 3, 1 + I, 2 + I, 3 + I, 1 + 2I, 2 + 2I, 3 + 3I, 1 + 3I, 2 + 3I, 2 + 3I\}$ be the mixed neutrosophic modulo integers 4.

We see o $(N(Z_4)) = 4^2 = 16$ that is by order of $N(Z_4)$ we mean the number of elements in $N(Z_4)$. We can construct 2 (15 $C_2$), that is 210 number of mixed neutrosophic groupoids using $N(Z_4)$.

In view of this we have the following results. We call a neutrosophic groupoid to be pure if it does not contain any element x where $x \neq 0$ and x has no I present with it. $(Z_nI, (x, y), *)$, $x, y \in Z_nI$ is a pure neutrosophic groupoid. We call $(N(Z_n), (x, y), *)$; $x, y \in N(Z_n)$ to be a neutrosophic groupoid which is not pure.

The following theorems can be easily proved by any reader.

**THEOREM 3.2.1:** *Let $(Z_nI, (x, y), *)$ be any pure neutrosophic groupoid of order n. There exists exactly $2C_{n-1}C_2$ number of groupoids of order n built using $Z_nI$.*

**THEOREM 3.2.2:** *Let $(N(Z_n), (x, y), *)$ be any neutrosophic groupoid (which is not pure) built using $N(Z)$. There exists $2_{n^2-1}C_2$) number of such groupoids of order $n^2$.*

Let $P = \{0, 2I\}$ be a left ideal of $Z_4I(u, v)$ where $u = 3$ and $v = 2$. Clearly $P = \{0, 2I\}$ is a right ideal of $Z_4I(v, u)$.

In view of this we have the following theorem the proof of which is straight forward.

**THEOREM 3.2.3:** *P is a left neutrosophic ideal of $Z_nI(u, v)$ if and only if P is a right neutrosophic ideal of $Z_nI(v, u)$; $v, u \in Z_nI \setminus \{0\}$.*
   *Suppose P is a left neutrosophic ideal of $(N(Z_n), (x, y), *)$; $x, y \in N(Z_n) \setminus \{0\}$ then P is a right neutrosophic ideal of $(N(Z_n), (y, x), *)$.*

In view of this we have the following theorem, the proof of which is left as an exercise to the reader.

**THEOREM 3.2.4:** *P is a neutrosophic left ideal of $(N(Z_n), (t, u), *)$ if and only if P is a neutrosophic right ideal of $(N(Z_n), (u, t), *)$ where $u, t \in N(Z_n)$.*

**THEOREM 3.2.5:** *No neutrosophic groupoid $(Z_nI, (t, u), *)$ $[N(Z_n), (t, u), *)]$ has $\{0\}$ to be an ideal.*

*Proof:* Follows by the very operation in $Z_nI$ or $(N(Z_n))$.

**THEOREM 3.2.6:** *The neutrosophic groupoid $Z_nI$ $(t, u)$ is an idempotent groupoid if and only if $t + u \equiv I \pmod{n}$.*

*Proof:* Recall a neutrosophic groupoid $Z_nI$ $(u, t)$ is an idempotent neutrosophic groupoid if and only if $a * a \equiv a \pmod{n I}$ for all $a \in Z_nI$ $(u, v)$.

Now $a * a = at + ua \equiv (t + u) a \equiv aI \pmod{n I}$ as $t + u \equiv I \pmod{nI}$. This is possible if and only if $t + u \equiv I \pmod{n}$.

We will illustrate this situation by an example.

*Example 3.2.4:* Let $Z_6I$ $(2I, 5I)$ be a pure neutrosophic groupoid.
Take $0 * 0 = 0$ (trivial). Consider

I * I       =    2I.
I * 5I × I  =    2I + 5I
            =    I (mod 6I).

2I * 2I     =    2I * 2I + 5I × 2I
            =    4I + 10I
            =    2I (mod 6I).
So
2I * 2I     =    2I
3I * 3I     =    3I × 2I + 3I × 5I
            =    6I + 15I
            =    3I (mod 6I)
4I * 4I     =    4I (easily proved)
5I * 5I     =    5I.

Thus $Z_6I$ $(2I, 5I)$ is a pure neutrosophic idempotent groupoid.
Infact there exists 4 neutrosophic idempotent groupoids constructed using $Z_6I$ given by $Z_6I$ $(2I, 5I)$, $Z_6I$ $(5I, 2I)$, $Z_6I$ $(3I, 4I)$ and $Z_6I$ $(4I, 3I)$. In case of $Z_9I$ we see there exists 7 neutrosophic idempotent groupoids built using $Z_9I$.

**THEOREM 3.2.7:** *Let $Z_nI$ (x, y) be a neutrosophic groupoid. If n is even there exists even number of neutrosophic idempotent groupoids. If n is odd then there exists an odd number of neutrosophic idempotent groupoids.*

This can be proved using simple number theoretic techniques.
Consider the neutrosophic groupoid $Z_6I$ (3I, 4I). Clearly $Z_6I$ (3I, 4I) is a neutrosophic semigroup.
For
(a * b) * c  =   (3Ia + 4Ib) * c
             =   (3Ia + 4Ib) 3I + 4I c
             =   3I a + 4Ic.
a * (b * c)  =   a * [3Ib + 4Ic]
             =   3aI + (3Ib + 4Ic) 4I
             =   3Ia + 4Ic.
So (a * b) * c = a * (b * c) for all a, b, c $\in Z_6I$.

Thus $Z_6I$ (3I, 4I) is a only neutrosophic semigroup and $Z_6I$ (4I, 3I) is also a neutrosophic semigroup and not a neutrosophic groupoid as the operation * is associative. We see 4I * 4I = 4I(mod 6I) and 3I * 3I = 3I(mod 6I).

In view of this we have in neutrosophic groupoids the following result.

**THEOREM 3.2.8:** *Let $Z_nI$ (u, t) be a neutrosophic groupoid. If $t^2 = t$ (mod nI) and $u^2 = u$ (mod nI) then $Z_nI$ (u, t) and $Z_nI$ (t, u) are neutrosophic semigroups.*

The proof is straight forward.

**THEOREM 3.2.9:** *Let $Z_nI$ (t, u) be a neutrosophic groupoid. If n is a prime number no groupoid built using with the operation x * y = ux + yt (mod nI); u, t $\in Z_nI$ is a neutrosophic semigroup.*

The proof of this theorem is also left as an exercise for the reader.

Recall as in case of groupoids. We say a neutrosophic groupoid G is normal if

  i. xG = G x
  ii. G (xy) = (Gx)y
  iii. x(yG) = xy G for all x, y ∈ G.

Let G be a neutrosophic groupoid H and K be any two proper neutrosophic subgroupoids of G with H ∩ K = ϕ. We say H is conjugate with K if there exists x ∈ G such that H = xK or Kx.

Several properties enjoyed by groupoids built using $Z_n(u, v)$ will be true for neutrosophic groupoids $Z_nI (u, v)$.

However for this new class of neutrosophic groupoids $(N(Z_n), (u, v), *)$ we do not know whether all properties are derivable. However every neutrosophic groupoid $Z_nI(u, v) \subseteq (N(Z_n) (u, v), *)$.

Now we will define Smarandache neutrosophic groupoid.

**DEFINITION 3.2.2:** *Let G be a neutrosophic groupoid under the operation \*. If H ⊆ G; H a proper subset of G is such that (H, \*) is a neutrosophic semigroup then we call G to be a Smarandache neutrosophic groupoid (S-neutrosophic groupoid).*

It is left as an exercise for the reader to prove that every neutrosophic subgroupoid of a Smarandache neutrosophic groupoid need not in general be a Smarandache neutrosophic subgroupoid.

We can build yet a new class of neutrosophic groupoids $(Z_nI, *)$ which is still restricted neutrosophic groupoid. To this end we say two neutrosophic number x, y ∈ $N(Z_n)$ or $Z_nI$ are relatively neutrosophic prime if (x, y) = I or (x, y) = 1. If n = 25; (8I, 9I) = I so 8I and 9I are relatively neutrosophic prime.

Take (2 + I, 7I) = 1 so 2 + I and 7I are relatively neutrosophic prime. We do not differentiate (x, y) = 1 or (x, y) = I; x, y ∈ $N(Z_n)$.

**DEFINITION 3.2.3:** *Let $Z_nI = \{0, I, 2I, …, (n-1)I\}$, $n > 3$, $n < \infty$. Define * a closed binary operation on $Z_nI$ as follows. For any $a, b \in Z_nI$ define $a * b = ta + ub$; $(t, u) = I$ $(t \neq u)$; $t, u \in Z_n I \setminus \{0\}$.*

*We see $(Z_nI, (t, u), *)$ is a neutrosophic groupoid defined as prime special neutrosophic groupoid.*

Take for example the set $Z_6I$.

***Example 3.2.5:*** *$(Z_6I, (3I, 5I), *)$ is a prime special neutrosophic groupoid. $(Z_6I, (2I, 4I), *)$ is not a prime special neutrosophic groupoid. $(Z_6I (I, 4I), *)$ is also a prime special neutrosophic groupoid $(Z_6I, (2I, 3I), *)$ is also a prime special neutrosophic groupoid. $(Z_6I, (2I, 5I), *)$ is a prime special neutrosophic groupoid. Thus we have using $Z_6I$, 18 prime special neutrosophic groupoids.*

Further we see the class of prime special neutrosophic groupoids built using $Z_nI$ is contained in the class of neutrosophic groupoids built using $Z_nI$.

It is an interesting and important problem to study the number of prime special neutrosophic groupoids built using $Z_nI$.

Now $\{Z_nI, (t, u), *\}$ is a non prime special neutrosophic groupoid if $(t, u) = rI$ if $t, u \in Z_nI$; $(t, u) = r$ if $t, u \in N(Z_n)$, $r \in Z_n$. We can also derive almost all properties for the new class of prime special neutrosophic groupoids.

***Example 3.2.6:*** *Let $(Z_8I, (2I, 4I), *)$ be a non prime special neutrosophic groupoids which is not a prime special neutrosophic groupoid. $(Z_8I, (4I, 4I), *)$ is also a neutrosophic groupoid which is non prime special neutrosophic groupoid.*

Based on this we define yet another new class of neutrosophic groupoids.

**DEFINITION 3.2.4:** *Let $\{Z_nI, *, (tI, tI)\}$ be a neutrosophic groupoid. This neutrosophic groupoid is known as the equal special neutrosophic groupoid.*

We see equal special neutrosophic groupoid is not a prime special neutrosophic groupoid or non prime special neutrosophic groupoid.

We will illustrate this by an example.

**Example 3.2.7:** Let $(Z_5I, (2I, 2I), *)$ be a equal special neutrosophic groupoid given by the following table.

| * | 0 | I | 2I | 3I | 4I |
|---|---|---|----|----|----|
| 0 | 0 | 2I | 4I | I | 3I |
| I | 2I | 4I | I | 3I | 0 |
| 2I | 4I | I | 3I | 0 | 2I |
| 3I | I | 3I | 0 | 2I | 4I |
| 4I | 3I | 0 | 2I | 4I | I |

Consider $(2I * 3I) * 4I = 3I$, $2I * (3I * 4I) = 2I (* 4I) = 2I$.

Thus $(Z_5I, (2I, 2I), *)$ is a neutrosophic groupoid. Thus we have 3 neutrosophic groupoids called the equal special neutrosophic groupoids.

**THEOREM 3.2.10:** *Let $(Z_nI, (tI, tI), *)$ be a equal special neutrosophic groupoid. There exists $(n - 2)$ equal special neutrosophic groupoids.*

Follows using number theoretic methods.

**THEOREM 3.2.11:** *The equal special neutrosophic groupoids are commutative neutrosophic groupoids.*

*Proof:* Given $(Z_nI, (tI, tI), *)$, $t < 3$ is a equal special neutrosophic groupoid. We have
$$a * b = tIa + tIb$$
$$= (a + b) tI,$$
$$b * a = tIb + tI a.$$
$$= (b + a) tI.$$
So $a * b = b * a$ for all $a, b \in Z_nI$. Thus $(Z_nI, (tI, tI), *)$ are commutative neutrosophic groupoids.

**THEOREM 3.2.12:** *($N(Z_n)$, $(tI, tI)$, \*) is a equal special neutrosophic groupoid. Clearly $(Z_nI, (tI, tI),$ \*) $\subseteq (N(Z_n) (tI, tI),$ \*). Thus $(Z_nI, (tI, tI),$ \*) is a equal special neutrosophic subgroupoid of the equal special neutrosophic groupoid.*

**THEOREM 3.2.13:** *The equal special neutrosophic groupoids $(Z_pI, (t, t),$ \*) are normal neutrosophic groupoids.*

*Proof:* Clearly the equal special neutrosophic groupoid $(Z_pI, (t, t),$ \*) are normal for we have a $(Z_pI, (t, t)) = (Z_pI (t, t))$ a for all a $\in Z_pI$.

Also it is easily proved $((Z_pI, (t, t)] x) y = (Z_pI (t, t)] (xy)$ and $(xy) (Z_pI(t, t)) = x (y (Z_pI, (t, t))$, p a prime.

Thus we have a new class of equal special neutrosophic normal groupoids.

Recall a neutrosophic groupoid G is said to be a P- groupoid if $(xy) x = x (yx)$ for all $x, y \in G$.

We will first illustrate it by an example.

*Example 3.2.8:* $Z_6I$ (4I, 4I) is a neutrosophic P – groupoid. For take a, b $\in Z_6I$ (4I 4I);

$$\begin{aligned} a * (b * a) &= a * (4Ib + 4aI) \\ &= 4aI + 4I (4Ib + 4aI) \\ &= 4aI + 16Ib + 16aI \\ &= 20aI + 16Ib. \end{aligned} \quad (1)$$

$$\begin{aligned} (a * b) * a &= (4aI + 4bI) * a \\ &= 16aI + 16bI + 4aI \\ &= 20aI + 16Ib. \end{aligned} \quad (2)$$

We see (1) and (2) are the same. Hence $Z_6I$ (4I, 4I) is a neutrosophic P-groupoid.

We have a class of neutrosophic P-groupoids which is evident from the following theorem.

**THEOREM 3.2.14:** *The neutrosophic groupoids $Z_nI$ (t, t) are neutrosophic P groupoids, $t \in Z_nI$.*

*Proof:* For every x, y ∈ $Z_nI$ (t, t) we have

$$x * (y * x) = x * (ty + tx)$$
$$= tx + t^2y + t^2x.$$
$$(x * y) * x = (tx + ty) * x$$
$$= t^2x + t^2y + tx.$$

Thus (x * y) * x = x * (y * x) for all x, y ∈ $Z_nI$ (t, t). Hence the claim.

We see we have infact a large class of neutrosophic P-groupoids which is evident from the following theorem.

**THEOREM 3.2.15:** *The neutrosophic groupoids ($N(Z_n)$, (t, t), *) are P- groupoids.*

The proof of the above theorem is left as an exercise for the reader to prove.
We say a neutrosophic groupoid G is said to be an alternative neutrosophic groupoid if (xy) y = x (yy) for all x, y ∈ G.

We first illustrate by an example the neutrosophic groupoid $Z_9I$ (5I, 5I) which is not an alternative groupoid.

*Example 3.2.9:* Let $Z_9I$ (5I, 5I) be a neutrosophic groupoid. To prove x * (y * y) ≠ (x * y) * y for x, y ∈ $Z_9I$ .
Consider
$$x * (y * y) = x * (5Iy + 5Iy)$$
$$= 5Ix + 25Iy + 25Iy$$
$$= 5x + 7Iy + 7Iy$$
$$= 5Ix + 5Iy.$$
Take
$$(x * y) * y = (5Ix + 5Iy) * y$$
$$= 25Ix + 25Iy + 5Iy$$
$$= 7Ix + 3Iy.$$

Thus $x * (y * y) \neq (x * y) * y$ for all $x, y \in Z_9I$ (5I, 5I). Hence $Z_9I$ (5I, 5I) is not an alternative neutrosophic groupoid.

In view of this we have the following theorem.

**THEOREM 3.2.16:** *$Z_pI$ (t, t) are not alternative neutrosophic groupoids $t \in Z_pI \setminus \{0\}$.*

*Proof:* Consider $x * (y * y)$ for any $x, y \in Z_pI$ (t, t).

$$\begin{aligned} x * (y * y) &= x * (ty + ty) \\ &= tx + t^2y + t^2y \quad (1) \\ (x * y) * y &= (tx + ty) * y \\ &= t^2x + t^2y + ty \quad (2) \end{aligned}$$

Clearly (1) and (2) are not equal for all t.
Hence the claim.

**THEOREM 3.2.17:** *All neutrosophic groupoids $Z_nI$ (t, t), n not a prime $t \in Z_nI \setminus \{0\}$ with $t^2 \equiv t(\mod nI)$ are alternative groupoids.*

*Proof:* Given $Z_nI$ (t, t) is a neutrosophic groupoid such that n is not a prime and $t^2 = t \pmod{nI}$. ($t \in Z_nI$)
Consider
$$\begin{aligned} x * (y * y) &= x * (ty + ty) \\ &= tx + t^2y + t^2y \\ &= tx + ty + ty \quad (1) \end{aligned}$$
Now
$$\begin{aligned} (x * y) * y &= (tx + ty) * y \\ &= t^2 x + t^2y + ty \\ &= tx + ty + ty \quad (2) \end{aligned}$$

From (1) and (2) we see $Z_nI$ (t, t) is a neutrosophic alternative groupoid.

The following theorem is left as an exercise for the reader to prove.

**THEOREM 3.2.18:** *Let $(N(Z_n), (t, t), *)$ be a neutrosophic groupoid where $t^2 = t \pmod{nI}$, $n$ a non prime. Then $(N(Z_n), (t, t), *)$ is an alternative neutrosophic groupoid.*

*Thus we have non trivial class of alternative neutrosophic groupoids as well as a class of neutrosophic groupoids which are not alternative.*

## 3.3 Neutrosophic Polynomial Groupoids

In this section we for the first time define the class of neutrosophic polynomial groupoids built using $Z_nI$ or $N(Z_n)$ or $ZI$ or $N(Z)$ or $QI$ or $N(Q)$ or $N(R)$ or $RI$ or $N(C)$ or $CI$.

$$Z_nI[x] = \left\{ \sum_{i=0}^{n} a_i x^i \,\middle|\, a_i \in Z_nI, n \leq \infty \right.$$

x a variable or an indeterminate$\}$.

$$N(Z_n)[x] = \left\{ \sum_{i=0}^{n} a_i x^i \,\middle|\, a_i \in N(Z_n), n \leq \infty \right\}.$$

Clearly $Z_nI[x] \subseteq N(Z_n)[x]$.

$$ZI[x] = \left\{ \sum_{i=0}^{\infty} a_i x \,\middle|\, a_i \in ZI, n < \infty \right\}$$

$$N(Z)[x] = \left\{ \sum_{i=0}^{\infty} a_i x^i \,\middle|\, a_i \in N(Z), n \leq \infty \right\}$$

and $ZI[x] \subseteq N(Z)[x]$. Similar notations for $RI$, $N(R)$ and so on.

Further $ZI[x] \subseteq QI[x] \subseteq RI[x] \subseteq CI[x]$ and $N(Z)[x] \subseteq N(Q)[x] \subseteq N(R)[x] \subseteq N(C)[x]$. Now we can build finite classes of neutrosophic groupoids of finite or infinite order using $Z_nI[x]$ and $N(Z_n)[x]$ and infinite classes of neutrosophic groupoid of infinite order using $QI[x]$ or $ZI[x]$ or $N(Q)[x]$ or $N(Z)[x]$ and so on.

Groupoids built using $N(Z_n)[x]$, $QI[x]$ and so on will in general be known as polynomial neutrosophic groupoids.

**DEFINITION 3.3.1:** *Let $G = \{Z_n I [x], *, (p, q); p\ q \in Z_n \setminus \{0\}$; p and q are primes\} be such that if*
$$a(x) = a_0 + a_1 x + \ldots + a_n x^n$$
*and*
$$b(x) = b_0 + b_1 x + \ldots + b_m x^m$$
*then*
$$\begin{aligned} a(x) * b(x) &= a_0 * b_0 + (a_1 * b_1)x + \ldots + (a_n * b_n)x^n \\ &\quad + (0 * b_{n+1}) x^{n+1} + \ldots + (0 * b_m) x^m \\ &= (pa_0 + qb_0)(mod\ nI) + (pa_1 + qb_1)(mod\ nI)x \\ &\quad + (pa_2 + qb_2)(mod\ nI) x^2 + \ldots + (pa_n + qb_n) \\ &\quad (mod\ nI)\ x^n + \ldots + (p0 + qb_m)\ (mod\ n\ I)\ x^m \end{aligned}$$
*where $a_i, b_i \in Z_n I$.*

*(If $m < n$ or $m = n$ one can easily define $*$ on $Z_n I [x]$). Thus $\{Z_n I[x], *, (p, q) / p, q \in Z_n \setminus \{0\}$, p and q primes\} = G is a neutrosophic polynomial groupoid of modulo integers of level one.*

We will illustrate this situation by some examples.

*Example 3.3.1:* Let $G = \{Z_{12}I[x], *, (3, 5)\}$ be the polynomial neutrosophic groupoid of modulo integer of level one.

*Example 3.3.2:* Let $G = \{Z_{19}I [x], *, (7, 11)\}$ be the polynomial neutrosophic groupoid of modulo integers of level one.
Clearly we see the cardinality of both the polynomial neutrosophic groupoids given examples 3.3.1 and 3.3.2 are of infinite order.

Now we give also finite order polynomial neutrosophic groupoids of level one.

*Example 3.3.3:* Let

$$G = \left\{ \sum_{i=0}^{9} a_i x^i \,\middle|\, a_i \in Z_4 I, (1,3) \right\}$$

be a neutrosophic polynomial groupoid of finite order.

*Example 3.3.4:* Let

$$G = \left\{ \sum_{i=0}^{4} a_i x^i \,\middle|\, a_i \in Z_{23} I, (5,13) \right\}$$

be the finite neutrosophic polynomial groupoid.

*Example 3.3.5:* Let

$$G = \left\{ \sum_{i=0}^{n} a_i x^i \,\middle|\, a_i \in N(Z_{16}); (7,11) \right\}$$

be an infinite neutrosophic polynomial groupoid as $n \leq \infty$.

*Example 3.3.6:* Let

$$G = \{ \sum_{i=0}^{n} a_i x^i \,;\, a_i \in N(Z_{43}), (23, 29); n \leq \infty \}$$

is an infinite neutrosophic polynomial groupoid.

*Example 3.3.7:* Let

$$G = \{ \sum_{i=0}^{n} a_i x^i \,;\, a_i \in N(Z_3); (1, 2) \}$$

be a finite neutrosophic polynomial groupoid.

*Example 3.3.8:* Let

$$G = \left\{ \sum_{i=0}^{5} a_i x^i \,\middle|\, a_i \in N(Z_6); (3,5) \right\}$$

is a finite neutrosophic polynomial groupoid built using $N(Z_6)$.

Now we can define neutrosophic polynomial integer (real or rational or complex) groupoid similar to the one given in the definition 3.3.1.

We will illustrate this only by examples.

*Example 3.3.9:* Let
$$G = \left\{ \sum_{i=0}^{n} a_i x^i \,\middle|\, n \leq \infty; a_i \in ZI; (13,17) \right\}$$

be a neutrosophic integer coefficient polynomial groupoid.

*Example 3.3.10:* Let
$$G = \left\{ \sum_{i=0}^{n} a_i x^i \,\middle|\, n \leq \infty; a_i \in 3ZI; (2,19) \right\}$$

be the neutrosophic integer coefficient polynomial groupoid.

*Example 3.3.11:* Let
$$G = \left\{ \sum_{i=0}^{n} a_i x^i \,\middle|\, n \leq \infty; a_i \in ZI; (23,53) \right\}$$

be a neutrosophic polynomial integer groupoid.

*Example 3.3.12:* Let $G = \{\Sigma a_i x^i; a_i \in N(Z); (3, 13)\}$ be a neutrosophic polynomial integer groupoid of infinite order.

*Example 3.3.13:* Let
$$G = \left\{ \sum_{i=0}^{n} a_i x^i \,\middle|\, n \leq \infty; a_i \in N(Z); (19,2) \right\}$$

be a neutrosophic polynomial integer groupoid of infinite order.

*Example 3.3.14:* Let
$$G = \left\{ \sum_{i=0}^{7} a_i x^i \,\middle|\, a_i \in N(Z); (3,17) \right\}$$

be an infinite neutrosophic polynomial integer groupoid.

*Example 3.3.15:* Let
$$P = \left\{\sum_{i=0}^{19} a_i x^i \,\middle|\, a_i \in ZI; (3, 29)\right\}$$

be an infinite neutrosophic polynomial integer groupoid.

**DEFINITION 3.3.2:** *Let*
$$G = \left\{\sum_{i=0}^{\infty} a_i x^i, *, a_i \in ZI; (p, q)\right\}$$

*be a neutrosophic polynomial integer groupoid. If $P \subseteq G$ is such that P is a neutrosophic polynomial integer groupoid under the operations of G then we call P to be a neutrosophic integer polynomial subgroupoid of G.*

It is to be noted in this definition we can replace ZI by N(Z) or $Z_nI$ or $N(Z_n)$ or N(Q) or N(R) or N(C) or QI or RI or CI.

We will illustrate however all these by some simple examples.

*Example 3.3.16:* Let
$$G = \left\{\sum_{i=0}^{n} a_i x^i, n \geq \infty, a_i \in N(Z_{12}), (3, 7)\right\}$$

be a neutrosophic polynomial groupoid.
Take
$$P = \left\{\sum_{i=0}^{10} a_i x^i, a_i \in N(Z_{12}), (3, 7)\right\} \subseteq G;$$

P is a neutrosophic polynomial subgroupoid of G.
    Consider
$$X = \left\{\sum_{i=0}^{20} a_i x^i, a_i \in (Z_{12}I), (3, 7)\right\} \subseteq G;$$

is also a neutrosophic polynomial subgroupoid of G.
Let
$$F = \left\{ \sum_{i=0}^{\infty} a_i x^i \,\middle|\, a_i \in Z_{12}I; (3,7) \right\} \subseteq G;$$

be a neutrosophic polynomial subgroupoid of G. We see P and X are finite neutrosophic polynomial subgroupoids of G where as F is an infinite neutrosophic polynomial subgroupoid of G.

*Example 3.3.17:* Let
$$G = \left\{ \sum_{i=0}^{n} a_i x^i \,\middle|\, a_i \in Z_7 I; n \leq \infty, (3,5) \right\}$$

be an infinite polynomial neutrosophic modulo integer groupoid.
Let
$$P = \left\{ \sum_{i=0}^{25} a_i x^i \,\middle|\, a_i \in Z_7 I; (3,5), * \right\} \subseteq G$$

be a finite polynomial neutrosophic subgroupoid of G.
Consider
$$X = \left\{ \sum_{i=0}^{\infty} a_i x^i \,\middle|\, a_i \in Z_7 I, (3,5), * \right\} \subseteq G;$$
X is an infinite polynomial subgroupoid of G.

*Example 3.3.18:* Let
$$G = \left\{ \sum_{i=0}^{\infty} a_i x^i \,\middle|\, a_i \in N(Z), *, (3,17) \right\}$$

be an infinite polynomial neutrosophic integer groupoid; G has no finite polynomial neutrosophic subgroupoid. However G has infinitely many polynomial neutrosophic subgroupoids.
Take
$$P = \left\{ \sum_{i=0}^{\infty} a_i x^i \,\middle|\, a_i \in N(mZ), *, (3,17) \right\} \subseteq G$$

are all infinite polynomial integer neutrosophic subgroupoids of G.

Since m = 2, 3, ..., ∞ we have infinitely many such subgroupoids but G has no finite polynomial neutrosophic subgroupoids.

***Example 3.3.19:*** Let

$$G = \left\{ \sum_{i=0}^{\infty} a_i x^i \,\middle|\, a_i \in QI, *, (29,3) \right\}$$

be a neutrosophic polynomial rational groupoid of infinite order. G has infinitely many polynomial neutrosophic subgroupoids of infinite order but has no neutrosophic polynomial subgroupoid of finite order.

Take

$$P = \left\{ \sum_{i=0}^{\infty} a_i x^{2i} \,\middle|\, a_i \in QI, (29,3) \right\} \subseteq G;$$

P is an infinite neutrosophic polynomial subgroupoid of infinite order.

***Example 3.3.20:*** Let

$$G = \left\{ \sum_{i=0}^{20} a_i x^i \,\middle|\, a_i \in N(Q), *, (-17, 23) \right\}$$

be a neutrosophic polynomial rational neutrosophic groupoid of infinite order.
Take

$$W = \left\{ \sum_{i=0}^{10} a_i x^i \,\middle|\, a_i \in N(Q), *, (-17, 23) \right\} \subseteq G$$

is also an infinite polynomial rational neutrosophic subgroupoid of G.

***Example 3.3.21:*** Let

$$G = \left\{ \sum_{i=0}^{n} a_i x^{2i} \,\middle|\, n \leq \infty, a_i \in N(R), *, (-2, 7) \right\}$$

be an infinite polynomial real neutrosophic groupoid.
Let
$$P = \left\{ \sum_{i=0}^{\infty} a_i x^i, *, (-2,7), a_i \in N(Q) \right\} \subseteq G,$$

P is an infinite polynomial real neutrosophic subgroupoid of G.

*Example 3.3.22:* Let
$$G = \left\{ \sum_{i=0}^{26} a_i x^i \,\bigg|\, a_i \in N(C), *, (7,17) \right\}$$

be an infinite complex polynomial neutrosophic groupoid.
Take
$$M = \left\{ \sum_{i=0}^{\infty} a_i x^i \,\bigg|\, a_i \in RI, *, (7,17) \right\} \subseteq G,$$

M is also an infinite complex polynomial subgroupoid of G.

Now having seen examples of these groupoids now we proceed onto state that all properties and definitions like ideals of a groupoid, normal groupoid, normal subgroupoid, conjugate subgroupoids and so on which can always be adopted to neutrosophic polynomial groupoids with appropriate simple modifications.

Only in case of homomorphism $\eta$ of neutrosophic groupoids G and G′; we demand $\eta(I) = I$ and $\eta$ should be a groupoid homomorphism. By no means I should be mapped on to any real number for the indeterminate can never be compensated, it should always continue to remain as an indeterminate. The reader is expected to give some examples of polynomial neutrosophic groupoid homomorphism.

We expect the reader to study the concept of substructures of polynomial neutrosophic groupoids.

We call a polynomial neutrosophic groupoid G to be a Smarandache polynomial neutrosophic groupoid if G has a proper subset S where S under the operations of G is a polynomial neutrosophic semigroup.

Now we define yet another new substructure in case of neutrosophic polynomial groupoids.

**DEFINITION 3.3.3:** *Let*

$$G = \left\{ \sum_{i=0}^{n} a_i x^i \,\middle|\, a_i \in N(Q), *, (p,q) \right\}$$

*be a neutrosophic polynomial groupoid of level one. Let*

$$P = \left\{ \sum a_i x^i \,\middle|\, a_i \in Q, *, (p,q) \right\};$$

*p, q are primes in Q ⊆ G. P is clearly a polynomial (sub)groupoid but is not a neutrosophic polynomial groupoid. We call P to be a pseudo neutrosophic polynomial subgroupoid of G. If G has no pseudo neutrosophic polynomial subgroupoid then we call G to be a pseudo simple neutrosophic polynomial groupoid.*

We will illustrate both these situations by some simple examples.

*Example 3.3.23:* Let

$$G = \left\{ \sum_{i=0}^{8} a_i x^i \,\middle|\, a_i \in N(Q), *, (3,7) \right\}$$

be a neutrosophic polynomial rational groupoid.
 Let

$$W = \left\{ \sum_{i=0}^{7} a_i x^i \,\middle|\, a_i \in 3Z, *, (3,7) \right\} \subseteq G$$

is a pseudo neutrosophic polynomial subgroupoid of G.

*Example 3.3.24:* Let

$$G = \left\{ \sum_{i=0}^{n} a_i x^i \,\middle|\, a_i \in NI, *, (3,13) \right\}$$

be a polynomial neutrosophic groupoid we see G has polynomial neutrosophic subgroupoids but G has no pseudo neutrosophic polynomial subgroupoid. Thus G is a pseudo simple neutrosophic polynomial groupoid.

*Example 3.3.25:* Let

$$G = \left\{ \sum_{i=0}^{\infty} a_i x^i \,\middle|\, a_i \in RI, *, (3, 23) \right\}$$

be a neutrosophic polynomial groupoid. It is easily verified G has no pseudo neutrosophic polynomial subgroupoids hence G is a pseudo simple neutrosophic polynomial groupoid.

In view of this we have the following theorem.

**THEOREM 3.3.1:** *Let*

$$G = \left\{ \sum_{i=0}^{n} a_i x^i \,\middle|\, a_i \in Z_n I \;(\; ZI \; or \; RI \; or \; QI \; or \; CI\;), \; *, \; (p, q);\right.$$

*p and q primes in $Z^+$} be a neutrosophic polynomial groupoid; clearly G is a pseudo simple neutrosophic groupoid.*

*Proof:* Since $Z_n$I or ZI or QI or RI or QI are pure neutrosophic sets we see the polynomial groupoids G built using them cannot have even a single non zero real coefficient. Hence G is always a pseudo simple neutrosophic groupoid in such cases. However G can have neutrosophic polynomial subgroupoids.

We give yet another theorem which gurantees the existence of pseudo neutrosophic polynomial subgroupoids.

**THEOREM 3.3.2:** *Let*

$$G = \left\{ \sum_{i=0}^{n} a_i x^i \,\middle|\, a_i \in N(Z_n\,)(\; or \; N(Z\,) \; or \; N(Q\,) \; or \; N(R\,) \; or \; N(C\,)),\right.$$

*\*, (p, q)}; be a neutrosophic polynomial groupoid. G has non trivial pseudo neutrosophic polynomial subgroupoid.*

*Proof:* Take
$$P = \left\{ \sum_{i=0}^{n} a_i x^i \,\middle|\, a_i \in Z_n \subseteq N(Z_n) \right.$$
(or $Z \subseteq N(Z)$, $Q \subseteq N(Q)$ or $R \subseteq N(R)$ or $C \subseteq N(C)$, *, (p, q)} $\subseteq$ G. P is clearly a pseudo neutrosophic polynomial subgroupoid of G.

Hence the theorem.

*Example 3.3.26:* Let
$$G = \left\{ \sum_{i=0}^{n} a_i x^i \,\middle|\, a_i \in N(Z_{12}), (7, 5), * \right\}$$
be a neutrosophic polynomial groupoid.
Take
$$W = \left\{ \sum_{i=0}^{n} a_i x^i \,\middle|\, a_i \in Z_{12}, (7, 5), * \right\} \subseteq G,$$
W is a pseudo neutrosophic polynomial subgroupoid of G.

*Example 3.3.27:* Let
$$G = \left\{ \sum_{i=0}^{n} a_i x^i \,\middle|\, a_i \in Z_{25} I *, (3, 17) \right\}$$
be a neutrosophic polynomial groupoid. Consider M = {$\Sigma a_i x^i$ | $a_i \in$ {0, 5I, 10I, 15I, 20I}, (3, 17)} $\subseteq$ G; M is a neutrosophic polynomial subgroupoid of G. However G has no pseudo neutrosophic polynomial subgroupoids. Thus G is a pseudo simple neutrosophic polynomial groupoid.

Note: If a neutrosophic polynomial groupoid has no neutrosophic polynomial subgroupoids as well as no pseudo neutrosophic polynomial subgroupoids then we define G to be a doubly simple neutrosophic polynomial groupoid.

*Example 3.3.28:* Let G = {$a_i x^2$ | $a_i \in Z_{23} I$, *, (3, 2)} be a neutrosophic polynomial groupoid. Clearly G is a doubly simple neutrosophic polynomial groupoid.

Next we proceed onto define level two neutrosophic polynomial groupoids.

**DEFINITION 3.3.4:** *Let $G = \{\Sigma a_i x^i \mid a_i \in Z_n$ (or $N(Z_n)$ or $ZI$ or $N(Z)$ or $RI$ or $QI$ or $N(Q)$ or $N(R))$, $*$, $(p, q)$ such that $p$ and $q$ are just numbers not necessarily prime but $(p, q) = 1\}$; $G$ is defined to be a polynomial neutrosophic groupoid of level two.*

We will illustrate by a few examples before we proceed onto define other levels of neutrosophic polynomial groupoids.

*Example 3.3.29:* Let
$$G = \left\{\sum_{i=0}^{\infty} a_i x^i \,\bigg|\, a_i \in ZI, *, (9, 16)\right\}$$

be the neutrosophic polynomial groupoid of level two.

*Example 3.3.30:* Let
$$G = \left\{\sum_{i=0}^{\infty} a_i x^i \,\bigg|\, a_i \in N(Q), *, (24, 49)\right\}$$

be the neutrosophic polynomial rational groupoid of level two.

*Example 3.3.31:* Let
$$G = \left\{\sum_{i=0}^{n} a_i x^i \,\bigg|\, a_i \in Z_{12}I, (9, 8), *\right\}$$

be a neutrosophic polynomial groupoid of level two.

*Example 3.3.32:* Let
$$G = \left\{\sum_{i=0}^{n} a_i x^i \,\bigg|\, a_i \in N(C), *, (24, 35)\right\}$$

be the neutrosophic polynomial complex groupoid of level two.

Substructures as in case of level one groupoids can be defined in case of level two neutrosophic polynomials groupoids also.

Further from the very context one can easily understand the level to which the polynomial neutrosophic groupoid belongs to. Now we define level three polynomial neutrosophic groupoids.

**DEFINITION 3.3.5:** *Let*

$$G = \left\{ \sum_{i=0}^{n} a_i x^i \,\middle|\, a_i \in N(Q) \right.$$

*(or QI, $Z_nI$ or $N(Z_n)$ or RI or $N(R)$ or $N(C)$ or CI) / $n \leq \infty$, \*, (p, q); (p, q) = d ≠ 1, p, q ∈ $Z^+$}, G is a polynomial neutrosophic groupoid defined as level three neutrosophic polynomial groupoid.*

We will illustrate this by some simple examples.

*Example 3.3.33:* Let

$$G = \left\{ \sum_{i=0}^{n} a_i x^i \,\middle|\, n \leq \infty, a_i \in Z_{17}I, *, (3, 15) \right\}$$

be a neutrosophic level three polynomial groupoid.

*Example 3.3.34:* Let

$$G = \left\{ \sum_{i=0}^{\infty} a_i x^i \,\middle|\, a_i \in ZI, *, (7, 147) \right\}$$

be a neutrosophic polynomial groupoid of level three.

*Example 3.3.35:* Let

$$G = \left\{ \sum_{i=0}^{\infty} a_i x^i \,\middle|\, a_i \in N(R), *, (3, 42) \right\}$$

be the neutrosophic polynomial groupoid of level three.

We can define substructure of level three neutrosophic polynomial groupoids as in case of level one neutrosophic polynomial groupoids.

Now we proceed onto define level four neutrosophic polynomial groupoids.

**DEFINITION 3.3.6:** *Let*
$$G = \left\{ \sum_{i=0}^{\infty} a_i x^i \,\middle|\, a_i \in N(Z_n) \right.$$

*(or $Z_n I$ or $N(Z)$ or $ZI$ or $RI$ or $N(R)$ or $N(Q)$ or $QI$ or $CI$ or $N(C)$), *, (t, t); $t \in Z^+$} be a neutrosophic polynomial groupoid. G is defined as the neutrosophic polynomial groupoid of level four.*

We will illustrate this situation by some simple examples.

*Examples 3.3.36:* Let
$$G = \left\{ \sum_{i=0}^{29} a_i x^i \,\middle|\, a_i \in QI, *, (5, 5) \right\}$$

be a neutrosophic polynomial groupoid of level four.

*Example 3.3.37:* Let
$$G = \left\{ \sum_{i=0}^{\infty} a_i x^i \,\middle|\, a_i \in N(R), *, (22, 120) \right\}$$

be a neutrosophic polynomial groupoid of level four.

*Example 3.3.38:* Let
$$G = \left\{ \sum_{i=0}^{5} a_i x^i \,\middle|\, a_i \in N(Z_{24}), *, (8, 18) \right\}$$

be a neutrosophic polynomial groupoid of level four.

Substructure of these groupoids can be constructed as in case of other level groupoids.

Finally we now proceed onto define level five groupoids.

**DEFINITION 3.3.7:** *Let*

$$G = \left\{ \sum_{i=0}^{\infty} a_i x^i \,\middle|\, a_i \in N(Z_n) \right.$$

*(or $Z_nI$, $ZI$, $N(Z)$ or $QI$ or $N(Q)$ or $N(R)$ or $RI$ or $CI$ or $N(C)$), \*, $(0,t)$ $t \neq 0$, $t \in Z^+$} be a neutrosophic polynomial groupoid. We define G to be a polynomial neutrosophic groupoid of level five.*

We will illustrate this by some examples.

*Example 3.3.39:* Let

$$G = \left\{ \sum_{i=0}^{\infty} a_i x^i \,\middle|\, a_i \in N(Q), *, (3,0) \right\}$$

be a polynomial neutrosophic groupoid of level five.

*Example 3.3.40:* Let

$$G = \left\{ \sum_{i=0}^{6} a_i x^i \,\middle|\, a_i \in N(Z_{12}), *, (0,8) \right\}$$

be a neutrosophic polynomial groupoid of level five.

*Example 3.3.41:* Let

$$G = \left\{ \sum_{i=0}^{\infty} a_i x^i \,\middle|\, a_i \in Z_{26}I, *, (0,9) \right\}$$

be a neutrosophic polynomial groupoid of level five.

*Example 3.3.42:* Let

$$G = \left\{ \sum_{i=0}^{9} a_i x^i \,\middle|\, a_i \in QI, *, (25, 0) \right\}$$

be a neutrosophic polynomial groupoid of level five.

Now we give results about these neutrosophic polynomials groupoids.

**THEOREM 3.3.3:** *The neutrosophic polynomial groupoids of level four are commutative groupoids.*

*Proof:* Follows from the very fact if G is any neutrosophic polynomial groupoid of level four then for a(x), b(x) ∈ G. We have a(x) * b(x) = ta(x) + tb(x) and b(x) * a(x) = tb(x) + tb(x), hence G is commutative groupoid. In case of building groupoids using $Z_nI$ we can choose t, u ∈ $Z_nI$.

**THEOREM 3.3.4:** *Let G be a neutrosophic polynomial groupoid of the form*

$$G = \left\{ \sum_{i=0}^{n} a_i x^i \,\middle|\, a_i \in Z_p I,\ p\ a\ prime,\ *,\ (tI, tI);\ tI < pI \right\}$$

*then G is a neutrosophic polynomial groupoid.*

The reader is expected to prove the theorem.

**THEOREM 3.3.5:** *Let*

$$G = \left\{ \sum_{i=0}^{n} a_i x^i\ ;\ n \leq \infty,\ a_i \in Z_n I, *,\ (tI, tI) \right\}$$

*be a neutrosophic polynomial groupoid of level four; then G is a neutrosophic polynomial P-groupoid.*

Proof is obvious

**THEOREM 3.3.6:** *Let*

$$G = \left\{ \sum_{i=0}^{n} a_i x^i \,\middle|\, n \leq \infty, a_i \in Z_p I\ ;\ p\ a\ prime,\ *,\ (tI, tI);\ tI < pI \right\}$$

*be a neutrosophic polynomial groupoid of level four. G is not a neutrosophic polynomial alternative groupoid.*

Proof is left as an exercise for the reader.

**THEOREM 3.3.7:** *Let*

$$G = \left\{ \sum_{i=0}^{\infty} a_i x^i \,\middle|\, a_i \in Z_n I \,;\, n \text{ is not a prime}, *, (tI, tI) \right\}$$

*is a neutrosophic alternative polynomial groupoid if and only if $(tI)^2 \equiv tI \pmod{n}$.*

*Proof:* Clearly this is proved using the alternative identity.

**THEOREM 3.3.8:** *Let*

$$G = \left\{ \sum_{i=0}^{\infty} a_i x^i \,\middle|\, a_i \in Z_n I, *, (0, tI) \right\}$$

*is a neutrosophic polynomial P-groupoid and neutrosophic polynomial alternative groupoid if and only if $(tI)^2 \equiv tI \pmod{n}$.*

Proof is left as an exercise for the reader.

*Example 3.3.43:* Let

$$G = \left\{ \sum_{i=0}^{8} a_i x^i \,\middle|\, a_i \in Z_6 I, *, (4I, 5I) \right\}$$

be a Smarandache polynomial neutrosophic groupoid.
For

$$A = \left\{ \sum_{i=0}^{8} 3I x^i \right\} \subseteq G$$

is a semigroup.
Also

$$P = \left\{ \sum_{i=0}^{8} a_i x^i \,\middle|\, a_i \in \{I, 3I, 5I\}, *, (4I, 5I) \right\} \subseteq G$$

is not Smarandache right ideal of G but P is a Smarandache left ideal of G.

*Example 3.3.44:* Let
$$G = \left\{ \sum_{i=0}^{9} a_i x^i \,\Big|\, a_i \in Z_8 I, *, (2I, 6I) \right\}$$
be a polynomial groupoid. Let
$$A = \left\{ \sum_{i=0}^{9} a_i x^i \,\Big|\, a_i \in \{0, 2I, 4I, 6I\}, *, (2I, 6I) \right\} \subseteq G$$

be a Smarandache subgroupoid. It is easily verified A is a Smarandache normal subgroupoid of G.

*Example 3.3.45:* Let
$$G = \left\{ \sum_{i=0}^{4} a_i x^i \,\Big|\, a_i \in Z_{10} I, *, (5I, 6I) \right\}$$

be a neutrosophic polynomial groupoid. It is easily verified G is a Smarandache strong Moufang groupoid.

*Example 3.3.46:* Let
$$G = \left\{ \sum_{i=0}^{8} a_i x^i \,\Big|\, a_i \in Z_{12} I, *, (3I, 9I) \right\}$$

be a neutrosophic polynomial groupoid. G is a Smarandache Moufang groupoid.

*Example 3.3.47:* Let
$$G = \left\{ \sum_{i=0}^{9} a_i x^i \,\Big|\, a_i \in Z_{12} I, *, (3I, 4I) \right\}$$

be a neutrosophic polynomial groupoid. It is easily verified that G is a Smarandache strong Bol polynomial neutrosophic groupoid.

*Example 3.3.48:* Let
$$G = \left\{ \sum_{i=0}^{\infty} a_i x^i \,\Big|\, a_i \in Z_4 I, 8, (2I, 3I) \right\}$$

be a neutrosophic polylnomial groupoid. It is easily verified G is a Smarandache Bol groupoid and is not a Smarandache strong Bol groupoid.

*Example 3.3.49:* Let
$$G = \left\{ \sum_{i=0}^{120} a_i x^i \,\middle|\, a_i \in Z_6 I, 8, (4I, 3I) \right\}$$
be a neutrosophic polynomial groupoid. G is a Smarandache strong neutrosophic polynomial P-groupoid.

*Example 3.3.50:* Let
$$G = \left\{ \sum_{i=0}^{27} a_i x^i \,\middle|\, a_i \in Z_4 I, *, (2I, 3I) \right\}$$

be a neutrosophic polynomial groupoid. G is a Smarandache strong neutrosophic polynomial P-groupoid.

*Example 3.3.51:* Let
$$P = \left\{ \sum_{i=0}^{\infty} a_i x^i \,\middle|\, a_i \in Z_{14} I, *, (7I, 8I) \right\}$$

be a neutrosophic polynomial groupoid. P is a Smaradanche strong alternative neutrosophic polynomial groupoid.

*Example 3.3.52:* Let
$$G = \left\{ \sum_{i=0}^{3} a_i x^i \,\middle|\, a_i \in Z_5 I, *, (I, 3I) \right\}$$
be a neutrosophic polynomial groupoid. Clearly G is not a Smarandache groupoid.

**THEOREM 3.3.9:** Let
$$G = \left\{ \sum_{i=0}^{7} a_i x^i \,\middle|\, a_i \in Z_n I, *, (tI, uI) \right\}$$
*be a polynomial neutrosophic groupoid. If $t + u \equiv I \pmod{nI}$ then G is a Smarandache neutrosophic idempotent groupoid.*

The proof is obvious by using simple number theoretic techniques.

**THEOREM 3.3.10:** *Let*

$$G = \left\{\sum_{i=0}^{\infty} a_i x^i \,\middle|\, a_i \in Z_n I, *, (tI, uI) / t + u = I \ (mod\ n)\right\}$$

*be Smarandache neutrosophic polynomial groupoid. G is a Smaradache P-groupoid if and only if $t^2 I \equiv tI\ (mod\ n)$ and $u^2 I = uI\ (mod\ n)$.*

Proof is left as an exercise for the reader.

**THEOREM 3.3.11:** *Let*

$$G = \left\{\sum_{i=0}^{\infty} a_i x^i \,\middle/\, a_i \in Z_n I, *, (tI, uI)\ with\ tI + uI \equiv I\ (mod\ nI)\right\}$$

*be a neutrosophic polynomial groupoid, G is a Smaradache strong Moufang groupoid if and only if $t^2 I \equiv tI\ (mod\ n)$ and $u^2 I = uI\ (mod\ n)$.*

Proof is left as an exercise for the reader.

Several results of this type can be derived for neutrosophic polynomial groupoids built using $Z_n I$.

## 3.4 Neutrosophic Matrix Groupoids

In this section we for the first time introduce a new class of neutrosophic matrix groupoids using $Z_n I$ or $N(Z_n)$ or $ZI$ or $N(Z)$, $QI$ or $N(Q)$, $RI$ or $N(R)$ or $N(C)$ or $CI$. We define some new classes and enumerate a few properties about them. Let $X = (a_1, a_2, \ldots, a_n)$ where $a_i \in Z_n I$ or $N(Z_n)$ or $QI$ or $RI$ or so on. We call X a neutrosophic row matrix with entries from $Z_n I$ or $N(Z_n)$ or $RI$ or $QI$ or so on.
 Likewise

$$P = \begin{bmatrix} y_1 \\ y_2 \\ \vdots \\ y_m \end{bmatrix}$$

where $y_i \in N(R)$ or $N(Z_n)$ or $Z_n I$ or $QI$ or $N(Q)$ or so on is a neutrosophic column matrix.

Now take

$$M_{n \times m} = \begin{bmatrix} m_{11} & \cdots & m_{1m} \\ m_{21} & \cdots & m_{2m} \\ \vdots & & \vdots \\ m_{n1} & \cdots & m_{nm} \end{bmatrix};$$

$m_{ij} \in N(Z_n)$ or $Z_n I$ or $N(R)$ or $RI$ or so on; we call $M_{n \times m}$ a neutrosophic n×m matrix. If n = m we call $M_{n \times m}$ to be a neutrosophic n×m square matrix.

With this convention we now proceed onto define several new classes of neutrosophic matrix groupoids.

**DEFINITION 3.4.1:** *Let $G = \{(x_1, \ldots, x_n) / x_i \in Z_n I, *, (p, q), p$ and $q$ primes$\}$ be a neutrosophic row matrix groupoid with the operation for *; $x = (x_1, \ldots, x_n)$ and $y = (y_1, \ldots, y_n) \in G$;*

$$\begin{aligned} x * y &= (x_1, \ldots, x_n) * (y_1, y_2, \ldots, y_n) \\ &= (px_1 + qy_1 \ (mod \ nI), \ldots, px_n + qy_n \ (mod \ nI)); \end{aligned}$$

*G is a neutrosophic row matrix groupoid of level one.*

We will illustrate this situation by some examples.

***Example 3.4.1:*** Let $G = \{(x_1, x_2, x_3, x_4, x_5, x_6) \mid x_i \in N(Q), *, (7, 13)\}$ be a neutrosophic row matrix groupoid of level one.

***Example 3.4.2:*** Let $G = \{(y_1, y_2, y_3, y_4, y_5, y_6, y_7, y_8, y_9, y_{10}) \mid y_i \in Z_9 I, *, (3, 7)\}$ be a neutrosophic row matrix groupoid of level one.

***Example 3.4.3:*** Let $P = \{(x_1, x_2, x_3) \mid x_i \in RI; *, (19, 3)\}$ be a neutrosophic row matrix groupoid of level one.

We see the row matrix neutrosophic groupoid given in examples 3.4.1 and 3.4.3 are of infinite order where as the groupoid in example 3.4.2 is of finite order.

**DEFINITION 3.4.2:** *Let*

$$G = \left\{ \begin{bmatrix} x_1 \\ x_2 \\ \vdots \\ x_m \end{bmatrix} \middle/ x_i \in ZI, *, (p, q); p \text{ and } q \text{ are distinct primes} \right\}.$$

*G is a column matrix neutrosophic groupoid of level one where * for any*

$$x = \begin{bmatrix} x_1 \\ x_2 \\ \vdots \\ x_m \end{bmatrix} \text{ and } y = \begin{bmatrix} y_1 \\ y_2 \\ \vdots \\ y_m \end{bmatrix}$$

*in G is defined by*

$$x*y = \begin{bmatrix} x_1 \\ x_2 \\ \vdots \\ x_m \end{bmatrix} * \begin{bmatrix} y_1 \\ y_2 \\ \vdots \\ y_m \end{bmatrix} = \begin{bmatrix} px_1 + qy_1 \\ px_2 + qy_2 \\ \vdots \\ px_m + qy_m \end{bmatrix}.$$

We will illustrate this situation by some simple examples.

***Example 3.4.4:*** Let

$$G = \left\{ \begin{bmatrix} x_1 \\ x_2 \\ x_3 \\ x_4 \\ x_5 \end{bmatrix} \middle| x_i \in N(Z), *, (3, 23) \right\}$$

be a column matrix neutrosophic groupoid of level one.

*Example 3.4.5:* Let

$$P = \left\{ \begin{bmatrix} x_1 \\ x_2 \\ x_3 \\ x_4 \\ x_5 \\ x_6 \\ x_7 \end{bmatrix} \mid x_i \in N(C), *, (19, 29) \right\}$$

be a column matrix neutrosophic groupoid of level one.

*Example 3.4.6:* Let

$$R = \left\{ \begin{bmatrix} x_1 \\ x_2 \\ x_3 \\ x_4 \end{bmatrix} \mid x_i \in Z_7 I, *, (3, 5) \right\}$$

be the column matrix neutrosophic groupoid of level one.

*Example 3.4.7:* Let

$$T = \left\{ \begin{bmatrix} x_1 \\ x_2 \\ x_3 \\ x_4 \\ x_5 \\ x_6 \\ x_7 \\ x_8 \end{bmatrix} \mid x_i \in N(Z_3), *, (1, 2) \right\}$$

be a neutrosophic column matrix groupoid of level one.

We see examples 3.4.4 and 3.4.5 groupoids are of infinite order where as groupoids in examples 3.4.6 and 3.4.7 are of finite order.

**DEFINITION 3.4.3:** *Let $M_{m \times n} = \{(m_{ij}) / m_{ij} \in N(Q); 1 \leq i \leq n, 1 \leq j \leq m, *, (p, q); p \text{ and } q \text{ are distinct primes}\}$ be the $n \times m$ matrix neutrosophic groupoid of level one.*

We will illustrate this situation by some examples.

*Example 3.4.8:* Let

$$G = \left\{ \begin{pmatrix} a_{11} & a_{12} & a_{13} & a_{14} \\ a_{21} & a_{22} & a_{23} & a_{24} \\ a_{31} & a_{32} & a_{33} & a_{34} \end{pmatrix} \mid a_{ij} \in ZI; \right.$$

$1 \leq i \leq 3; 1 \leq j \leq 4; *, (13, 43)\}$ be the $3 \times 4$ neutrosophic matrix groupoid of level one.

*Example 3.4.9:* Let

$$P = \left\{ \begin{pmatrix} a_{11} & a_{12} \\ a_{21} & a_{22} \\ a_{31} & a_{32} \\ a_{41} & a_{42} \\ a_{51} & a_{52} \\ a_{61} & a_{62} \\ a_{71} & a_{72} \end{pmatrix} \mid a_{ij} \in N(Z_5), *, (3, 2); 1 \leq i \leq 6, 1 \leq j \leq 2 \right\}$$

be a $6 \times 2$ neutrosophic matrix groupoid of level one. Clearly P is a finite groupoid.

*Example 3.4.10:* Let

$$W = \left\{ \begin{pmatrix} a_1 & a_2 & a_3 & a_4 & a_5 \\ a_6 & a_7 & a_8 & a_9 & a_{10} \\ a_{11} & a_{12} & a_{13} & a_{14} & a_{15} \\ a_{16} & a_{17} & a_{18} & a_{19} & a_{20} \end{pmatrix} \right| $$

$a_i \in N(Q)$, *, (31, 43), $1 \leq i \leq 20$} be the 4×5 neutrosophic matrix groupoid of level one of infinite order.

When n = m in the definition 3.4.3 we get square neutrosophic matrix groupoids of level one.

We will illustrate this situation by some simple examples.

*Example 3.4.11:* Let

$$P = \left\{ \begin{pmatrix} a_1 & a_2 & a_3 \\ a_4 & a_5 & a_6 \\ a_7 & a_8 & a_9 \end{pmatrix} \mid a_i \in N(R), *, (3, 17)\ 1 \leq i \leq 9 \right\}$$

be the $3 \times 3$ neutrosophic matrix groupoids of level one. Clearly P is an infinite groupoid.

*Example 3.4.12:* Let

$$M = \left\{ \begin{pmatrix} a & b \\ c & d \end{pmatrix} \mid a, b, c, d \in ZI, *, (23I, 53I) \right\}$$

be the $2 \times 2$ neutrosophic matrix groupoid of level one of infinite order.

*Example 3.4.13:* Let

$$W = \left\{ \begin{pmatrix} a_1 & a_2 & a_3 & a_4 \\ 0 & a_5 & a_6 & a_7 \\ 0 & 0 & a_8 & a_9 \\ 0 & 0 & 0 & a_{10} \end{pmatrix} \mid a_i \in N(Z_{19}I), *, (17, 2); \right.$$

$1 \leq i \leq 10$} be the $4 \times 4$ square matrix neutrosophic groupoid of finite order.

*Example 3.4.14:* Let R = {$20 \times 20$ matrices with entries from $Z_{12}I$, *, (3I, 11I)} be the $20 \times 20$ neutrosophic matrix groupoid of level one of finite order.

Now having seen level one groupoid we proceed on to define level two neutrosophic matrix groupoids.

**DEFINITION 3.4.4:** *Let $G = \{m \times n$ matrices with entries from $Z_nI$ or $N(Z_n)$ or $QI$ or $N(Q)$ or $ZI$ or $N(Z)$ or so on, *, $(p, q)$ such that $(p, q) = 1$; $p$ and $q$ need not necessarily be primes$\}$. $G$ is defined to be the $m \times n$ neutrosophic matrix groupoid of level two. If $m = 1$ then $G$ is a $1 \times n$ row neutrosophic matrix groupoid of level two. If $n = 1$ we get a $m \times 1$ column neutrosophic matrix groupoid of level two.*

*If $m = n$ then we get a square neutrosophic matrix groupoid of level two. If $m \neq n$ then $G$ is a rectangular neutrosophic matrix groupoid of level two.*

We will illustrate the definition by some examples.

*Example 3.4.15:* Let

$$G = \left\{ \begin{bmatrix} x_1 \\ x_2 \\ x_3 \\ x_4 \\ x_5 \\ x_6 \end{bmatrix} \mid x_i \in Z_7I, 1 \leq i \leq 6, *, (3, 4) \right\}$$

be a column neutrosophic matrix groupoid of level two. Clearly G is of finite order.

*Example 3.4.16:* Let $W = \{(y_1, y_2, y_3) \mid y_i \in N(R), 1 \leq i \leq 3, *, (8, 9)\}$ be the row neutrosophic matrix groupoid of level two of infinite order.

*Example 3.4.17:* Let

$$P = \left\{ \begin{pmatrix} a & b \\ c & d \end{pmatrix} \mid a, b, c, d \in N(C), *, (12, 25) \right\}$$

be the square neutrosophic matrix groupoid of infinite order.

*Example 3.4.18:* Let

$$M = \left\{ \begin{bmatrix} a_1 & a_2 & a_3 & a_4 & a_5 & a_6 \\ a_7 & a_8 & a_9 & a_{10} & a_{11} & a_{12} \end{bmatrix} \mid a_i \in N(Z_{29}); \right.$$

$1 \le i \le 12, *, (24, 25)\}$ be the rectangular neutrosophic matrix groupoid of finite order.

Now having seen level two groupoids we will just give the method by which level three, level four and level five matrix groupoids are built.

In the definition 3.4.4 if we take p, q such that $(p, q) = d \ne 1$ then we call G to be a neutrosophic matrix groupoid of level three. If in the definition $p = q = t$ is taken then we define the neutrosophic matrix groupoid to be a level four groupoid.

If instead of (p, q) we take one of p or q to be zero then we define those matrix neutrosophic groupoids to be level five matrix groupoids.

We will illustrate this situation by some examples.

*Example 3.4.19:* $G = \{(x_1, x_2, x_3, x_4) \mid x_i \in N(Z), *, (25, 35)\}$ be the row matrix neutrosophic groupoid of level three. Clearly G is of infinite order.

*Example 3.4.20:* Let

$$G = \left\{ \begin{bmatrix} x_1 \\ x_2 \\ x_3 \\ x_4 \\ x_5 \\ x_6 \\ x_7 \end{bmatrix} \mid x_i \in N(Z_{12}), *, (8, 6) \right\}$$

be the column matrix neutrosophic groupoid of level three of finite order.

*Example 3.4.21:* Let

$$G = \left\{ \begin{pmatrix} a & b \\ c & d \end{pmatrix} \middle| \; a, b, c, d \in N(R), *, (27, 45) \right\}$$

be the square matrix neutrosophic groupoid of level three of infinite order.

*Example 3.4.22:* Let P = {all 10 × 9 neutrosophic matrices with entries from $Z_{11}I$, *, (9, 6)} be the 10 × 9 neutrosophic matrix groupoid of level three of finite order.

*Example 3.4.23:* Let

$$G = \left\{ \begin{bmatrix} x_1 \\ x_2 \\ x_3 \\ x_4 \\ x_5 \\ x_6 \\ x_7 \\ x_8 \end{bmatrix} \middle| \; x_i \in N(Z), *, (5, 5); 1 \leq i \leq 8 \right\}$$

be a 8 × 1 neutrosophic column matrix groupoid of level four of infinite order.

*Example 3.4.24:* Let G = {[$x_1$, $x_2$, $x_3$]| $x_i \in N(Q)$, $1 \leq i \leq 3$, *, (17, 17)} be a 1 × 3 neutrosophic row matrix groupoid of level four. Clearly G is a groupoid of infinite order.

*Example 3.4.25:* Let G = {all 12 × 15 neutrosophic matrices with entries from $Z_6$, *, (3, 3)} be a rectangular neutrosophic matrix groupoid of level four of finite order.

*Example 3.4.26:* Let G = {8 × 8 neutrosophic matrices with entries from ZI, *, (19, 19)} be a square neutrosophic matrix groupoid of level four of infinite order.

Now we proceed onto define neutrosophic groupoids of level five.

*Example 3.4.27:* Let

$$G = \left\{ \begin{bmatrix} x_1 \\ x_2 \\ x_3 \\ \vdots \\ x_{25} \end{bmatrix} \mid x_i \in N(Z_4), *, (3, 0); 1 \leq i \leq 25 \right\}$$

be a neutrosophic column groupoid matrix of level five of finite order.

*Example 3.4.28:* Let G = {[$x_1, x_2, \ldots, x_7$] such that $x_i \in$ ZI, *, (0, 8), $1 \leq i \leq 7$} be a neutrosophic row groupoid of level five of infinite order.

*Example 3.4.29:* Let V = {3 × 3 square matrices with entries from N(C), *, (0, 7)} be a square neutrosophic complex groupoid of infinite order of level five.

*Example 3.4.30:* Let P = {20 × 5 rectangular matrices with entries from N($Z_8$), *, (6, 0)} be a rectangular neutrosophic groupoid of level five of finite order.
    Here we wish to state from now on wards we will not state to which level the groupoid belongs, it will be known from the context and further the properties defined will hold good for all levels of groupoids.

We will illustrate these by suitable examples.

**DEFINITION 3.4.5:** *Let (G, *, (p, q)) be a matrix neutrosophic groupoid. Suppose P $\subseteq$ G and if (P, *, (p, q)) is again a neutrosophic matrix groupoid then we call (P, *, (p, q)) to be a* *neutrosophic matrix subgroupoid. If G has no neutrosophic matrix subgroupoids then we define G to be a simple neutrosophic matrix groupoid.*

We will first illustrate this situation by some simple examples.

**Example 3.4.31:** Let G = {($x_1$, $x_2$, $x_3$), *, (3, 5), $x_i \in Z_7I$; $1 \le i \le 3$} be a neutrosophic row matrix groupoid of level one of finite order.
Choose P = {(x, x, x) | x ∈ $Z_7I$, *, (3, 5)} ⊆ G, P is a neutrosophic row matrix subgroupoid of G of level one.

**Example 3.4.32:** Let G = {set of all 3 × 3 neutrosophic matrices with entries from N(Z), *, (7, 8)} be a neutrosophic matrix groupoid of level two of infinite order. Choose P = {all 3 × 3 neutrosophic matrices with entries from 3ZI, *, (7, 8)} ⊆ G; P is a neutrosophic matrix subgroupoid of level two of infinite order.

**Example 3.4.33:** Let G = {all 2 × 9 neutrosophic matrices with entries from N(C), *, (9, 13)} be a neutrosophic matrix groupoid of level two.
Take P = {all 2 × 9 neutrosophic matrices with entries from N(R); *, (9, 13)} ⊆ G; P is a neutrosophic matrix subgroupoid of G of infinite order and of level two.

**Example 3.4.34:** Let G = {all 7 × 2 neutrosophic matrices with entries from RI, *, (2, 15)} be the neutrosophic matrix groupoid of level two. Take W = {all 7 × 2 neutrosophic matrices with entries from 3ZI, *, (2, 15)} ⊆ G, W is a neutrosophic matrix subgroupoid of G.

**Example 3.4.35:** Let G = {all 2 × 5 neutrosophic matrices with entries from ZI, *, (3, 15)} be a neutrosophic matrix groupoid of level three. Let W = {all 2 × 5 neutrosophic matrices with entries from 3ZI, *, (3, 15)} ⊆ G, W is a neutrosophic matrix subgroupoid of G.

**Example 3.4.36:** Let G = {all 1 × 5 row matrices with entries from $Z_{30}I$, *, (6, 15)} be a neutrosophic row matrix groupoid. Take W = {(a, a, a, a, a)| a ∈ $Z_{30}I$, *, (6, 15)} ⊆ G; W is a neutrosophic row matrix subgroupoid of G of level three.

**Example 3.4.37:** Let M = {all 10 × 8 neutrosophic matrices from QI, *, (17, 34)} be a neutrosophic matrix groupoid of level three. Take P = {all 10 × 8 neutrosophic matrices with entries from 13ZI, *, (17, 34)} ⊆ M; P is a neutrosophic matrix subgroupoid of M of level three.

**Example 3.4.38:** Let

$$G = \left\{ \begin{pmatrix} a & a & a \\ a & a & a \\ a & a & a \end{pmatrix} \mid a \in Z_7I, *, (3, 6) \right\}$$

be a neutrosophic matrix groupoid of level three; clearly G has no proper subgroupoids. Thus G is a simple neutrosophic matrix groupoid.

**Example 3.4.39:** Let G = {(a a a a a) where a ∈ $Z_{17}I$, (16, 14), *} be a neutrosophic matrix groupoid, clearly G is a simple matrix neutrosophic groupoid of level three.

**Example 3.4.40:** Let

$$G = \left\{ \begin{bmatrix} a & a & a \\ a & a & a \\ a & a & a \\ a & a & a \\ a & a & a \end{bmatrix} \text{ where } a \in ZI, *, (12, 38) \right\}$$

be a neutrosophic matrix groupoid of level three. G is not simple, G has infinitely many neutrosophic subgroupoids.

*Example 3.4.41:* Let

$$G = \left\{ \begin{bmatrix} a \\ a \\ a \\ a \\ a \\ a \\ a \\ a \end{bmatrix} \mid a \in QI, *, (9, 48) \right\}$$

be a neutrosophic matrix groupoid of level three. G has infinitely many neutrosophic matrix subgroupoids.

*Example 3.4.42:* Let

$$G = \left\{ \begin{bmatrix} a \\ a \\ a \\ a \\ a \\ a \\ a \\ a \end{bmatrix} \mid a \in Z_{43}I, *, (24, 40) \right\}$$

be a neutrosophic matrix groupoid of level three. G is a simple neutrosophic matrix groupoids.

*Example 3.4.43:* Let G = {all 5 × 5 neutrosophic matrices with entries from ZI, *, (13, 13)} be a neutrosophic matrix groupoid of level four.

Take W = {all 5 × 5 upper triangular neutrosophic matrices with entries from mZI, *, (13, 13)} $\subseteq$ G is a neutrosophic matrix subgroupoid of G for every m $\in Z^+ \setminus \{1\}$.

*Example 3.4.44:* Let G = [(a, a, a, a) | a $\in Z_5I$, *, (3, 3)} be a neutrosophic matrix groupoid of level four. G is a simple neutrosophic matrix groupoid of level four.

***Example 3.4.45:*** Let $G = \{(a, a, a, a) \mid a \in Z_6I, *, (4, 4)\}$ be a neutrosophic matrix groupoid of level four. Let $V = \{(a, a, a, a) \mid a \in \{0, 2I, 4I\}, *, (4, 4)\} \subseteq G$, V is a neutrosophic matrix subgroupoid of G.

***Example 3.4.46:*** Let $G = \left\{ \begin{pmatrix} a & a & a & a \\ a & a & a & a \end{pmatrix} \mid a \in Z_{23}I, *, (5, 5) \right\}$ be a neutrosophic matrix groupoid G is a simple neutrosophic matrix groupoid.

**THEOREM 3.4.1:** *Let*

$$G = \left\{ \begin{pmatrix} a & ... & a \\ a & ... & a \\ \vdots & & \vdots \\ a & ... & a \end{pmatrix} \mid a \in Z_pI, p \text{ a prime, } *, (t, t); 0 \leq t < p \right\}$$

*be any neutrosophic matrix groupoid. G is a simple neutrosophic matrix groupoid of level four.*

The proof is left as an exercise for the reader.

**THEOREM 3.4.2:** *Let*

$$G = \left\{ \begin{pmatrix} a & a & ... & a \\ a & a & ... & a \\ \vdots & \vdots & & \vdots \\ a & a & ... & a \end{pmatrix} \right/$$

$a \in Z_nI$ *or RI or QI or CI or ZI, *, (t, t); $0 \leq t < n$ if $a \in Z_nI$ otherwise $t \in Z^+ \setminus \{1\}$/ be a neutrosophic matrix groupoid. G is not a simple neutrosophic matrix groupoid of level four.*

The proof is left as an exercise for the reader.

**Example 3.4.47:** Let G = {all 7 × 7 square matrices with entries from $Z_{12}I$, *, (3, 0)} is a neutrosophic matrix groupoid of level five. Let P = {All 7 × 7 square matrices with entries from {0, 3I, 6I, 9I}, *, (3, 0)} ⊆ G, P is a neutrosophic matrix subgroupoid of level five.

**Example 3.4.48:** Let

$$G = \left\{ \begin{bmatrix} a \\ b \\ c \\ d \\ e \\ f \end{bmatrix} \mid a, b, c, d, e, f \in ZI, *, (0, 8) \right\}$$

be a neutrosophic matrix groupoid of level five.
Take

$$P = \left\{ \begin{bmatrix} a \\ a \\ a \\ a \\ a \\ a \\ a \end{bmatrix} \mid a \in 5ZI, *, (0, 8) \right\} \subseteq G,$$

P is a neutrosophic matrix subgroupoid of level five.

**Example 3.4.49:** Let G = {(a, a, a, a, a, a, a, a, a), a ∈ $Z_{43}$ I, *, (0, 24)} be a neutrosophic matrix groupoid of level five. G is a simple neutrosophic matrix groupoid of level five.

**Example 3.4.50:** Let G = {All 8 × 8 neutrosophic matrices with entries from N(Q), *, (21, 0)} be a neutrosophic matrix groupoid of level 5.
P = {all 8 × 8 neutrosophic matrices with entries from QI, *, (21, 0)} ⊆ G is a neutrosophic matrix subgroupoid of G of level five.

***Example 3.4.51:*** Let

$$G = \left\{ \begin{pmatrix} a & a & a & a \\ a & a & a & a \\ a & a & a & a \end{pmatrix} \mid a \in Z_{11}I, *, (90, 0) \right\}$$

be a neutrosophic matrix groupoid of level five. G is a simple neutrosophic matrix groupoid.

In view of this we have the following theorem the proof of which is left as an exercise for the reader.

**THEOREM 3.4.3:** *Let*

$$G = \left\{ \begin{pmatrix} a & a & ... & a \\ a & a & ... & a \\ \vdots & \vdots & & \vdots \\ a & a & ... & a \end{pmatrix} \right.$$

*such that $a \in Z_pI$, p a prime, *, (0, t) $0 \leq t < p$} be a neutrosophic matrix groupoid of level five. G is a simple neutrosophic matrix groupoid of level five.*

This proof is also left as an exercise for the reader.

Now we give yet another theorem which states as follows.

**THEOREM 3.4.4:** *Let $G = \{m \times n$ matrices with entries from $Z_nI$, *, (t, 0)} (the entries can be from $N(Z_n)$ or ZI or N(z) or N(Q) or N(R) or N(C) or QI or RI or CI). G is not a simple neutrosophic matrix groupoid of level five.*

The proof is left as an exercise for the reader.

As in case of other groupoids we can in case of neutrosophic matrix groupoids also define the notion of ideal (right and left), normal subgroupoids, normal groupoids, Smarandache groupoids, groupoids satisfying special identities. However we give some interesting properties satisfied by these neutrosophic matrix groupoids and a few examples to substantiate them. For more refer [20].

**THEOREM 3.4.5:** *Let $G = \{m \times n$ matrices with all entries from $Z_nI$, *, $(t, u)$; $t$ and $u$ primes with $t + u = 1$ or $I \pmod{n}\}$. Then $G$ is a simple neutrosophic matrix groupoid.*

*Proof:* Let

$$P = \left\{ \begin{pmatrix} a & a & \ldots & a \\ a & a & \ldots & a \\ \vdots & \vdots & & \vdots \\ a & a & \ldots & a \end{pmatrix} \mid a \in Z_nI, *, (t, u), t + u = 1 \text{ or } I \bmod n \right\}.$$

It is easily verified using simple number theoretic techniques.

**COROLLARY 3.4.1:** *Let $G$ be as in theorem 3.4.5 then $G$ is a neutrosophic groupoid which is an idempotent groupoid.*

*Example 3.4.52:* Let

$$G = \left\{ \begin{pmatrix} a & a & a \\ a & a & a \\ a & a & a \end{pmatrix} \mid a \in Z_{10}I, *, (3, 7) \right\};$$

$G$ is an idempotent (neutrosophic matrix) groupoid.

It is further verified that $G$ has no left or right ideals.

**THEOREM 3.4.6:** *Let $G = \{$all $m \times n$ neutrosophic matrices with entries from $Z_nI$, *, $(t, t)$; $0 \leq t < n\}$ be a neutrosophic matrix groupoid. $G$ is a neutrosophic matrix P-groupoid.*

*Proof:* Let $A = (a_{ij})$ and $B = (b_{ij})$ be in $G$. To show
$$(A * B) * A = A * (B * A).$$
Consider
$$\begin{aligned}(A * B) * A &= (ta_{ij} + tb_{ij}) * A \\ &= (t^2 a_{ij} + t^2 b_{ij} + ta_{ij}) \end{aligned} \quad \text{I}$$
Consider
$$\begin{aligned}A * (B * A) &= (a_{ij}) * (tb_{ij} + ta_{ij}) \\ &= (ta_{ij} + t^2 b_{ij} + t^2 a_{ij}) \end{aligned} \quad \text{II}$$

I and II are identical hence G is a neutrosophic matrix P-groupoid.

*Example 3.4.53:* Let

$$G = \left\{ \begin{pmatrix} a & b \\ c & d \end{pmatrix} \mid a, b, c, d \in Z_6I, *, (5, 5) \right\}$$

be a neutrosophic matrix groupoid. Clearly G is a neutrosophic matrix P-groupoid.

Let

$$A = \begin{pmatrix} a & b \\ c & d \end{pmatrix}$$

and

$$B = \begin{pmatrix} x & y \\ z & w \end{pmatrix}$$

be in G.
To show

$(A * B) * A = A * (B * A)$

$$= \left[ \begin{pmatrix} a & b \\ c & d \end{pmatrix} * \begin{pmatrix} x & y \\ z & w \end{pmatrix} \right] * \begin{pmatrix} a & b \\ c & d \end{pmatrix}$$

$$= \begin{pmatrix} 5a + 5x (\mod 6) & 5b + 5y (\mod 6) \\ 5c + 5z (\mod 6) & 5d + 5w (\mod 6) \end{pmatrix} * \begin{pmatrix} a & b \\ c & d \end{pmatrix}$$

$$= \begin{pmatrix} 7a + 7x + 5a (\mod 6) & 7b + 7y + 5b (\mod 6) \\ 7c + 7z + 5c (\mod 6) & 7d + 7w + 5d (\mod 6) \end{pmatrix}$$

$$= \begin{pmatrix} x & y \\ z & w \end{pmatrix} \qquad\qquad \text{I}$$

$$A * \begin{pmatrix} x & y \\ z & w \end{pmatrix} * \begin{pmatrix} a & b \\ c & d \end{pmatrix}$$

$$= A * \begin{pmatrix} 5x + 5a(\bmod 6) & 5y + 5b(\bmod 6) \\ 5z + 5c(\bmod 6) & 5w + 5d(\bmod 6) \end{pmatrix}$$

$$= \begin{pmatrix} 5a + 7x + 7a(\bmod 6) & 5b + 7y + 7b(\bmod 6) \\ 5c + 7z + 7c(\bmod 6) & 5d + 7w + 7d(\bmod 6) \end{pmatrix}$$

$$= \begin{pmatrix} x & y \\ z & w \end{pmatrix} \qquad \text{II}$$

I and II are equal to $\begin{pmatrix} x & y \\ z & w \end{pmatrix}$. Thus G is a neutrosophic matrix P- groupoid.

**THEOREM 3.4.7:** *Let $G = \{m \times n$ neutrosophic matrices with entries from $Z_n$ I; n not a prime, *, (t, t); $1 \leq t < n\}$ be a neutrosophic matrix groupoid. G is an alternative neutrosophic matrix groupoid if and only if $t^2 \equiv t \pmod{n}$.*

*Proof:* To prove the theorem it is sufficient if we show for $A = (a_{ij})$ and $B = (b_{ij})$ in G
$$(A * B) * B = A * (B * B).$$
Consider
$(A * B) * B \quad = \quad (tb_{ij} + tb_{ij})*B$
$\qquad\qquad\qquad = \quad (ta_{ij} + t^2 b_{ij} + tb_{ij}) \pmod{n}$
$\qquad\qquad\qquad = \quad (ta_{ij} + tb_{ij} + tb_{ij}) \qquad\qquad\qquad$ I
Take
$A * (B * B) \quad = \quad (a_{ij})* (tb_{ij} + tb_{ij})$
$\qquad\qquad\qquad = \quad (ta_{ij} + t^2 b_{ij} + t^2 b_{ij})$
$\qquad\qquad\qquad = \quad (ta_{ij} + tb_{ij} + tb_{ij}) \qquad\qquad\qquad$ II

I and II are identical if and only if $t^2 \equiv t \pmod{n}$.
    Thus G is a neutrosophic alternative matrix groupoid.

Note if n is a prime then $t^2 = t$ (mod n) is impossible. In view of this we can say the neutrosophic matrix groupoid described in the above theorem is not an alternative groupoid if n is a prime.

**THEOREM 3.4.8:** *The neutrosophic matrix groupoid $G = \{m \times p$ matrices with entries from $Z_nI$; *, (0, t) (or (t, 0)); $1 \leq t < n\}$ is an alternative neutrosophic matrix groupoid and a neutrosophic matrix P-groupoid if and only if $t^2 \equiv t$ (mod n).*

The reader is expected to prove the above theorem.

*Example 3.4.54:* Let $G = \{5 \times 8$ neutrosophic matrices with entries from $Z_6I$, *, (4, 0)$\}$ be a neutrosophic matrix groupoid. G is easily verified to be both a neutrosophic matrix P-groupoid as well as neutrosophic matrix alternative groupoid.

We see also $G = \{m \times n$, neutrosophic matrices with entries from $Z_6I$ with (3, 0)$\}$ (or (0, 3) or (0, 4) or (4, 0)) are all both P-neutrosophic matrix groupoids as well as neutrosophic matrix alternative groupoids.

*Example 3.4.55:* Let $G = \{8 \times 9$ neutrosophic matrices with entries from $Z_{10}I$, *, (1, 5)$\}$ be a neutrosophic matrix groupoid. It is easy to verify G is a Smarandache neutrosophic matrix groupoid.

*Example 3.4.56:* Let $G = \{$all $3 \times 4$ neutrosophic matrices with entries from $Z_8I$, *, (2, 6)$\}$ be a neutrosophic matrix groupoid. $P = \{$all $3 \times 4$ neutrosophic matrices with entries from $\{0, 2I, 4I, 6I\} \subseteq Z_8I$, *, (2, 6)$\} \subseteq G$; P is a Smarandache normal matrix groupoid of G.

*Example 3.4.57:* Let $G = \{$all $2 \times 2$ square matrices with entries from $Z_5I$, *, (3, 3)$\}$ be the neutrosophic matrix groupoid. Clearly G is Smarandache inner commutative groupoid.

*Example 3.4.58:* Let $G = \{$all $1 \times 10$ neutrosophic matrices with entries from $Z_{12}I$, *, (3, 9)$\}$ be a neutrosophic matrix groupoid.

It is easily verified that G is a Smarandache Moufang groupoid and is not a Smarandache strong Moufang groupoid.

**Example 3.4.59:** Let G = {8 × 9 neutrosophic matrices with entries from $Z_{12}I$, *, (3, 4)} be a neutrosophic matrix groupoid. G is a Smarandache strong neutrosophic Bol groupoid.

**Example 3.4.60:** Let G = {all 18 × 1 neutrosophic matrices with entries from $Z_4I$, *, (2, 3)} be a neutrosophic column matrix groupoid.
It is easily verified G is a Smarandache Bol matrix neutrosophic groupoid and is not a Smarandache strong Bol neutrosophic matrix groupoid.

**Example 3.4.61:** Let G = {all 10×6 neutrosophic matrices with entries from $Z_nI$, *, (0, 2); n = 2m} G is a Smarandache neutrosophic matrix groupoid of level five.

Several results true in case of general groupoids built using $Z_n$ can be derived for these neutrosophic matrix groupoids built using $Z_nI$.

## 3.5 Neutrosophic Interval Groupoids

In this section we just introduce a new class of groupoids called neutrosophic interval groupoids.
We will give the basic notations.
$$o(Z_nI) = \{[0, a_i] \mid a_i \in Z_nI\}$$
$$o(N(Z_n)) = \{[0, a + bI] \mid a, b \in Z_n\}$$
$$o(Q^+I) = \{[0, a] \mid a \in Q^+I; a = 0$$
is also allowed
$$o(N(Q^+ \cup \{0\})) = \{[0, a + bI] \mid a, b \in Q^+ \cup \{0\}\}$$
$$o(Z^+I) = \{[0, a] \mid a \in Z^+I \cup \{0\}\}$$
$$o\{N(Z^+I)) = \{[0, a + bI] \mid a, b \in Z^+ \cup \{0\}\}$$
$$o(R^+I) = \{[0, a] \mid a \in R^+I \cup \{0\}\}$$
$$o(N(R^+ \cup \{0\})) = \{[0, a + bI] \mid a, b \in R^+ \cup \{0\}\}$$
$$o(C^+I) = \{[0, a] \mid a \in C^+I \cup \{0\}\}$$
and

$$o\{N(C^+)) = \{a + bI \,/\, a, b \in C^+ \cup \{0\}\}$$

are special neutrosophic intervals. We only take positive values and always the least value is zero.

**DEFINITION 3.5.1:** *Let $o(Z_nI)$ be the collection of intervals. Define * for any two intervals $[0, a]$, $[0, b]$ in $o(Z_nI)$ as $[0, a] * [0, b] = [0, ta + pb(\bmod n)]$ where $t, p \in Z_n \setminus \{0\}$ such that $t$ and $p$ are primes with $(t, p) = 1$. We see $G = \{o(Z_n, I), *, (t, p)\}$ is a interval neutrosophic groupoid of level one.*

We will illustrate this by some examples.

***Example 3.5.1:*** Let $G = \{[0, a] \mid a \in Z_pI, *, (7, 3)\}$ be a neutrosophic interval groupoid of level one.

In the definition 3.5.1 we can replace $o(Z_nI)$ by $o(N(Z_n))$ or $o(Z^+I)$ or $o(N(Z^+))$ or $o(R^+I)$ or $o(N(R^+))$ or $o(N(Q^+I))$ or so on. Still we get only neutrosophic interval groupoids.

We will illustrate all this situation only by some examples.

***Example 3.5.2:*** Let $G = \{[0, a] \mid a \in Z^+I, *, (19, 17)\}$ be the neutrosophic interval groupoid of level one.

***Example 3.5.3:*** Let $G = [0, a + bI] \mid a, b \in Q^+ \cup \{0\}, *, (3, 23)\}$ be the neutrosophic interval groupoid of level one.

***Example 3.5.4:*** Let $P = \{[0, a + bI] \mid a, b \in R^+ \cup \{0\}, *, (13, 11)\}$ be a neutrosophic interval groupoid of level one.

***Example 3.5.5:*** Let $W = \{[0, a] \mid a \in C^+I, *, (3, 29)\}$ be the complex neutrosophic interval groupoid of level one.

Now we just indicate how neutrosophic interval groupoids of different levels are defined.
In the definition 3.5.1 if we take instead of $t, u$; $t$ and $u$ prime, $t, u$ such that $(t, u) = 1$ but $t$ and $u$ need not be primes then we get neutrosophic interval groupoids of level two.

We will first illustrate level two neutrosophic interval groupoids.

*Example 3.5.6:* Let K = {[0, a] | a ∈ Q$^+$I, *, (15, 17)} be the neutrosophic interval groupoid of level two.

*Example 3.5.7:* Let T = {[0, a + bI] | a, b ∈ R$^+$ ∪ {0}, *, (18, 25)} be the neutrosophic interval groupoid of level two.

*Example 3.5.8:* Let S = {[0, a] | a ∈ Z$^+$I, *, (27, 40)} be the neutrosophic interval groupoid of level two.

Now if in the definition 3.5.1 we take t, u such that (t, u) = d ≠ 1 then we define those groupoids to be neutrosophic interval groupoids of level three.

We will illustrate this by some examples.

*Example 3.5.9:* Let T = {[0, a] | a ∈ $Z_{24}$I; *, (6, 14)} be a neutrosophic interval groupoid of level three.

*Example 3.5.10:* Let M = {[0, a + bI] | a, b ∈ R$^+$ ∪ {0}, *, (24, 38)} be the neutrosophic interval groupoid of level three.

*Example 3.5.11:* Let P = {[0, a] | a ∈ (C$^+$I), *, (25, 45)} be the neutrosophic interval groupoid of level three.

*Example 3.5.12:* Let B = {[0, a + bI] | a, b ∈ C$^+$ ∪ {0}} *, (8, 18)} be a neutrosophic matrix groupoid of level three.

If in the definition 3.5.1 we take t = u then we call the interval neutrosophic matrix groupoid to be a level four interval neutrosophic matrix groupoid.

We will illustrate this situation by some examples.

*Example 3.5.13:* Let T = {[0, a] | a ∈ $Z_{28}$I, *, (15, 15)} be a neutrosophic interval groupoid of level four.

*Example 3.5.14:* Let $T = \{[0, a + bI] \mid a, b \in R^+ \cup \{0\}, (7, 7)\}$ be a neutrosophic interval groupoid of level four.

*Example 3.5.15:* Let $W = \{[0, a] \mid a \in C^+I; (9, 9)\}$ be the neutrosophic interval groupoid of level four.

*Example 3.5.16:* Let $P = \{[0, a + bI] \mid a, b \in C^+ \cup \{0\}, (20, 20)\}$ be the neutrosophic interval groupoid of level four.

Now if in the definition 3.5.1 we take $u = 0$ or $t = 0$ then we get neutrosophic interval groupoid of level five.

We will illustrate this by examples.

*Example 3.5.17:* Let $S = \{[0, a] \mid a \in Z_9I; *, (0, 8)\}$ be the neutrosophic interval groupoid of level five.

*Example 3.5.18:* Let $S = \{[0, a + bI] \mid a, b \in N(Z_{28}), *, (9, 0)\}$ be the neutrosophic interval groupoid of level five.

We will not distinguish between the levels of neutrosophic interval groupoids or the neutrosophic sets which are being used to build these groupoids as it clear by the structure.

**DEFINITION 3.5.2:** *Let $G = \{[0, a] \mid [0, a] \in o(Z_nI), *, (p, q)\}$ be a neutrosophic interval groupoid.*
*Let $P \subseteq o(Z_nI)$ such that $(P, * (p, q))$ is a neutrosophic interval groupoid. We call P to be a neutrosophic interval subgroupoid of G.*

We will illustrate this by some examples.

*Example 3.5.19:* Let $G = \{[0, a] \mid a \in Z_8I, *, (3, 0)\}$ be a neutrosophic interval groupoid of finite order. $P = \{[0, a] \mid a \in \{0, 2I, 4I, 6I\}, *, (3, 0)\}$, G is a neutrosophic interval subgroupoid of G.

*Example 3.5.20:* Let T = {[0, a + bI] | a + bI ∈ N($Q^+$), *, (3, 8)} be a neutrosophic interval groupoid. P = {[0, bI] | bI ∈ $Q^+$I, *, (3, 8)} ⊆ T; is a neutrosophic interval subgroupoid of G.

*Example 3.5.21:* Let G = {[0, a + bI] | a, bI ∈ N($Z^+$I), *, (20, 0)} be a neutrosophic interval groupoid. Choose M = {[0, a] | a ∈ $Z^+$I, *, (20, 0)} ⊆ G, M is a neutrosophic interval subgroupoid of G.

*Example 3.5.22:* Let G = {[0, a] | a ∈ $R^+$I ∪ {0}, *, (18, 24)} be a neutrosophic interval groupoid. Let V = {[0, a] | a ∈ $Q^+$I ∪ {0}, *, (18, 24)} ⊆ G; V is a neutrosophic interval subgroupoid of G.

Now it may so happen G be a neutrosophic interval groupoid, P may be a proper subset of G which is just an interval groupoid and not a neutrosophic interval groupoid. We call P to be a pseudo neutrosophic interval subgroupoid of G.

If G has no pseudo neutrosophic interval subgroupoids we call G to be a pseudo simple neutrosophic interval groupoid.

We will illustrate this situation by some examples.

*Example 3.5.23:* Let G = {[0, a + bI], a, b ∈ N($Z_{12}$), *, (3, 4)} be a neutrosophic interval groupoid. Take W = {[0, a] | a ∈ $Z_{12}$, *, (3, 4)} ⊆ G is a real interval groupoid; W is called as pseudo neutrosophic interval subgroupoid of G.

*Example 3.5.24:* Let G = {[0, a] | a ∈ $Z^+$I ∪ {0}, *, (3, 24)} be a neutrosophic interval groupoid we see G has neutrosophic interval subgroupoids but has no pseudo neutrosophic interval subgroupoid. Thus G is a simple pseudo neutrosophic groupoid.

The following theorem gives a very large class of pseudo simple neutrosophic groupoids.

**THEOREM 3.5.1:** *Let G = {[0, a] | a ∈ $Z_n$I or $Z^+$I or $Q^+$I or $R^+$I or $C^+$I, *, (t, u)} be a neutrosophic interval groupoid. G is a pseudo neutrosophic interval groupoid.*

The proof is left as an exercise to the reader.

**THEOREM 3.5.2:** *Let $G = \{[0, a + bI] / a, b \in N(Q^+)$ or $N(Z_n)$ or $N(R^+)$ or $N(Z^+)$ or $N(C^+)$, *, $(t, u)\}$ be a neutrosophic interval groupoid. G has pseudo neutrosophic interval subgroupoids that G is not a simple pseudo neutrosophic interval groupoid.*

This proof is also left as an exercise to the reader.

We call a neutrosophic interval groupoid G to be doubly simple neutrosophic interval groupoid if G has no neutrosophic interval subgroupoid as well as G has no pseudo neutrosophic interval subgroupoid.

*Example 3.5.25:* Let $G = \{[0, a] \mid a \in Z_5I, *, (3, 4)\}$ be a neutrosophic interval groupoid. Clearly G has no neutrosophic interval subgroupoids. Further G has no pseudo neutrosophic interval subgroupoid. Thus G is a doubly simple neutrosophic interval groupoid.

We will illustrate this situation by a theorem which guarantees the existence of a large class of doubly simple neutrosophic interval groupoids.

**THEOREM 3.5.3:** *Let $G = \{[0, a] \mid a \in Z_pI, p$ a prime, *, $(t, u); 0 < t, u \leq p - 1\}$ be a neutrosophic interval groupoid G is a doubly simple neutrosophic interval groupoid.*

The proof is simple and is left as an exercise to the reader.
    As in case of general groupoids we can in case of neutrosophic interval groupoids also define the notion of right ideal, left ideal, ideal, normal, subgroupoids, normal neutrosophic interval groupoids, Smarandache neutrosophic interval groupoids and all concepts like neutrosophic interval groupoids satisfying special identities like Moufang, Bol, alternative. Idempotent neutrosophic interval groupoids, neutrosophic interval P-groupoids and their Smarandache analogue and Smarandache strong analogue.

We will however give some examples of these concepts and a few interesting theorems.

***Example 3.5.26:*** Let G = {[0, a] | a ∈ $Z_4I$, *, (2, 3)} be a neutrosophic interval groupoid.

This groupoid has left ideals given by P = {[0, a] | a ∈ {0, 2I}, *, (2, 3)} ⊆ G and T = {[0, a] | a ∈ {I, 3I}, *, (2, 3)} ⊆ G. Clearly P and T are not neutrosophic interval right ideals of G.

However P and T right ideals of G′ = {[0, a] | a ∈ $Z_4I$, *, (3, 2)}.

In view of this we have the following theorem.

**THEOREM 3.5.4:** *Let G = {[0, a] | a ∈ $Z_nI$, *, (t, u)} and G′ = {[0, a] | a ∈ $Z_nI$, *, (u, t)} be neutrosophic interval groupoids. P is a left neutrosophic interval ideal of G if and only if P is a right neutrosophic interval ideal of G.*

The proof is got by simple number theoretic computations and the simple translation of rows to columns. The reader is expected to prove the same.

***Example 3.5.27:*** Let G = {[0, a] | a ∈ $Z_{10}I$, *, (3, 7)} be a neutrosophic interval groupoid G has no left or right ideals.

**THEOREM 3.5.5:** *Let G = {[0, a] | a ∈ $Z_pI$; p a prime, (t, t); t < p, *} be a neutrosophic interval groupoid. G is a normal groupoid.*

Proof is left as an exercise for the reader.

**THEOREM 3.5.6:** *Let G = {[0, a] | a ∈ $Z_nI$, *, (t, t)} be a neutrosophic interval groupoid. G is a P-groupoid.*

This proof is also left for the reader as a simple number theoretic exercise.

**THEOREM 3.5.7:** *Let $G = \{[0, a] \mid a \in Z_pI$, $p$ a prime, $*$, $(t, t)$; $1 \leq t < p\}$ is a neutrosophic groupoid. $G$ is not an alternative interval groupoid.*

*Proof:* We have to prove for any $x, y, \in G$;
$$(x * y) * y \neq x * (y * y).$$
Take $x = [0, a]$ and $y = [0, b]$. Consider

$$\begin{aligned}(x * y) * y &= \{[0, a] * [0, b]) * [0, b] \\ &= [0, ta + tb] * [0, b] \\ &= [0, t^2a + t^2b + tb \text{ (mod } p)]\end{aligned} \quad \text{I}$$

Consider
$$\begin{aligned}x * (y * y) &= [0, a] * ([0, b] * [0, b]) \\ &= [0, a] * ([0, tb + tb]) \\ &= [0, at + t^2b + t^2b \text{ (mod } p)]\end{aligned} \quad \text{II}$$

Clearly I and II are different. Hence $G$ is not an alternative neutrosophic interval groupoid.

**THEOREM 3.5.8:** *Let $G = \{[0, a] \mid a \in Z_nI$, $(t, t)$; $1 < t < n$, $n$ not a prime, $*\}$ be a neutrosophic interval groupoid. $G$ is an alternative neutrosophic interval groupoid if and only if $t^2 \equiv t \text{(mod } n)$.*

The proof is left to the reader.

Now we will illustrate these situations by some simple examples.

***Example 3.5.28:*** Let $G = \{[0, a] \mid a \in Z_{24}I, *, (11, 11)\}$ be a neutrosophic interval groupoid which is alternative.
    Suppose in $G$ we take instead of $(11, 11)$; $(10, 10)$ then $G = \{[0, a] \mid a \in Z_{24}I, *, (10, 10)\}$ is not an alternative groupoid.

***Example 3.5.29:*** Let $G = \{[0, a] \mid a \in Z_{20}I, *, (5, 5)\}$ be a neutrosophic interval groupoid, $G$ is clearly a neutrosophic alternative interval groupoid.

**THEOREM 3.5.9:** *Let $G = \{[0, a] \mid a \in Z_nI, *, (t, 0)\}$ be a neutrosophic interval groupoid. G is a P-groupoid and alternative interval groupoid if and only if $t^2 = t(\bmod n)$.*

The proof is left as a simple exercise.

We will illustrate this by some simple examples.

*Example 3.5.30:* Let $G = \{[0, a] \mid a \in Z_{40}, *, (16, 0)\}$ be a neutrosophic interval groupoid G is an alternative interval groupoid and G is also a neutrosophic interval P-groupoid.

*Example 3.5.31:* Let $G = \{[0, a] \mid a \in Z_{29}, *, (t, 0)\}$ be an interval neutrosophic groupoid. G is not a P-groupoid for any t; $1 < t < 29$.

*Example 3.5.32:* Let $G = \{[0, a] \mid a \in Z_{10}I, *, (5, 6)\}$ be a neutrosophic interval groupoid. G is a Smarandache strong neutrosophic interval Moufang groupoid.

*Example 3.5.33:* Let $G = \{[0, a] \mid a \in Z_{12}I, (3, 9), *\}$ be a neutrosophic interval groupoid; G is a Smarandache Moufang neutrosophic interval groupoid which is not a Smarandache strong Moufang neutrosophic interval groupoid.

*Example 3.5.34:* Let $G = \{[0, a], *, (3, 4), a \in Z_{12}I\}$ be a neutrosophic interval groupoid. G is a Smaradache strong Bol neutrosophic interval groupoid.

*Example 3.5.35:* Let $G = \{[0, a] \mid a \in Z_4I, *, (2, 3)\}$ be a neutrosophic interval groupoid. G is a Smaradanche neutrosophic interval Bol groupoid but is not a Smarandahce strong neutrosophic interval Bol groupoid.

*Example 3.5.36:* Let $G = \{[0, a] \mid a \in Z_6I, (4, 3), *\}$ be a neutrosophic interval groupoid. G is a Smarandache strong neutrosophic interval P-groupoid.

*Example 3.5.37:* Let G = {[0, a], *, a ∈ $Z_{12}I$, (5, 10)} be a neutrosophic interval groupoid. G is only a Smarandache neutrosophic interval P-groupoid and is not a Smarandache strong neutrosophic interval P-groupoid

*Example 3.5.38:* Let G = {[0, a] | a ∈ $Z_{14}I$, *, (97, 8)} be a neutrosophic interval groupoid. G is a Smarandache neutrosophic strong interval alternative groupoid.

Several properties enjoyed by interval groupoids can also be derived for neutrosophic interval groupoids with appropriate modifications.

When we define homomorphism of neutrosophic interval groupoids it is essential that I is mapped only on to I for during homomorphisms it is impossible for the indeterminate I to be changed to real, this indeterminate must remain as an indeterminate only.

### 3.6 Neutrosophic Interval Matrix Groupoids

In this section we will define for the first time the new notion of neutrosophic interval matrix groupoids and describe a few of the properties associated with them. We will be using only the notations given in section 3.5.

We will first give some essential notations.

Let $(a_1, a_2, \ldots, a_n) = X$ be such that $a_i \in o(Z_nI)$ or $o(N(Z_n))$, $o(Z^+I)$, $o(N(Q^+ \cup \{0\}\} \; o(Z^+I)$, $o(N(Z^+I))$ and $o(R^+I)$ so on. X is known as the neutrosophic interval row matrix.

$$Y = \begin{bmatrix} x_1 \\ x_2 \\ \vdots \\ x_m \end{bmatrix}$$

where $x_i$'s are neutrosophic intervals from $o(Z_nI)$, $o(N(Z_n))$, $o(Q^+I)$, $o(N(Q^+ \cup \{0\}))$ and so on. Y will in general be known as the neutrosophic interval column matrix.

Let

$$M = \begin{bmatrix} a_{11} & \cdots & a_{1m} \\ \vdots & & \vdots \\ a_{n1} & \cdots & a_{nm} \end{bmatrix}$$

be the neutrosophic interval matrix. The entries in M are from $o(Z_nI)$, $o(N(Z_n))$, $o(Q^+I)$ and so on.

If m = n then we call M to be a neutrosophic interval square matrix.

Now we will define groupoids using these neutrosophic interval matrices. We can have five levels of interval matrix groupoids we will define them without mentioning the levels for from the very context one can easily understand to which level the groupoid belongs to. We will now make a formal definition and give examples of them.

**DEFINITION 3.6.1:** *Let $G = \{(a_1, a_2, ..., a_n) \mid a_i = [0, x_i] \in o(Z_nI\}$ or $o(N(Z_n))$, $o(Z^+I)$, $o(N(Z^+ \cup 0))$ and so on, *, (p, q); q, p $\in Z_n$ or $Z^+ \cup \{0\}\}$ where for $X = (a_1, ..., a_n)$ and $Y = (b_1, ..., b_n)$ in G. $X * Y = ([0, x_1], ..., [0, x_n]) * ([0, y_1], ..., [0, y_n]) = ([0, px_1 + qy_1], [0, px_2 + qy_2], ..., [0, px_n + qy_n]) \in G$.*

It is easily verified G is a groupoid and G is defined as a row neutrosophic interval matrix groupoid.

We will illustrate this by some examples.

*Example 3.6.1:* Let $G = \{([0, a_1], [0, a_2], ..., [0, a_8]) \mid a_i \in Z_9I, *, (3, 8)\}$ be a neutrosophic $1 \times 8$ row interval matrix groupoid built using $o(Z_9I)$.

*Example 3.6.2:* Let $G = \{([0, a_1], [0, a_2], ..., [0, a_{25}]) \mid a_i \in N(Q), * (24, 0)\}$ be a neutrosophic row interval matrix groupoid of level five.

***Example 3.6.3:*** Let $G = \{([0, a_1], [0, a_2], …, [0, a_{14}]) \mid a_i \in N(Z^+), i = 1, 2, …, 14; *, (9, 17)\}$ be a neutrosophic row interval matrix groupoid of level two.

***Example 3.6.4:*** Let $G = \{([0, a_1], [0, a_2], [0, a_3]) \mid a_i \in N(Q^+), (5, 19), *\}$ be a neutrosophic row interval matrix groupoid of level one.

***Example 3.6.5:*** Let $G = \{([0, a_1], [0, a_2], …, [0, a_{140}]), *, a_i \in N(R^+), (4, 28), *\}$ be a neutrosophic row interval matrix groupoid of level three.

***Example 3.6.6:*** Let $G = \{([0, a_1], [0, a_2], …, [0, a_{48}]) \mid a_i \in Z^+I, (3, 3), *\}$ be a neutrosophic row interval matrix groupoid of level four.

***Example 3.6.7:*** Let $G = \{([0, a_1], [0, a_2], …, [0, a_{18}]) \mid a_i \in Z_{25}I, 1 \leq i \leq 18, *, (24, 8)\}$ be a neutrosophic row interval matrix groupoid.

**DEFINITION 3.6.2:** *Let*

$$G = \left\{ \begin{bmatrix} x_1 \\ x_2 \\ \vdots \\ x_n \end{bmatrix} \right/$$

$x_i = [0, a_i], a_i \in o(Z_nI)$ *or* $o(N(Z_n))$ *or* $o(Z^+I)$ *or* $o(Q^+I)$ *and so on.* $*, (p, q); p, q \in Z_n$ *or* $Z^+ \cup \{0\}\}$ *be such that for any*

$$x = \begin{bmatrix} x_1 \\ x_2 \\ \vdots \\ x_n \end{bmatrix} \text{ and } y = \begin{bmatrix} y_1 \\ y_2 \\ \vdots \\ y_n \end{bmatrix}$$

*where* $x_i = [0, a_i]$ *and* $y_i = [0, b_i]; a_i, b_i \in o(Z_nI)$ *or* $o(N(Z_n))$ *or so on.*

$$x * y = \begin{bmatrix} x_1 \\ x_2 \\ \vdots \\ x_n \end{bmatrix} * \begin{bmatrix} y_1 \\ y_2 \\ \vdots \\ y_n \end{bmatrix} = \begin{bmatrix} 0, pa_1 + qb_1 \\ 0, pa_2 + qb_2 \\ \vdots \\ 0, pa_n + qb_n \end{bmatrix}.$$

It is easily verified G is a neutrosophic column matrix interval groupoid.

We will illustrate this situation by some examples.

*Example 3.6.8:* Let

$$G = \left\{ \begin{bmatrix} x_1 \\ x_2 \\ x_3 \\ x_4 \\ x_5 \end{bmatrix} \mid x_i = [0, a_i], a_i \in Z_8 I; 1 \leq i \leq 5; *, (7, 1) \right\}$$

be a neutrosophic interval column matrix groupoid of level one.

*Example 3.6.9:* Let

$$G = \left\{ \begin{bmatrix} x_1 \\ x_2 \\ x_3 \\ x_4 \\ x_5 \\ x_6 \\ x_7 \end{bmatrix} \mid x_i = [0, a_i]; a_i \in N(Q^+), *, (8, 92), 1 \leq i \leq 7 \right\}$$

be a neutrosophic interval column matrix groupoid level three.

*Example 3.6.10:* Let

$$G = \left\{ \begin{bmatrix} y_1 \\ y_2 \\ y_3 \end{bmatrix} \mid y_i \in [0, m_i]; m_i \in R^+ I; 1 \leq i \leq 3, *, (9, 28) \right\}$$

be a neutrosophic column interval matrix groupoid of level two.

*Example 3.6.11:* Let

$$G = \left\{ \begin{bmatrix} x_1 \\ x_2 \\ x_3 \\ x_4 \\ \vdots \\ x_{27} \end{bmatrix} \mid x_i = [0, t_i]; t_i \in N(Z_{28}); 1 \le i \le 27; *, (13, 13) \right\}$$

be a neutrosophic column matrix interval groupoid of level four.

*Example 3.6.12:* Let

$$G = \left\{ \begin{bmatrix} x_1 \\ x_2 \\ x_3 \\ x_4 \end{bmatrix} \mid x_i = [0, n_i], n_i \in Z_5 I, *, (3, 0) \right\}$$

be a neutrosophic column interval matrix groupoid of level five. Clearly G is a groupoid of finite order.

*Example 3.6.13:* Let

$$G = \left\{ \begin{bmatrix} x_1 \\ x_2 \end{bmatrix} \mid x_i = [0, a_i] \in Z^+ I, *, (8, 26), 1 \le i \le 2 \right\}$$

be a neutrosophic column interval matrix groupoid of level three. Clearly G is of infinite order.

**DEFINITION 3.6.3:** *Let $G = \{M = (m_{ij}) \mid m_{ij} = [0, a_{ij}]; 1 \le i \le n, 1 \le j \le m; a_i \in N(Q^+), *, (t, u)\}$ be such that, for any $M = (m_{ij})$ and $N = (n_{ij}); M * N = ([0, a_{ij}]), * ([0, b_{ij}])$ (where $n_{ij} = [0, b_{ij}]$) $= ([0, ta_{ij} + ub_{ij}])$.*

$G$ is a neutrosophic $n \times m$ interval matrix groupoid when $n = m$ we call $G$ to be a neutrosophic interval square matrix groupoid.

We will illustrate these situation by some examples.

*Example 3.6.14:* Let

$$G = \left\{ \begin{pmatrix} x_1 & x_2 \\ x_3 & x_4 \end{pmatrix} \mid x_i = [0, a_i]; a_i \in Z_{81}I, *, (9, 25) \right\}$$

be a neutrosophic interval square matrix groupoid of finite order.

*Example 3.6.15:* Let

$$G = \left\{ \begin{pmatrix} y_1 & y_2 & y_3 & y_4 \\ y_5 & y_6 & y_7 & y_8 \end{pmatrix} \mid y_i \in [0, a_i]; a_i \in Z^+I; *, (3, 11) \right\}$$

be a neutrosophic interval rectangular matrix groupoid.

*Example 3.6.16:* Let

$$G = \left\{ \begin{bmatrix} y_1 & y_2 & y_3 \\ y_4 & y_5 & y_6 \\ y_7 & y_8 & y_9 \\ y_{10} & y_{11} & y_{12} \end{bmatrix} \mid y_i = [0, m_i] \ m_i \in N(R^+), *, (0, 7) \right\}$$

be a neutrosophic interval groupoid of level five.

*Example 3.6.17:* Let

$$G = \left\{ \begin{bmatrix} a_1 & a_2 & a_3 \\ a_4 & a_5 & a_6 \\ a_7 & a_8 & a_9 \\ a_{10} & a_{11} & a_{12} \\ a_{13} & a_{14} & a_{15} \end{bmatrix} \mid a_i = [0, m_i], m_i \in N(Q^+), *, (14, 26) \right\}$$

be a neutrosophic interval matrix groupoid of level three.

*Example 3.6.18:* Let

$$G = \left\{ \begin{bmatrix} a_1 & a_2 & a_3 & a_4 \\ a_5 & a_6 & a_7 & a_8 \\ a_9 & a_{10} & a_{11} & a_{12} \\ a_{13} & a_{14} & a_{15} & a_{16} \end{bmatrix} \right| $$

$a_i = [0, d_i]$; $d_i \in R^+I$; $1 \leq i \leq 16$, *, $(0, 26)\}$ be a neutrosophic interval matrix groupoid of infinite order.

*Example 3.6.19:* Let

$$G = \left\{ \begin{bmatrix} a_1 & a_2 \\ a_3 & a_4 \\ \vdots & \vdots \\ a_{49} & a_{50} \end{bmatrix} \middle| a_i = [0, n_i]; n_i \in Z_4I, (2, 2), * \right\}$$

be a neutrosophic interval matrix groupoid of finite order.

Now having defined neutrosophic interval matrix groupoids we now proceed onto define substructures in them and then illustrate it by examples.

**DEFINITION 3.6.4:** *Let $G = \{M = (m_{ij}) \mid m_{ij} = [0, a_{ij}]; a_{ij} \in N(Z^+)$ or $N(Z_n)$ or $Z_nI$ or $N(R^+)$ and so on $1 \leq i \leq m; 1 \leq j \leq t$, *, $(p, q)\}$ be a neutrosophic interval matrix groupoid. Let $P \subseteq G$; P is proper subset of G such that P itself is a neutrosophic interval matrix groupoid. We define P to be a neutrosophic interval matrix subgroupoid of G.*

We will illustrate this by some simple examples.

*Example 3.6.20:* Let

$$G = \left\{ \begin{bmatrix} a_1 & a_2 \\ a_3 & a_4 \\ a_5 & a_6 \end{bmatrix} \middle| a_i = [0, m_i]; m_i \in Z_9I; 1 \leq i \leq 6, *, (3, 6) \right\}$$

be a neutrosophic interval matrix groupoid of level three.

Take

$$P = \left\{ \begin{bmatrix} a & a \\ a & a \\ b & b \end{bmatrix} \mid a = [0, m], b = [0, n]; m, n \in Z_9I; *, (3, 6) \right\} \subseteq G.$$

P is a neutrosophic interval matrix subgroupoid of G of level three.

*Example 3.6.21:* Let G = {all 9 × 9 interval neutrosophic matrices with entries from $o(Z^+I)$, *, (11, 11)} be a neutrosophic matrix interval groupoid of level four.

Take W = {all 9 × 9 upper triangular interval neutrosophic matrices, with entries from $o(Z^+I)$, *, (11, 11)} ⊆ G; T is a neutrosophic matrix interval subgroupoid of G of level four.

*Example 3.6.22:* Let

$$G = \left\{ \begin{bmatrix} a_1 & a_2 \\ a_3 & a_4 \\ \vdots & \vdots \\ a_{19} & a_{20} \end{bmatrix} \mid \right.$$

$a_i = [0, n_i]$, $n_i \in N(R^+ \cup \{0\}))$; $1 \leq i \leq 20$; *, (0, 16)} be a neutrosophic matrix interval groupoid of level five.

Take

$$P = \left\{ \begin{bmatrix} a_1 & a_1 \\ a_2 & a_2 \\ \vdots & \vdots \\ a_{10} & a_{10} \end{bmatrix} \mid \right.$$

$a_i = [0, m_i]$; $m_i \in N(R^+ \cup \{0\}))$, $1 \leq i \leq 10$, *, (0, 16)} ⊆ G. P is a neutrosophic matrix interval subgroupoid of G of level five.

We can define as in case of other groupoids the notion of ideals, left ideal, right ideal, normal groupoid, normal subgroupoid, Smarandache groupoids, groupoids satisfying special identities and so on.

We will give some examples, a few important properties about these structures.

***Example 3.6.23:*** Let

$$G = \left\{ \begin{bmatrix} a_1 & a_2 \\ a_3 & a_4 \end{bmatrix} \mid a_i = [0, m_i];\ m_i \in Z_4I,\ 1 \leq i \leq 4, *, (2, 3) \right\}$$

be a neutrosophic interval matrix groupoid.
  Take

$$P = \left\{ \begin{bmatrix} a_1 & a_2 \\ a_3 & a_4 \end{bmatrix} \mid a_i = [0, m_i];\ m_i \in \{0, 2I\}, *, (2, 3) \right\} \subseteq G,$$

T is also a neutrosophic interval matrix left ideal of G.

**THEOREM 3.6.1:** *Let $G = \{m \times p$ neutrosophic interval matrices from $o(Z_nI)$, *, $(t, u)\}$ be a neutrosophic matrix interval groupoid. G is an idempotent neutrosophic interval matrix groupoid or neutrosophic interval matrix idempotent groupoid if and only if $t + u \equiv 1 (mod\ n)$.*

*Proof:* Let $M = (m_{ij}) \in G$ with $m_{ij} = [0, a_{ij}]$, $a_{ij} \in Z_nI$; $1 \leq i \leq m$, $1 \leq j \leq p$. To show $M * M = M$.
Consider
$$\begin{aligned} M * M &= (m_{ij}) * (m_{ij}) \\ &= [0, a_{ij}] *, [0, a_{ij}] \\ &= [0, ta_{ij} + ua_{ij}\ (mod\ n)] \\ &= [0, a_{ij}(t + u)\ (mod\ n)] \\ &= [0, a_{ij}] \\ &= M \end{aligned}$$

if and only if $t + u = 1\ (mod\ n)$. Hence the claim.

**THEOREM 3.6.2:** *Let $G = \{$all $m \times p$ neutrosophic interval matrices with entries from $o(Z_nI)$, *, $(t, u)\}$ and $G' = \{$all $m \times p$ neutrosophic interval matrices with entries from $o(Z_nI)$, *, $(u, t)\}$ be two neutrosophic interval matrix groupoids. P is a left ideal of G if and only if P is a right ideal of G.*

*Proof:* Simple number theoretic computations will yield the result.

The reader is expected to prove.

*Example 3.6.24:* Let G = {set of all 3 × 1 column interval neutrosophic matrices with entries from $o(Z_{10}I)$, *, (3, 7)} be a neutrosophic column interval matrix groupoid. G has no left or right ideals.

**THEOREM 3.6.3:** *Let G = {collection of all 1 × m neutrosophic interval row matrices with entries from $o(Z_nI)$; *, (t, u)} be a neutrosophic interval row matrix groupoid. If n = t + u; where t and u are primes then G has no left or right ideals.*

The reader can prove this theorem using some simple number theoretic techniques.

**COROLLARY 3.6.1:** *If in the above theorem n = p, p a prime and if (t, u) = 1 with t + u = p then also G has no left or right ideals.*

*Example 3.6.25:* Let
$$G = \left\{ \begin{bmatrix} a & a & a \\ a & a & a \end{bmatrix} \mid a \in o(Z_5I), *, (2, 2) \right\}$$
be a neutrosophic matrix interval groupoid. G has no subgroupoids.

*Example 3.6.26:* Let
$$G = \left\{ \begin{bmatrix} a & a \\ a & a \\ a & a \\ a & a \\ a & a \\ a & a \end{bmatrix} \mid a \in o(Z_7I), *, (3, 3) \right\}$$

be a neutrosophic interval matrix groupoid. G has no subgroupoids.

*Example 3.6.27:* Let

$$G = \left\{ \begin{bmatrix} a & a \\ a & a \end{bmatrix} \mid a \in o\ (Z_6I),\ *,\ (2,\ 2) \right\}$$

be a neutrosophic interval matrix groupoid. G has subgroupoid

$$P = \left\{ \begin{bmatrix} a & a \\ a & a \end{bmatrix} \mid a \in o(0,\ 2I,\ 4I),\ *,\ (2,\ 2) \right\} \subseteq G$$

is a neutrosophic interval matrix subgroupoid of G.

**THEOREM 3.6.4:** *Let*

$$G = \left\{ \begin{bmatrix} a & a \\ a & a \end{bmatrix} \mid a \in o(Z_pI);\ p\ a\ prime,\ *,\ (t,\ t),\ t < p \right\}$$

*be the neutrosophic interval matrix groupoid of level four. G is a normal groupoid.*

The proof is left as an exercise for the reader.

**THEOREM 3.6.5:** *Let*

$$G = \left\{ \begin{pmatrix} a & a & a \\ a & a & a \\ a & a & a \end{pmatrix} \mid a \in o(Z_nI);\ *,\ (t,\ t) \right\}$$

*be a neutrosophic interval matrix groupoid of level four. G is a P-groupoid.*

The reader is requested to supply the proof.

**THEOREM 3.6.6:** *Let*

$$G = \left\{ \begin{bmatrix} a & a & a & a \\ a & a & a & a \\ a & a & a & a \\ a & a & a & a \\ a & a & a & a \end{bmatrix} \;\middle|\; a \in o(Z_nI),\, *,\, (t, t),\, 1 < t < n \right\}$$

*be a neutrosophic interval matrix groupoid. If n is a prime then G is not an alternative groupoid.*

*Proof:* We see for any $X = (a)$ and $Y = (b)$ in G we have to show that

$$(X * Y) * Y \neq X * (Y * Y).$$

Consider
$(X * Y) * Y \quad = \quad (ta + tb) * Y$
$\quad\quad\quad\quad\quad = \quad (t2a + t2b + tb) \pmod{n})$      I

Consider
$X * (Y * Y) \quad = \quad X * (tb + tb)$
$\quad\quad\quad\quad\quad = \quad ((tb + t^2 a + t^2 b) \pmod{n})$      II

I = II if and only if $t^2 \equiv t \pmod{n}$ this is possible if and only if n is not a prime number.

**THEOREM 3.6.7:** *Let*

$$G = \left\{ \begin{bmatrix} a & a & a & a \\ a & a & a & a \\ a & a & a & a \\ a & a & a & a \end{bmatrix} \;\middle|\; a \in o(Z_nI),\, *,\, (0, t) \right\}$$

*be a neutrosophic matrix interval groupoid. G is a P-groupoid and an alternative groupoid if and only if $t^2 \equiv t \pmod{n}$.*

Proof is simple and hence is left as an exercise for the reader.

*Example 3.6.28:* Let G = {(a, a, a, a, a, a, a, a) | a ∈ o($Z_{10}$I), *, (91, 5)} be a neutrosophic matrix interval groupoid. G is a Smarandache neutrosophic matrix interval groupoid.

*Example 3.6.29:* Let

$$G = \left\{ \begin{bmatrix} a & a \\ a & a \end{bmatrix} \mid a \in o\ (Z_6 I),\ *,\ (4, 5) \right\}$$

be a neutrosophic matrix interval groupoid.

$$P = \left\{ \begin{bmatrix} a & a \\ a & a \end{bmatrix} \mid a \in \{[0, 0], [0, 2I], [0, 4I]\},\ *,\ (4, 5) \right\}$$

is a neutrosophic interval matrix subgroupoid of G but is not a Smarandache subgroupoid of G.

*Example 3.6.30:* Let

$$G = \left\{ \begin{bmatrix} a \\ a \\ a \\ a \end{bmatrix} \mid a \in o\ (Z_6 I),\ *,\ (4, 5) \right\} \text{ be}$$

a neutrosophic matrix interval groupoid.

$$A = \left\{ \begin{bmatrix} a \\ a \\ a \\ a \end{bmatrix} / a \in \{[0, 0], [0, I], [0, 3I], [0, 5I]\},\ *,\ (4, 5) \right\} \subseteq G$$

is a Smarandache left ideal of G which is clearly not a Smarandache right ideal of G.

*Example 3.6.31:* Let

$$G = \left\{ \begin{bmatrix} a & a \\ a & a \\ a & a \end{bmatrix} / a \in \{[0\ 0], [0\ 2I], [0, 4I]\},\ *,\ (2, 4) \right\}$$

is an ideal which is not a Smarandache ideal or even a Smarandache groupoid of G.

*Example 3.6.32:* Let

$$G = \left\{ \begin{bmatrix} a & a & a \\ a & a & a \\ a & a & a \\ a & a & a \end{bmatrix} \Big/ a \in o(Z_8I), *, (2, 6) \right\}$$

be the neutrosophic matrix interval groupoid.

Let

$$A = \left\{ \begin{bmatrix} a & a & a \\ a & a & a \\ a & a & a \\ a & a & a \end{bmatrix} \Big| a \in \{0, [0, 2I], [0, 4I], [0, 6I]\}, *, (2, 6) \right\} \subseteq G;$$

A is clearly a Smarandache neutrosophic matrix interval normal groupoid of G.

*Example 3.6.33:* Let $G = \{(a, a, a, \ldots, a) | a \in o(Z_{10}I), *, (5, 6)\}$ be a $1 \times 26$ neutrosophic row matrix interval groupoid. G is a Smarandache strong Moufang groupoid.

*Example 3.6.34:* Let

$$G = \left\{ \begin{bmatrix} a & a & a & a & a & a \\ a & a & a & a & a & a \\ a & a & a & a & a & a \end{bmatrix} \Big| a \in o(Z_{12}I); *, (3, 9) \right\}$$

be a neutrosophic matrix interval groupoid. G is only Smarandache Moufang groupoid but is not a Smarandache strong Moufang groupoid.

***Example 3.6.35:*** Let

$$G = \left\{ \begin{bmatrix} a & a & a \\ a & a & a \\ a & a & a \\ a & a & a \\ a & a & a \\ a & a & a \\ a & a & a \\ a & a & a \\ a & a & a \end{bmatrix} \mid a \in o(Z_{12}I), *, (3, 4) \right\}$$

be a neutrosophic matrix interval groupoid. G is a Smarandache strong Bol groupoid.

***Example 3.6.36:*** Let

$$G = \left\{ \begin{bmatrix} a & a & a \\ a & a & a \\ a & a & a \\ a & a & a \\ a & a & a \\ a & a & a \\ a & a & a \\ a & a & a \\ a & a & a \\ a & a & a \\ a & a & a \\ a & a & a \end{bmatrix} \mid a \in o(Z_6I), *, (4, 3) \right\}$$

be a neutrosophic matrix interval groupoid. G is a Smarandache strong P-groupoid.

*Example 3.6.37:* Let

$$G = \left\{ \begin{bmatrix} a & a & a & a \\ a & a & a & a \\ a & a & a & a \\ a & a & a & a \end{bmatrix} \mid a \in o(Z_4I), *, (2, 3) \right\}$$

be a neutrosophic matrix interval groupoid. G is only a Smarandache Bol groupoid but G is not a Smarandache strong Bol groupoid.

*Example 3.6.38:* Let

$$G = \left\{ \begin{bmatrix} a & a & a & a & a \\ a & a & a & a & a \end{bmatrix} \mid a \in o(Z_4I), *, (2, 3) \right\}$$

be a neutrosophic matrix interval groupoid. G is a Smarandache strong P-groupoid.

*Example 3.6.39:* Let

$$G = \left\{ \begin{bmatrix} a & a & a & a \\ a & a & a & a \\ a & a & a & a \\ a & a & a & a \\ a & a & a & a \\ a & a & a & a \end{bmatrix} \mid a \in o(Z_6I), *, (3, 5) \right\}$$

be a neutrosophic matrix interval groupoid. Clearly G is only a Smarandache strong P-groupoid.

*Example 3.6.40:* Let

$$G = \left\{ \begin{bmatrix} a & a & a \\ a & a & a \\ a & a & a \end{bmatrix} \mid a \in o(Z_{12}I), *, (5, 10) \right\}$$

be a neutrosophic matrix interval groupoid. G is a Smarandache P-groupoid and is not a Smarandache strong P-groupoid.

*Example 3.6.41:* Let

$$G = \left\{ \begin{bmatrix} a & a \\ a & a \\ a & a \\ a & a \\ a & a \end{bmatrix} \mid a \in o(Z_{14}I), *, (7, 8) \right\}$$

be a neutrosophic matrix interval groupoid. G is a Smarandache, strong alternative groupoid.

*Example 3.6.42:* Let

$$G = \left\{ \begin{bmatrix} a & a & a & a & a \\ a & a & a & a & a \\ a & a & a & a & a \end{bmatrix} \mid a \in o(Z_{12}I), *, (1, 6) \right\}$$

be a neutrosophic matrix interval groupoid. G is only a Smarandache strong alternative groupoid.

*Example 3.6.43:* Let

$$G = \left\{ \begin{bmatrix} a & a \\ a & a \end{bmatrix} \mid a \in o(Z_9I), *, (4, 4) \right\}$$

be a neutrosophic matrix interval groupoid. Clearly G is not a Smarandache groupoid.

**THEOREM 3.6.8:** *Let*

$$G = \left\{ \begin{bmatrix} a & a & a \\ a & a & a \\ a & a & a \end{bmatrix} \mid a \in o(Z_nI), *, (0, 2) \text{ and } n = 2m \right\}$$

*be a neutrosophic interval matrix groupoid. Clearly G is a Smarandache groupoid.*

Several other results not mentioned in this section but satisfied by general groupoids are true in case of Smarandache neutrosophic interval matrix groupoids with appropriate modifications.

### 3.7 Neutrosophic Interval Polynomial Groupoids

In this section we will be defining and discussing about neutrosophic interval polynomial groupoids or neutrosophic polynomial interval coefficients groupoids.

The notions given in section 3.5 will be used in this section. Here also we can define five levels of neutrosophic interval polynomial groupoids using $o(Z_nI)$, $o(N(Z_n))$, $o(Z^+I)$, $o(N(Z^+))$, $o(Q^+I)$, $o(N(Q^+))$ and so on.

From the context one can easily understand to which level the neutrosophic interval groupoid belongs to.

**DEFINITION 3.7.1:** *Let*

$$G = \left\{ \sum_{i=0}^{\infty} [0, a_i] x^i \ / \ a_i \in Z_nI, *, (p, q). \ p, q \in Z_n \right\}$$

*where * for any two interval neutrosophic polynomials.*

$$x = \sum_{i=0}^{\infty} [0, a_i] x^i$$

*and*

$$y = \sum_{i=0}^{\infty} [0, b_i] x^i$$

*is defined as follows.*

$$\begin{aligned} x * y &= \left( \sum_{i=0}^{\infty} [0, a_i] x^i \right) * \left( \sum_{i=0}^{\infty} [0, b_i] x^i \right) \\ &= \sum_{i=0}^{\infty} ([0, a_i] * [0, b_i]) x^i \end{aligned}$$

$$= \sum_{i=0}^{\infty}([0, pa_i + qb_i(\bmod n\,)]x^i.$$

*G is easily verified to be a groupoid which we call as the neutrosophic polynomial interval groupoid.*

*If p and q are two distinct primes we call G to be a level one groupoid.*

*If p and q are distinct such that (p, q) = 1 we call the groupoid G to be a level two groupoid. If p and q are such that (p, q) = d ≠ 1 we call G to be a level five groupoid.*

In definition 3.7.1 we can take instead of $Z_nI$. $N(Z_n)$ or $Z^+I$ or $N(Z^+)$ or $N(R^+)$ or $N(Q^+)$ or $Q^{+7}I$ or so on.

We will give some illustrations before we proceed on to prove some results.

*Example 3.7.1:* Let

$$G = \left\{ \sum_{i=0}^{25}[0, a_i]x^i \mid a_i \in Q^+I \cup \{0\}, *, (13, 41) \right\}$$

be a neutrosophic polynomial interval groupoid of level one.

*Example 3.7.2:* Let

$$V = \left\{ \sum_{i=0}^{45}[0, a_i]x^i \mid a_i \in N(R^+ \cup \{0\}), *, (12, 25) \right\}$$

be a neutrosophic interval polynomial groupoid of level two.

*Example 3.7.3:* Let

$$V = \left\{ \sum_{i=0}^{150}[0, a_i]x^i \mid a_i \in Z_{20}I, *, (15, 10) \right\}$$

be a neutrosophic polynomial interval groupoid of level three.

*Example 3.7.4:* Let

$$W = \left\{\sum_{i=0}^{\infty}[0,a_i]x^i \mid a_i \in C^+I \cup \{0\}, *, (3, 3)\right\}$$

be a neutrosophic polynomial interval groupoid of level four.

*Example 3.7.5:* Let

$$P = \left\{\sum_{i=0}^{\infty}[0,a_i]x^i \mid a_i \in Z^+I \cup \{0\}, *, (0, 22)\right\}$$

be a neutrosophic polynomial interval groupoid of level five.

Now we have seen different levels of neutrosophic polynomial interval groupoids. Now we will define the notion of subgroupoids.

**DEFINITION 3.7.2:** *Let*

$$G = \left\{\sum_{i=0}^{\infty}[0,a_i]x^i \,/\, a_i \in Z_nI \text{ or } Z^+I \text{ or so on}, *, (p, q)\right\}$$

*be a neutrosophic polynomial interval groupoid.*
*Suppose $P \subseteq G$ be a proper subset of G if P under the operations of G is a neutrosophic polynomial interval groupoid, then we define P to be a neutrosophic polynomial interval subgroupoid of G.*

We will illustrate this situation by some examples.

*Example 3.7.6:* Let

$$G = \left\{\sum_{i=0}^{\infty}[0,a_i]x^i \mid a_i \in Z_{40}I, *, (23, 17)\right\}$$

be a neutrosophic polynomial interval groupoid of level one.

Take

$$P = \sum_{i=0}^{5}[0,a_i]x^i \mid a_i \in Z_{40}I, *, (23, 17)\} \subseteq G$$

is a neutrosophic polynomial interval subgroupoid of G.

Suppose

$$S = \left\{\sum_{i=0}^{\infty}[0,a_i]x^i \mid \right.$$

$a_i \in \{0, 2I, 4I, 6I, 8I, 10I, 1I, \ldots, 36I, 38I\}, *, (23, 17)\} \subseteq G$, S is also a neutrosophic polynomial interval subgroupoid of G. We see P is a finite order where as S is of infinite order.

*Example 3.7.7:* Let

$$G = \left\{\sum_{i=0}^{25}[0,a_i]x^i \mid a_i \in Q^+I \cup \{0\}, *, (24, 35)\right\}$$

be a neutrosophic polynomial interval groupoid of level two.

Take

$$W = \left\{\sum_{i=0}^{25}[0,a_i]x^i \mid a_i \in Z^+I \cup \{0\}, *, (24, 35)\right\} \subseteq G$$

is a neutrosophic polynomial interval subgroupoid of G of level two.

Consider

$$P = \left\{\sum_{i=0}^{10}[0,a_i]x^i \mid a_i \in 4Z^+I \cup \{0\}, *, (24, 35)\right\} \subseteq G,$$

P is also a neutrosophic polynomial interval subgroupoid of G.

*Example 3.7.8:* Let

$$G = \left\{ \sum_{i=0}^{\infty}[0,a_i]x^i \mid a_i \in N(R^+ \cup \{0\}), *, (38, 90)\right\}$$

be a neutrosophic polynomial interval groupoid.
Take

$$P = \left\{ \sum_{i=0}^{25}[0,a_i]x^i \mid a_i \in 15; Z^+I \cup \{0\}, *, (38, 90)\right\} \subseteq G$$

is a neutrosophic polynomial interval subgroupoid of G of level three.

*Example 3.7.9:* Let

$$G = \left\{ \sum_{i=0}^{40}[0,a_i]x^i \mid a_i \in Z^+I \cup \{0\}, *, (10, 10)\right\}$$

be a neutrosophic polynomial interval groupoid of G of level four.
Consider

$$W = \left\{ \sum_{i=0}^{20}[0,a_i]x^i \mid a_i \in 20\, Z^+I \cup \{0\}, *, (10, 10)\right\} \subseteq G$$

is a neutrosophic polynomial interval subgroupoid of G of level four.

*Example 3.7.10:* Let

$$G = \left\{ \sum_{i=0}^{\infty}[0,a_i]x^i \mid a_i \in N(C^+) \cup \{0\}, *, (26, 0)\right\}$$

be a neutrosophic polynomial interval groupoid of level five.

$$P = \left\{ \sum_{i=0}^{\infty} [0, a_i] x^i \mid a_i \in C^+I \cup \{0\}, *, (26, 0) \right\} \subseteq G$$

is a neutrosophic polynomial interval subgroupoid of G.

We can as in case of other groupoids define the notion of left ideal, right ideal, ideal, normal groupoid, normal subgroupoid, Smarandache groupoids and Smarandache groupoids satisfying special identities.

**Chapter Four**

# APPLICATION OF THESE NEW CLASSES OF GROUPOIDS AND NEUTROSOPHIC GROUPOIDS

In this chapter we give some of the application of these new classes of groupoids. In the first place these polynomial groupoids can be used in theory of cryptography as the operations used by the user cannot be easily understood by the intruder so the information cannot be broken unless he has a good knowledge of groupoids and their properties.

Secondly these matrix groupoids can serve as storage systems in a special way when the information stored is a confidential one.

One of the nice applications would be using interval groupoids in the construction of semiautomaton and Smarandache automaton. Interested reader can unravel the probable applications of matrix interval groupoids and neutrosophic matrix interval groupoids when the field of study demands a non-associate structure with closure operations and interval solutions.

Smarandache groupoids can be used in computing or information science.

Smarandache groupoids can be used in biology to describe certain aspects in the crossing of organisms in genetics and in considerations of metabolisms.

**Chapter Five**

# SUGGESTED PROBLEMS

Here we suggest over 200 problems. Solving these will enable the reader to understand the concepts described in this book.

1. Find subgroupoids of $\{Z^+, *, (7, 9)\}$.

2. Find subgroupoids of $\{Q^+, *, (3/7, 2/11)\}$.

3. Let $P_{2\times 2} = \{$all $2 \times 2$ matrices with entries from $Z^+ \cup \{0\}\}$; define * on $P_{2\times 2}$ by $A * B = 5A + 7B$ for all $A, B \in P_{2\times 2}$. $P = \{P_{2\times 2}, *, (5, 7)\}$ is a matrix groupoid.
   a. Is P a normal matrix groupoid?
   b. Can P have normal matrix subgroupoids?
   c. Find atleast three matrix subgroupoids of P.
   d. Is P a Bol matrix groupoid?
   e. Can P be an alternative matrix groupoid?

Justify your answers.

4. Let S = {Collection of all 2 × 4 matrices with entries from $Q^+ \cup \{0\}$, *, (7, 3)} be a 2 × 4 matrix groupoid.
   a. Find atleast 3 subgroupoids.
   b. Can S have normal matrix subgroupiods?
   c. Is S a simple matrix groupoid?
   d. Is S a commutative matrix groupoid?

5. Let T = {all 1 × 5 row matrices with entries from $Z_{12}$, *, (3, 5)}. That is for A, B ∈ T; A * B = 3A + 5B, '+' mod 12 be the matrix groupoid.
   a. How many elements are in T?
   b. Find subgroupoids in T.
   c. Is T a normal matrix groupoid?
   d. Can T be a Bol groupoid?
   e. Is T a Moufang groupoid?
   Please justify your answers.

6. Let P = {7 × 1 matrices of the form
$$\left\{ \begin{bmatrix} a_1 \\ a_2 \\ a_3 \\ a_4 \\ a_5 \\ a_6 \\ a_7 \end{bmatrix} \middle| a_i \in Z_7, *, (3,2), 1 \le i \le 7 \right\}$$
be a 7 × 1 column matrix groupoid.
   a. Does P have subgroupoids?
   b. Can P have normal subgroupoids?
   c. Is P a S-groupoid?
   d. Is P a left alternative groupoid? Justify your claim.

7. Let W = $\left\{ \begin{bmatrix} a_1 & a_2 & a_3 \\ a_4 & a_5 & a_6 \\ a_7 & a_8 & a_9 \end{bmatrix} \middle| a_i \in Z_{11}; 1 \le i \le 9, *, (5, 6) \right\}$ be a 3 × 3 matrix groupoid. Find subgroupoids of W. What is the cardinality of W?

8. Let R = $\left\{ \begin{bmatrix} a & a \\ a & a \end{bmatrix} \middle| a \in Z_{13}; * (7, 5) \right\}$ be a 2 × 2 matrix groupoid.

   Find the number of elements in R. Is R simple? Justify your claim.

9. Let T = $\left\{ \begin{bmatrix} a & a & a & a \\ 0 & a & a & a \\ 0 & 0 & a & a \\ 0 & 0 & 0 & a \end{bmatrix} \middle| a \in Z_{24}; *, (2, 19) \right\}$ be a 4 × 4 matrix groupoid.
   
   a. Find the cardinality of T.
   b. Is T simple?
   c. Can T have normal subgroupoids?
   d. Is T a normal groupoid?
   e. Can T be a Bol groupoid?
   f. Is T a Moufang groupoid?
   g. Can T be an alternative groupoid? Justify your claim.

10. Obtain some interesting results about matrix groupoids built using $Q^+ \cup \{0\}$ using $(p, q) = (3, 5)$.

11. Can the matrix groupoid P = $\{A_{5\times 3} = (a_{ij}); 1 \le i \le 5, 1 \le j \le 3; a_{ij} \in R^+ \cup \{0\}, (23, 41)\}$ be Moufang? Justify your claim.

12. Can P in problem (11) be an alternative groupoid? Substantiate your answer.

13. Obtain some interesting properties about the matrix groupoids constructed using $C^+ \cup \{0\}$.

14. Give an example of a matrix groupoid which is simple.

15. Give an example of a matrix groupoid which is normal.

16. Obtain conditions for a matrix groupoid to be a Smarandache matrix groupoid.

17. Give an example of a matrix groupoid which is not a S-matrix groupoid.

18. Does there exist a Bol matrix groupoid constructed using $Z^+ \cup \{0\}$?

19. Can there be a matrix groupoid which is Moufang constructed using $Q^+ \cup \{0\}$?

20. Can any matrix groupoid constructed using $R^+ \cup \{0\}$ be right or left alternative? Justify your claim.

21. Find subgroupoids of the matrix groupoid

    $$G = \left\{ \begin{bmatrix} a_1 & a_2 & a_3 \\ a_4 & a_5 & a_6 \\ a_7 & a_8 & a_9 \\ a_{10} & a_{11} & a_{12} \end{bmatrix} \text{ where } a_i \in Z_{19};\ 1 \le i \le 12,\ *,\ (3, 7) \right\}.$$

    a. Is G a Bol matrix groupoid?
    b. Can G be right alternative?
    c. Does G have normal subgroupoids?
    d. Can $H = \left\{ \begin{bmatrix} a & a & a \\ a & a & a \\ a & a & a \\ a & a & a \end{bmatrix} \middle| a \in Z_{19},\ *,\ (3, 7) \right\} \subseteq G$ be normal?

    Justify your claims.

22. Let $H = \left\{ \begin{bmatrix} a_1 & 0 & b_1 & 0 \\ 0 & a_2 & 0 & b_2 \\ 0 & 0 & a_3 & b_3 \\ a_4 & 0 & 0 & b_4 \\ a_5 & b_5 & 0 & 0 \\ 0 & a_6 & b_6 & 0 \end{bmatrix} \middle| a_i, b_i \in Z_{17}, 1 \le i \le 6; *, (8, 9) \right\}$

be a matrix groupoid.
   a. Find atleast 2 proper subgroupoids of H.
   b. What is the cardinality of H?
   c. Can H satisfy any one of the special identities?

23. Let $P = \left\{ \begin{bmatrix} a_1 & a_2 & a_3 \\ 0 & a_4 & a_5 \\ 0 & 0 & a_6 \end{bmatrix} \middle| a_i \in Z_{240}; 1 \le i \le 6, (2, 3) \right\}$ be a

matrix groupoid find;
   a. Atleast 2 subgroupoids of P.
   b. Is P simple?
   c. Can P be a S-groupoid?
   d. Is P a P-groupoid?
   e. Can P be a strong Bol groupoid?
   f. Can P be a strong Moufang groupoid? Justify your claim.
      Replace (2, 3) by (20, 30) and answer all the six problems (a) to (f).

24. Let $V = \{(a_1, a_2, \ldots, a_9) | a_i \in Z_{24}; 1 \le i \le 9; (16, 9)\}$ be a matrix groupoid.
   a. Is V a P-groupoid?
   b. Is V an idempotent groupoid?
   c. Does V contain normal subgroupoids?
   d. Is V an alternative groupoid? Justify your claim.

   If $W = \{(a_1, \ldots, a_9) | a_i \in Z_{24}; 1 \le i \le 9, (9, 16)\}$ be a matrix groupoid. Is V isomorphic with W? Is $|V| = |W|$? Justify all your claims.

25. Prove G = {(x₁, …, x₂₀)| xᵢ ∈ $Z_{20}$; 1 ≤ i ≤ 20; *, (16, 5)} is a S-idempotent matrix groupoid. Is G a Smarandache P-groupoid? Justify your claim.

26. Prove T = $\left\{ \begin{bmatrix} a_1 & a_2 \\ a_3 & a_4 \\ a_5 & a_6 \\ a_7 & a_8 \\ a_9 & a_{10} \end{bmatrix} \middle| a_i \in Z_{12}; 1 \le i \le 10, *, (4, 9) \right\}$ is a matrix S-P-groupoid.

27. Consider T = $\left\{ \begin{bmatrix} a_1 & a_2 & a_3 \\ a_4 & a_5 & a_6 \\ a_7 & a_8 & a_9 \\ a_{10} & a_{11} & a_{12} \end{bmatrix} \middle| a_i \in Z_5, * (3, 2) \right\}$ be a matrix groupoid. How many groupoids are in the class of 4 × 3 matrix groupoids C($Z_5$, (4 × 3)) constructed using $Z_5$?
    a. How many in that class are S-matrix P-groupoids?
    b. Can this class of groupoids contain matrix S-strong Bol groupoids? Justify your claim.

28. Prove the 1 × 12 row matrix groupoid V = {{(a₁, a₂, …, a₁₂) | aᵢ ∈ $Z^+ \cup \{0\}$}, *, (11, 15)} is of infinite order.
    a. Is V a S-matrix groupoid?
    b. Find two subgroupoids of V.
    c. Can V have normal subgroupoids? Justify your claim.
    d. Is the groupoid W = {(a, a, …, a) / a ∈ $Z^+ \cup \{0\}$, *, (11, 15)} normal in V? Substantiate your claim.

29. Suppose in problem (28) set $Z^+ \cup \{0\}$ is replaced by $Z_{25}$ so that V is of finite order, can V be normal? Does atleast V contain normal subgroupoids? What are the identities satisfied by V? Is V an alternative groupoid? Is V an idempotent groupoid? Justify your claim.

30. Let T = { $Z_{22}^{3\times 8}$, *, (1, 11)} be a 3 × 8 matrix groupoid built using $Z_{22}$. Is T a S-groupoid? Justify your answer.

31. Let W = { $Z_{46}^{7\times 2}$, *, (1, 23)} be a 7 × 2 matrix groupoid built using $Z_{46}$.
   a. Is W a S-idempotent matrix groupoid? Prove your claim.
   b. Can W be a S-strong Bol matrix groupoid? Justify.
   c. Prove W is to a S-P-matrix groupoid.

32. Let M = { $Z_{25}^{9\times 3}$, *, (p, q), p, q, $\in Z_{25}$} be a 9 × 3 matrix groupoid.
   a. Does there exists for a suitable p, q $\in Z_{25}$, so that M is a;
      i. S-P-matrix groupoid.
      ii. S-idempotent matrix groupoid.
      iii. S-alternative matrix groupoid.
      iv. S-strong Bol matrix groupoid.

33. Let G = {$C^+$, *, 10 × 8, (9, 19)} be the set of all 10×8 complex matrix groupoid. Obtain some interesting properties about G.

34. Let P = $\left\{ \sum_{i=0}^{10} a_i x_i \mid a_i \in Z_{17}, *, (3, 11) \right\}$ be a polynomial groupoid built using $Z_{17}$. Find subgroupoids of P. Is P a S-groupoid? Prove or disprove!

35. Obtain some interesting properties about polynomial groupoids built using $Z_n$ or $Q^+ \cup \{0\}$ or $R^+ \cup \{0\}$ or $Z^+ \cup \{0\}$.

36. Give an example of a simple polynomial groupoid built using $Z_n$.

37. Obtain some interesting properties about polynomial S-groupoids.

38. Give an example of a S-strong Bol polynomial groupoid.

39. Does their exist a S-strong Moufang polynomial groupoid built using $Z^+ \cup \{0\}$? Justify your claim!

40. Construct an example of a polynomial groupoid which is not a S-polynomial groupoid.

41. Obtain a necessary and sufficient condition for a S-polynomial groupoid to be a S-alternative polynomial groupoid built using $Z_n$.

42. Can an alternative polynomial groupoid $G = \left\{ \sum_{i=0}^{n} a_i x^i \mid a_i \in Z^+ \cup \{0\}, *, p, q \in Z^+ (p, q) \right\}$ be built? Justify your claim. Find atleast 3 polynomial subgroupoids of G.

43. When is $T = \left\{ \sum_{i=0}^{24} a_i x^i \mid a_i \in Z_{28}, *, (p, q), p, q \in Z_{28} \right\}$;
    a. S-strong polynomial Bol groupoid (obtain condition on p and q)?
    b. S-strong alternative polynomial groupoid?
    c. S-strong P-groupoid?
    d. S-groupoid?
    e. S-idempotent groupoid?
    Obtain conditions on p and q.

44. Can $V = \left\{ \sum_{i=0}^{15} a_i x^i \mid a_i \in Z_{19}; *, (p, q), p, q, \in Z_{19} \right\}$ ever be a
    a. S-idempotent groupoid?
    b. S-P-groupoid?
    c. S-alternative groupoid (Here V is a polynomial groupoid)?

45. Obtain some interesting results about polynomial groupoids built using $Z_n$, n not a prime.

46. Prove there exist infinite number of polynomial groupoids built using $Z_n$.

47. Can the polynomial groupoid $J = \left\{ \sum_{i=0}^{29} a_i x^i \mid a_i \in Z_{24}; *, (20, 3) \right\}$ have identity element? Justify your claim.
    a. Is J a S-groupoid?
    b. Is J a S-strong Bol polynomial groupoid? Justify your claim.

48. Let $P = \left\{ \sum_{i=0}^{10} a_i x^i \mid a_i \in Z_{12}, *, (3, 9) \right\}$ be a polynomial groupoid built using $Z_{12}$.
    a. Find subgroupoids of P.
    b. Can P have a subgroupoid of order 4?
    c. What is the cardinality of P?
    d. Is P a S-groupoid?
    e. Does P have normal subgroupoids?
    Justify your answers.

49. Can any polynomial groupoid built using $R^+ \cup \{0\}$ be normal?

50. Obtain conditions for two polynomial subgroupoids built using $Z_{45}$ to be S-semi conjugate subgroupoids.

51. Can $G = \left\{ \sum_{i=0}^{5} a_i x^i \mid a_i \in Z_8, *, (2, 8) \right\}$ the polynomial groupoid have S-semi conjugate subgroupoids?

52. Define the notion of S-conjugate subgroupoids of a polynomial groupoid G built using $Z^+ \cup \{0\}$ or $Z_n$ or $R^+ \cup \{0\}$ or $Q^+ \cup \{0\}$.

53. Give examples of S-conjugate subgroupoids built using $Z_{144}$.

54. Let $T = \left\{ \sum_{i=0}^{20} a_i x^i \mid a_i \in Z_{12}, *, (1, 3) \right\}$ be a polynomial groupoid built using $Z_{12}$. Can T have S-conjugate polynomial subgroupoids? Is $P = \left\{ \sum_{i=0}^{20} a_i x^i , a_i \in \{0, 3\ 6, 9\}, *, (3, 1) \right\}$ conjugate with $S = \left\{ \sum_{i=0}^{20} a_i x^i , a_i \in \{2, 5, 8, 11\} \right.$ where P and S polynomial subgroupoids of T? Prove T is a S-polynomial groupoid.

55. Define S-inner commutative polynomial groupoid and give some examples.

56. Is $G = \left\{ \sum_{i=0}^{8} a_i x^i \mid a_i \in Z_5, *, (3, 3) \right\}$ a S-inner commutative polynomial groupoid? Justify your claim.

57. Prove every S-commutative polynomial groupoid is a S-inner commutative polynomial groupoid and not conversely.

58. Give an example of a S-Moufang polynomial groupoid built using $Z_{20}$.

59. Is $G = \left\{ \sum_{i=0}^{45} a_i x^i \mid a_i \in Z_{10}, *, (5, 6) \right\}$, the polynomial groupoid a S-strong Moufang groupoid? Justify your claim.

60. Is $P = \left\{ \sum_{i=0}^{20} a_i x^i \mid a_i \in 12\ (3, 9) \right\}$ the polynomial groupoid a S-strong Moufang Groupoid? Justify your claim.

61. Prove every S-strong Moufang groupoid is a S-Moufang groupoid and not conversely.

62. Is $W = \left\{ \sum_{i=0}^{\infty} a_i x^i \mid a_i \in Z_4, (2, 3) \right\}$ the polynomial groupoid a S-strong Bol groupoid? Can W be just a S-Bol groupoid? Justify your claim.

63. Give an example of S-strong P-groupoid built using $Z_6$.

64. Let $B = \left\{ \sum_{i=0}^{25} a_i x^i \mid a_i \in Z_6, *, (4, 3) \right\}$ be a polynomial groupoid built using $Z_6$. Is B a S-strong P-groupoid? Prove your claim.

65. Let $C = \left\{ \sum_{i=0}^{14} a_i x^i \mid a_i \in Z_4, *, (2, 3) \right\}$ be a polynomial groupoid. Is C a S-strong P-groupoid? Justify your claim.

66. Prove if G a S-polynomial groupoid built using $Z_n$(n even) and if G is a S-strong polynomial P-groupoid then every pair in G need not satisfy the P-groupoid identity.

67. Every S-strong polynomial P-groupoid is a S-polynomial P-groupoid and not conversely.

68. Define the notion of S-strong right alternative groupoid and illustrate it by an example.

69. Is $G = \left\{ \sum_{i=0}^{20} a_i x^i \mid a_i \in Z_{14}, *, (7, 8) \right\}$ the polynomial groupoid a S-strong alternative groupoid?

70. Let $G = \left\{ \sum_{i=0}^{\infty} a_i x^i \mid a_i \in Z_{12}, *, (1, 6) \right\}$ be a polynomial groupoid built using $Z_{12}$. Is G a S-strong polynomial alternative groupoid?

71. Give an example of a S-strong polynomial right alternative groupoid which is not a S-strong polynomial left alternative groupoid.

72. Give an example of S-alternative polynomial groupoid which is not a S-strong alternative polynomial groupoid.

73. Can $G = \left\{\sum_{i=0}^{19} a_i x^i \mid a_i \in Z_{19}, *, (3, 12)\right\}$ the polynomial groupoid be a strong Bol groupoid?

74. Let $G = \left\{\sum_{i=0}^{\infty} a_i x^i \mid a_i \in Z_{19}, *, (4, 15)\right\}$ be a polynomial groupoid. Can G be a S-Bol groupoid? Justify your claim.

75. Let $T = \left\{\sum_{i=0}^{27} a_i x^i \mid a_i \in Z_8, *, (3, 6)\right\}$ be a polynomial groupoid. Can T be a S-strong P-groupoid?

76. Let $L = \left\{\sum_{i=0}^{\infty} a_i x^i \mid a_i \in Z_9, *, (2, 4)\right\}$ be a polynomial groupoid. Can L be a S-strong right alternative groupoid?

77. Find S-right ideals of L given in problem (76). Does L in problem (76) have S-ideals? Justify.

78. Let G = {set of all 3 × 5 interval matrices built using $Z_I^+$, *, (3, 7)} be an interval matrix groupoid. Find subgroupoids of G. Is G a S-interval matrix groupoid?

79. Obtain some interesting properties about row matrix interval groupoid built using $Z_n^I$.

80. Let P = {all 3 × 1 interval matrices constructed using $Z_5^I$, *, (2, 3)} be the 3 × 1 interval matrix groupoid built on $Z_5^I$. Is P a S-interval matrix groupoid?

81. Let T = {{set of all 3 × 3 interval matrix built using $Z_{19}^I$ }, *, (15, 5)} be the matrix interval groupoid. What are the identities satisfied by T? Is T finite? Can T be a S-strong Bol interval matrix groupoid? Justify.

82. Give an example of a 3 × 8 matrix interval groupoid built using $Z_{15}^I$, which is a S-groupoid.

83. Give an example of a 6 × 3 matrix interval groupoid built using $Z_7^I$ which is not a S-groupoid.

84. Give an example of a 5 × 5 matrix interval groupoid built using $Z_{12}^I$ which is a S-strong Moufang groupoid.

85. Give an example of a S-row matrix interval groupoid built using $Z_{24}^I$ which is a S-alternative row matrix interval groupoid.

86. Give an example of a 5 × 8 matrix interval groupoid built using $Z_{23}^I$ which is simple.

87. Does their exists a S-inner commutative matrix interval groupoid constructed using $Z_{43}^I$?

88. Obtain some interesting properties about square matrix interval groupoids built using $Z_p^I$; p an odd prime.

89. Construct a class of simple matrix interval groupoids using $Z_{41}^I$.

90. Give an example of a 10 × 2 matrix interval normal groupoid built using $Z_7^I$.

91. Give an example of a 1 × 8 matrix interval groupoid constructed using $Z_{25}^I$ which has normal matrix interval subgroupoids but is not a normal matrix interval groupoid.

92. Give an example of a 9 × 9 matrix interval groupoid which is a not a S-strong right alternative 9 × 9 matrix groupoid but a S-strong left alternative 9 × 9 matrix groupoid.

93. Can a 4 × 10 matrix interval groupoid built using $R_1^+$ be a S-groupoid? Justify your claim.

94. Let G = { $Z_9^I$ [x], *, (3, 7)} be an interval polynomial groupoid built using the intervals in $Z_9^I$,
    a. Does G have subgroupoids?
    b. Is G a S-groupoid?
    c. Is G a S-strong Bol groupoid?
    d. Can G be normal?
    e. Can G have alteast normal subgroupoids?
    f. Can G be a S-strong alternative groupoid? Justify your answers.

95. Let P = $\left\{ \sum_{i=0}^{4} a_i x^i \mid a_i \in Z_4^I \right.$ = {all intervals built using 0, 1, 2, 3}, *, (3, 1)} be a polynomial interval coefficient groupoid.
    a. What is the cardinality of P?
    b. Find subgroupoids of P.
    c. Does P satisfy any one of the standard identities?

96. Let W = $\left\{ \sum_{i=0}^{\infty} a_i x^i \mid a_i \in Z_{12}^I \right.$ = all intervals built using $Z_{12}$, *, (1, 3)} be an infinite polynomial interval coefficient groupoids built using $Z_{12}^I$.
    a. Prove W is an infinite S-groupoid.
    b. Is W a S-strong Moufang groupoid? Prove your answer.
    c. Can W be a S-Bol groupoid? Justify.
    d. Find two S-subgroupoids of W which are S-conjugate subgroupoids.
    e. Is W a S-idempotent groupoid? Justify your claim.

97. Let B = $\left\{ \sum_{i=0}^{5} a_i x^i \mid a_i \in Z_5^I \right.$ = {all intervals built using $Z_5$ = {0, 1, 2, 3, 4}, *, (3, 3)} be a polynomial interval groupoid built using $Z_5^I$.

a. Find subgroupoids of B.
b. Prove B is S-groupoid.
c. Is B a commutative groupoid?
d. Find the cardinality of B.
e. Is B inner commutative?

98. Let $K = \left\{ \sum_{i=0}^{20} a_i x^i \mid a_i \in Z_8^I = \{\text{all intervals built using } Z_8 = \{0, 1, 2, \ldots, 7\}\}, *, (2, 8) \right\}$ be a polynomial interval groupoid built using $Z_8^I$.
   a. Prove K is a S-groupoid.
   b. Prove K is of finite cardinality.
   c. Prove K is a S-normal groupoid.

99. Let $P = \left\{ \sum_{i=0}^{17} a_i x^i \mid a_i \in Z_6^I = \{\text{all intervals built using } \{0, 1, 2, 3, 4, 5\}\}, *, (0, 2) \right\}$ be a interval polynomial groupoid built using $Z_6^I$ in the variable x.
   a. Prove P is a non commutative groupoid.
   b. Is P a P-groupoid? Prove your claim.
   c. If (0, 2) is replaced by (0, 3) or (3, 0) can P be a P-groupoid and an alternative groupoid?

100. Let $G = \left\{ \sum_{i=0}^{\infty} a_i x^i \mid a_i \in \{\text{all intervals built using } Z_{24} = \{0, 1, 2, \ldots, 24\}\}, *, (0, t) \right\}$ be a polynomial interval groupoid in the variable x with coefficients from $Z_{24}^I$. Prove there exists a $t \in Z_{24}$ such that for that t, G is both a P-groupoid and an alternative groupoid.

101. Let $W = \left\{ \sum_{i=0}^{25} a_i x^i \mid a_i \in Z_{19}^I = \{\text{all intervals built using } \{0, 1, 2, \ldots, 18\}\}, *, (0, t) \right\}$ be a interval polynomial groupoid built using $Z_{19}^I$. Prove for no $t \in Z_{19} \setminus \{0, 1\}$. W can be a P-groupoid or a alternative groupoid.

102. Obtain some interesting results on polynomial interval groupoids built using $Z_I^+$.

103. Let $S = \left\{ \sum_{i=0}^{9} a_i x^i \,/\, a_i \in Z_{11}^I = \{\text{all intervals built using } \{0, 1, 2, \ldots, 10\}\}, *, (3, 0) \right\}$ be a interval polynomial groupoid. Can S satisfy the Moufang identity?

104. Can the interval polynomial groupoid $F = \left\{ \sum_{i=0}^{\infty} a_i x^i \,/\, a_i \in \{\text{all intervals built using } Z_n = \{0, 1, 2, \ldots, n-1\} = \{Z_n^I, *, (t, 0)\} \right.$ be a Bol groupoid?

105. Will $T = \left\{ \sum_{i=0}^{4} a_i x^i \mid a_i \in Z_{20}^I = \{\text{all intervals built using } Z_{20}\}, *, (0, 4) \right\}$ be an alternative interval polynomial groupoid?

106. Find all subgroupoids of $P = \left\{ \sum_{i=0}^{3} a_i x^i \mid a_i \in Z_4^I = \{\text{all intervals built using } Z_4 = \{0, 1, 2, 3\}, *, (0, 3)\} \right\}$ the interval polynomial groupoid built using intervals from $Z_4^I$. Find the order of P. Does the order of the subgroupoids divide the order of P?

107. Find all subgroupoids and ideals if any of the interval polynomial groupoid $T = \left\{ \sum_{i=0}^{3} a_i x^i \mid a_i \in Z_{12}^I = \{\text{all intervals built using } Z_{12} = \{0, 1, 2, \ldots, 1\}, *, (4, 0)\} \right.$. Does there exists in T a subgroupoid which is not an ideal? What is the order of T?

108. Let $S = \left\{ \sum_{i=0}^{\infty} a_i x^i \mid a_i \in Z_p^I = \{\text{all intervals built using } Z_p = \{0, 1, 2, \ldots, p-1\} \text{ where p is a prime}, *, (t, t); t < p \right\}$ be a interval

polynomial groupoid built using coefficients from $Z_p^I$. Prove S is a normal groupoid.

109. Prove groupoid S in problem (108) is a P-groupoid.

110. Prove groupoid S in problem (108) is a commutative groupoid.

111. Can groupoid P in problem (108) be an alternative groupoid?

112. Let F be the collection of all interval polynomial groupoids G where $G = \left\{ \sum_{i=0}^{\infty} a_i x^i \mid a_i \in Z_n^I \right\}$ = {all intervals using the modulo integers {0, 1, 2, …, n – 1}, n is not a prime, *, (t, t) where $t \in Z_n \setminus \{0, 1\}$ and t is such that $t^2 = t \pmod{n}$. Prove F contains only alternative groupoids. Find the number of such groupoids when n = 128.

113. Let $S = \left\{ \sum_{i=0}^{7} a_i x^i \mid a_i \in Z_{18}^I, *, (9, 9) \right\}$ be a interval polynomial groupoid built using the intervals from $Z_{18}^I$. Does S contain any normal subgroupoid?

114. Is the interval polynomial groupoid $T = \left\{ \sum_{i=0}^{5} a_i x^i \mid a_i \in Z_{25}^I; 0 \le i \le 5, *, (17, 17) \right\}$ a normal groupoid. Is T a S-groupoid? How many subgroupoids does T have? What is the cardinality of T?

115. Can the polynomial interval groupoid $L = \left\{ \sum_{i=0}^{3} a_i x^i \mid a_i \in Z_{21}^I, *, (7, 7) \right\}$ be a P-groupoid?
   a. Find all the subgroupoids of L.
   b. What is the cardinality of L?
   c. Does the order of each subgroupoid divide the order of L?
   d. Will L be a S-groupoid?

- e. Can L be a normal groupoid?
- f. Does L contain any normal subgroupoids?
- g. Will L be an idempotent groupoid?

116. Does the polynomial interval groupoid $V = \left\{ \sum_{i=0}^{3} a_i x^i \mid a_i \in Z_{11}^I; *, (4, 7) \right\}$ satisfy the Bol identity?

117. Can the groupoid V given in problem (116) be a Moufang groupoid or an idempotent groupoid? Justify your claim.

118. Does the polynomial interval groupoid $M = \left\{ \sum_{i=0}^{\infty} a_i x^i \mid a_i \in Z_{21}^I; *, (3, 7) \right\}$ have ideals? What is the order of M? Is M a S groupoid? Find all right ideals of M?

119. How many semigroups are in the interval polynomial groupoid $K = \left\{ \sum_{i=0}^{6} a_i x^i \mid a_i \in Z_7^I; *, (3, 4) \right\}$? Is K a S-groupoid? What is the order of K?

120. Construct a polynomial interval groupoid using $Z_{27}^I$ such that it is not a S-groupoid.

121. Construct a polynomial interval groupoid T using $Z_{292}^I$ which is an alternative groupoid.

122. Does there exist a simple polynomial interval groupoid using the interval set $Z_{27}^I$ ?

123. Find the number of subgroupoids of the polynomial interval groupoid $X = \left\{ \sum_{i=0}^{5} a_i x^i \mid a_i \in Z_3^I; *, (1, 2) \right\}$.

124. Is X given in problem 123 a S-groupoid?

125. What is the order of X given in problem (123)?

126. Let $G = \left\{ \sum_{i=0}^{\infty} a_i x^i \mid a_i \in Z_{17}^I ; *, (11, 7) \right\}$ be a polynomial interval groupoid.
    a. Is G simple?
    b. Can G contain semigroups?

127. Obtain some interesting properties about interval matrix groupoids built using intervals from $Z_I^+$.

128. Let P = {6 × 7 interval matrices built using intervals from $Z_3^I$} be a interval matrix groupoid. Is P simple? Find all subgroupoids of P. What is the cardinality of P?

129. Let Y = {All 3 × 3 interval matrices with entries from $Z_7^I$}. P = {Y, *, (3, 5)} be the interval matrix groupoid using $Z_7^I$. What is the cardinality of P? Find all subgroupoids of P. Is P a Moufang groupoid? Can P be a normal groupoid? Does P contain ideals? Find all right ideals of P. If (3, 5) is replaced by (4, 4) in P, will P be simple?

130. Let B = {L, *, (3, 9)} where L = {all 7 × 2 interval matrices with entries from $Z_I^+$} be a matrix interval groupoid.
    a. Can B be a S-groupoid?
    b. Find at least 5 subgroupoids of B.
    c. Does B contain right ideal?
    d. Can B have left ideals?
    e. Will B have an ideal?
    f. Can B be a normal groupoid?
    g. Can B satisfy any standard identities?

131. Obtain some interesting properties about the groupoid F = {the 10 × 10 interval matrix groupoid F built using $R_I^+$, *, (3, 11)}. Can F be simple? Find atleast 3 proper subgroupoids of F. Can F be a S-groupoid?

132. Let D = {All 2 × 17 interval matrices using $Z_4^I$; *, (3, 1)} be a matrix interval groupoid,
    a. What is the order of D?
    b. How many subgroupoids exists in D?
    c. Is D simple?
    d. Is D normal?
    e. Is D a Moufang groupoid?
    f. Can D have ideals?
    g. Is every left ideal of D a right ideal of D? Justify your answer.
    h. Can we say order the subgroupoids of D will divide the order the groupoid D?

133. Let G = {all 2 × 2 interval matrix built using $Z_5^I$, *, (3, 3)} be the interval matrix groupoid built using $Z_5^I$.
    a. Find the number of elements in G.
    b. Is G a commutative groupoid?
    c. Is G a normal groupoid?
    d. Can G have left ideals which are not right ideals and vice versa?
    e. How many subgroupoids does G contain?

134. If in problem (133) $Z_5$I is replaced by $Z_n$I; n a composite number. Find the solution from a to e of problem (133).

135. Let G = {1 × 7 interval matrices with entries from $Z_{13}^I$, *, (3, 11)} be a interval matrix groupoid.
    a. Find subgroupoids of G.
    b. Is G a normal groupoid?
    c. Can G be a Bol groupoid?
    d. Will G be a S-groupoid?
    e. Find the cardinality of G
    f. Find right ideals of G.
    g. Prove in general right ideal in G is not a left ideal of G and vice versa.

136. Let T = {all 3 × 9 interval matrices with entries from $Z_8^I$, *, (0, 2)} be the interval matrix groupoid.
    a. Prove T is a S-groupoid.
    b. Find subgroupoids of G
    c. Find the cardinality of T.
    d. Does T have S-subgroupoids?
    e. Find ideals in T.
    f. Prove in T a left ideal in general is not a right ideal.

137. Let X = {all 2 × 2 interval matrices with entries from $Z_{12}^I$, *, (6, 6)}. Derive the important features enjoyed by X.

138. Let Y = {5 × 5 interval matrices with entries from $Z_6^I$, *, (3, 4)} be a matrix interval groupoid built using $Z_6^I$.
    a. Is Y a S-groupoid?
    b. Prove Y is of finite order.
    c. Is Y a S-strong Bol-groupoid?
    d. Find subgroupoids of Y.
    e. Will Y be a S-strong Moufang groupoid?
    f. Can Y contain normal subgroupoids?

139. Construct a S-Moufang groupoid with 3 × 2 interval matrices built using $Z_{22}^I$.

140. Find all S-P- 3 × 3 matrix interval groupoids built using $Z_9^I$.

141. Let G = {3 × 3 interval matrices with entries from $Z_6^I$, *, (2, 3)} be an interval matrix groupoid. Find all the special identities satisfied by G.

142. Let T = {3 × 2 interval matrices with entries from $Z_6^I$, *, (2, 3)} be the interval matrix groupoid built using $Z_6^I$. S = {3 × 2 interval matrices with entries from $Z_{12}^I$, *, (2, 3)} be the interval matrix groupoid built using $Z_{12}^I$. Construct a SG homomorphism between T and S.

143. Let W = {2 × 5 interval matrices with entries from $Z_{19}^I$, *, (3, 7)} be a matrix interval groupoid built using $Z_{19}^I$. V = {2 × 5 interval matrices with entries from $Z_{19}^I$, *, (3, 6)} be a matrix interval groupoid built using $Z_{19}^I$. Give a SG homomorphism between W and V. Find all the common properties enjoyed by W and V.

144. T = {3 × 2 interval matrices with entries from $Z_{19}^I$, *, (3, 12)} be a interval matrix groupoid built using $Z_{19}^I$.
   a. Is T a S-Bol groupoid ?
   b. Find the order of T.
   c. Find all subgroupoids of T.
   d. Is T a S-groupoid?
   e. Is T a normal groupoid?
   f. Can T have normal subgroupoids?
   g. Obtain any other interesting property enjoyed by T.

145. Let B = {All 1 × 8 interval matrices with entries from $Z_{12}^I$, *, (3, 6)} be an interval matrix groupoid.
   a. Find the cardinality of B
   b. Find all subgroupoids of B
   c. Is B a S-strong Moufang groupoid?
   d. Can B be a S-strong right alternative groupoid?
   e. Can B have normal subgroupoids?
   f. Find an ideal in B.

146. Let W = {3 × 1 interval matrices with entries from $Z_{12}^I$, *, (1, 3)} be the matrix interval groupoid built using $Z_{12}^I$.
   a. Find two S-subgroupoids of W which are S-conjugate.
   b. Is W a S-groupoid?
   c. Is W a normal groupoid?
   d. How many semigroups does W contain?

147. Let V = {3 × 3 interval matrices with entries from $Z_5^I$, *, (3, 3)} be a matrix interval groupoid built using $Z_5^I$.

      a. Is V a S-inner commutative groupoid?
      b. Find cardinality of V.
      c. Find all subgroupoids of V.
      d. Does V have normal subgroupoids?
      e. Will V satisfy any one of the special identities?

148. Let G = {4 × 5 interval matrices with entries from $Z_6^I$, *, (4, 5)} be a matrix interval groupoid. Prove every S-left ideal of G is not a S-right ideal of G.

149. Define the notion of S-seminormal matrix interval groupoid. Illustrate it by some examples.

150. Let G = {all 6 × 6 interval matrices with entries from $Z_8^I$, *, (2, 6)} be a matrix interval groupoid.
      a. Prove G is a S-groupoid.
      b. Does G have a S-normal subgroupoid?
      c. Find all S-subgroupoids of G.
      d. What is the cardinality of G?
      e. Can G be a S-strong Moufang groupoid?

151. Let T = {Set of all 3 × 6 interval matrices with entries from $Z_{12}^I$, *, (3, 9)} be a interval matrix groupoid. Prove T is a not a S-strong Moufang groupoid. Prove T is only a S-Moufang groupoid.

152. Prove every S-strong Moufang groupoid is a S-Moufang groupoid and not conversely.

153. W = {set of all 2 × 8 interval matrices with entries from $Z_{10}^I$, *, (5, 6)} be an interval matrix groupoid. Prove W is a S-strong Moufang groupoid.

154. Let P = {All (4 × 4) interval matrices with entries from $Z_{25}^I$, *, (17, 17)} be an interval matrix groupoid. Is P a normal groupoid?

155. Let S = {all 7 × 7 interval matrices with entries from $Z_{21}^I$, *, (7, 7)} be an interval matrix groupoid. Is S a P-groupoid? Justify your claim.

156. Let V = {3 × 3 interval matrices with entries from $Z_{18}^I$, *, (9, 9)} be a interval matrix groupoid;
    a. Is V a normal groupoid?
    b. Does V contain normal subgroupoids?
    c. Is V an idempotent groupoid?
    d. Find all subgroupoids of V.
    e. Is V a S-groupoid?
    f. Find all S-subgroupoids of V. (Provided V is a S-groupoid).

157. Let G = {all 3 × 12 interval matrices with entries from $Z_4^I$, *, (2, 3)} be a matrix interval groupoid.
    a. Can {0} be an ideal of G?
    b. Prove left ideals of G are not right ideals of G.
    c. Find the cardinality of G.
    d. Find all subgroupoids of G.
    e. Is G a S-groupoid?
    f. Is G a normal groupoid?

158. Let V = {all 3 × 5 interval matrices with entries from $Z_6^I$, *, (3, 5)} be a matrix interval groupoid.
    a. Is V a S-strong P-groupoid?
    b. Is V a S-P-groupoid?
    c. Is V a S-groupoid?
    d. What is the cardinality of V?
    e. Find all subgroupoids of G.
    f. Find all S-subgroupoids of G.

159. Let B = {2 × 1 interval matrices with entries from $Z_4^I$, *, (2, 3)} be a matrix interval groupoid.
    a. Prove B is a S-strong P-groupoid.
    b. Find all subgroupoids of B.
    c. What is the cardinality of B?
    d. Find all S-subgroupoids of B.

e. Can B have ideals?
    f. Can a right ideal in B in general be a left ideal?

160. Let R = {All 2 × 2 interval matrices with entries from $R_I^+$, *, (3, 13)} be an interval matrix groupoid.
    a. Find atleast 3 subgroupoids of R.
    b. Find 3 left ideals of R.
    c. Find 3 right ideals of R.
    d. Find 3 ideals of R.
    e. Can R satisfy special identities?
    f. Can R have normal subgroupoids?

161. Obtain some interesting results about the matrix interval groupoid G = {all 7 × 7 interval matrices with entries from $Q_I^+$, *, (1, 2)}.

162. Let G = {$Z_9I$, *, (3, 2)} be a neutrosophic groupoid. Find whether G is a S-neutrosophic groupoid. What is the cardinality of G?

163. Let G = {$Z_{11}I$, *, (3, 4)} be a neutrosophic groupoid. Can G be written as a partition of conjugate groupoids? Justify your answer.

164. Find the conjugate groupoids of G = {$Z_{15}I$, *, (9, 8)}.

165. How many neutrosophic groupoids of order 7 exist?

166. Find all right ideals of the neutrosophic groupoid G = {$Z_{21}I$, *, (3, 7)}.

167. Let G = {$N(Z_{25})$, *, (3, 8)} be a neutrosophic groupoid. What are the important properties enjoyed by G? How many groupoids can be built using $N(Z_{25})$?

168. Is P = {$N(Z_5)$, *, (2I+ 1, 4I)} a neutrosophic groupoid? Justify your claim.

169. Is R = {$N(Z_7)$, *, (4 + 3I, 2)} a neutrosophic groupoid?

170. Let W = {N(Z$^+$), *, (1+3I, 7+2I)}. Is W a neutrosophic groupoid? What are the special properties enjoyed by W?

171. Let S = {N(Z$^+$∪ {0}), *, {3I, 2 + 2I)} be a special neutrosophic groupoid. Find subgroupoids of S. Does S satisfy any one of the special identities? Justify your answer.

172. Does there exist a P-neutrosophic groupoid of order 15?

173. Give an example of a P-neutrosophic groupoid built using $Z_{18}I$.

174. Find all the ideals of the neutrosophic groupoid W = {N(Z$_{21}$), *, (7, 14)}. Find atleast 3 neutrosophic subgroupoids. Is W an-idempotent neutrosophic groupoid?

175. Does there exist an infinite neutrosophic groupoid which is a P-groupoid?

176. Give an example of a S- neutrosophic groupoid of order 32.

177. Suppose G is a neutrosophic groupoid of finite order. Will every neutrosophic subgroupoid of G divide the order of G? Justify your claim.

178. What is the order of the neutrosophic groupoid G built using $Z_{21}$ where G = {N($Z_{21}$), *, (2, 3)}?

179. Obtain some interesting results about neutrosophic interval groupoid.

180. Give an example of a neutrosophic interval matrix Bol groupoid.

181. Show there exist neutrosophic interval matrix groupoids which do not satisfy any special identities.

182. Does there exists a neutrosophic polynomial interval groupoid of finite order?

183. Let $G = \left\{ \sum_{i=0}^{12} [0, a_i]x^i \mid a_i \in x_4I, *, (3,2) \right\}$ be a neutrosophic polynomial interval groupoid.
   a. Does G satisfy any of the special identities?
   b. Is G normal?
   c. Is G a Smarandache groupoid?
   d. Find atleast 2 subgroupoids of G.
   e. Does G have Smarandache ideals?

184. Let $G = \left\{ \sum_{i=0}^{\infty} [0, a_i]x^i \mid a_i \in Z_{15}I, *, (3,13) \right\}$ be a neutrosophic polynomial interval groupoid of infinite order.
   a. Find subgroupoids in G.
   b. Does there exists right ideals in G which are not left ideals?
   c. Is G a Smarandache groupoid?
   d. Does G satisfy any special identities?

185. Obtain some interesting results about neutrosophic interval groupoids of finite order.

186. Let G = {[0, $a_i$]| $a_i \in Z_{26}I$, *, (3, 6)} be a neutrosophic interval groupoid;
   a. Find the order of G.
   b. Does G have subgroupoids?
   c. Is G a Smarandache groupoid? Justify your claim.

187. Let G = {[$a_1, a_2, \ldots, a_9$]| $a_i = [0, x_i]$; $x_i \in Z_{25}I$, *, (12, 13)} be a neutrosophic row matrix interval groupoid;
   a. Does G satisfy any special identities?
   b. What is the order of G?
   c. Is G a normal groupoid?

188. Let $P = \left\{ \begin{bmatrix} a_1 & a_6 \\ a_2 & a_7 \\ a_3 & a_8 \\ a_4 & a_9 \\ a_5 & a_{10} \end{bmatrix} \middle| \begin{array}{l} a_i = [0, x_i]; \\ x_i \in Z_{12}I, *, (9,9); \\ 1 \le i \le 10 \end{array} \right\}$ be a neutrosophic matrix interval groupoid.

   a. Is P a finite order groupoid?
   b. Is P an alternative groupoid?
   c. Does P satisfy the Moufang identity?
   d. Find some subgroupoids and ideals in P.
   e. Is P a Smarandache groupoid?

189. Obtain some interesting results about neutrosophic polynomial interval groupoids built using $o(NC^+)$.

190. Does there exists a neutrosophic polynomial interval groupoid which satisfies both the Bol identity and Moufang identity?

191. Give an example of a neutrosophic polynomial interval groupoid which is a P-groupoid and an alternative groupoid.

192. Obtain some interesting properties about neutrosophic matrix interval groupoids built using $o(Z_{23}^I)$.

193. Obtain some interesting results about neutrosophic polynomial interval groupoids built using $o(N(Z_p))$, p a prime.

194. Enumerate the properties enjoyed by neutrosophic polynomial interval groupoids $G = \left\{ \sum_{i=0}^{10} [0, a_i] x^i \middle| \begin{array}{l} a_i \in Z_n I, *, (t,t); \\ 1 < t < n \end{array} \right\}$.

195. Does there exist a neutrosophic row matrix interval groupoid of order 26? Justify your claim.

196. Is it always true that in case of finite neutrosophic matrix groupoids the order of the subgroupoid divides the order of the groupoid? Justify your claim.

197. Let $G = \left\{ \sum_{i=0}^{n} [0,a_i]x^i \mid a_i \in Z^+I \cup \{0\}, *, (3,7) \right\}$ be the neutrosophic polynomial interval groupoid.
   a. Find subgroupoids of G.
   b. Is G a Smarandache groupoid? Justify.

198. Let $G = \left\{ \begin{bmatrix} a & a \\ a & a \end{bmatrix} \mid a \in o(Z_8I), *, (3,5) \right\}$ be a neutrosophic matrix interval groupoid.
   a. What is the order of G?
   b. Find subgroupoids of G.
   c. Does G satisfy any special identities?

199. Let $G = \left\{ \begin{bmatrix} a & a \\ a & a \\ a & a \\ a & a \end{bmatrix} \mid a \in o(N(Z_{11})), *, (6,5) \right\}$ be a neutrosophic matrix interval groupoid.
   a. Is G a finite order groupoid?
   b. Enumerate all the special identities satisfied by G.

200. Let $P = \left\{ \sum_{i=0}^{22} [0,a_i]x^i \mid a_i \in Z_5I, *, (3,2) \right\}$ be a neutrosophic polynomial interval groupoid.
   a. Find the order of P.
   b. Is G Smarandache strong Bol groupoid?
   c. Does G have Smarandache ideals?

201. Obtain some interesting properties about finite neutrosophic polynomial interval groupoids built using $Z_n$. $n < \infty$.

202. Let $G = \left\{ \sum_{i=0}^{12} [0, a_i] \mid a_i \in Z_{12}I, *, (9,4) \right\}$ be a neutrosophic polynomial interval groupoid and $G' = \left\{ \begin{bmatrix} a & a & a & a \\ a & a & a & a \end{bmatrix} \mid a \in Z_{24}I, *, (12,4) \right\}$ be a neutrosophic matrix interval groupoid. Does there exists a groupoid homomorphism of G to G'?

203. Does there exists a neutrosophic matrix interval groupoid built using $o(Z_nI)$ which has Smarandache conjugate subgroupoids?

204. Give an example of a pseudo simple neutrosophic matrix interval groupoid.

205. Give an example of a doubly simple neutrosophic polynomial interval groupoid.

# FURTHER READING

# INDEX

# ABOUT THE AUTHORS

**Dr.W.B.Vasantha Kandasamy** is an Associate Professor in the Department of Mathematics, Indian Institute of Technology Madras, Chennai. In the past decade she has guided 13 Ph.D. scholars in the different fields of non-associative algebras, algebraic coding theory, transportation theory, fuzzy groups, and applications of fuzzy theory of the problems faced in chemical industries and cement industries. She has to her credit 646 research papers. She has guided over 68 M.Sc. and M.Tech. projects. She has worked in collaboration projects with the Indian Space Research Organization and with the Tamil Nadu State AIDS Control Society. She is presently working on a research project funded by the Board of Research in Nuclear Sciences, Government of India. This is her $50^{th}$ book.

On India's 60th Independence Day, Dr.Vasantha was conferred the Kalpana Chawla Award for Courage and Daring Enterprise by the State Government of Tamil Nadu in recognition of her sustained fight for social justice in the Indian Institute of Technology (IIT) Madras and for her contribution to mathematics. The award, instituted in the memory of Indian-American astronaut Kalpana Chawla who died aboard Space Shuttle Columbia, carried a cash prize of five lakh rupees (the highest prize-money for any Indian award) and a gold medal.
She can be contacted at vasanthakandasamy@gmail.com
Web Site: http://mat.iitm.ac.in/home/wbv/public_html/

**Dr. Florentin Smarandache** is a Professor of Mathematics at the University of New Mexico in USA. He published over 75 books and 150 articles and notes in mathematics, physics, philosophy, psychology, rebus, literature.

In mathematics his research is in number theory, non-Euclidean geometry, synthetic geometry, algebraic structures, statistics, neutrosophic logic and set (generalizations of fuzzy logic and set respectively), neutrosophic probability (generalization of classical and imprecise probability). Also, small contributions to nuclear and particle physics, information fusion, neutrosophy (a generalization of dialectics), law of sensations and stimuli, etc. He can be contacted at smarand@unm.edu

**Dr. Moon Kumar Chetry** is presently working as a scientist in Defence Research and Development Organisation (DRDO). He has obtained his doctor of philosophy in Groupoids from department of mathematics, Indian Institute of Technology, Madras.